\documentclass[11pt]{article}

\usepackage{amsmath,amsthm}
\usepackage {latexsym}
\usepackage{amssymb}

\newcommand{\les}{\lesssim}

\newcommand{\const}{\mbox{\rm const}}
\newcommand{\beeq}{\begin{equation}}
\newcommand{\eneq}{\end{equation}}

\newcommand{\vecc}{{\underline v}}

\newcommand{\bear}{\begin{eqnarray}}
\newcommand{\eear}{\end{eqnarray}}
\newcommand{\beq}{\begin{equation}}
\newcommand{\eeq}{\end{equation}}

\newcommand{\supp}{\mbox{\rm supp}}
\newcommand{\spec}{{\rm spec}}

\newcommand{\half}{\frac{1}{2}}

\newcommand{\eps}{{\varepsilon}}
\newcommand{\R}{{\mathbb R}}

\newcommand{\Compl}{{\mathbb C}}

\newcommand{\calg}{\,{\mathfrak g}}
\newcommand{\calG}{{\mathcal G}}
\newcommand{\calF}{{\mathcal F}}
\newcommand{\calL}{{\mathcal L}}
\newcommand{\calS}{{\mathcal S}}
\newcommand{\calB}{{\mathcal B}}
\newcommand{\calJ}{{\mathcal J}}

\newcommand{\calM}{{\mathcal M}}

\newcommand{\calN}{{\mathcal N}}
\newcommand{\calP}{{\mathcal P}}

\newcommand{\sign}{\mbox{sign}}

\newcommand{\Laplace}{\triangle}
\newcommand{\Lapl}{\frac12{\triangle}}

\newcommand{\kato}{{\mathcal K}}
\newcommand{\trip}{|\!|\!|}
\newcommand{\la}{\langle}
\newcommand{\ra}{\rangle}

\renewcommand{\div}{\mbox{\rm div}}

\def\nn{\nonumber}
\def\calge1{\calg_{\vec{e_1}}}

\def\bm{\left( \begin{array}{cc}}
\def\endm{\end{array}\right)}
\def\ker{{\rm ker}}
\def\Ran{{\rm Ran}}
\def\Hil{{\mathcal H}}
\def\Dom{{\rm Dom}}

\newtheorem{theorem}{Theorem}
\newtheorem{lemma}[theorem]{Lemma}
\newtheorem{defi}[theorem]{Definition}
\newtheorem{cor}[theorem]{Corollary}
\newtheorem{prop}[theorem]{Proposition}

\theoremstyle{remark}
\newtheorem{remark}[theorem]{Remark}

\def\gatil{\tilde{\gamma}}

\def\il{\int\limits}

\def\y{{y}}
\def\x{{x}}

\def\Bp{B^+}
\def\Btp{\tilde{B}^+}

\def\norm[#1][#2]{\Vert #1 \Vert_{#2}}

\def\pizer{\pi^{(0)}}
\def\Uzer{U^{(0)}}
\def\Zzer{Z^{(0)}}
\def\Zone{Z^{(1)}}
\def\azer{\alpha^{(0)}}
\def\pione{\pi^{(1)}}
\def\Uone{U^{(1)}}
\def\aone{\alpha^{(1)}}
\def\pitwo{\pi^{(2)}}
\def\Utwo{U^{(2)}}
\def\Ztwo{Z^{(2)}}
\def\Zthree{Z^{(3)}}

\def\pithree{\pi^{(3)}}
\def\Uthree{U^{(3)}}

\def\Udis{U_{\rm dis}}
\def\Utildis{\tilde{U}_{\rm dis}}
\def\Uroot{U_{\rm root}}
\def\Uhyp{U_{\rm hyp}}
\def\Utilhyp{\tilde{U}_{\rm hyp}}
\def\Proot{P_{\rm root}}
\def\Pim{P_{\rm im}}
\def\xizer{\xi^{(0)}}
\def\etazer{\eta^{(0)}}
\def\Lin{{\mathcal L}}
\def\trip{|\!|\!|}

\def\Util3{\tilde{U}^{(3)}}

\def\xitil{\tilde{\xi}}
\def\etatil{\tilde{\eta}}
\def\calTzer{{\mathcal T}^{(0)}}

\begin{document}

\title {Stable manifolds for an orbitally unstable NLS}
\author{W. Schlag\thanks{The author was partially supported by the NSF grant
DMS-0300081 and a Sloan Fellowship}}

\maketitle

\section{Introduction}
\label{sec:intro}

We consider the cubic NLS in $\R^3$
\beeq
\label{eq:NLS}
 i\partial_t \psi + \Laplace \psi = -|\psi|^2 \psi.
\eneq
This equation is locally well-posed in $H^1(\R^3)=W^{1,2}(\R^3)$.
Let $\phi=\phi(\cdot,\alpha)$ be the ground state of
\beeq
\label{eq:ground}
-\Laplace \phi+\alpha^2\phi=\phi^3.
\eneq
By this we mean that $\phi>0$ and $\phi\in C^2(\R^3)$. It is a classical
fact (see Coffman~\cite{Cof}) that such solutions exist and are unique for the cubic
nonlinearity.  Moreover, they are radial and smooth.
Similar facts are known for more general nonlinearities,
see e.g.,~Berestycki and Lions~\cite{BerLio} for existence
and Kwon~\cite{Kwo} for uniqueness in greater generality.

Clearly,
$\psi=e^{it\alpha^2}\phi$ solves~\eqref{eq:NLS}.
We seek an $H^1$-solution $\psi$ of the form $\psi= W+R$  where
\begin{align}
W(t,x) &= e^{i\theta(t,x)} \phi(x-y(t),\alpha(t)) \label{eq:W}\\
\theta(t,x) &= v(t)\cdot x-\int_0^t (|v(s)|^2-\alpha^2(s))\, ds + \gamma(t) \label{eq:theta} \\
y(t) & = 2\int_0^t v(s)\,ds + D(t) \label{eq:y}
\end{align}
is the usual soliton with a moving set of parameters $\pi(t):=(\gamma(t),v(t),D(t),\alpha(t))$,
and $R$ is a small perturbation.
Performing a Galilei transform, we may assume that $W(0,x)=\phi(x,\alpha)=\alpha\phi(\alpha x,1)$
for some $\alpha>0$. The final equality holds because of the cubic nonlinearity.

With each $\alpha>0$ we associate the matrix operator
\beeq
\label{eq:Hilintro}
\Hil = \Hil(\alpha)=\bm -\Laplace + \alpha^2 -2\phi^2(\cdot,\alpha) & -\phi^2(\cdot,\alpha) \\
\phi^2(\cdot,\alpha) & \Laplace - \alpha^2 + 2\phi^2(\cdot,\alpha)
\endm. \eneq This operator is closed on the domain
$W^{2,2}(\R^3)\times W^{2,2}(\R^3)$ and its spectrum is known to
be located on $\R\cup i\R$ with essential spectrum equal to
$(-\infty,-\alpha^2]\cup [\alpha^2,\infty)$. As proved by
Weinstein~\cite{Wei1} and~\cite{Wei2}, $\Hil(\alpha)$ has a
rootspace of dimension eight at zero, and
$\ker(\Hil^3)=\ker(\Hil^2) \ne \ker(\Hil)$. In fact, $\dim \ker
(\Hil) =4$. On the other hand, any discrete spectrum different
from zero is known to consist entirely of eigenvalues whose
geometric and algebraic multiplicities coincide. Due to the $L^2$
supercritical nature of the problem, $\Hil(\alpha)$ {\em does}
have purely imaginary eigenvalues,  see Section~\ref{sec:imspec}
below. Moreover, due to the standard symmetries of the spectrum
(which follow from the commutation properties of $\Hil$ with the
Pauli matrices) we know that these purely imaginary eigenvalues
are symmetric with respect to the real axis, together with their
multiplicities.
 It is  well-known that the {\em supercritical equation}~\eqref{eq:NLS} is
 {\em orbitally unstable},
see Berestycki and Cazenave~\cite{BerCaz}. This is in contrast to
the {\em orbital stability} of the {\em subcritical equations}
that was proved by Cazenave and Lions~\cite{CazLio} and
Weinstein~\cite{Wei1}, \cite{Wei2}. In addition, let us mention
the general framework developed by Shatah, Strauss, and Grillakis
which gives a criterion for the orbital stability/instability of
solitary waves for large classes of Hamiltonian PDE.
See~\cite{Sha}, \cite{ShaStr}, \cite{Grill}, \cite{GSS1},
\cite{GSS2}. The recent book by Sulem and Sulem~\cite{SulSul}
contains an exposition of this work.

\medskip\noindent
For our main theorem we need to impose the following {\bf spectral
conditions :}

\smallskip
{\em $\Hil(\alpha)$ does not have any embedded eigenvalues in the
essential spectrum and the edges $\pm\alpha^2$ of the essential
spectrum are not resonances.}

\smallskip The property that the edges $\pm\alpha^2$
of the essential spectrum are not resonances\footnote{Strictly
speaking, this means that $\pm\alpha^2$ are neither eigenvalues
nor resonances, but we will show that eigenvalues cannot occur
there.} means that $(\Hil\pm\alpha^2)^{-1}$ exists as a bounded
operator between the weighted spaces
\[L^{2,\sigma}(\R^3)\times L^{2,\sigma}(\R^3)\to L^{2,-\sigma}(\R^3)\times L^{2,-\sigma}(\R^3)\]
where $\sigma>1$.

\begin{theorem}
\label{thm:main} Impose\footnote{By scaling invariance, if they
hold for one $\alpha>0$, then they hold for all $\alpha>0$. This
is due to the monomial nonlinearity.} the spectral conditions
 for all $\alpha>0$ and fix any $\alpha_0>0$. Then there exist a real-linear subspace
$\calS\subset L^2(\R^3)$ of co-dimension nine and a small
$\delta>0$ with the following properties: Let \beeq
\label{eq:calB} \calB:= \Big\{R_0\in L^2(\R^3)\:|\: \trip
R_0\trip:= \|R_0\|_{W^{1,2}\cap W^{1,1}}<\delta\Big\} \eneq and
let $\Sigma:=\{f\in L^2(\R^3)\:|\: \trip f\trip <\infty\}$. Then
there exists a  map $\Phi:\calB\cap \calS\to \Sigma$ with the
properties
\begin{align}
 \trip \Phi(R_0)\trip &\les \trip R_0\trip^2 \qquad \forall R_0\in
 \calB\cap \calS\label{eq:r0til} \\
\trip \Phi(R_0)-\Phi(\tilde{R}_0) \trip &\les  \delta\trip R_0-
\tilde{R}_0\trip \qquad\forall R_0,\tilde{R}_0\in  \calB\cap \calS
\label{eq:lip}
\end{align}
and so that for any $R_0\in~\calB\cap \calS$ the NLS~\eqref{eq:NLS}
has a global $H^1$ solution $\psi(t)$ for $t\ge0$
with initial condition $\psi(0)=\phi(\cdot,\alpha_0)+R_0+\Phi(R_0)$. Moreover,
\[ \psi(t)= W(t,\cdot)+R(t) \]
where $W$ as in~\eqref{eq:W} is governed by a path $\pi(t)$ of parameters which converges to some
terminal vector $\pi(\infty)$ such that $\sup_{t\ge0}|\pi(t)-\pi(\infty)|\les \delta^2$
and so that
\[ \|R(t)\|_{W^{1,2}} \les \delta,\qquad \|R(t)\|_\infty \les \delta t^{-\frac32} \]
for all $t>0$. Finally, there is scattering:
\[ R(t)= e^{it\Laplace} f_0 + o_{L^2}(1) \text{\ \ as $t\to\infty$ }\]
for some $f_0\in L^2(\R^3)$.
\end{theorem}

Concerning the spectral conditions, we remark that it is
well-known that imbedded eigenvalues and resonances are unstable
under perturbations. See the recent work by Cuccagna, Pelinovsky,
and Vougalter \cite{CPV}, \cite{CP} for precise statements to this
effect for matrix Schr\"odinger operators, as well as Costin,
Soffer~\cite{CosSof} for the scalar case. The proof of
Theorem~\ref{thm:main} does not rely on the specific structure of
the cubic nonlinearity and applies to other supercritical
nonlinearities as well. Hence it should be possible to formulate
Theorem~\ref{thm:main} without any spectral conditions for a
"generic perturbation" (in a suitable sense) of the cubic
nonlinearity. However, we have chosen to present
Theorem~\ref{thm:main} as stated in order not to obscure the main
ideas and we plan to return to the issue of generic perturbations
elsewhere. Moreover, in a forthcoming paper with
Krieger~\cite{KS}, the analogue of Theorem~\ref{thm:main} in one
dimension will be proved without any spectral conditions for
supercritical monomial nonlinearities. Rather, in that case these
spectral conditions can be proved by adapting some arguments of
Perelman~\cite{Pe2}. It is to be expected that the spectral
conditions therefore also hold in $\R^3$, although this yet needs
to be proved. On the other hand, it may be worth noting that we
reduce the question of threshold resonances for the matrix
operator $\Hil(\alpha)$ to the question whether or not the scalar
Schr\"odinger operator $L_-=-\Laplace+\alpha^2-\phi^2$ has a
resonance at $\alpha^2$ (more precisely, we show that if it does
not have an eigenvalue or resonance, then neither does the matrix
operator). The latter issue should in principle be accessible, say
by numerics.

The method of proof of Theorem~\ref{thm:main} also extends to the
case of more derivatives, i.e., $R_0\in W^{k,1}\cap W^{k,2}$, for
$k\ge2$, but we do not elaborate on this here. We will refer to
the fact that $R_0$ needs to satisfy
\beeq \label{eq:init_small}
 \|R_0\|_{W^{1,2}}+\|R_0\|_{W^{1,1}}<\delta
\eneq
as the {\em smallness condition.} It is not hard to see from
a close inspection of our proof that $W^{1,1}$ can be improved to
some $W^{1,p}$, $1<p<\frac43$, but we choose $p=1$ for simplicity.

To understand the origin of~$\calS$, we  need to introduce the
Riesz projections $P_s(\alpha)$ and
${\rm Id}-P_s(\alpha)$ (the index $s$ here stands for {\em stable}).
They are invariant under $\Hil(\alpha)$ and
\[ \spec(\Hil(\alpha)P_s(\alpha))=(-\infty,-\alpha^2]\cup [\alpha^2,\infty), \qquad
\spec(\Hil(\alpha)({\rm Id}-P_s(\alpha))= \{\pm i\sigma\}\cup
\{0\}.\] Here $\pm i\sigma$ are precisely the unique pair of
simple,  purely imaginary eigenvalues of $\Hil(\alpha)$,
$\sigma>0$. Finally, let $P^+_{u}(\alpha)$ be the Riesz projection
such that
\[ \spec(\Hil(\alpha) P_u^+(\alpha))= \{0\}\cup\{i\sigma\}.\]
The notation $P_u^+$ is meant to indicate the unstable modes as $t\to+\infty$.
The real-linear, finite-codimensional subspace $\calS$ above is
precisely the set of $R_0$ so that
\beeq
\label{eq:orth}
 P_u^+(\alpha_0) \binom{R_0}{\bar{R}_0}=0.
\eneq
The codimension of $\calS$ is simply the number of unstable
(or non-decaying) modes of the linearization: eight in the root
space and one exponentially unstable mode. The stable manifold
$\calM$ is  the surface described by the parameterization
$R_0\mapsto R_0+\Phi(R_0)$ where $R_0$ belongs to a small ball $\calB\cap\calS$ inside
of~$\calS$. The inequality~\eqref{eq:r0til} means that $\calS$ is
the tangent space to~$\calM$ at zero, whereas~\eqref{eq:lip}
expresses that~$\calM$  is given in terms of a  Lipschitz
parameterization. It is easy to see that it also the graph of a
Lipschitz map $\tilde{\Phi}: \calS\cap\calB\to\Sigma$. Indeed,
define $\tilde{\Phi}$ as
\[ R_0+P_{\calS}\Phi(R_0) \mapsto R_0 + \Phi(R_0), \]
where $\calS$ is the projection onto $\calS$ which is induced by the Riesz-projection
$I-P_u^+(\alpha_0)$ (the latter operates on $L^2\times L^2$, whereas we need only the
first coordinate of this projection, see Remark~\ref{rem:Jinv} below for the details of this).
The left-hand side is clearly in $\calS$. Moreover,
to see that this map is well-defined as well as Lipschitz, note that \eqref{eq:lip} implies that
\begin{align*}
& (1-C\delta)\trip R_1-R_0\trip
\le \trip R_1-R_0 + \Phi(R_1)-\Phi(R_0) \trip \le (1+C\delta) \trip R_1-R_0\trip \\
& (1-C\delta)\trip R_1-R_0\trip
\le \trip R_1-R_0 + P_{\calS}\Phi(R_1)-P_{\calS}\Phi(R_0) \trip \le (1+C\delta) \trip R_1-R_0\trip.
\end{align*}

Theorem~\ref{thm:main} should be understood as follows: The
instability result of Berestycki and Cazenave \cite{BerCaz} shows
that one can have finite time blow-up for initial data
$\psi_0=\phi(\cdot,\alpha)+R_0$ where $R_0$ can be made
arbitrarily small in any reasonable norm. On the other hand, one
may ask what the obstruction for global existence and even
stronger, for asymptotic stability, is in the orbitally unstable
case. Naturally, the first guess is the unstable subspace of the
linearized evolution~$e^{-it\Hil}$ with~$\Hil$ as
in~\eqref{eq:Hilintro}. This refers to the finite-dimensional
subspaces of those $f\in L^2(\R^3)\times L^2(\R^3)$ for which
$e^{-it\Hil}f$ does not decay locally as $t\to\infty$. Clearly,
this subspace contains all the (generalized) eigenspaces of all
eigenvalues of~$\Hil(\alpha)$ that lie on~$i\R^+\cup\{0\}$.
Conversely, we show in Sections~\ref{sec:lin_L2}
and~\ref{sec:lin_dis} that (for much more general systems
than~\eqref{eq:Hilintro}) \beeq \label{eq:L2intro} \sup_{t\ge0}
\big\|e^{-it\Hil(\alpha)}(I-P_u^+(\alpha))\big\|_{2\to2} < \infty.
\eneq While this bound was proved by Weinstein~\cite{Wei1}
and~\cite{Wei2} in the {\em subcritical case}, in which
$I-P_u^+=P_s$ only projects out the rootspace at zero, we are not
aware of a reference for~\eqref{eq:L2intro}. Moreover, adapting
the method of proof from~\cite{GolSch} from the scalar case
considered there to the matrix case needed here, we show that
\beeq \label{eq:Linftyintro}
\big\|e^{-it\Hil(\alpha)}(I-P_u^+(\alpha))\big\|_{1\to\infty} \les
t^{-\frac32} \eneq for all $t>0$. In view of~\eqref{eq:L2intro}
and \eqref{eq:Linftyintro}, it is conceivable that at least in
first approximation, the condition~\eqref{eq:orth} should ensure
stability. On the other hand, since it is based on linearization,
one needs to expect quadratic corrections. This is precisely the
content of~\eqref{eq:r0til}. So the statement of the theorem is
that after quadratic corrections, \eqref{eq:orth} gives the
desired asymptotic stability.

In the subcritical (monomial, say) case the linearized operator
\eqref{eq:Hilintro} has a root space of dimension eight and no
imaginary eigenvalues. Since there is asymptotic stability in this
case, one would naturally expect that the root space should not
contribute to the "bad directions" in the supercritical case,
i.e., that the codimension of the true stable manifold  should
really be one for~\eqref{eq:NLS}. This can indeed be achieved by
letting all symmetries of the NLS act on the manifold $\calM$ from
Theorem~\ref{thm:main}. In this way we regain eight dimensions
(six from the Galilei transforms, one from modulation, and one
from scaling) provided we show that these symmetries act
transversally to~$\calM$. This is done in the following theorem.

\begin{theorem}
\label{thm:main2} Impose the spectral conditions for all
$\alpha>0$ and fix any $\alpha_0>0$. Then there exist a small
$\delta>0$ and a Lipschitz manifold $\calN$ of size\footnote{This
means that $\calN$ is the graph of a Lipschitz map $\Psi$ with
domain $\calB\cap\tilde\calS$ where $\tilde\calS$ is a subspace of
codimension one and with $\calB$ as in \eqref{eq:calB}.} $\delta$
and codimension one so that $\phi(\cdot,\alpha_0)\in\calN$ with
the following property: for any choice of initial data
$\psi(0)\in\calN$ the NLS~\eqref{eq:NLS} has a global $H^1$
solution $\psi(t)$ for $t\ge0$.  Moreover,
\[ \psi(t)= W(t,\cdot)+R(t) \]
where $W$ as in~\eqref{eq:W} is governed by a path $\pi(t)$ of
parameters so that $|\pi(0)-(0,0,0,\alpha_0)|\les\delta$ and which
converges to some terminal vector $\pi(\infty)$ such that
$\sup_{t\ge0}|\pi(t)-\pi(\infty)|\les \delta^2$. Finally,
\[ \|R(t)\|_{W^{1,2}} \les \delta,\qquad \|R(t)\|_\infty \les \delta t^{-\frac32} \]
for all $t>0$ and there is scattering:
\[ R(t)= e^{it\Laplace} f_0 + o_{L^2}(1) \text{\ \ as $t\to\infty$ }\]
for some $f_0\in L^2(\R^3)$.
\end{theorem}

\noindent   This result raises the question of deciding the
behavior of solutions with initial data $\phi(0)\in \calB\setminus
\calN$. One reasonable possibility would be that any initial data
$\phi(0)\in\calB\setminus\calN$ leads to blow up either as time goes to $+\infty$
or $-\infty$.

The motivation for studying these question was two-fold. First, it
is natural to seek stable manifolds for unstable problems. There
is a substantial ODE literature in this context, but the PDE case
is much less studied. One example where stable manifolds have been
constructed for Schr\"odinger equations is the recent paper by
Tsai and Yau~\cite{TY}, although the problem they consider is
rather different from the one studied here.

Second, there is a large literature concerning {\em asymptotic
stability} questions for {\em subcritical equations}, see the
papers by Soffer, Weinstein~\cite{SofWei1}, \cite{SofWei2},
Pillet, Wayne~\cite{PW}, and Buslaev, Perelman~\cite{BP1},
\cite{BP2}, as well as Cuccagna~\cite{Cuc}. Recently, there has
also been some work on the multi-soliton case, see Rodnianski,
Soffer, and the author~\cite{RSS1}, \cite{RSS2}, as well as
Perelman~\cite{Pe3}. Most of these papers are based on a
Lyapunov-Schmidt reduction, i.e., on splitting the evolution into
a finite dimensional part and a complementary part on which the
linearized evolution needs to be dispersive. In the subcritical
case the finite dimensional part is exactly the root space,
assuming as one usually does, that there are no other eigenvalues
than zero. The dimension of this finite dimensional part then
coincides with the number of parameters in our soliton (namely
eight). This is natural, since both are intimately related to the
family of symmetries of the NLS~\eqref{eq:NLS}, see~\cite{Wei1}.
This fact allows one to write down a system of ODEs for the
parameter paths called the modulation equations which ensure that
the finite dimensional part is  not present at all. In the context
of asymptotic stability of solitons this method was first
implemented by Soffer and Weinstein~\cite{SofWei1},
\cite{SofWei2}. Our second motivation for Theorem~\ref{thm:main}
was the question to what extent these asymptotic stability methods
also apply to the $L^2$ supercritical case, which is orbitally
unstable. As explained before, for the case of supercritical
monomial nonlinearities the time evolution of the linearized
problem has exponentially growing solutions. Needless to say,
these modes cannot be controlled by the modulation equations.
Rather, they need to be eliminated by a different mechanism. To
first order, the unstable modes of the linearization need to be
removed from the initial conditions. This is the origin
of~\eqref{eq:orth}. However, this is only an approximation and
quadratic corrections need to be made.

\section{The linearization, Galilei transforms, and $\calJ$-invariance}
\label{sec:ansatz}

As in \cite{RSS2} we require the soliton paths in~\eqref{eq:W} to be {\em admissible}.
The constant $\delta$ which appears in the following definition is the same small
constant as in Theorem~\ref{thm:main}. It will be specified later.

\begin{defi}
\label{def:adm}
We say that a path $\pi:[0,\infty)\to\R^8$ with $\pi(t):=(\gamma(t),{v}(t),D(t),\alpha(t))$ is {\em admissible} provided it belongs to ${\rm Lip}([0,\infty),\R^8)$, the limit
\[ \lim_{t\to\infty} (\gamma(t),{v}(t),D(t),\alpha(t)) =: (\gamma(\infty),{v}(\infty),D(\infty),
\alpha(\infty)) \]
exists, and such that the entire path lies within a $\delta$-neighborhood of those limiting values for
all times $t\ge0$.
Moreover, we assume that
\begin{align*}
& |{v}(t)-{v}(\infty)| = o(t^{-1}) \text{\ \ as\ \ }t\to\infty\\
& \int_0^\infty \int_s^\infty (|\dot{{v}}(\sigma)|+|\dot{\alpha}(\sigma)|)\,d\sigma\, ds < \infty.
\end{align*}
Under these conditions,
define a constant parameter vector $\pi_\infty=(\gamma_\infty,{v}_\infty,D_\infty,\alpha_\infty)$ as
\begin{align*}
\gamma_\infty &:= \gamma(\infty) + 2\int_0^\infty \int_s^\infty ({v}(\sigma)\cdot \dot{{v}}(\sigma)
-\alpha(\sigma)\dot{\alpha}(\sigma))\,d\sigma\,ds \\
{v}_\infty &:= {v}(\infty)\\
D_\infty &:= D(\infty) -2\int_0^\infty \int_s^\infty \dot{{v}}(\sigma)\, d\sigma\,ds \\
\alpha_\infty &:= \alpha(\infty) \\
\end{align*}
With these parameters, define the usual {\em Galilei transform} to be
\beeq
\label{eq:gal}
 \calg_\infty(t) = e^{i(\gamma_\infty+{v}_\infty\cdot x-t|{v}_\infty|^2)}
\, e^{-i(2t{v}_\infty+D_\infty)\cdot \vec{p}},
\eneq
where $\vec{p}:=-i\nabla$.
\end{defi}

The action of $\calg_\infty(t)$ on functions is
\[ (\calg_\infty(t)f)(x)= e^{i(\gamma_\infty+{v}_\infty\cdot x-t|{v}_\infty|^2)} f(x-2t{v}_\infty-D_\infty),\]
they are unitary on $L^2$, isometries on all $L^p$,
and the commutation property $e^{it\Laplace}\calg_\infty(0)=\calg_\infty(t)e^{it\Laplace}$ holds.
The inverse of $\calg_\infty(t)$ is
\[ \calg_\infty(t)^{-1} = e^{i(2t{v}_\infty+D_\infty)\cdot \vec{p}}\,e^{-i(\gamma_\infty+{v}_\infty\cdot x-t|{v}_\infty|^2)}=e^{-i(\gamma_\infty+{v}_\infty\cdot D_\infty+{v}_\infty\cdot x+t|{v}_\infty|^2)}\,
e^{i(2t{v}_\infty+D_\infty)\cdot \vec{p}}.
\]
Moreover, the Galilei transform~\eqref{eq:gal} generates an eight-parameter family of solitons:
Let $\phi(\cdot,\alpha_\infty)$ be the ground state of~\eqref{eq:ground} with $\alpha=\alpha_\infty$.
Then
\beeq
\label{eq:Winfty}
 W_\infty(t,\cdot) = \calg_\infty(t) [e^{it\alpha_\infty^2}\phi(\cdot,\alpha_\infty)]
\eneq
solves \eqref{eq:NLS}, where $W_\infty$ is a soliton as introduced in~\eqref{eq:W}
but with the constant parameter path $\pi_\infty$.
For future reference, let
\beeq
\label{eq:inf_path}
y_\infty(t) := 2t{v}_\infty + D_\infty,\; \theta_\infty(t,x) := {v}_\infty\cdot x - t(|{v}_\infty|^2-\alpha_\infty^2)+\gamma_\infty.
\eneq
We will use repeatedly that
\beeq
\label{eq:ti}
\theta_\infty(t,x+y_\infty)=t(|v_\infty|^2+\alpha_\infty^2)+
v_\infty\cdot(x+D_\infty)+\gamma_\infty.
\eneq
With these notations, $W_\infty$ in~\eqref{eq:Winfty} takes the form
\[ W_\infty(t,x) = e^{i\theta_\infty(t,x)}\phi(x-y_\infty(t),\alpha_\infty). \]

\begin{lemma}
\label{lem:rho_infty}
Suppose $\pi$ is an admissible path and
let $\theta,y$ and $\theta_\infty,y_\infty$ be as in \eqref{eq:theta}, \eqref{eq:y},
 and \eqref{eq:inf_path}, respectively.
Define
\beeq
\label{eq:rhoinf}
 \rho_\infty(t,x) := \theta(t,x+y_\infty)-\theta_\infty(t,x+y_\infty).
\eneq
Then
\[ \rho_\infty(t,x) = (1+|x|)o(1), \quad y(t)-y_\infty(t)=o(1)\]
as $t\to\infty$.
\end{lemma}
\begin{proof}
In view of the definition of $\pi_\infty$,
\begin{align}
\theta(t,x+y_\infty)-\theta_\infty(t,x+y_\infty) &= v(t)\cdot(x+2tv_\infty+D_\infty) -
\int_0^t (|v(s)|^2-\alpha^2(s))\,ds + \gamma(t) \nn \\
&\quad  - v_\infty\cdot(x+2tv_\infty+D_\infty) + t(|v_\infty|^2-\alpha_\infty^2)-\gamma_\infty \nn \\
& = ({v}(t)-{v}_\infty)\cdot(x+2t{v}_\infty+D_\infty) +2\int_0^\infty \int_s^\infty ({v}\cdot\dot{{v}}-\alpha\dot{\alpha})(\sigma)\,d\sigma ds \nn \\
&\quad -\gamma_\infty + \gamma(t) -2\int_t^\infty\int_s^\infty ({v}\cdot\dot{{v}}-\alpha\dot{\alpha})(\sigma)\,d\sigma ds \nn\\
& =  ({v}(t)-{v}_\infty)\cdot(x+2t{v}_\infty+D_\infty) -2\int_t^\infty\int_s^\infty ({v}\cdot\dot{{v}}-\alpha\dot{\alpha})(\sigma)\,d\sigma ds  \nn \\
&\quad -\gamma(\infty) + \gamma(t).
\label{eq:rhorep}
\end{align}
Using Definition~\ref{def:adm}  implies the desired
bound on $\rho_\infty$. As for $y(t)-y_\infty(t)$, the definition of $D_\infty$ implies that
\beeq
\label{eq:yrep}
y_\infty(t)-y(t) = 2tv_\infty+D_\infty-2\int_0^t v(s)\,ds - D(t)=D(\infty)-D(t)-2\int_t^\infty\int_s^\infty \dot{v}(\sigma)\,d\sigma\,ds,
\eneq
which goes to zero as $t\to\infty$.
\end{proof}

Recall from Section~\ref{sec:intro} that we seek an $H^1$ solution $\psi(t)$ of the cubic
NLS~\eqref{eq:NLS} of the
form $\psi=W+R$. The following standard lemma derives the equation for $R$, or rather for
the vector $\binom{R}{\bar{R}}$.

\begin{lemma}
\label{lem:Zeq}
Assume that $\pi(t)=(\gamma(t),{v}(t),D(t),\alpha(t))$ is admissible, see Definition~\ref{def:adm}, and
let $W=W(t,x)$ be as in~\eqref{eq:W}. Let $0<T\le \infty$.
Then $\psi\in C([0,T),H^1(\R^3))\cap C^1([0,T),H^{-1}(\R^3))$ solves~\eqref{eq:NLS}
with $\psi=W+R$ iff $Z=\binom{R}{\bar{R}}$ solves the equation
\begin{align}
 i\partial_t Z + \bm \Laplace + 2|W|^2 & W^2 \\ -\bar{W}^2 & -\Laplace -2|W|^2 \endm Z
 & = \dot{v} \binom{-xe^{i\theta} \phi(\cdot-y,\alpha)}{xe^{-i\theta} \phi(\cdot-y,\alpha)} +
\dot{\gamma} \binom{-e^{i\theta} \phi(\cdot-y,\alpha)}{e^{-i\theta} \phi(\cdot-y,\alpha)}\nn \\
&\quad + i\dot{\alpha}
\binom{e^{i\theta}\partial_\alpha \phi(\cdot-y,\alpha)}{e^{-i\theta}\partial_\alpha \phi(\cdot-y,\alpha)}
 + i\dot{D}
\binom{-e^{i\theta}\nabla \phi(\cdot-y,\alpha)}{-e^{-i\theta}\nabla \phi(\cdot-y,\alpha)} \nn \\
&\quad + \binom{-2|R|^2W-\bar{W}R^2-|R|^2R}{2|R|^2\bar{W}+W\bar{R}^2+|R|^2\bar{R}} \label{eq:Zsys}
\end{align}
in the sense of $C([0,T),H^1(\R^3)\times H^1(\R^3))\cap C^1([0,T),H^{-1}(\R^3)\times H^{-1}(\R^3))$.
Here $y$ and $\theta$ are the functions from~\eqref{eq:y} and~\eqref{eq:theta}, and $\alpha=\alpha(t)$.
For future reference, we denote the matrix operator on the left-hand side
of~\eqref{eq:Zsys} by $-\Hil(\pi(t))$, i.e.,
\beeq
\label{eq:Hpidef}
\Hil(\pi(t)) := \bm - \Laplace - 2|W|^2 & -W^2 \\ \bar{W}^2 & \Laplace+ 2|W|^2 \endm.
\eneq
\end{lemma}
\begin{proof}
Let $\phi=\phi(\cdot,\alpha(t))$ for ease of notation.
Direct differentiation shows that $W(t,x)$ satisfies
\[ i\partial_tW+\Laplace W = -|W|^2W -W(\dot{v}x+\dot{\gamma})-ie^{i\theta}\nabla \phi\cdot \dot{D}+
ie^{i\theta}\dot{\alpha}\partial_\alpha \phi.\]
Hence $W+R$ is a solution of \eqref{eq:NLS} iff
\[ i\partial_t R + \Laplace R = -2|W|^2 R - 2|R|^2 W -|R|^2R - W^2\bar{R} -\bar{W}R^2
-e^{i\theta}\phi(\dot{v}x+\dot{\gamma})
- ie^{i\theta}\nabla \phi\cdot \dot{D} +
ie^{i\theta}\dot{\alpha}\partial_\alpha \phi. \]
Joining this equation with its conjugate leads to the system~\eqref{eq:Zsys}.
Conversely, if $Z(0)$ is of the form
\[ Z(0) = \binom{Z_1(0)}{\overline{Z_1(0)}}, \]
and $Z(t)$ solves~\eqref{eq:Zsys}, then $Z(t)$ remains of this form for all times.
This is simply the statement that the system~\eqref{eq:Zsys} is invariant under the transformation
\beeq
\label{eq:Jdef} \calJ:f\mapsto \overline{Jf}, \qquad J=\bm 0&1 \\1&0 \endm, \quad f=\binom{f_1}{\overline{f_1}},
\eneq
which can be checked by direct verification. This fact allows us to go back from the system
to the scalar equation.
\end{proof}

The issue of $\calJ$-invariance is of great importance. The $\calJ$-invariant vectors
in $L^2(\R^3)\times L^2(\R^3)$
form a real-linear subspace, namely
\[ \Big\{ \binom{f}{\bar{f}}\:\Big|\: f\in L^2(\R^3)\Big\}. \]
Writing $f=f_1+if_2$ it can be seen to be isomorphic to the subspace
\[ \Big\{ \binom{f_1}{f_2}\:\Big|\: f_1,f_2\in L^2(\R^3), \; f_1,f_2 \text{\ are real-valued}\Big\}, \]
which is clearly linear, but only over $\R$.
We need to insure that all vectorial solutions we construct
belong to this subspace. Only then is it possible to revert to the scalar NLS~\eqref{eq:NLS}.

As usual, it will be convenient to transform~\eqref{eq:Zsys} to a stationary frame.
In addition, a modulation will be performed. The details are as follows.

\begin{lemma} \label{lem:UPDE}
Let $\pi(t)$ be an admissible path and let $\pi_\infty$ be the constant vector associated
with it as in Definition~\ref{def:adm}. Given a vector $Z=\binom{Z_1}{Z_2}$, introduce
$U$, as well as $M(t), \calG_\infty (t)$ by means of
\beeq
\label{eq:Udef}
 U(t)=\bm e^{i\omega(t)} & 0 \\ 0 & e^{-i\omega(t)} \endm \binom{\calg_\infty(t)^{-1} Z_1(t)}
{\overline{\calg_\infty(t)^{-1} \overline{Z_2(t)}}}
= M(t) \calG_\infty (t) Z(t),
\eneq
where $\omega(t)=-t\alpha_\infty^2$.
Then $Z(t)$ solves \eqref{eq:Zsys} in the $H^1$ sense iff
 $U=\binom{U_1}{U_2}$ as in~\eqref{eq:Udef} satisfies the following PDE in the $H^1$ sense
(with $\phi_\infty=\phi(\cdot,\alpha_\infty)$):
\beeq
\label{eq:UPDE}
i\dot{U}(t) + \bm \Laplace +2\phi_\infty^2-\alpha_\infty^2 & \phi_\infty^2 \\ -\phi_\infty^2 & -\Laplace -2\phi_\infty^2 +\alpha_\infty^2 \endm U = \dot{\pi}\partial_\pi W(\pi) +N(U,\pi) + VU
\eneq
where we are using the formal notations
\begin{align}
 V = V(t) &:= \bm 2(\phi_\infty^2(x)-\phi^2(x+y_\infty-y)) & \phi_\infty^2(x)-e^{2i\rho_\infty} \phi^2(x+y_\infty-y) \\
-\phi_\infty^2(x)+ e^{-2i\rho_\infty} \phi^2(x+y_\infty-y) & -2(\phi_\infty^2(x)-\phi^2(x+y_\infty-y)) \endm \label{eq:V}\\
 \dot{\pi}\partial_\pi W(\pi) &:= \dot{{v}} \binom{-(x+y_\infty)e^{i\rho_\infty} \phi(x+y_\infty-y)}{(x+y_\infty) e^{-i\rho_\infty} \phi(x+y_\infty-y)} + \dot{\gamma} \binom{-e^{i\rho_\infty}\phi(x+y_\infty-y)}{e^{-i\rho_\infty}\phi(x+y_\infty-y)} \label{eq:vdotetc} \\
 & \quad +i\dot{\alpha} \binom{e^{i\rho_\infty}\partial_\alpha\phi(x+y_\infty-y)}{e^{-i\rho_\infty}\partial_\alpha\phi(x+y_\infty-y)}+i\dot{D} \binom{e^{i\rho_\infty}\nabla\phi(x+y_\infty-y)}
{e^{-i\rho_\infty}\nabla\phi(x+y_\infty-y)} \nn \\
N(U,\pi) &:= \binom{-2|U_1|^2e^{i\rho_\infty}\phi(x+y_\infty-y)-e^{-i\rho_\infty}\phi(x+y_\infty-y)U_1^2-|U_1|^2U_1}{2|U_2|^2e^{-i\rho_\infty}\phi(x+y_\infty-y)+e^{i\rho_\infty}\phi(x+y_\infty-y)U_2^2+|U_2|^2U_2}. \label{eq:NUpi}
\end{align}
Here $\rho_\infty=\rho_\infty(t,x)$ is as in Lemma~\ref{lem:rho_infty}, $\phi(x+y_\infty-y)=\phi(x+y_\infty(t)-y(t),\alpha(t))$, and $\omega=\omega(t)$ is as in~\eqref{eq:Udef}.
Finally, $Z$ is $\calJ$-invariant iff $U$ is $\calJ$-invariant, and $U$ is $\calJ$-invariant
iff $U(0)$ is $\calJ$-invariant.
\end{lemma}
\begin{proof}
Throughout this proof we will adhere to the convention that
$\phi=\phi(\cdot,\alpha(t))$ whereas $\phi_\infty=\phi(\cdot,\alpha_\infty)$).
Write the equation~\eqref{eq:Zsys} for $Z$ in the form
\beeq
\label{eq:Zeq}
 i\partial_t Z -\Hil_\infty Z = F + (\Hil(\pi(t))-\Hil_\infty)Z
\eneq
where
\beeq
\label{eq:Hinfty}
 \Hil_\infty = \bm -\Laplace - 2|W_\infty|^2 & -W_\infty^2 \\ \bar{W}_\infty^2 & \Laplace +2|W_\infty|^2 \endm,
\eneq
see \eqref{eq:Winfty} and~\eqref{eq:inf_path}.
With $\calG_\infty(t)$ defined as in~\eqref{eq:Udef}, and with $p=-i\nabla$,
\begin{align}
i\frac{d}{dt} \calG_\infty(t)f &= \binom{i\dot{\calg}_\infty(t)^{-1} f_1}{i\overline{\dot{\calg}_\infty(t)^{-1}
\overline{f_2}}} = \binom{-(2{v}_\infty\cdot p + |{v}_\infty|^2)\calg_\infty(t)^{-1} f_1}{-(2{v}_\infty\cdot p - |{v}_\infty|^2)\overline{\calg_\infty(t)^{-1} \overline{f_2(t)}}} \nn \\
& = \bm -(2{v}_\infty\cdot p + |{v}_\infty|^2) & 0 \\
0 & -(2{v}_\infty\cdot p - |{v}_\infty|^2) \endm \calG_\infty(t) f
\label{eq:gal_der}
\end{align}
for any $f=\binom{f_1}{f_2}$.
Furthermore,
\begin{align}
& M(t)\calG_\infty(t) \Hil_\infty \binom{f_1}{f_2} \label{eq:comm_GH} \\
 &= \bm e^{i\omega(t)} & 0 \\ 0 & e^{-i\omega(t)} \endm
\binom{-(\Laplace +2\phi_\infty^2)\calg_\infty(t)^{-1}f_1 -
\phi_\infty^2e^{2i\theta_\infty(t,x+y_\infty)}\calg_\infty(t)^{-1} f_2}{\phi_\infty^2e^{-2i\theta_\infty(t,x+y_\infty)}
\overline{\calg_\infty(t)^{-1}\overline{f_1}} + (\Laplace +2\phi_\infty^2)\overline{\calg_\infty(t)^{-1}\overline{f_2}}}
\nn \\
& \quad - \bm e^{i\omega(t)} & 0 \\ 0 & e^{-i\omega(t)} \endm
\bm -|{v}_\infty|^2+2i{v}_\infty\cdot\nabla & 0 \\ 0 & |{v}_\infty|^2+2i{v}_\infty\cdot\nabla \endm \calG_\infty(t) \binom{f_1}{f_2} \nn \\
&= \bm e^{i\omega(t)} & 0 \\ 0 & e^{-i\omega(t)} \endm
\binom{-(\Laplace +2\phi_\infty^2)\calg_\infty(t)^{-1}f_1 -
\phi_\infty^2
e^{2i[\theta_\infty(t,x+y_\infty)-(t|{v}_\infty|^2+{v}_\infty\cdot( x+D_\infty)+\gamma_\infty)]}\overline{\calg_\infty(t)^{-1}\overline{f_2}}}{\phi_\infty^2e^{-2i[\theta_\infty(t,x+y_\infty)-
(t|{v}_\infty|^2+{v}_\infty\cdot( x+D_\infty)+\gamma_\infty)]}\calg_\infty(t)^{-1}f_1 + (\Laplace +2\phi_\infty^2)\overline{\calg_\infty(t)^{-1}\overline{f_2}}}
\nn \\
& \quad - \bm e^{i\omega(t)} & 0 \\ 0 & e^{-i\omega(t)} \endm
\bm -|{v}_\infty|^2+2i{v}_\infty\cdot\nabla & 0 \\ 0 & |{v}_\infty|^2+2i{v}_\infty\cdot\nabla \endm \calG_\infty(t) \binom{f_1}{f_2}. \nn
\end{align}
Now
\[ \theta_\infty(t,x+y_\infty)-(t|v_\infty|^2+v_\infty\cdot( x+D_\infty)+\gamma_\infty) =
t\alpha_\infty^2,\]
see \eqref{eq:ti}.
Hence, by the definition of $\omega(t)$ (and dropping the argument $t$ from $M$ and $\calG_\infty$ for
simplicity),
\begin{align}
\eqref{eq:comm_GH} &= \bm e^{i\omega(t)} & 0 \\ 0 & e^{-i\omega(t)} \endm
\bm -(\Laplace +2\phi_\infty^2) & -
\phi_\infty^2 e^{2it\alpha_\infty^2} \\
\phi_\infty^2 e^{-2it\alpha_\infty^2} & \Laplace +2\phi_\infty^2 \endm \calG_\infty
f \nn \\
& \quad - \bm e^{i\omega(t)} & 0 \\ 0 & e^{-i\omega(t)} \endm
\bm -|{v}_\infty|^2-2{v}_\infty\cdot p & 0 \\ 0 & |{v}_\infty|^2-2{v}_\infty\cdot p \endm \calG_\infty f  \nn \\
&  =  \bm -\Laplace - 2\phi_\infty^2 & -\phi_\infty^2 \\
\phi_\infty^2  & \Laplace +2\phi_\infty^2 \endm M\calG_\infty f - \bm -|{v}_\infty|^2-2{v}_\infty\cdot p & 0 \\ 0 & |{v}_\infty|^2-2{v}_\infty\cdot p \endm \calG_\infty f. \label{eq:H_phi}
\end{align}
Denote the first matrix operator in~\eqref{eq:H_phi} by $\Hil_\phi$.
Hence, in combination with \eqref{eq:gal_der} one concludes from~\eqref{eq:Zeq} that
\begin{align*}
i\dot{U} &= i\dot{M}\calG_\infty Z + iM\dot{\calG}_\infty Z + M\calG_\infty \Hil_\infty +
M\calG_\infty(F+(\Hil(\pi(t))-\Hil_\infty)Z) \nn \\
&= \bm -\dot{\omega} & 0 \\ 0 & \dot{\omega} \endm M\calG_\infty Z +
\bm -(2{v}_\infty\cdot p + |{v}_\infty|^2) & 0 \\
0 & -(2{v}_\infty\cdot p - |{v}_\infty|^2) \endm M\calG_\infty Z + \Hil_\phi M\calG_\infty Z \nn\\
&\quad + \bm |{v}_\infty|^2+2{v}_\infty\cdot p & 0 \\ 0 & -|{v}_\infty|^2+2{v}_\infty\cdot p \endm M\calG_\infty Z + M\calG_\infty(F+(\Hil(\pi(t))-\Hil_\infty)Z)  \nn \\
& = \bm -\Laplace +\alpha_\infty^2 - 2\phi_\infty^2 & -\phi_\infty^2 \\
\phi_\infty^2  & \Laplace -\alpha_\infty^2 +2\phi_\infty^2 \endm U(t) + M\calG_\infty(F+(\Hil(\pi(t))-\Hil_\infty)\calG_\infty^{-1}M^{-1} U).
\end{align*}
It remains to compute the terms
\begin{align}
 \dot{\pi}\partial_\pi W(\pi)+N(U,\pi) &= M(t)\calG_\infty(t) F(t) \label{eq:Fdef} \\
 V &= M(t)\calG_\infty(t)(\Hil(\pi(t))-\Hil_\infty)\calG_\infty(t)^{-1}M(t)^{-1} \label{eq:Vdef}
\end{align}
In view of \eqref{eq:Zsys}, one has
\begin{align*}
 F &= \dot{v} \binom{-xe^{i\theta} \phi(x-y)}{xe^{-i\theta} \phi(x-y)} +
\dot{\gamma} \binom{-e^{i\theta} \phi(x-y)}{e^{-i\theta} \phi(x-y)} \\
&\quad + i\dot{\alpha} \binom{e^{i\theta}\partial_\alpha \phi(x-y)}{e^{-i\theta}\partial_\alpha \phi(x-y)}
 + i\dot{D} \binom{-e^{i\theta}\nabla \phi(x-y)}{-e^{-i\theta}\nabla \phi(x-y)}
 + \binom{-2|R|^2W-\bar{W}R^2-|R|^2R}{2|R|^2\bar{W}+W\bar{R}^2+|R|^2\bar{R}}.
\end{align*}
Now
\[ \theta(t,x+y_\infty)-(\alpha_\infty^2t+v_\infty\cdot(x+D_\infty)+t|v_\infty|^2+\gamma_\infty)
=\theta(t,x+y_\infty)-\theta_\infty(t,x+y_\infty)=\rho_\infty(t,x),\]
see \eqref{eq:ti} and Lemma~\ref{lem:rho_infty}.
Thus, the first term of $M\calG F$ is
\begin{align*}
 & \dot{v} \bm e^{i\omega} & 0\\ 0 & e^{-i\omega} \endm \binom{-(x+y_\infty)e^{i\theta(t,x+y_\infty)}\, e^{-i(t|v_\infty|^2+v_\infty\cdot (x+D_\infty)+\gamma_\infty)} \phi(x+y_\infty-y)}{(x+y_\infty)
e^{-i\theta(t,x+y_\infty)}\, e^{i(t|v_\infty|^2+v_\infty\cdot(x+D_\infty)+\gamma_\infty)} \phi(x+y_\infty-y)} \\
&= \dot{v}  \binom{-(x+y_\infty)e^{i\rho_\infty(t,x)} \phi(x+y_\infty-y)}{(x+y_\infty)e^{-i\rho_\infty(t,x)} \phi(x+y_\infty-y)}.
\end{align*}
This gives the $\dot{v}$ term in~\eqref{eq:vdotetc}. The other terms involving $\dot{\alpha},
\dot{\gamma}$, and $\dot{D}$ are treated similarly, and we skip the details. The cubic term in~\eqref{eq:Zsys} is also easily transformed, and it leads to the nonlinear
term $N(U,\pi)$in~\eqref{eq:NUpi}. We skip that calculation as well.
Finally, it remains to transform $\Hil(\pi(t))-\Hil_\infty$. One has
\[
 \Hil(\pi(t))-\Hil_\infty =
 \bm 2(\phi_\infty^2(\cdot-y_\infty)-\phi^2(\cdot-y))
& e^{2i\theta_\infty}\phi_\infty^2(\cdot-y_\infty)  - e^{2i\theta} \phi^2(\cdot-y)  \\
e^{-2i\theta} \phi^2(\cdot-y) - e^{-2i\theta_\infty} \phi_\infty^2(\cdot-y_\infty) & - 2(\phi_\infty^2(\cdot-y_\infty)-\phi^2(\cdot-y))
\endm
\]
where $\phi_\infty=\phi(\cdot-y_\infty(t),\alpha_\infty)$, $\phi=\phi(\cdot-y(t),\alpha(t))$
for simplicity.
It is easy to check that
\[
\calG_\infty(t) (\Hil(\pi(t))-\Hil_\infty) = \bm 2(\phi_\infty^2(x)-\phi^2(x+y_\infty-y)) & * \\
-e^{-2it\alpha_\infty^2} (\phi_\infty^2(x)-e^{-2i\rho_\infty} \phi^2(\cdot+y_\infty-y)) & * \endm \calG_\infty(t).
\]
After conjugation by the matrix $M(t)$ this takes the desired form~\eqref{eq:V} and we are done.
For the final statements concerning $\calJ$-invariance, observe first that
the transformation~\eqref{eq:Udef} from $Z$ to $U$ preserves $\calJ$-invariance. Second,  the equation~\eqref{eq:UPDE} is $\calJ$-invariant, which shows that it suffices to assume the
$\calJ$-invariance of $U(0)$ to
guarantee it for all $t\ge0$.  To check the $\calJ$-invariance of~\eqref{eq:UPDE}, note that the
right-hand side of~\eqref{eq:UPDE} transforms like
\[ \calJ[\dot{\pi}\partial_\pi W(\pi) +N(U,\pi) + VU] =
- [\dot{\pi}\partial_\pi W(\pi) +N(\calJ U,\pi) + V\cal J U],
\]
while the left-hand side transforms as follows:
\begin{align*}
& \calJ [ i\dot{U}(t) + \bm \Laplace +2\phi_\infty^2-\alpha_\infty^2 & \phi_\infty^2 \\ -\phi_\infty^2 & -\Laplace -2\phi_\infty^2 +\alpha_\infty^2 \endm U ] \\
& = -i\dot{\calJ U}(t) - \bm \Laplace +2\phi_\infty^2-\alpha_\infty^2 & \phi_\infty^2 \\ -\phi_\infty^2 & -\Laplace -2\phi_\infty^2 +\alpha_\infty^2 \endm \calJ U
\end{align*}
Combining these statements yields the desired $\calJ$-invariance of~\eqref{eq:UPDE}.
\end{proof}

In what follows, we will need to bound the nonlinear term $N(U,\pi)$ in various norms.
For future reference we therefore include the following lemma.

\begin{lemma}
\label{lem:NUbd}
Let $\pi$ be an admissible path and let $N(U,\pi)$ be as in~\eqref{eq:NUpi}.
Then
\begin{align}
\|N(U,\pi)\|_1 &\les \min(\|U\|_\infty^2, \|U\|_2^2) + \|U\|_3^3 \label{eq:NUbd1}\\
 \|N(U,\pi)\|_2 &\les \min(\|U\|_\infty^2, \|U\|_4^2) + \|U\|_6^3 \label{eq:NUbd}\\
 \|\nabla N(U,\pi)\|_2 &\les \min(\|U\|_\infty^2, \|U\|_4^2) + \|U\nabla U\|_2 + \||U|^2 \nabla U\|_2
\label{eq:NabNUbd2} \\
\|\nabla N(U,\pi)\|_1 &\les \min(\|U\|_\infty^2, \|U\|_2^2) + \|U\nabla U\|_1 + \||U|^2 \nabla U\|_1
\label{eq:NabNUbd1}
\end{align}
\end{lemma}
\begin{proof}
Direct estimation of the terms of the matrix on the right-hand side of~\eqref{eq:NUpi}.
\end{proof}

\section{The linearized problem and the root spaces at zero}
\label{sec:linroot}

Recall that  $\phi=\phi(\cdot,\alpha)$ is the ground state of $-\Laplace \phi+\alpha^2\phi=\phi^3$. Define
\beeq
\label{eq:Hilal}
\Hil(\alpha) :=   \bm -\Laplace -2\phi^2+\alpha^2 & -\phi^2 \\ \phi^2 & \Laplace +2\phi^2
-\alpha^2 \endm.
\eneq
Hence the matrix operator on the left-hand side of~\eqref{eq:UPDE}
is equal to~$-\Hil(\alpha_\infty)$, i.e., \eqref{eq:UPDE} can be rewritten as
\[ i\partial_tU -\Hil(\alpha_\infty) U = \dot{\pi}\partial_\pi W(\pi)+N(U,\pi)+VU \text{\ \ or\ \ } i\partial_tU -\Hil(t) U = \dot{\pi}\partial_\pi W(\pi)+N(U,\pi), \]
where $\Hil(t):=\Hil(\alpha_\infty)+V(t)$.
The main goal of this section and the following one is to characterize the entire discrete
spectrum of $\Hil(\alpha)$. More precisely, we will show that zero is the only real eigenvalue in $[-\alpha^2,\alpha^2]$,
and that there is a unique pair of imaginary eigenvalues $\pm i\sigma$, each of which is simple.
As mentioned in Section~\ref{sec:intro}, we will need to impose the following spectral conditions on the essential spectrum:

\begin{defi}
\label{def:spec_ass}
For all $\alpha>0$
we assume that $\Hil(\alpha)$ has no embedded eigenvalues in its essential spectrum
$(-\infty,-\alpha^2]\cup[\alpha^2,\infty)$, and that
the edges $\pm\alpha^2$ are not resonances.
\end{defi}

Now fix some $\alpha>0$.
Because of the aforementioned properties of $\Hil(\alpha)$, one has the representation
\[ L^2(\R^3)\times L^2(\R^3) = \calN + \calL + (\calN^*+\calL^*)^\perp \]
where $\calL,\calL^*$ are the sum of the eigenspaces corresponding to the purely imaginary eigenvalues
of $\Hil(\alpha)$ and~$\Hil(\alpha)^*$, respectively, and $\calN, \calN^*$ are the root spaces
of $\Hil(\alpha)$ and~$\Hil(\alpha)^*$, respectively, i.e.,
\[ \calN = \bigcup_{n=1}^\infty \ker(\Hil(\alpha)^n), \quad
\calN^* = \bigcup_{n=1}^\infty \ker((\Hil(\alpha)^*)^n).\]
The sum here is direct but not orthogonal. In particular, this representation shows that
\beeq
\label{eq:RanPs}
 \Ran(P_s(\alpha))= (\calN^*+\calL^*)^\perp
\eneq
where $I-P_s(\alpha)$ is the Riesz projection corresponding to the discrete spectrum of $\Hil(\alpha)$,
see Lemma~\ref{lem:L2split} below.
In \cite{Wei1}, Weinstein showed that the root spaces $\calN(\alpha)$ and $\calN^*(\alpha)$
of $\Hil(\alpha)$ and $\Hil^*(\alpha)$, respectively, are (with $\phi=\phi(\cdot,\alpha)$)
\begin{align}
 \calN=\calN(\alpha) &= \span \Big\{ \binom{i\phi}{-i\phi}, \binom{\partial_\alpha\phi}{\partial_\alpha\phi}, \binom{\partial_j\phi}{\partial_j\phi}, \binom{ix_j\phi}{-ix_j\phi}\:\Big|\: 1\le j\le 3\Big\} \label{eq:N}\\
 \calN^*=\calN(\alpha)^* &= \span \Big\{ \binom{\phi}{\phi}, \binom{i\partial_\alpha\phi}{-i\partial_\alpha\phi}, \binom{i\partial_j\phi}{-i\partial_j\phi}, \binom{x_j\phi}{x_j\phi}\:\Big|\: 1\le j\le 3\Big\}. \label{eq:N*}
\end{align}
Showing that the root spaces contain the sets on the right-hand side is just a
matter of direct computation. The difficulty lies with showing the equality.
Moreover, Weinstein showed that $\ker(\Hil^2(\alpha))=\ker(\Hil^3(\alpha))$
(his argument only applies to certain nonlinearities, which include the cubic NLS in $\R^3$).

In order to apply the linear dispersive $L^1(\R^3)\to L^\infty(\R^3)$
estimates from Section~\ref{sec:lin_L2} and~\ref{sec:lin_dis} to~\eqref{eq:UPDE}, one needs to
project $U$ onto $\Ran(P_s)$. Following common practice, see Soffer, Weinstein~\cite{SofWei1},
\cite{SofWei2}, and Buslaev, Perelman~\cite{BP1},
we will make an appropriate choice of the path~$\pi(t)$ in order to insure that $U(t)$
is perpendicular to $\calN^*$.
However, for technical reasons it is advantageous to impose an orthogonality condition
onto a {\em time-dependent} family of functions rather than $\calN^*$ itself. We introduce this
family in the following definition. In view of Lemma~\ref{lem:rho_infty}, it approaches
$\calN^*$ in the limit $t\to\infty$.

\begin{defi}
\label{def:rootspace}
Assume that $\pi$ is an admissible path and
let $y,\theta$ be as in~\eqref{eq:y}, \eqref{eq:theta},  $y_\infty,\theta_\infty$ as in~\eqref{eq:inf_path}, and $\rho_\infty$ as in~\eqref{eq:rhoinf}. With these functions, define
\begin{align*}
\xi_1(t) &:= \binom{e^{i\rho_\infty}\phi(\cdot+y_\infty-y,\alpha(t))}{e^{-i\rho_\infty}\phi(\cdot+y_\infty-y,\alpha(t))},\qquad
\xi_2(t) := \binom{ie^{i\rho_\infty}\partial_\alpha\phi(\cdot+y_\infty-y,\alpha(t))}{-ie^{-i\rho_\infty}\partial_\alpha\phi(\cdot+y_\infty-y,\alpha(t))} \\
\xi_{2+\ell}(t) &:= \binom{e^{i\rho_\infty}x_\ell\phi(\cdot+y_\infty-y,\alpha(t))}{e^{-i\rho_\infty}x_\ell\phi(\cdot+y_\infty-y,\alpha(t))},\qquad
\xi_{5+\ell}(t) := \binom{ie^{i\rho_\infty}\partial_\ell\phi(\cdot+y_\infty-y,\alpha(t))}{-ie^{-i\rho_\infty}
\partial_\ell\phi(\cdot+y_\infty-y,\alpha(t))}
\end{align*}
for $\ell=1,2,3$.
We also introduce another family $\{\eta_j\}_{j=1}^8$ by
\beeq
\label{eq:etadef}
 \eta_j = \bm -i & 0\\ 0 & i \endm \xi_j \text{\ \ for any\ \ }1\le j\le 8.
\eneq
\end{defi}

By inspection, $\calJ \xi_j=\xi_j$ for $1\le j\le 8$ and we chose $\eta_j$ in such a way
that $\calJ \eta_j=\eta_j$ for each $j$.
Clearly, while the $\xi_j$ correspond to $\Hil^*$, the $\eta_j$ correspond to~$\Hil$,
cf.~\eqref{eq:N} and~\eqref{eq:N*}.
Let $U$ be as in Lemma~\ref{lem:UPDE}. We refer to the condition that
\beeq
\label{eq:OC}
 \la U(t),\xi_j(t) \ra =0
\eneq
for all $t\ge 0$, $1\le j\le 8$ as the~{\em orthogonality condition}.
As usual, 
the orthogonality condition~\eqref{eq:OC} leads to an ODE for the path~$\pi(t)$.
Following \cite{Cuc}, 
we first modify the $\gamma$ parameter.

\begin{lemma}
\label{lem:sigma_tild}
Let $\pi(t)$ be an admissible path as in Definition~\ref{def:adm}. Set
\beeq
\label{eq:dotga}
\dot{\gatil}(t):=\dot{\gamma}(t)+\dot{v}(t)\cdot y_\infty(t)
\eneq
 and $\gatil(\infty):=0$,i.e.,
\[  \gatil(t):= -\int_t^\infty \Big[\dot{\gamma}(s)+\dot{v}(s)\cdot y_\infty(s)\Big]\,ds. \]
Then the function $\dot{\pi}\partial_\pi W(\pi)$ on the right-hand side of~\eqref{eq:UPDE}
satisfies
\[
\dot{\pi}\partial_\pi W(\pi)= i\Big[  \sum_{\ell=1}^3 (\dot{D_\ell}\eta_{5+\ell}-\dot{v}_\ell\eta_{2+\ell}) + \dot{\alpha}\eta_2
- \dot{\gatil}\eta_1 \Big]
\]
where the functions $\{\eta_j\}_{j=1}^8$ are as in~\eqref{eq:etadef}.
\end{lemma}
\begin{proof} By inspection. \end{proof}

The following lemma records some useful facts about the two families in
Definition~\ref{def:rootspace}.

\begin{lemma}
\label{lem:orth}
Let $\phi=\phi(\cdot,\alpha(t))$ be the ground state of~\eqref{eq:ground} and
let $\{\xi_j\}_{j=1}^8$ and $\{\eta_j\}_{j=1}^8$ be as in Definition~\ref{def:rootspace}.
Then
\begin{align*}
\la\xi_1,\eta_j \ra &= 2 \la\partial_\alpha\phi,\phi\ra \text{\ \ if $j=2$ and $=0$ else}, \\
\la\xi_2,\eta_j \ra &= -2 \la\partial_\alpha\phi,\phi\ra \text{\ \ if $j=1$ and $=0$ else}, \\
\la\xi_{2+\ell},\eta_j \ra &= -2\la\phi,\phi \ra \text{\ \ if $j=5+\ell$ and $=0$ else}, \\
\la \xi_{5+\ell},\eta_j \ra &= 2\la\phi,\phi \ra \text{\ \ if $j=2+\ell$ and $=0$ else}.
\end{align*}
Here $\partial_\alpha \la\phi,\phi\ra = 2\la\partial_\alpha\phi,\phi\ra=-\alpha^{-1}\|\phi\|_2^2$.
Moreover, let $\sigma_3=\bm 1&0\\0&-1\endm$ and
\beeq
\label{eq:Edef} E := [ -|v(t)-v_\infty|^2 -2i (v(t)-v_\infty)\cdot\nabla +\alpha_\infty^2-\alpha(t)^2]\sigma_3
\eneq
and $\Hil(t)=\Hil(\alpha_\infty)+V(t)$, see~\eqref{eq:Hilal} and~\eqref{eq:V}. Then
\begin{align*}
\Hil(t)^* \xi_1 &= E\xi_1 \\
\Hil(t)^* \xi_2 &= -2i\alpha(t) \xi_1+ E\xi_2 \\
\Hil(t)^* \xi_{2+\ell} &= -2i\xi_{5+\ell} + E\xi_{2+\ell} \\
\Hil(t)^* \xi_{5+\ell} &= E\xi_{5+\ell}
\end{align*}
for any $\ell=1,2,3$.
\end{lemma}
\begin{proof}
The statements about the scalar products are checked by direct verification.
That $\partial_\alpha \la\phi,\phi\ra=-\alpha^{-1}\|\phi\|_2^2$
follows from the fact that the ground states $\phi(\cdot,\alpha)$ of~\eqref{eq:ground}
satisfy $\phi(\cdot,\alpha)=\alpha\phi(\alpha x,1)$.
For the action of
$\Hil(t)^*$ on $\{\xi_j(t)\}_{j=1}^8$ it is convenient to introduce
\[ R^{-1} = \bm e^{-i\rho_\infty}  & 0 \\ 0 & e^{i\rho_\infty} \endm\]
and to set $\xi_j^{(0)}=R^{-1}\xi_j$ for $1\le j\le 8$. Then direct computation yields
\begin{align*}
 \Hil(t)^*\xi_1^{(0)} &=(\alpha_\infty^2-\alpha(t)^2)\sigma_3\xi_1^{(0)}\\
 \Hil(t)^*\xi_2^{(0)} &=-2\alpha\xi_1^{(0)}+(\alpha_\infty^2-\alpha(t)^2)\sigma_3\xi_2^{(0)} \\
\Hil(t)^*\xi_{2+\ell}^{(0)} &=2\xi_{5+\ell}^{(0)}+(\alpha_\infty^2-\alpha(t)^2)\sigma_3\xi_{2+\ell}^{(0)}\\
 \Hil(t)^*\xi^{(0)}_{5+\ell} &=
(\alpha_\infty^2-\alpha(t)^2)\sigma_3\xi_{5+\ell}^{(0)},
\end{align*}
which relates the $\xi_j^{(0)}$ to the root space of $\Hil(t)^*$.
Since
\[ \Hil(t)^*\xi_j = R\Hil(t)^*\xi_j^{(0)} + R\Big[\bm \Laplace&0\\0&-\Laplace\endm,R^{-1}\Big]\xi_j, \]
and
\begin{align*}
E -(\alpha_\infty^2-\alpha(t)^2)\sigma_3 &=R\Big[\bm \Laplace&0\\0&-\Laplace\endm,R^{-1}\Big] \\
& =
\bm -i\Laplace \rho_\infty-|\nabla\rho_\infty|^2-2i\nabla\rho_\infty\cdot\nabla & 0\\
0 & -i\Laplace\rho_\infty+|\nabla\rho_\infty|^2-2i\nabla\rho_\infty\cdot\nabla \endm,
\end{align*}
the lemma follows from~\eqref{eq:rhorep}.
\end{proof}

We can now derive the usual modulation equations for the admissible path $\pi$
under the orthogonality assumption~\eqref{eq:OC}.

\begin{lemma}
\label{lem:modul}
Assume that $\pi$ is an admissible path and that $U$ is an $H^1$ solution of~\eqref{eq:UPDE}
with an initial condition $U(0)$
which satisfies the  orthogonality assumptions~\eqref{eq:OC}
at time $t=0$.
Then $U$ satisfies the orthogonality assumptions~\eqref{eq:OC} for all times iff $\pi$ satisfies the
modulation equations (with $E$ as in \eqref{eq:Edef} and with $\phi=\phi(\cdot,\alpha(t))$)
\begin{align*}
\dot{\alpha}\alpha^{-1}\|\phi\|_2^2 &= \la U,\dot{\xi}_1 \ra -i\la U, E\xi_1 \ra
-i \la N(U,\pi),\xi_1 \ra \\
\dot{\gatil}\alpha^{-1}\|\phi\|_2^2 &= \la U,\dot{\xi}_2 \ra -i \la U, E\xi_2 \ra
-i \la N(U,\pi),\xi_2 \ra \\
2\dot{D_\ell} \|\phi\|_2^2 &= \la U,\dot{\xi}_{2+\ell} \ra -i \la U, E\xi_{2+\ell} \ra
-i \la N(U,\pi),\xi_{2+\ell} \ra \\
2\dot{v}_\ell \|\phi\|_2^2 &= \la U,\dot{\xi}_{5+\ell} \ra - i \la U, E\xi_{5+\ell} \ra
-i \la N(U,\pi),\xi_{5+\ell} \ra
\end{align*}
for all $1\le\ell\le 3$.
\end{lemma}
\begin{proof}
Clearly,  for any $1\le j\le 8$,
\[  \la U(t),\xi_j(t)\ra=0 \text{\ \ for all\ \ }t\ge0 \]
is equivalent with
\[ \la U(0),\xi_j(0)\ra=0 \text{\ \ and\ \ } \la \partial_t U,\xi_j \ra = -\la U,\dot{\xi}_j \ra \text{\ \ for all\ \ }t\ge0. \]
Starting from
\[ i\partial_t U -\Hil(t) U = \dot{\pi}\partial_\pi W(\pi)+N(U,\pi),\]
the modulation equations now follow from the previous two lemmas.
\end{proof}

Later it will be important to have a family of functions that plays
the same role for $Z(t)$ as $\{\xi_j\}_{j=1}^8$ does for $U(t)$.
The following lemma introduces this family and establishes some statements
for it  analogous to the ones we just obtained for $\{\xi_j\}_{j=1}^8$.

\begin{lemma}
\label{lem:xiZ}
Fix an admissible path $\pi$ and
let $\theta,y$ be as in \eqref{eq:theta} and \eqref{eq:y}, respectively. Define
\begin{align*}
\xitil_1(t,x) &:= \binom{e^{i\theta(t,x)}\phi(x-y(t),\alpha(t))}{e^{-i\theta(t,x)}\phi(x-y(t),\alpha(t))},\quad
\xitil_2(t,x) := \binom{ie^{i\theta(t,x)}\partial_\alpha\phi(x-y(t),\alpha(t))}{-ie^{-i\theta(t,x)}\partial_\alpha\phi(x-y(t),\alpha(t))} \\
\xitil_{2+\ell}(t) &:= \binom{e^{i\theta(t,x)}(x_\ell-y_\ell(t))\phi(x-y(t),\alpha(t))}{e^{-i\theta(t,x)}(x_\ell-y_\ell(t))\phi(x-y(t),\alpha(t))},\quad
\xitil_{5+\ell}(t) := \binom{ie^{i\theta(t,x)}\partial_\ell\phi(x-y(t),\alpha(t))}{-ie^{-i\theta(t,x)}\partial_\ell\phi(x-y(t),\alpha(t))}
\end{align*}
for $\ell=1,2,3$. Then
\[ \xi_j(t)= M(t)\calG_\infty(t) \xitil_j(t) \text{\ \ provided\ \ }j\ne 3,4,5 \]
and
\[ \xi_{2+\ell}(t)= M(t)\calG_\infty(t) \xitil_{2+\ell}(t) - (y_\infty-y)_\ell(t)\xi_1(t) \text{\ \ provided\ \ }\ell = 1,2,3.\]
Also, let $U$ and $Z$ be related by \eqref{eq:Udef}. Then  $U$ satisfies the orthogonality
condition~\eqref{eq:OC} iff $Z(t)$ satisfies
\[  \la Z(t),\xitil_j(t) \ra =0 \text{\ for all\ \ }1\le j\le 8,\; t\ge0.\]
Finally, introduce $\{\etatil_j(t)\}_{j=1}^8$  as in~\eqref{eq:etadef}, i.e.,
\beeq
\nn
 \etatil_j = \bm -i & 0\\ 0 & i \endm \xitil_j \text{\ \ for any\ \ }1\le j\le 8.
\eneq
Then the same scalar product relations hold as in Lemma~\ref{lem:orth}. Indeed,
\begin{align*}
\la\xitil_1,\etatil_j \ra &= 2 \la\partial_\alpha\phi,\phi\ra \text{\ \ if $j=2$ and $=0$ else}, \\
\la\xitil_2,\etatil_j \ra &= -2 \la\partial_\alpha\phi,\phi\ra \text{\ \ if $j=1$ and $=0$ else}, \\
\la\xitil_{2+\ell},\etatil_j \ra &= -2\la\phi,\phi \ra \text{\ \ if $j=5+\ell$ and $=0$ else}, \\
\la \xitil_{5+\ell},\etatil_j \ra &= 2\la\phi,\phi \ra \text{\ \ if $j=2+\ell$ and $=0$ else}.
\end{align*}
Finally, with $\Hil(\pi(t))$ as in~\eqref{eq:Hpidef} we have the relations
\beeq
\label{eq:calSj}
i\dot{\xitil}_j(t) - \Hil^*(\pi(t)) \xitil_j(t) = \dot{\pi}_* \calS_j(t)
\eneq
for all $1\le j\le8$. Here $\dot{\pi}_*:=(\dot{\alpha},\dot{v},\dot D,\dot{\gamma}_*)$,
$\dot{\gamma}_*(t):=\dot\gamma(t)+\dot v(t)\cdot y(t)$, and $\dot{\pi}_* \calS_j$
is  a short-hand notation for a linear combination of $\{i\xitil_j\}_{j=1}^8$ as well
as derivatives of these functions, and
with coefficients from the vector $\pm \dot{\pi}_*$.
\end{lemma}
\begin{proof}
This can be read off from the definitions of $M(t)$, $\calG_\infty(t)$.
\end{proof}

We now make two remarks. The first one concerns how to insure the
orthogonality condition for the transformed solution $U(t)$ (which depends
on some path) at time $t=0$ by a condition which is path independent.
The second one concerns the $\calJ$-invariance of eigenfunctions.

\begin{remark}
\label{rem:init}
Let $R_0\in L^2(\R^3)$ 
be such that $Z_0=\binom{R_0}{\bar{R}_0}$  satisfies $Z_0\in {\calN^*}^\perp$.
According to Lemma~\eqref{eq:UPDE} the transformed 
 initial condition is
\[ U(0)=M(0)\calG_\infty(0)Z_0.\]
We claim that then $\la U(0),\xi_j(0)\ra=0$ for all $1\le j\le 8$, which is
precisely the condition of the previous lemma.  Here
$\{\xi_j\}_{j=1}^8$ are the functions from Definition~\ref{def:rootspace}
defined relative to {\em any  admissible path}~$\pi$ as long as it starts at
$\pi(0)=(\alpha_0,0,0,0)$ as required by Theorem~\ref{thm:main} (this is
of course no restriction, since the initial soliton in the theorem is
as good as any other modulo a Galilei transform).
We verify this claim for $\xi_1(0)$, the other cases being similar.
First, since $\pi(0)=(\alpha_0,0,0,0)$, one checks directly from the definitions that
\begin{align*}
 \xi_1(0) &= \binom{e^{-i(v_\infty\cdot(x+D_\infty)+\gamma_\infty)}\phi(\cdot+D_\infty,\alpha_0)}
{e^{i(v_\infty\cdot(x+D_\infty)+\gamma_\infty)}\phi(\cdot+D_\infty,\alpha_0)} \\
& = \calG_\infty(0) \binom{\phi(\cdot,\alpha_0)}{\phi(\cdot,\alpha_0)}.
\end{align*}
Since also $M(0)={\rm Id}$, it therefore follows that
\[
\la U(0),\xi_1(0) \ra =
\Big\la \calG_\infty(0)M(0)Z_0, \calG_\infty(0)\binom{\phi(\cdot,\alpha_0)}{\phi(\cdot,\alpha_0)} \Big\ra
= \Big\la Z_0, \binom{\phi(\cdot,\alpha_0)}{\phi(\cdot,\alpha_0)} \Big\ra = 0,
\]
by unitarity of $\calG_\infty(0)$ and the assumption on $Z_0$.
\end{remark}

\begin{remark}
\label{rem:Jinv}
By inspection, all root spaces in this section are $\calJ$-invariant. This is a general fact.
Indeed, one checks easily that $J\Hil(\alpha)J=-\Hil(\alpha)$. Therefore, if
$\Hil(\alpha)f=i\sigma f$ with $\sigma\in\R$, it follows that $\Hil(\alpha)Jf=-i\sigma J f$
where as usual $J=\bm 0&1\\1&0\endm$. Hence
\[ \calJ \ker (\Hil-i\sigma I) = \overline{J \ker (\Hil-i\sigma I)} = \ker (\Hil-i\sigma I)\]
for any $\sigma\in\R$. A similar argument shows that the root spaces at zero are also $\calJ$-invariant.
In particular, one concludes from this that the Riesz projections
$P_s, P_{\rm root}, P_{\rm im}$ preserve the
space of $\calJ$-invariant functions in~$L^2(\R^3)\times L^2(\R^3)$.
This can also easily be seen directly: Let $P$ be any Riesz projection corresponding to
an eigenvalue of $\Hil(\alpha)$ on~$i\R$, i.e.,
\[ P = \frac{1}{2\pi i}\oint_\gamma (zI-\Hil(\alpha))^{-1}\, dz \]
where $\gamma$ is a small positively oriented circle centered at  that eigenvalue.
Since $J\Hil(\alpha)J=-\Hil(\alpha)$, one concludes that
\[
JPJ = \frac{1}{2\pi i}\oint_\gamma J(zI-\Hil(\alpha))^{-1}J\, dz = \frac{1}{2\pi i}\oint_\gamma
(\Hil(\alpha)+zI)^{-1}\, dz.
\]
Thus, if $F=\binom{F_1}{F_2}$, then $-\bar{\gamma}=\gamma$ (in the sense of oriented curves)
implies that
\[
\overline{JPF} = -\frac{1}{2\pi i} \oint_\gamma (\Hil(\alpha)+\bar{z}I)^{-1}\, d\bar{z}\; \overline{JF}
= \frac{1}{2\pi i}\oint_\gamma (zI-\Hil(\alpha))^{-1}\, dz\; \overline{JF}= P\calJ F,
\]
so $\calJ\circ P =P\circ \calJ$, as claimed.

Another important issue related to $\calJ$-invariance is whether or not a solution $\pi(t)$
of the system of modulation equations is real-valued. This is of course crucial, and is indeed
the case if $U(t)$ is $\calJ$-invariant for all $t\ge0$ and if
$\pi(0)\in\R^8$. However, we will not take up this issue here, but rather when we start
solving the modulation equations by means of a contraction scheme, see Lemma~\ref{lem:linexist} below.
\end{remark}

\section{The linearized problem and the discrete spectrum}
\label{sec:imspec}

In this section we describe the entire discrete spectrum of the linearized Hamiltonian obtained from the cubic
NLS~\eqref{eq:NLS}. Recall that the nonlinearity $|\psi|^{2\beta}\psi$ has two scalar elliptic operators associated with
it, namely,
\[ L_-:= -\Laplace+\alpha^2-\phi(\cdot,\alpha)^{2\beta}, \quad L_+:= -\Laplace+\alpha^2-(2\beta+1)\phi(\cdot,\alpha)^{2\beta}\]
where $\phi(\cdot,\alpha)$ is a ground state of the equation
\[ -\Laplace\phi + \alpha^2\phi = \phi^{2\beta}.\]
The meaning of $L_-,L_+$ is that the linearized operator of the
NLS \[ i\partial_t \psi + \Laplace \psi +|\psi|^{2\beta} \psi=0\]
takes the form
\[ \bm 0 & iL_-\\-iL_+ & 0\endm \]
provided the perturbation is written as $R=u+iv$ and this matrix
acts on $\binom{u}{v}$ (in contrast, \eqref{eq:Hilal} acts on
$\binom{R}{\bar{R}}$).
 We are interested in the range
$\frac{2}{3}<\beta\le 1$, which is supercritical. The restriction
$\beta\le 1$ has to do with Weinstein's work~\cite{Wei1}, where it
is imposed. We recall that it is know that $L_-$ has zero as
lowest eigenvalue (with $\phi$ as ground state), whereas $L_+$ has
a unique negative eigenvalue $E_0$, and a kernel spanned by
$\partial_j\phi$, $1\le j\le 3$. In one dimension, it is known
that $L_-$ and $L_+$ do not have eigenvalues in $(0,\alpha^2)$,
see Perelman~\cite{Pe2}. In dimension three it can be checked
numerically that this property also holds. This is accomplished by
showing that the associated Birman-Schwinger kernels
\begin{align*}
K_-(x,y):=\frac{\phi(x)^\beta\phi(y)^\beta}{4\pi|x-y|} \text{\ \
for the case of\ \ }L_-
\\
K_+(x,y):=\frac{(2\beta+1)\phi(x)^\beta\phi(y)^\beta}{4\pi|x-y|}
\text{\ \ for the case of\ \ }L_+
\end{align*}
have the corresponding number of eigenvalues in $(1,\infty)$: One
for $K_-$ and four for $K_+$. We will use this property of
$L_-,L_+$ in what follows -- details of the numerical work will
appear elsewhere.

The arguments in this section apply to general nonlinearities as
above, but we will present most proofs only for $\beta=1$. We
start with the existence of the imaginary spectrum. Similar
 statements were proved earlier by Grillakis~\cite{Grill}.

\begin{prop}
\label{prop:imag}
Let $\Hil(\alpha)$ be as in~\eqref{eq:Hilal} with $\alpha>0$.
Then there exist a positive integer $N$,
as well as an increasing finite sequence $\{\sigma_j(\alpha)\}_{j=1}^N$ of positive real numbers
so that for all $1\le j\le N$,
\[ \Hil(\alpha) f_j^{\pm}(\alpha) = \pm i\sigma_j(\alpha)  f_j^{\pm}(\alpha) \]
where $f_j^\pm(\alpha)$ are exponentially decreasing, $C^\infty$ functions
with $\|f_j^\pm(\alpha)\|_2=1$.
\end{prop}
\begin{proof}
The main point here is to show that $N>0$. It will be convenient
to apply a change of variables as follows: \beeq \label{eq:cov}
\bm 1 & i \\1 & -i\endm^{-1}  \bm -\Laplace -2\phi^2+\alpha^2 &
-\phi^2 \\ \phi^2 & \Laplace +2\phi^2 -\alpha^2 \endm \bm 1 & i
\\1 & -i\endm = \bm 0 & i L_{-} \\ -iL_{+} & 0 \endm \eneq where
$L_{-}=-\Laplace+\alpha^2-\phi^2$ and $L_+ =
-\Laplace+\alpha^2-3\phi^2$. Since $L_{-}\phi=0$ and $\phi>0$, the
function $\phi$ is the ground state of $L_{-}$, which implies that
$L_{-}\ge0$, $\ker{L_-}={\rm span} \{\phi\}$, and that for all
$f\perp\phi$ $\|L_{-}f\|_2\ge \|f\|_2$ since there is no
eigenvalue in $(0,\alpha^2)$. As far as $L_{+}$ is concerned,
clearly
\[\ker{L_{+}}\supset{\rm span} \{\partial_j\phi\:|\:1\le j\le 3\}.\]
In fact, Weinstein~\cite{Wei1} (see Proposition~2.8 and the proof
on page 483 of~\cite{Wei1} which applies to the cubic NLS
in~$\R^3$) showed that equality holds here and that $L_{+}$ has
exactly one negative eigenvalue $E_0$, see  (E.10) on page~489
in~\cite{Wei1}. Since
\[
 \bm 1 & 0\\ 0 & -1\endm \bm 0 & i L_{-} \\ -iL_{+} & 0 \endm \bm 1 & 0\\ 0 & -1\endm^{-1}
= - \bm 0 & i L_{-} \\ -iL_{+} & 0 \endm,
\]
we conclude that any purely imaginary spectrum has to come in conjugate pairs, as claimed above.
In fact, it is known that the spectrum can only be real or imaginary, see~\cite{BP1} or
Section~\ref{sec:lin_L2} below. Moreover, $N<\infty$ follows from Fredholm's alternative.
Finally, to show that $N>0$, consider the function
\[ g(\lambda) = \la (L_+-\lambda)^{-1}\phi,\phi\ra,\]
which is analytic on $(E_0,\alpha^2)$. At $\lambda=0$ it is well-defined because $\phi\perp \ker(L_+)$.
Moreover, $g'(\lambda)>0$ on that interval, and $g(\lambda)\to-\infty$ as $\lambda\to E_0+$.
Finally, it is clear that $L_+(\partial_\alpha \phi)=-2\alpha\phi$.
Since $\phi(x,\alpha)=\alpha\phi(\alpha x,1)$, one has
\[ 2\la \partial_\alpha \phi(\cdot,\alpha),\phi(\cdot,\alpha)\ra = \partial_\alpha \|\phi(\cdot,\alpha)\|_2^2 = -\alpha^{-2} \|\phi(\cdot,1)\|_2^2<0. \]
Therefore,
\[ g(0)= \frac{1}{-2\alpha} \la \partial_\alpha\phi,\phi\ra >0.\]
We conclude that $g(\lambda_0)=0$ for some $E_0<\lambda_0<0$. Hence,
\[ (L_+-\lambda_0)^{-1} \phi = \chi\ne0, \;\chi\perp \phi,\]
which implies that $\chi=\sqrt{L_-}\chi_1$ for some $\chi_1\perp\phi$.
Moreover, $\chi\in C^\infty$ by elliptic regularity (using that $\phi$, and thus also
the coefficients of $L_+$, are
smooth). It is also true that $\chi_1$ is smooth. In order to see this, use complex interpolation
as well as elliptic regularity on $L_-$ to conclude that $L_{-}^{-\half}\chi\in H^k$ for  all $k\ge1$.
This implies that $\chi_1$ is smooth.
The conclusion is that
\[ \phi = (L_+-\lambda_0)\sqrt{L_-}\chi_1, \text{\ \ or\ \ }0=\sqrt{L_-}(L_+-\lambda_0)\sqrt{L_-}\chi_1,\]
which in turn yields that
\[ \la \sqrt{L_-}L_+\sqrt{L_-}\chi_1,\chi_1\ra = \lambda_0 \la L_-\chi_1,\chi_1\ra < 0.\]
Hence, the self-adjoint operator $\sqrt{L_-}L_+\sqrt{L_-}$ with domain $H^4(\R^3)$ has negative
spectrum (for a detailed proof of the fact that $\sqrt{L_-}L_+\sqrt{L_-}$ is self-adjoint
with this domain see Lemma~11.10 in~\cite{RSS2}).
We would like to say that it therefore must have a negative eigenvalue.
In order to see this, we consider the variational problem
\beeq
\label{eq:var}
 \inf_{\substack{\|f\|_2=1,f\in H^2\\ f\perp \phi}} \la \sqrt{L_-}L_+\sqrt{L_-} f,f \ra =: \lambda_1<0
\eneq
and we need to show that it has a minimizer, say $f_\infty\in H^2, \|f_\infty\|_2=1, f_\infty\perp \phi$.
Note that~\eqref{eq:var} is well-defined on $H^2$ rather than $H^4$ as long as the scalar product is
 interpreted as
\beeq
\label{eq:ener}
 \int_{\R^3}[|\nabla\sqrt{L_-}f|^2+\alpha^2|\sqrt{L_-}f|^2 - 3\phi^2|\sqrt{L_-}f|^2]\,dx.
\eneq
If the minimizer $f_\infty$ exists, then by Lagrange multipliers
\[ \sqrt{L_-}L_+\sqrt{L_-}f_\infty = \lambda f_\infty + \mu \phi \]
for some $\lambda,\mu\in \R$ (in particular, $f_\infty\in H^4$).
It is clear that necessarily $\lambda=\lambda_1$ and $\mu=0$ as desired.
It therefore remains to find this minimizer. Pick a minimizing sequence $f_n\in H^2$ of functions
which satisfy both constraints.
In view of~\eqref{eq:ener}, $\|f_n\|_2=1$, as well as $\|\sqrt{L_-}f\|_2\asymp \|f\|_{H^1}$
for all $f\perp\phi$, one concludes that $\sup_n\|f_n\|_{H^2}<\infty$. Hence,
without loss of generality, $f_n \rightharpoonup f_\infty$
in the weak sense in~$H^2$. By Rellich's compactness lemma, $f_n \to f_\infty$ in
$H^1_{\rm loc}$.
Therefore, the constraint $f_\infty\perp\phi$ holds.
The difficulty lies with showing the other constraint $\|f_\infty\|_2=1$.
We claim first that
\beeq
\label{eq:lower}
 \la \sqrt{L_-}L_+\sqrt{L_-}f_\infty,f_\infty\ra \le \lambda_1.
\eneq
By definition,
\[ \lambda_1=\lim_{n\to\infty}\la \sqrt{L_-}L_+\sqrt{L_-}f_n,f_n\ra = \lim_{n\to\infty}\int_{\R^3}[|\nabla\sqrt{L_-}f_n|^2+\alpha^2|\sqrt{L_-}f_n|^2 - 3\phi^2|\sqrt{L_-}f_n|^2]\,dx. \]
Clearly, because of weak convergence in $H^2$,
\[ \liminf_{n\to\infty}\int_{\R^3}[|\nabla\sqrt{L_-}f_n|^2+\alpha^2|\sqrt{L_-}f_n|^2]\, dx \ge \int_{\R^3}[|\nabla\sqrt{L_-}f_\infty|^2+\alpha^2|\sqrt{L_-}f_\infty|^2]\, dx
\]
and by the decay of $\phi$ at infinity, also
\[
\lim_{n\to\infty}\int_{\R^3} 3\phi^2|\sqrt{L_-}f_n|^2\,dx = \int_{\R^3} 3\phi^2|\sqrt{L_-}f_\infty|^2\,dx.
\]
Combining these statements yields \eqref{eq:lower}.
Note that because of $\lambda_1<0$ one concludes from~\eqref{eq:lower} in particular
that $f_\infty\ne0$.
Assume now that $0<\|f_\infty\|_2<1$. Then $g:=\frac{f_\infty}{\|f_\infty\|_2}$ satisfies
$\|g\|_2=1$, $g\perp\phi$, and
\[  \la \sqrt{L_-}L_+\sqrt{L_-}g,g\ra \le \|f_\infty\|_2^{-2}\lambda_1< \lambda_1<0, \]
which contradicts the definition of $\lambda_1$. Hence indeed $\|f_\infty\|_2=1$, and
thus $f_\infty$ is a minimizer and thus also an eigenfunction
\[  \sqrt{L_-}L_+\sqrt{L_-} f_\infty = \lambda_1 f_\infty.\]
This implies, with $\lambda_1=:-\sigma^2, \sigma>0$ and $v:=\sqrt{L_-}f_\infty$,
\[ L_-L_+ v = -\sigma^2 v, \quad v\ne0,\;v\perp \phi. \]
We can find $u\in H^2$ so that $L_-u=-\sigma v$, and $u$ is unique up to a multiple of $\phi$.
Then
\[
 L_{-}(u+c\phi) = -\sigma v,\quad L_{-}(L_+v - \sigma (u+c\phi))=0
\]
for any $c\in\R$. Since $\sigma\ne0$, we can choose $c_0\in\R$ in such a way that in fact
\[ L_+v - \sigma (u+c_0\phi)=0.\]
Renaming $u+c_0\phi$ into $u$, we obtain the system
\[ \bm 0 & iL_- \\ -iL_+ & 0\endm \binom{v}{u}=-i\sigma \binom{v}{u},\]
which shows that $\binom{v}{u}$ is an eigenvector with eigenvalue $-i\sigma$.
Similarly, one finds an eigenfunction with eigenvalue $i\sigma$.
Finally, the statements concerning the smoothness and exponential decay of the eigenfunctions
follow from elliptic regularity, and Agmon's argument in Lemma~\ref{lem:agmon} below, respectively.
\end{proof}

\begin{remark}
The previous proof is in some sense the converse of some arguments in Weinstein's papers~\cite{Wei1}
and~\cite{Wei2}. Indeed, there one uses that $\partial_\alpha \|\phi(\cdot,\alpha)\|_2^2 >0$
for the {\em subcritical case} to show that $\la L_+ f,f\ra\ge0$ for all $f\perp \phi$. Note that
we have exploited the opposite effect, namely that $\la L_+ f,f\ra<0$ for some $f\perp \phi$
which we derived from  $\partial_\alpha \|\phi(\cdot,\alpha)\|_2^2 < 0$
for the {\em supercritical case}.
\end{remark}

Next we show that zero is the only point of the discrete spectrum
on the real axis. Since any such point would have to be an
eigenvalue, we just need to show that zero is the only eigenvalue
in the interval $(-\alpha^2,\alpha^2)$. The following lemma is
somewhat stronger, since it proves this for the closed interval
$[-\alpha^2,\alpha^2]$. The argument is an adaptation of
Proposition~2.1.2 in Perelman~\cite{Pe2}. It is based on the fact
that $L_+$ does not have any eigenvalues in $(0,\alpha^2)$.

\begin{lemma}
\label{lem:realev}
The only eigenvalue of $\Hil(\alpha)$ in the interval $[-\alpha^2,\alpha^2]$ is
zero.
\end{lemma}
\begin{proof}
Suppose not.  Then $\Hil(\alpha)^2$ has an eigenvalue $E\in (0,\alpha^4]$. For simplicity and without
loss of generality, let us choose $\alpha=1$. Then there is $\psi\in L^2(\R^3)$, $\psi\ne0$,  such that
\[ L_{-}L_+\psi = E\psi \]
with $0<E\le 1$. Clearly, $\psi\perp\phi$ and $\psi\in H^4_{loc}(\R^3)$ by elliptic regularity.
Define $A:= PL_+P$ where $P$ is the projection orthogonal to $\phi$. Since
\[ \ker(L_+) = {\rm span} \{\partial_j \phi\:|\:1\le j\le 3\},
\]
and $\la \phi,\partial_\alpha\phi\ra\ne0$,
we conclude that
\[
\ker(A) = {\rm span} \{\partial_j\phi,\,\phi\:|\:1\le j\le 3\}.
\]
Moreover, let $E_0<0$ be the unique negative eigenvalue of $L_+$. Then
consider (as before) the function
\[ g(\lambda):= \la (L_+-\lambda)^{-1}\phi,\phi\ra \]
which is differentiable on the  interval $(E_0, 1)$  due to the orthogonality
of $\phi$ to the kernel of $L_+$. Moreover,
\[ g'(\lambda)=\la (L_+-\lambda)^{-2}\phi,\phi\ra>0,\; g(0)= -\half \la \phi,\partial_\alpha\phi\ra>0. \]
The final inequality here is due to the supercritical nature of our problem.
Since also $g(\lambda)\to -\infty$ as $\lambda\to E_0$, it follows that $g(\lambda_1)=0$ for some $E_0<\lambda_1<0$.
Moreover, this is the only zero of $g(\lambda)$ with $E_0<\lambda<1$.
If we set
\[ \eta:= (L_+-\lambda_1)^{-1} \phi,\]
then
\[ A\eta=\lambda_1\eta, \qquad \la \eta,\phi\ra=0.\]
Conversely, if
\[ Af=\lambda f\]
for some $-\infty<\lambda<1$, $\lambda\ne0$,  and $f\in
L^2(\R^3)$, then $f\perp \phi$ and
\[ (PL_+P-\lambda)f=(A-\lambda)f=0.\]
Since also
\[ E_0\la f,f\ra\le \la L_+f,f\ra = \lambda\la f,f\ra \]
it follows that $\lambda\ge E_0$. If $\lambda=E_0$, then $f$ would necessarily
have to be the ground state of $L_+$ and thus of definite sign. But
then $\la f,\phi\ra\ne0$, which is impossible. Hence $E_0<\lambda<1$.
But then $g(\lambda)=0$ implies that $\lambda=\lambda_1$ is unique.
In summary, $A$ has eigenvalues $\lambda_1$ and $0$ in $(-\infty,1)$, with $\lambda_1$
being a simple eigenvalue and $0$ being an eigenvalue of multiplicity four.
Now define
\[ \calF:= {\rm span} \{ \psi,\,\eta,\, \partial_j\phi,\,\phi\:|\:1\le j\le 3\}.
 \]
We claim that
\beeq
\label{eq:dimclaim}
 \dim(\calF)=6.
\eneq
Since $\phi$ is perpendicular to the other functions, it suffices to show that
\[ c_1\psi+c_2\eta + \sum_{j=3}^5 c_j\partial_j\phi=0 \]
can only be the trivial linear combination.
Apply $L_+$. Then
\[ c_1 L_+\psi+c_5 L_+\eta =0\]
and therefore
\begin{align*}
c_1 \la L_+\psi,\psi\ra + c_2 \la L_+\eta,\psi\ra &=0 \\
c_1 \la L_+\psi,\eta\ra + c_2 \la L_+\eta,\eta\ra &=0.
\end{align*}
This is the same as
\begin{align*}
c_1 E\la L_{-}^{-1} \psi,\psi\ra + c_2 \lambda_1\la \eta,\psi\ra &=0 \\
c_1 \lambda_1\la \psi,\eta\ra + c_2 \lambda_1\la \eta,\eta\ra &=0.
\end{align*}
The determinant of this system is
\[ E\lambda_1\la L_{-}^{-1} \psi,\psi\ra \la \eta,\eta\ra - \lambda_1^2 |\la \eta,\psi\ra|^2 <0.\]
Hence $c_1=c_2=0$ and therefore also $c_3=c_4=c_5=0$, as desired.
Thus, \eqref{eq:dimclaim} holds.
Finally, we claim that
\beeq
\label{eq:qform}
 \sup_{\|f\|_2=1,\,f\in\calF} \la Af,f\ra <1.
\eneq
If this is true, then by the min-max principle and \eqref{eq:dimclaim} we would obtain
that the number of eigenvalues of $A$ in the interval $(-\infty,1)$ (counted with multiplicity)
would have to be at least six. On the other, we showed before that this number is exactly five,
leading to a contradiction. Hence, the lemma will follow once we verify~\eqref{eq:qform}.
Since $\la PL_{-}^{-1}Pf,f\ra <\la f,f\ra$ for all $f\ne0$, and since $E\le 1$ by assumption,
 this in turn follows from the stronger claim that
\beeq
\label{eq:qform'}
 \la Af,f\ra \le E\la PL_{-}^{-1}P f,f\ra
\eneq
for all $f=a\psi+b\phi+\vec{c}\cdot\nabla \phi + d\eta$. Clearly, we can take $b=0$. Then the
left-hand side of \eqref{eq:qform'} is equal to
\begin{align}
&\la L_+(a\psi),a\psi+\vec{c}\cdot\nabla\phi+d\eta\ra +\la L_+(\vec{c}\cdot\nabla\phi+d\eta),
a\psi+\vec{c}\cdot\nabla\phi+d\eta\ra \nn\\
&= E\la L_-^{-1}(a\psi),a\psi+\vec{c}\cdot\nabla\phi+d\eta\ra +E\la \vec{c}\cdot\nabla\phi+d\eta,
L_{-}^{-1}(a\psi)\ra + \la L_+(d\eta),d\eta\ra \nn\\
&= E\la L_-^{-1}(a\psi),a\psi+\vec{c}\cdot\nabla\phi+d\eta\ra +E\la \vec{c}\cdot\nabla\phi+d\eta,
L_{-}^{-1}(a\psi)\ra + \lambda_1\|d\eta\|_2^2,
\label{eq:lhs1}
\end{align}
whereas the right-hand side of \eqref{eq:qform'} is
\beeq
= E\la L_-^{-1}(a\psi),a\psi+\vec{c}\cdot\nabla\phi+d\eta\ra +E\la \vec{c}\cdot\nabla\phi+d\eta,
L_{-}^{-1}(a\psi)\ra + E\la L_{-}^{-1}(\vec{c}\cdot\nabla\phi+d\eta),\vec{c}\cdot\nabla\phi+d\eta\ra.
\label{eq:rhs1}
\eneq
Since
\[ \lambda_1\|d\eta\|_2^2\le 0,\quad E\la L_{-}^{-1}(\vec{c}\cdot\nabla\phi+d\eta),\vec{c}\cdot\nabla\phi+d\eta\ra
\ge0,
\]
we see that \eqref{eq:rhs1} does indeed dominate \eqref{eq:lhs1},
and \eqref{eq:qform'} follows.
\end{proof}

It will be important for us that the proof of
Lemma~\ref{lem:realev}  applies to all critical and supercritical
nonlinearities $|\psi|^{2\beta}\psi$, $\frac23\le \beta\le 1$ (the
latter restriction has to do with Weinstein's paper~\cite{Wei1},
where it is needed). In the critical case, $\lambda_1=0$ and
$\eta=-\half\partial_\alpha\phi$ and $0$ is the only eigenvalue of
$A$ in $(E_0,1)$. Otherwise, the argument is the same.

In the subcritical case the proof of Lemma~\ref{lem:realev} breaks
down. In fact, the statement is false: The proof of the following
lemma shows that there has to be at least one pair $\pm \lambda$
of real eigenvalues with $0<\lambda<\alpha^2$ in the subcritical
case. It is reasonable to expect that there should be exactly one
such pair, but we do not address that here.

\begin{lemma}
\label{eq:disc_spec}
For any nonlinearity $|\psi|^{2\beta}\psi$ with $\frac23<\beta\le1$ the discrete spectrum of the linearized operator
$\Hil(\alpha)$ consists
of zero and a single pair of imaginary eigenvalues $\pm i\sigma$, $\sigma>0$, each of which is simple (i.e.,
$N=1$ in Proposition~\ref{prop:imag}).
\end{lemma}
\begin{proof}
Let $\gamma$ be a sufficiently large contour (say an ellipse) that contains zero and does not
intersect $(-\infty,-\alpha^2]\cup[\alpha^2,\infty)$. Use this contour to define the Riesz projection
associated with the discrete spectrum. Since the Riesz projection depends continuously on
the parameter $\beta$ determining the nonlinearity, we conclude that its rank is constant.
For $\beta=\frac23$ (the critical case) Weinstein showed that the dimension of the root space at zero
 is ten. Moreover, it is known
that the spectrum has to be real in that case, and we just showed that zero is the only point
of the discrete spectrum. Hence the rank is ten for all $\frac23\le \beta\le1$. Moreover, in
the supercritical case Weinstein showed that the root space at zero has dimension eight (in fact, it
is enough for us now to know that it has rank at least eight, which is obvious), and we know
that there are no other real eigenvalues and at least one pair of purely imaginary ones
(in fact, the existence would now also follow). So there has to be a unique pair, each of which is simple.
We are done.
\end{proof}

Next, we turn to the issue of resonances of $\Hil(\alpha)$ at the
edges of the essential spectrum. A ``resonance'' at $\alpha^2$ (or
$-\alpha^2$) here refers to the existence of a solution $f$ to
$\Hil(\alpha)f=\alpha^2 f$ (or $=-\alpha^2 f$) so that $f\not\in
L^2(\R^3)$, but such that \beeq \label{eq:wl2}
 \int_{\R^3} |f(x)|^2(1+|x|)^{-2\gamma}\, dx <\infty
\eneq for all $\gamma>\half$. As will be explained in
Section~\ref{sec:lin_L2}, if $\pm\alpha^2$ are neither resonances
nor eigenvalues (we have already excluded the latter), then
$(\Hil(\alpha)\mp\alpha^2)^{-1}$ is bounded on suitable weighted
$L^2$ spaces. This will be important in order to establish the
dispersive estimates for $e^{it\Hil(\alpha)}$. The proof of the
following lemma is similar to that of Lemma~\ref{lem:realev}, and
is an adaptation of the argument in Appendix~1 of Perelman's
paper~\cite{Pe2}. It shows that if the {\em scalar} operator $L_-$
does not have a resonance at $\alpha^2$ (the edge of its
continuous spectrum), then the {\em matrix} operator
$\Hil(\alpha)$ does not have a resonance at $\pm\alpha^2$. While
it is easy to verify this in one dimension (say numerically, see
also \cite{Pe2} where this fact is mentioned but not proved), it
is not clear to the author how to accomplish this in three
dimensions (even numerically). Hence we leave it as a conditional
statement.

\begin{lemma}
\label{lem:res} Suppose that $L_-$ has neither an eigenvalue nor a
resonance at $\alpha^2$. Then the edges $\pm\alpha^2$ are not
resonances of $\Hil(\alpha)$, i.e., there do not exist solutions
$f$ of $\Hil(\alpha)f= \pm\alpha^2 f$ which satisfy~\eqref{eq:wl2}
but are not in $L^2$.
\end{lemma}
\begin{proof}
We again set $\alpha=1$. By symmetry, it suffices to consider the
right edge $\alpha^2$.  Suppose then there is such a solution $f$
with $\Hil(\alpha)f=f$. Write $f=\binom{\psi}{\tilde{\psi}}$. Then
$iL_{-}\tilde{\psi}=\psi$ and $-iL_+\psi=\tilde\psi$. In
particular, $L_-L_+\psi =\psi$ and
\[ \int_{\R^3} |\psi(x)|^2\, dx=\infty, \quad \int_{\R^3} |\psi(x)|^2(1+|x|)^{-2\gamma}\, dx <\infty \]
for all $\gamma>\half$. Clearly, $\la \psi,\phi\ra =0$, the latter
inner product being well-defined because of the rapid decay of
$\phi$ and~\eqref{eq:wl2}. Furthermore, $\psi\in H^4_{loc}(\R^3)$
by elliptic regularity. Pick a smooth cut-off $\chi\ge0$ which is
constant $=1$ around zero, and compactly supported. Define for any
$0<\eps<1$
\[ \psi^\eps := \psi \chi(\eps\cdot)+\mu(\eps)\phi, \quad \mu(\eps):= - \frac{\la \psi \chi(\eps\cdot),\phi\ra}{\la \phi,\phi\ra}.
\]
Clearly, $\la \psi^\eps,\phi\ra=0$ and $|\mu(\eps)|=o(1)$ as
$\eps\to0$ (in fact, like $e^{-C/\eps}$). It follows that
\[ \|\psi^\eps\|_2^2 = M_0(\eps)+o(1),\quad M_0(\eps):= \int_{\R^3}|\psi(x)|^2\chi(\eps x)^2\, dx\]
with $M_0(\eps)\to\infty$ as $\eps\to0$. We now claim that \beeq
\label{eq:L+claim}
 \la  L_+\psi^\eps,\psi^\eps \ra = \|\psi^\eps\|_2^2 + \la (L_+-1)\psi,\psi\ra +o(1)
\eneq as $\eps\to0$. We will need to justify that
\[ M(\eps):=\la
(L_+-1)\psi,\psi\ra\] is a finite expression. We first show that
this justification, as well as~\eqref{eq:L+claim} can be reduced
to showing that $\nabla \psi\in L^2(\R^3)$. Write
$L_{-}=-\Laplace+1+V_1$ and $L_+=-\Laplace+1+V_2$, with Schwartz
functions $V_1,V_2$ (they are of course explicitly given in terms
of $\phi$, but we are not going to use that now). We start from
the evident expression \beeq \nonumber
 \la  L_+\psi^\eps,\psi^\eps \ra = \|\psi^\eps\|_2^2 + \la (L_+-1)\psi^\eps,\psi^\eps\ra
= \|\psi^\eps\|_2^2 + \la (-\Laplace+V_2)\psi^\eps,\psi^\eps\ra.
\eneq By the rapid decay of $V_2$ and~\eqref{eq:wl2},
\[ \la (-\Laplace+V_2)\psi^\eps,\psi^\eps\ra = \int_{\R^3} |\nabla \psi^\eps(x)|^2\, dx
+ \int_{\R^3} V_2(x)|\psi(x)|^2\, dx + o(1).
\]
Assuming $\nabla\psi\in L^2$, we calculate further that
\begin{align*}
\int_{\R^3} |\nabla \psi^\eps(x)|^2\, dx &= \int_{\R^3} \Big| \nabla\psi(x)\chi(\eps x)+ \eps\psi(x)\nabla\chi(\eps x)\Big|^2\, dx \\
&= \int_{\R^3} |\nabla\psi(x)|^2\, dx + \int_{\R^3} |\nabla\psi(x)|^2(\chi(\eps x)^2-1)\, dx \\
& \quad + 2\eps\int_{\R^3} \psi(x)\chi(\eps
x)\nabla\psi(x)\cdot\nabla\chi(\eps x)\, dx
+ \eps^2\int_{\R^3}\psi(x)^2|\nabla\chi(\eps x)|^2\, dx \\
&=  \int_{\R^3} |\nabla\psi(x)|^2\, dx + o(1).
\end{align*}
To pass to the last line, estimate the error terms using
$\nabla\psi\in L^2$ and~\eqref{eq:wl2} (any $\gamma<1$ works
here). This proves~\eqref{eq:L+claim} provided we interpret $\la
(L_+-1)\psi,\psi\ra$ as
\[ \int_{\R^3} [|\nabla\psi(x)|^2+V_2(x)|\psi(x)|^2]\,dx.\]
To prove $\nabla\psi\in L^2$, we start from the definition, i.e.,
\[ (-\Laplace+1+V_1)(-\Laplace+1+V_2)\psi=\psi\]
which can be written as \beeq \label{eq:def_eq}
 (\Laplace^2-2\Laplace)\psi + (-\Laplace+1)V_2\psi+V_1(-\Laplace+1)\psi +V_1V_2\psi =0.
\eneq At least formally, integrating by parts against $\psi$
yields that
\[
\|\Laplace\psi\|_2^2+2\|\nabla\psi\|_2^2 \le \int
(|V_1|+|V_2|+|V_1V_2|)|\psi(x)|^2\, dx + \|V_2\psi\|_2^2 +
\|V_1\psi\|_2^2 + \half\|\Laplace\psi\|_2^2,
\]
and thus, in particular, $\nabla\psi\in L^2$. To make this
precise, we write the defining equations as $-iL_+\psi=\tilde\psi$
and $iL_-\tilde\psi=\psi$. Then
\[ -\Laplace \psi = i\tilde\psi-\psi-V_2\psi \]
which implies that (via \eqref{eq:wl2}) \beeq \label{eq:D2w}
\int_{\R^3} |\Laplace\psi(x)|^2\la x\ra^{-2\gamma}\, dx < \infty
\qquad \forall\,\gamma>\half. \eneq It is now a simple matter to
deduce from this  and \eqref{eq:wl2} that \beeq \label{eq:nablaw}
 \int_{\R^3} |\nabla\psi(x)|^2\la x\ra^{-2\gamma}\, dx < \infty \qquad \forall\,\gamma>\half.
\eneq This can be done in various ways. For example, it follow
from Gauss' integral theorem applied to $\div(u\nabla u\la
x\ra^{-2\gamma})$ that  for all $u\in H^2$ with compact support
\[ \int_{\R^3} |\nabla u(x)|^2\la x\ra^{-2\gamma}\, dx \le C\int_{\R^3}(|u(x)|^2+|\Laplace u(x)|^2)\la x\ra^{-2\gamma}\, dx.\]
Setting $u=\psi^\eps$, and letting $\eps\to0$ yields the desired
inequality~\eqref{eq:nablaw} for $\nabla u$. To make our heuristic
argument leading to $\nabla\psi\in L^2$ precise, we
pair~\eqref{eq:def_eq} with $\chi(\eps x)^4\psi$ and integrate by
parts. This yields
\begin{align*}
\label{eq:schag1}
&\la \Laplace\psi,\Laplace(\chi(\eps\cdot)^4\psi)\ra + 2\la \nabla\psi,\nabla(\chi(\eps\cdot)^4\psi)\ra \\
& = -\la V_2\psi,(-\Laplace+1)\chi(\eps\cdot)^4\psi\ra -\la
(-\Laplace+1)\psi,V_1\chi(\eps\cdot)^4\psi\ra -\la V_1V_2
\psi,\chi(\eps\cdot)^4\psi\ra.
\end{align*}
The terms on the right-hand side are uniformly bounded as
$\eps\to0$ due to \eqref{eq:wl2}, \eqref{eq:D2w}
and~\eqref{eq:nablaw} and the rapid decay of $V_1,V_2$ and their
derivatives. The second term on the left-hand side satisfies
\begin{align*}
&\int_{\R^3} \nabla\psi(x)(4\eps\nabla\chi(\eps x)\chi(\eps x)^3\psi(x)+\chi(\eps x)^4\nabla\psi(x)) \, dx \\
&\ge \half\int_{\R^3} |\nabla\psi(x)|^2\chi(\eps x)^4\, dx -
C\eps^2\int_{\R^3} |\psi(x)|^2|\nabla\chi(\eps x)|^2
\chi(\eps x)^2\, dx \\
&\ge \half\int_{\R^3} |\nabla\psi(x)|^2\chi(\eps x)^4\, dx - O(1)
\end{align*}
by \eqref{eq:wl2}. Similarly,
\begin{align*}
&\la \Laplace\psi,\Laplace(\chi(\eps\cdot)^4\psi)\ra \\
&=\int_{\R^3} \Laplace\psi(x) \Big((4\eps^2\Laplace\chi(\eps x)\chi(\eps x)^3+ 12\eps^2|\nabla\chi(\eps x)|^2\chi(\eps x)^2)\psi(x)\\
&\quad +8\eps\nabla\chi(\eps x)\cdot\nabla\psi(x)\chi(\eps x)^3+\chi(\eps x)^4\Laplace\psi(x)\Big)\, dx \\
&\ge \half\int_{\R^3} |\Laplace\psi(x)|^2\chi(\eps x)^4\, dx -
C\eps^2\int_{\R^3} |\nabla\psi(x)|^2\chi(\eps x)^2
\, dx \\
&\quad -C\eps^4\int_{\R^3} |\psi(x)|^2[|\nabla\chi(\eps x)|^4
+\chi(\eps x)^2]\, dx\\
&\ge \half\int_{\R^3} |\Laplace\psi(x)|^2\chi(\eps x)^4\, dx -
O(1)
\end{align*}
by \eqref{eq:wl2} and \eqref{eq:nablaw}. Combining these estimates
and invoking the monotone convergence theorem yields
$\Laplace\psi\in L^2$ and $\nabla\psi\in L^2$. It is easy to see
that the previous argument allows a better conclusion than $L^2$,
namely that $\la x\ra^b \Laplace \psi\in L^2$ for any $b<\half$
and similarly for $\nabla \psi$. In fact, a much stronger
conclusion is possible for $\Laplace\psi$: recall that
$g_1:=\psi+i\tilde\psi$ and $g_2:=\psi-i\tilde\psi$ satisfy
\[ \bm -\Laplace+1+W_1 & W_2 \\ -W_2 & \Laplace-1-W_1 \endm
\binom{g_1}{g_2}=\binom{g_1}{g_2}\] where $W_1,W_2$ are again
exponentially decaying potentials. This implies that
\begin{align*}
\Laplace g_1 = W_1 g_1+W_2 g_2 \\
\Laplace g_2 = 2g_2 + W_1g_2 + W_2 g_1.
\end{align*}
Hence $\la x\ra^b \Laplace g_1\in L^2$ for all $b>0$. Similarly,
\[ g_2=(\Laplace-2)^{-1}[W_1g_2 + W_2 g_1]\]
is exponentially decaying, which implies that  $\la x\ra^b
\Laplace g_2\in L^2$ for all $b>0$. Consequently, the same
estimate holds for $\Laplace \psi$ as well as for
$f:=(L_+-1)\psi$. Hence $\la (L_+-1)\psi,\psi\ra=\la f,\psi\ra$ is
well-defined as a usual scalar product. Moreover, one has
\[ L_-f=-(L_--1)\psi \text{\ \ or\ \
}\psi=-(L_--1)^{-1}L_-f=-f-(L_--1)^{-1}f.\] We conclude that \beeq
\label{eq:fpsi} \la f+\psi,f\ra = -\la (L_--1)^{-1}f,f\ra <0,
\eneq where the final inequality follows from $L_-\ge1$ on
$\{\phi\}^\perp$, as well as our assumption that $L_-$ has neither
an eigenvalue nor a resonance at $\alpha^2=1$. Recall that this
insures that for any $\tau>0$
\[ \| \la x\ra^{-1-\tau}(L_--1)^{-1}h\|_2 \le C\|\la
x\ra^{1+\tau}h\|_2 \] for all $h$ for which the right-hand side is
finite, in particular for $h=f$. The inequality~\eqref{eq:fpsi}
will play a crucial role in estimating a quadratic form as in
Lemma~\ref{lem:realev}. To see this, let
\[ \calF_\eps:={\mathrm
span}\{\psi^\eps,\partial_1\phi,\partial_2\phi,\partial_3\phi,\eta,\phi\}.
\]
As in the proof of Lemma~\ref{lem:realev} one shows that
$\dim\calF_\eps=6$, as least if $\eps>0$ is sufficiently small
(use that $\la L_+\psi^\eps,\psi^\eps\ra\to\infty$ as $\eps\to0$).
It remains to show that for small $\eps>0$ \beeq \label{eq:tschud}
\max_{f\in\calF_\eps} \frac{\la PL_+P f,f\ra}{\la f,f\ra}<1 \eneq
where $P$ is the projection orthogonal to $\phi$. If so, then this
would imply that $A=PL_+P$ has at least six eigenvalues (with
multiplicity) in $(-\infty,1)$. However, we have shown in the
proof of Lemma~\ref{lem:realev} that there are exactly five such
eigenvalues. To prove~\eqref{eq:tschud}, it suffices to consider
the case $f\perp \phi$. Compute
\begin{align*}
& \frac{\la
L_+(a\psi^\eps+\vecc\cdot\nabla\phi+d\eta),a\psi^\eps+\vecc\cdot\nabla\phi+d\eta\ra}
{\|a\psi^\eps+\vecc\cdot\nabla\phi+d\eta\|^2} \\
& =\frac{|a|^2(\|\psi^\eps\|_2^2 + M(\eps)+o(1))+2\Re\lambda_1\la
a\psi^\eps,d\eta\ra +
\lambda_1\|d\eta\|_2^2}{|a|^2\|\psi^\eps\|_2^2 + 2\Re\la
a\psi^\eps, \vecc\cdot\nabla\phi+d\eta\ra +
\|\vecc\cdot\nabla\phi+d\eta\|_2^2} \\
& \le \max_{x\in\Compl^5} \frac{|x_1|^2(1+\delta^2
M(\eps)+o(\delta^2))+2\delta \lambda_1\Re\la x_1\psi^\eps,x_5 e_4
\ra + \lambda_1|x_5|^2}{|x_1|^2+2\delta\Re\la
x_1\psi^\eps,x_2e_1+x_3e_2+x_4e_3+x_5e_4\ra +
\|x_2e_1+x_3e_2+x_4e_3+x_5e_4\|_2^2}
\end{align*}
where we have set $\delta^2:=\|\psi^\eps\|_2^2$ and
\[ e_1=\frac{\partial_1\phi}{\|\partial_1\phi\|_2},
\quad  e_2=\frac{\partial_2\phi}{\|\partial_2\phi\|_2}, \quad
e_3=\frac{\partial_3\phi}{\|\partial_3\phi\|_2}, \quad
e_4=\frac{\eta}{\|\eta\|_2}.\] Note that $\eta$ is a radial
function, since it is given by $(L_+-\lambda_1)^{-1}\phi$ and both
$\phi$ and the kernel of $(L_+-1)^{-1}$ are radial. Hence
$e_j\perp e_4$ for $1\le j\le 3$. Set
\[ b_j^\eps:=  \la \psi^\eps,e_j\ra \text{\ \ for\ \ }1\le j\le
4.\] Then $b_j^\eps\to b_j^0:= \la \psi,e_j\ra$ as $\eps\to0$ by
the exponential decay of the $e_j$. Let $B^\eps,C^\eps$ (which
depend on $\eps$) be $5\times 5$ Hermitian matrices so that
\[ C^\eps_{11}:=1+\delta^2 M(\eps)+o(\delta^2), \;C^\eps_{15}=C^\eps_{51}:=\lambda_1
\delta b^\eps_4,\; C^\eps_{55}:= \lambda_1\] and $C^\eps_{ij}=0$
else,
\[ B^\eps_{1j}=B^\eps_{j1}:=\delta b^\eps_{j-1} \text{\ \ for\ \ } 2\le j\le 5\]
and $B^\eps_{ij}=0$ else.  In view of the preceding,
\[
\max_{f\in\calF_\eps} \frac{\la PL_+P f,f\ra}{\la f,f\ra} \le
\max_{x\in\Compl^5} \frac{\la Cx,x\ra}{\la (I+B)x,x\ra}.
\]
Clearly, the right-hand side equals the largest eigenvalue of the
Hermitian matrix \[ (I+B^\eps)^{-\half}C^\eps(I+B^\eps)^{-\half}=
C-\half(BC+CB)+\frac38(B^2C+CB^2)+\frac14BCB+O(\delta^3),
\]
where we have dropped the $\eps$ in the notation on the right-hand
side. With some patience one can check that the right-hand side
equals the matrix $D$ which is given by (dropping $\eps$ from the
notation)
\[\left[
\begin{array}{lllll} 1+\delta^2 M_1 & -\frac{\delta}{2} b_1 & -\frac{\delta}{2} b_2 &
-\frac{\delta}{2} b_3 & \frac{\delta}{2}(\lambda_1-1)b_4\\
-\frac{\delta}{2} b_1 & \frac{\delta^2}{4} b_1^2 &
\frac{\delta^2}{4} b_1b_2 & \frac{\delta^2}{4} b_1b_3 &
\frac{\delta^2}{4}(1-\half\lambda_1)b_1b_4 \\
-\frac{\delta}{2} b_2 & \frac{\delta^2}{4} b_1b_2 &
\frac{\delta^2}{4} b_2^2 & \frac{\delta^2}{4} b_2b_3 &
\frac{\delta^2}{4}(1-\half\lambda_1)b_2b_4 \\
-\frac{\delta}{2} b_3 & \frac{\delta^2}{4} b_1b_3 &
\frac{\delta^2}{4} b_2b_3 & \frac{\delta^2}{4} b_3^2 &
\frac{\delta^2}{4}(1-\half\lambda_1)b_3b_4 \\
\frac{\delta}{2}(\lambda_1-1)b_4 &
\frac{\delta^2}{4}(1-\half\lambda_1)b_1b_4 &
\frac{\delta^2}{4}(1-\half\lambda_1)b_2b_4 &
\frac{\delta^2}{4}(1-\half\lambda_1)b_3b_4 & \lambda_1
+\frac{\delta^2}{4}(1-\lambda_1) b_4^2
\end{array}\right] + o(\delta^2)
\]
where $M_1:=M(\eps)-\frac34\lambda_1
b_4^2+\frac34(b_1^2+b_2^2+b_3^2+b_4^2)$. When $\delta=0$, this
matrix has eigenvalues $1$, $0$, $\lambda_1<0$, and $0$ has
multiplicity three. When $\delta\ne0$ but very small, the largest
eigenvalue will be close to one, of the form $1+x$ with $x$ small.
We need to see that $x<0$. Collecting powers\footnote{this was
done by means of Maple} of $x$ in $ \det(D-(1+x)I)$ we arrive at
the condition \begin{align*}
 (1-\lambda_1)x &= \delta^2
[M(\eps)(1-\lambda_1)+(b_1^2+b_2^2+b_3^2)(1-\lambda_1)+b_4^2(1-\lambda_1)^2]+o(\delta^2)\\
&=\delta^2
[M(\eps)(1-\lambda_1)+((b_1^0)^2+(b_2^0)^2+(b_3^0)^2)(1-\lambda_1)+(b_4^0)^2(1-\lambda_1)^2]+o(\delta^2).
\end{align*}
We have
\[ b_j^0 = \la \psi,e_j\ra = -\la (L_+-1)\psi,e_j\ra = -\la
f,e_j\ra \text{\ \ for\ \ }1\le j\le 3.\] On the other hand,
\[ b_4^0 = \la \psi,e_4\ra = -\la f,e_4\ra + \la\psi,L_+e_4\ra =
-\la f,e_4\ra + \lambda_1 b_4^0\] and thus,
\[ b_4^0 = -(1-\lambda_1)^{-1}\la f,e_4\ra.\]
Since $\lambda_1<0$ in the supercritical case, we obtain that
\begin{align*}
 (1-\lambda_1)x &\le (1-\lambda_1)\delta^2[M(\eps)+\sum_{j=1}^4
\la f,e_j\ra^2] + o(\delta^2)  \le (1-\lambda_1)\delta^2[M(\eps)+
\la f,f\ra ] + o(\delta^2) \\
& = (1-\lambda_1)\delta^2 \la f+\psi,f\ra + o(\delta^2) =
-(1-\lambda_1)\delta^2\la (L_--1)^{-1}f,f\ra + o(\delta^2)
\end{align*}
which yields that $x<0$ for $\delta$ small. But $\eps>0$ small
implies that  $\delta$ is  small and we are done.
\end{proof}

We now present a simple continuity statement which will be important in the following
two sections.

\begin{cor}
\label{cor:imag_proj}
We can choose the $f^{\pm}(\alpha)$ in Proposition~\ref{prop:imag} to be $\calJ$-invariant, i.e.,
$\calJ f^{\pm}(\alpha)=f^{\pm}(\alpha)$. Since $\|f^{\pm}(\alpha)\|_2=1$, they
are therefore unique up to a sign. Choose this sign consistently, i.e., so that $f^{\pm}(\alpha)$
varies continuously with $\alpha$. In that case there is the bound
\beeq
\label{eq:fal_cont}
|\sigma(\alpha_1)-\sigma(\alpha_2)|+\|f^{\pm}(\alpha_1)-f^{\pm}(\alpha_2)\|_2 \le C(\alpha_1)|\alpha_1-\alpha_2|
\eneq
for all $\alpha_1,\alpha_2>0$ which are sufficiently close.
Let $\Pim^{\pm}(\alpha)$ denote the Riesz projection onto $f^{\pm}(\alpha)$,
respectively. Then one has, relative to the operator norm on $L^2\times L^2$,
\beeq
\label{eq:Riesz_cont}
\| \Pim^{\pm}(\alpha_1) - \Pim^{\pm}(\alpha_2)\| \le C(\alpha_1)|\alpha_1-\alpha_2|
\eneq
for all $\alpha_1,\alpha_2$ as above. Moreover, the Riesz projections admit the
explicit representation
\beeq
\label{eq:Riesz_rep}
 \Pim^{\pm}(\alpha) = f^{\pm}(\alpha) \la \cdot, \tilde{f}^{\pm}(\alpha)\ra,
\eneq
where $\Hil(\alpha)^* \tilde{f}^{\pm}(\alpha)=\mp i\sigma \tilde{f}^{\pm}(\alpha)$,
and $\|\tilde{f}^{\pm}(\alpha)\|_2=1.$
\end{cor}
\begin{proof}
By Remark~\ref{rem:Jinv}, $\ker(\Hil(\alpha)\mp i\sigma)$ is $\calJ$-invariant. Thus,
$\calJ f^\pm(\alpha)=\lambda f^\pm(\alpha)$ for some $\lambda\in\Compl$. It is easy to see that
this requires that $|\lambda|^2=1$. Let $e^{2i\beta}=\lambda$. It follows
that $\calJ (e^{i\beta} f^\pm(\alpha)) = e^{i\beta} f^\pm(\alpha)$, leading to our
choice of $\calJ$-invariant eigenfunction. Using the well-known fact that
\[ \ker[\Hil(\alpha)\mp i\sigma(\alpha)]= \ker[(\Hil(\alpha)\mp i\sigma(\alpha))^2], \]
see Lemma~\ref{lem:spec} below,
one easily obtains (by means of the Riesz projections) that
\[
\|(\Hil(\alpha)-z)^{-1}\| \les |z\mp i\sigma(\alpha)|^{-1} \text{\ \ provided\ \ }|z\mp i\sigma(\alpha)|<r_0(\alpha).
\]
In conjunction with the resolvent identity, this yields
\[ |\sigma(\alpha_1)-\sigma(\alpha_2)| \le C(\alpha_1)|\alpha_1-\alpha_2|, \]
as well as~\eqref{eq:Riesz_cont}. However, the latter
clearly implies the remaining bound in~\eqref{eq:fal_cont}.
Finally, by the Riesz representation theorem, we necessarily have
that~\eqref{eq:Riesz_rep} holds with some choice of $\tilde{f}^{\pm}(\alpha)\in L^2\times L^2$.
Since $\Pim^{\pm}(\alpha)^2=\Pim^{\pm}(\alpha)$, one checks that
\[ \Pim^{\pm}(\alpha)^*\tilde{f}^{\pm}(\alpha)=\tilde{f}^{\pm}(\alpha).\]
However, writing down $\Pim^{\pm}(\alpha)$ explicitly shows that
\begin{align*}
 \Pim^{+}(\alpha)^* &= \Big( \frac{1}{2\pi i}\oint_\gamma (-\Hil(\alpha)+zI)^{-1}\, dz \Big)^*
= -\frac{1}{2\pi i}\oint_\gamma (-\Hil(\alpha)^*+\bar{z}I)^{-1}\, d\bar{z} \\
&= \frac{1}{2\pi i}\oint_{-\bar{\gamma}} (-\Hil(\alpha)^*+zI)^{-1}\, dz
\end{align*}
which is equal to the Riesz projection corresponding to the eigenvalue $-i\sigma$ of~$\Hil(\alpha)^*$.
Here $\gamma$ is a small, positively oriented, circle around $i\sigma$. A similar calculation applies to
$\Pim^{-}(\alpha)$. Hence $\Hil(\alpha)^*\tilde{f}^{\pm}(\alpha)=\mp i\sigma(\alpha) \tilde{f}^{\pm}(\alpha)$, as claimed.
In view of~\eqref{eq:Riesz_rep},
\[ \|\tilde{f}^{\pm}(\alpha)\|_2^2 = \|\Pim^{+}(\alpha) \tilde{f}^{\pm}(\alpha)\|_2 \le
\|\tilde{f}^{\pm}(\alpha)\|_2.\]
which implies that $\|\tilde{f}^{\pm}(\alpha)\|_2\le1$. On the other hand,
\[ 1=\|f^{\pm}(\alpha)\|_2 = \|\Pim^{+}(\alpha) f^{\pm}(\alpha)\|_2 \le \|f^{\pm}(\alpha)\|_2
\|\tilde{f}^{\pm}(\alpha)\|_2=\|\tilde{f}^{\pm}(\alpha)\|_2,\]
and we are done.
\end{proof}

\section{The contraction scheme: part I}
\label{sec:contract1}

We now set up the contraction map that will lead to a proof of Theorem~\ref{thm:main}.
According to Lemmas~\ref{lem:Zeq} and~\ref{lem:UPDE}, in order to solve the cubic NLS~\eqref{eq:NLS}
with $\psi(t)=W(t)+R(t)$, we need to find an admissible path
$\pi(t)$ as well as a function
 \[ Z\in C([0,\infty),H^1(\R^3)\times H^1(\R^3))\cap C^1([0,\infty),H^{-1}(\R^3)\times H^{-1}(\R^3))\]
so that $Z(t)$ is $\calJ$-invariant and such that $(\pi(t),Z(t))$
together satisfy~\eqref{eq:Zsys}. As initial conditions we will choose,
with $R_0$ satisfying~\eqref{eq:init_small} and~\eqref{eq:orth} as well as with
some $\alpha=\alpha(R_0)$,
\beeq
\label{eq:init_cond}
 \pi(0):=(\alpha_0,0,0,0),\qquad Z(0):=\binom{R_0}{\bar{R}_0}+hf^+(\alpha)
+\sum_{j=1}^8 a_j \eta_j(\alpha)
\eneq
where  $h\in\R$,
$f^+(\alpha)$ is an eigenvector  of~$\Hil(\alpha)$ with eigenvalue $i\sigma$, $a_j\in \R$, and
$\calN(\alpha)=\{\eta_j(\alpha)\}_{j=1}^8$ is the rootspace of $\Hil(\alpha)^*$.
The contraction argument will be set in the following space.
The parameter $\alpha_0>0$ is fixed.

\begin{defi}
\label{def:boot}
Let $q>2$ be large and fixed. For any sufficiently small $\delta>0$
define
\begin{align*} X_\delta &:=\Big\{ (\pi,Z)\in  {\rm Lip}([0,\infty),\R^8)\times
\big[L^\infty((0,\infty),(H^1(\R^3))^2)\cap L^\infty_{\rm loc}((0,\infty),(Y_q(\R^3))^2)\big]  
\:\big|\:\\
& \text{conditions \eqref{eq:boot1}--\eqref{eq:boot4} are valid}\Big\}
\end{align*}
where $Y_q=\{f\in H^1(\R^3)\:|\: \nabla f\in L^4+L^q\}$ and for a.e.~$t\ge0$,
\begin{align}
|\dot{\alpha}(t)|+|\dot{v}(t)|+|\dot{\gatil}(t)|+|\dot{D}(t)| &\le \delta^2\la t\ra^{-3}
\label{eq:boot1}\\
 \|Z(t)\|_2 +  \|\nabla Z(t)\|_2 &\le c_0\,\delta \label{eq:boot2} \\
t^{\frac32}\|Z(t)\|_\infty &\le c_0\,\delta   \label{eq:boot3} \\
t^{\frac34}\|\nabla Z(t)\|_{L^4+L^q} &\le \delta. \label{eq:boot4}
\end{align}
Here $\la t\ra=(1+t^2)^{\half}$.
We also require that $\pi(0)=(\alpha_0,0,0,0)$.
Here $c_0$ is a sufficiently small constant and $\gatil$ is defined as in Lemma~\ref{lem:sigma_tild}.
Finally, we require that for a.e.~$t\ge0$
\beeq
\label{eq:U_space}
 Z(t) =  \calJ Z(t) = \overline{J Z(t)}
\eneq
where $J=\bm 0&1\\1&0\endm$.
\end{defi}

Note that any path in $X_\delta$ is admissible for small $\delta$.
In \eqref{eq:boot4} one would like to take $q=\infty$, but for technical reasons it is better
to take finite but very large~$q$. We will assume that some such large~$q$ was chosen and it will
be kept fixed. Note that $Y_q\hookrightarrow L^\infty(\R^3)$, so \eqref{eq:boot3} is
meaningful.

In what follows, we will need to deal with several paths simultaneously.
Therefore, our notations will need to indicate relative to which paths
Galilei transforms, root spaces, etc.~are defined. For example, $\calG_\infty(\pi)(t)$
will mean the (vector) Galilei transform from~\eqref{eq:Udef} defined in terms of~$\pi$,
and $\{\xi_j(\pi)(t)\}_{j=1}^8$ will be the set of functions from Definition~\ref{def:rootspace}
which are obtained from~$\pi$.

The contraction scheme is based on the linearized equation~\eqref{eq:Zsys}. Indeed, given
$(\pizer,\Zzer)\in X_\delta$ with $\Zzer=\binom{R^{(0)}}{\bar{R}^{(0)}}$, we solve for
\begin{align}
& i\partial_t Z(t) + \bm \Laplace + 2|W(\pizer)|^2 & W^2(\pizer) \\ -\bar{W}^2(\pizer) & -\Laplace -2|W(\pizer)|^2 \endm Z(t)  \nn\\
 & = \dot{v} \binom{-xe^{i\theta(\pizer)(t)} \phi(\cdot-y(\pizer)(t),\azer(t))}
{xe^{-i\theta(\pizer)(t)} \phi(\cdot-y(\pizer)(t),\azer(t))} +
\dot{\gamma} \binom{-e^{i\theta(\pizer)(t)} \phi(\cdot-y(\pizer)(t),\azer(t))}{e^{-i\theta(\pizer)(t)} \phi(\cdot-y(\pizer)(t),\azer(t))}\nn \\
&\quad + i\dot{\alpha}
\binom{e^{i\theta(\pizer)(t)}\partial_\alpha \phi(\cdot-y(\pizer)(t),\azer(t))}{e^{-i\theta(\pizer)(t)}\partial_\alpha \phi(\cdot-y(\pizer)(t),\azer(t))}
 + i\dot{D}
\binom{-e^{i\theta(\pizer)(t)}\nabla \phi(\cdot-y(\pizer)(t),\azer(t))}{-e^{-i\theta(\pizer)(t)}\nabla \phi(\cdot-y(\pizer)(t),\azer(t))} \nn \\
&\quad + \binom{-2|R^{(0)}|^2W(\pizer)(t)-\bar{W}(\pizer)(t)(R^{(0)})^2-|R^{(0)}|^2R^{(0)}}{2|R^{(0)}|^2\bar{W}(\pizer)(t)+W(\pizer)(t)(\bar{R}^{(0)})^2+|R^{(0)}|^2\bar{R}^{(0)}} \label{eq:Zsyszer}
\end{align}
with initial condition~\eqref{eq:init_cond}. The vector $\dot{\pi}$ here will be determined by
means of the orthogonality condition $\la Z(t),\dot{\xitil}_j\ra =0$ for all $1\le j\le8$, $t\ge0$,
cf.~Lemma~\ref{lem:xiZ}. In this section it will be convenient to work on the level of the
transformed solutions $\Uzer,U$ and the following definition makes this precise. The reader should
note that \eqref{eq:ZtoU}--\eqref{eq:Unew} are the same as~\eqref{eq:Zsyszer}, whereas~\eqref{eq:pinew}
is related to the aforementioned
orthogonality condition on~$Z(t)$.

\begin{defi}
\label{def:contract}
Suppose $(\pizer,\Zzer)\in X_\delta$ and set
\[ \Uzer(t):= M(\pizer)(t)\calG_\infty(t)\Zzer(t),\]
cf.~\eqref{eq:Udef}.
Let $\pizer_\infty$ be the constant vector associated with
\[ \pizer(t)=(\azer(t),v^{(0)}(t),D^{(0)}(t),\gamma^{(0)}(t)) \]
as in Definition~\ref{def:adm}. Let $(\pi,Z)$ be defined as the solutions of the linear problems
\begin{align}
Z(t) &:= \calG_\infty(t)^{-1}M(\pizer)(t)^{-1}U(t) \label{eq:ZtoU} \\
i\partial_t U - \Hil(\azer_\infty)U &= \dot{\pi} \partial_\pi W(\pizer) + N(\Uzer,\pizer)+V(\pizer)U
\label{eq:Unew} \\
\la \dot{\pi} \partial_\pi W(\pizer),\xi_j(\pizer) \ra &= -i\la U,\dot{\xi}_j(\pizer) \ra
- \la U, E(\pizer)\xi_j(\pizer) \ra - \la N(\Uzer,\pizer),\xi_j(\pizer) \ra \label{eq:pinew}
\end{align}
for $1\le j\le8$.
The notation on the right-hand side of~\eqref{eq:Unew} is analogous to that in
\eqref{eq:V}, \eqref{eq:vdotetc}, \eqref{eq:NUpi}, and the matrix
operators $E(\pizer)$ are those from~\eqref{eq:Edef}. The initial conditions are,
with $R_0$ satisfying the smallness condition~\eqref{eq:init_small},
\begin{align}
U(0) &= \calG_\infty(\pizer)(0)Z(0)=\calG_\infty(\pizer)(0)\big[\binom{R_0}{\bar{R}_0}+hf^+(\azer_\infty)+\sum_{j=1}^8 a_j \eta_j(\azer_\infty)\big]\label{eq:Unewinit} \\
\pi(0) &= (\alpha_0,0,0,0) \label{eq:pinewinit}
\end{align}
where $h, \{a_j\}_{j=1}^8\in\Compl$ are constants
(later we will make a unique choice of these constants in terms of the data $(\pizer,\Uzer)$, and in fact
they will be chosen real-valued).
Here, for any $\alpha>0$ we set $\calN(\alpha)=\{\eta_j(\alpha)\}_{j=1}^8$, and we define
$f^{\pm}(\alpha)$ via
\[ \Hil(\alpha)f^{\pm}(\alpha) = \pm i\sigma f^{\pm}(\alpha), \qquad \sigma>0. \]
We are assuming for simplicity that there is a unique pair $\{f^{\pm}(\alpha)\}$
of simple eigenvectors of $\Hil(\alpha)$ with imaginary eigenvalues.
\end{defi}

The main point of this section as well as the next is to show that
the map
\beeq
\label{eq:Psidef}
\Psi:(\pizer,\Zzer)\mapsto (\pi,Z),
\eneq
as given by \eqref{eq:ZtoU}--\eqref{eq:pinew},
defines a contraction on $X_\delta$ relative to a suitable norm
provided the parameters $h, \{a_j\}_{j=1}^8$ are chosen correctly.
As a first step, we show in this section
that $\Psi:X_\delta\to X_\delta$ for $\delta>0$ small and
provided $h$ is chosen properly. Before doing so,
we add some clarifying remarks on Definition~\ref{def:contract}. In particular,
we need to prove the existence of solutions to~\eqref{eq:Unew}, \eqref{eq:pinew}.

We start with a simple technical statement that improves on
Lemma~\ref{lem:rho_infty} by means of the stronger assumptions~\eqref{eq:boot1}.

\begin{lemma}
\label{lem:rho_infty2}
Let $\theta,y$ and $\theta_\infty,y_\infty$ be as in \eqref{eq:theta}, \eqref{eq:y},
 and \eqref{eq:inf_path}, respectively. Let $\rho_\infty$ be as in Lemma~\ref{lem:rho_infty}.
Under the conditions of Definition~\ref{def:boot} the bounds
\[ |\rho_\infty(t,x)| \les \delta^2 (1+|x|)\la t\ra^{-1}, \quad |y(t)-y_\infty(t)|\les \delta^2\la t\ra^{-1}\]
hold for all $t\geq 0$. Moreover,
\beeq
\label{eq:dott}
 |\dot{\rho}_\infty(t,x)| \les \delta^2 (1+|x|)\la t\ra^{-2}, \quad |\dot{y}(t)-\dot{y}_\infty(t)|\les \delta^2\la t
\ra^{-2}
\eneq
for all $t\ge0$. In particular, one has the bounds
\[ \|V(t)\|_{L^1\cap L^\infty}\les \delta^2\la t\ra^{-1},\quad \|\dot{\xi}_j(t)\|_{L^1\cap L^\infty}\les \delta^2 \la t\ra^{-2}, \quad \|E(t)\xi_j(t)\|_{L^1\cap L^\infty}\les \delta^2 \la t\ra^{-2},\]
where $V(t)$ is the matrix from~\eqref{eq:V}, and $E(t)$ is the matrix operator from~\eqref{eq:Edef}.
\end{lemma}
\begin{proof}
In view of the definitions,
\begin{align*}
\theta(t,x+y_\infty)-\theta_\infty(t,x+y_\infty) &= ({v}(t)-{v}_\infty)\cdot(x+2t{v}_\infty+D_\infty) + \gamma(t)-\gamma(\infty) \\
& \quad  -2\int_t^\infty\int_s^\infty ({v}\cdot\dot{{v}}-\alpha\dot{\alpha})(\sigma)\,d\sigma ds.
\end{align*}
Now $|\dot{\gamma}(t)|\le \delta^2\la t\ra^{-2}$ because of~\eqref{eq:dotga}.
Using Definition~\ref{def:boot}  therefore implies the desired
bound on $\rho_\infty$. As for $y(t)-y_\infty(t)$, the definition of $D_\infty$ implies that
\[ y_\infty(t)-y(t) = 2tv_\infty+D_\infty-2\int_0^t v(s)\,ds - D(t)=D(\infty)-D(t)-2\int_t^\infty\int_s^\infty \dot{v}(\sigma)\,d\sigma\,ds,\]
which is no larger than $C\delta^2\la t\ra^{-1}$, as claimed.
\end{proof}

We will make frequent use of the following simple observation: If $U(t)$ and $Z(t)$ are
related by~\eqref{eq:ZtoU}, then $U(t)$ satisfies \eqref{eq:boot2}--\eqref{eq:boot4} iff
$Z(t)$ does (possibly at the loss of a small multiplicative constant).

\begin{lemma}
\label{lem:linexist}
Let $R_0\in H^1(\R^3)$ satisfy $\binom{R_0}{\bar{R}_0}\perp \calN(\alpha_0)^*$.
Given $(\pizer,\Zzer)\in X_\delta$ and any $h, \{a_j\}_{j=1}^8\in \Compl$, there exist unique solutions
\[ (\pi,Z)\in {\rm Lip}([0,\infty),\Compl^8)\times \big[ C([0,\infty),H^1(\R^3)\times H^1(\R^3))
\cap C^1([0,\infty), H^{-1}(\R^3)\times H^{-1}(\R^3)) \big]
\]
of \eqref{eq:ZtoU}--\eqref{eq:pinew} with initial conditions~\eqref{eq:Unewinit}, \eqref{eq:pinewinit}.
Moreover, if $\delta>0$ is sufficiently small, then
for any value of $h\in\Compl$, there is a unique choice of $\{a_j\}_{j=1}^8=\{a_j(h)\}_{j=1}^8\in\Compl^8$ so that
$U(t)$ satisfies the orthogonality conditions
\beeq
\label{eq:UnewOC}
 \la U(t),\xi_j(\pizer)(t)\ra =0 \text{\ \ for all\ \ }t\ge0, \; 1\le j\le 8,
\eneq
cf.~\eqref{eq:OC}. Moreover, if $h\in\R$ then also $\{a_j(h)\}_{j=1}^8\in\R^8$ and in that case
$U(t)$, and therefore also $Z(t)$, is $\calJ$-invariant and $\pi(t)\in \R^8$ for all $t\ge0$.
\end{lemma}
\begin{proof}
For the existence statement, solve \eqref{eq:pinew} for $\dot{\pi}$, which
can be done as in Lemma~\ref{lem:modul}. Plugging the result into~\eqref{eq:Unew}
leads to a linear equation for $U$, which takes the following form:
\beeq
\label{eq:Ualone}
i\partial_t U - \Hil(\azer_\infty)U = \Lin(U,\pizer) + \tilde{N}(\Uzer,\pizer)+V(\pizer)U.
\eneq
Here $\Lin(U,\pizer)$ is the linear term which is obtained by replacing $\dot{\pi}\partial_\pi W(\pizer)$
on the right-hand side of~\eqref{eq:Unew} with the expressions that result by solving~\eqref{eq:pinew}
for~$\dot{\pi}$. See Lemmas~\ref{lem:sigma_tild} and~\ref{lem:modul} for the details of this process.
Moreover, in this way one picks up the final term in~\eqref{eq:pinew} which leads
to the modified nonlinear term $\tilde{N}(\Uzer,\pizer)$ in~\eqref{eq:Ualone}.
We will need to bound this nonlinear term. For this purpose, we record the estimate
\beeq
\label{eq:NUtilbd}
\|\tilde{N}(\Uzer,\pizer)\|_{W^{k,p}} \les \|N(\Uzer,\pizer)\|_{W^{k,p}} +
\min(\|N(\Uzer,\pizer)\|_1,\|N(\Uzer,\pizer)\|_2).
\eneq
Viewed as a linear operator in $U$, $\Lin(\cdot,\pizer)$ has finite rank and co-rank.
In fact, both its range and co-kernel are spanned by eight exponentially decreasing, smooth functions
(which depend on time).
Moreover, by Lemma~\ref{lem:rho_infty2} it satisfies the bound
\beeq
\label{eq:Linbd}
 \|\Lin(U,\pizer)\|_{W^{k,p}} \le C_{k,p}\, \delta^2 \la t\ra^{-2}\|U\|_2
\eneq
for any integer $k\ge0$, and $1\le p\le\infty$.
The equations \eqref{eq:pinew} are chosen precisely in order to ensure that
\[ \frac{d}{dt} \la U(t),\xi_j(\pizer)(t)\ra =0 \text{\ \ for all\ \ }t\ge0, \;1\le j\le8.\]
On the other hand, in Remark~\ref{rem:init} we showed that
\[
\la U(0),\xi_1(\pizer)(0)\ra =0
\]
is the same as (with $\calN(\alpha_0)^*=\{\xi_j(\alpha_0)\}_{j=1}^8$, see~\eqref{eq:N*})
\begin{align*}
0 &= \Big\la \calG_\infty(\pizer)(0)\Big[ \binom{R_0}{\bar{R}_0}+ h f^+(\azer_\infty)+\sum_{j=1}^8 a_j \eta_j(\azer_\infty)\Big] ,
\calG_\infty(\pizer)(0) \binom{\phi(\cdot,\alpha_0)}{\phi(\cdot,\alpha_0)}\Big\ra \\
& =\Big\la \binom{R_0}{\bar{R}_0}, \binom{\phi(\cdot,\alpha_0)}{\phi(\cdot,\alpha_0)}\Big\ra
+  \Big\la h f^+(\azer_\infty)+\sum_{j=1}^8 a_j \eta_j(\azer_\infty) ,
\binom{\phi(\cdot,\alpha_0)}{\phi(\cdot,\alpha_0)}\Big\ra \\
& = h \la f^+(\azer_\infty) , \xi_1(\alpha_0) \ra
+ \sum_{j=1}^8 a_j \la \eta_j(\azer_\infty), \xi_1(\alpha_0) \ra
\end{align*}
and similarly for $\{\xi_k(\pizer)(0)\}_{k=2}^8$. Here we used that $\binom{R_0}{\bar{R}_0}\in \calN^*(\alpha_0)^\perp$ by assumption, as well as  that $\calG_\infty(\pizer)(0) $ is unitary.
Hence~\eqref{eq:UnewOC} holds for all times iff for any~$h\in\Compl$
we can find $\{a_j\}_{j=1}^8\in\Compl^8$ such that
\beeq
\label{eq:ajdt}
 0 = h \la f^+(\azer_\infty) , \xi_\ell(\alpha_0)\ra
+ \sum_{j=1}^8 a_j \la \eta_j(\azer_\infty), \xi_\ell(\alpha_0)\ra \text{\ \ for all\ \ }1\le\ell\le 8.
\eneq
However, $\|\eta_j(\azer_\infty)-\eta_j(\alpha_0)\|\les\delta^2$ because
$|\azer_\infty-\alpha_0|\les \delta^2$. Together with Lemma~\ref{lem:orth}
this implies that the matrix
\[ \calM:=\{\la \eta_j(\azer_\infty), \xi_\ell(\alpha_0)\ra \}_{j,\ell=1}^8 \]
is invertible with norm $\les 1$. Hence $\{a_j\}_{j=1}^8=\{a_j(h)\}_{j=1}^8\in\Compl^8$ is
uniquely determined for any~$h\in\Compl$. For future reference, we note the estimate
\beeq
\label{eq:aj_bd}
\sum_{j=1}^8 |a_j(h)| \les \delta^2\,|h|,
\eneq
which follows from the fact that $\la f^+(\azer_\infty) ,\xi_\ell(\azer_\infty)\ra =0$ and thus
\begin{align*}
&  |\la f^+(\azer_\infty) , \xi_\ell(\alpha_0)\ra|
 = |\la f^+(\azer_\infty) , \xi_\ell(\alpha_0)-\xi_\ell(\azer_\infty)\ra|
 \les \delta^2.
\end{align*}
It is important to realize that the assumption $\binom{R_0}{\bar{R}_0}\perp \calN(\alpha_0)^*$
is precisely used in~\eqref{eq:aj_bd};
if we drop this assumption, then $a_j(h)\not\to0$ as $h\to0$.
Finally, we note that for any $\calJ$-invariant functions $f,g$ one has $\la f,g\ra\in\R$. Hence
for $h\in\R$ both the matrix $\calM$ as well as the vector
\[ \{h \la f^+(\azer_\infty) , \xi_\ell(\alpha_0)\ra\}_{\ell=1}^8\]
 are real-valued so that
in fact $\{a_j(h)\}_{j=1}^8\in\R^8$ (recall that $\calG_\infty(\pizer)(0)$ preserves $\calJ$-invariance).
In view of the preceding, any solution of~\eqref{eq:Ualone} with initial condition~\eqref{eq:Unewinit}
and this choice of $a_j(h)$ will satisfy the orthogonality condition~\eqref{eq:UnewOC} on its interval of existence.

To prove global existence of solutions to the linear problem~\eqref{eq:Ualone},
we perform a contraction argument in $C([0,T],L^2)$ on some finite
time interval $[0,T]$ (one can take $T=1$). Given
any initial condition $U(0)\in H^1\times H^1$, and any $\tilde{U}\in C([0,T],H^1(\R^3))$, we solve
\[
i\partial_t U - \Hil(\azer_\infty)U = \Lin(\tilde{U},\pizer) + \tilde{N}(\Uzer,\pizer)+V(\pizer)\tilde{U}
\]
via the evolution $e^{-it\Hil(\azer_\infty)}$. I.e., write the solution as
\beeq
\label{eq:DuhU}
 U(t)= e^{-it\Hil(\azer_\infty)}U(0)-i\int_0^t e^{-i(t-s)\Hil(\azer_\infty)}
[\Lin(\tilde{U},\pizer) + \tilde{N}(\Uzer,\pizer)+V(\pizer)\tilde{U} ](s)\, ds
\eneq for all times $t\ge0$.
In addition to the bounds in~\eqref{eq:Linbd}, we note the following two bounds: First,
\beeq
\label{eq:NUL2est} \sup_{s\ge0}\|\tilde{N}(\Uzer,\pizer)(s)\|_{L^2} \les \delta^2,
\eneq
which follows from \eqref{eq:NUtilbd}, \eqref{eq:NUbd} of Lemma~\ref{lem:NUbd}, and \eqref{eq:boot2}
applied to~$\Uzer$, and second
\beeq
\label{eq:VUest}
 \sup_{s\ge0}\|V(\pizer)(t)\tilde{U}(s)\|_{H^1} \les \delta^3,
\eneq
which follows from Lemma~\ref{lem:rho_infty2} and again~\eqref{eq:boot2}.
Apply the  linear estimate on the evolution $e^{-it\Hil(\azer_\infty)}$
given by Theorem~\ref{thm:L2stable}. 
Note that in contrast to these estimates,
here we are including the entire discrete spectrum, which possibly leads to exponential growth. However,
on a time interval of length $T=1$, say, we can always achieve that the map $\tilde{U}\mapsto U$
is a contraction in the norm of $C([0,T],L^2\times L^2)$ for small $\delta$. Since the size of
this $\delta$ does not depend on the size of the initial condition, we can restart this procedure
and thus obtain a global solution of~\eqref{eq:DuhU} that belongs to $C([0,\infty),L^2\times L^2)$.
Typically, this solution will grow exponentially.
Next, we wish to estimate the first derivative of~\eqref{eq:DuhU}
by means of the $L^2$ bound in
Theorems~\ref{thm:L2stable} which will lead to the improved statement that
\[ U\in C([0,\infty),H^1(\R^3)\times H^1(\R^3))
\cap C^1([0,\infty), H^{-1}(\R^3)\times H^{-1}(\R^3))\]
solves \eqref{eq:Ualone} in the strong sense.
Inserting this solution $U$ into  equation~\eqref{eq:pinew} then yields the path~$\pi$.
Indeed, simply integrate in time using the initial condition~\eqref{eq:pinewinit}.
It remains to show that for $T=1$,
\beeq
\label{eq:str1}
 \sup_{0\le t\le T} \Big\|\nabla \int_0^t e^{-i(t-s)\Hil(\azer_\infty)} \,\tilde{N}(\Uzer,\pizer) \, ds \Big\|_2
 \les \delta^2.
\eneq
Here we omitted the other terms on the right-hand side of \eqref{eq:Ualone}, i.e., $\Lin(U,\pizer)$
and~$V(\pizer)\tilde{U}$, since they satisfy the bounds~\eqref{eq:Linbd} and~\eqref{eq:VUest},
respectively,
and thus yield the desired $L^2$ estimate on the derivative
(for the issue of interchanging the evolution with a gradient,
see Corollary~\ref{cor:strich_der} below).
In view of~\eqref{eq:NUtilbd} and~\eqref{eq:NabNUbd2} one has
\begin{align}
\|\nabla \tilde{N}(\Uzer,\pizer)(t)\|_2 &\les \min(\|\Uzer(t)\|_\infty^2, \|\Uzer(t)\|_4^2) + \|\Uzer\nabla \Uzer(t)\|_2 + \||\Uzer|^2 \nabla \Uzer(t)\|_2 \nn \\
&\les \delta^2 + \|\Uzer(t)\|_{L^4\cap L^{\frac{2q}{q-2}}} \|\nabla \Uzer(t)\|_{4+q}+\|\,|\Uzer|^2 \nabla \Uzer(t)\|_2 \nn\\
&\les \delta^2 + \delta^2 t^{-\frac34} + \|\,|\Uzer|^2 \nabla \Uzer(t)\|_2 \label{eq:lang1}
\end{align}
The first two terms in~\eqref{eq:lang1} contribute a finite amount
to~\eqref{eq:str1}, as desired. The final term
in~\eqref{eq:lang1}, however, is too singular at $t=0$ and we
therefore need to invoke the Strichartz estimates from
Corollary~\ref{cor:strich_der} to control it. More precisely, we
split $|U|^2 U(s)$ into $P_s(\azer_\infty)|U|^2 U(s)$ and
$(I-P_s(\azer_\infty))|U|^2 U(s)$. The latter does not present a
problem, since the range of $I-P_s(\azer_\infty)$ is spanned by
finitely many Schwartz functions. Thus,
\begin{align*}
& \sup_{0\le t\le T} \Big\|\nabla \int_0^t e^{-i(t-s)\Hil(\azer_\infty)} (I-P_s(\azer_\infty))|\Uzer|^2(s) \Uzer(s) \, ds \Big\|_2  \\
&\les e^{CT} \sup_{s\ge0} \|\nabla(I-P_s(\azer_\infty))|\Uzer|^2(s) \Uzer(s)\|_2 \les e^{CT} \sup_{s\ge0} \|\Uzer(s)\|_3^3 \les e^{CT} \delta^3,
\end{align*}
as desired. For  $P_s(\azer_\infty)|U|^2 U (s)$ we use the following Strichartz estimate:
\begin{align*}
\sup_{0\le t\le T} \Big\| \nabla \int_0^t e^{-i(t-s)\Hil(\azer_\infty)} \, P_s(\azer_\infty) |\Uzer|^2\Uzer(s)\,ds \Big\|_2 & \les \Big( \int_0^T \| \, |\Uzer|^2\Uzer(s) \|_{L^{\frac43}(\R^3)}^{\frac85}\, ds \Big)^{\frac58} \\
& \qquad +\Big( \int_0^T \|\,  |\Uzer|^2 \nabla \Uzer(s) \|_{L^{\frac43}(\R^3)}^{\frac85}\, ds \Big)^{\frac58}.
\end{align*}
It will suffice to deal with the term on the right-hand side containing $\nabla \Uzer$, since the one
without any derivatives is easier.
The corresponding integrand is estimated in terms of~\eqref{eq:boot2} and~\eqref{eq:boot4} as follows:
For all $0<s\le T$,
\begin{align}
 \||\Uzer|^2 \nabla \Uzer(s) \|_{L^{\frac43}(\R^3)}  &\les \|\nabla \Uzer(s)\|_{2} \|\Uzer(s)\|_8^2 \les
 \delta \|\Uzer(s)\|_6^{\frac32} \|\Uzer\|_\infty^{\frac12} \nn \\
&\les \delta\, \delta^{\frac32} (\delta s^{-\frac34}+\delta)^{\frac12} \le C(T) \delta^3 s^{-\frac38},
\label{eq:38}
\end{align}
where we used the Sobolev embedding bound
\[ \|\Uzer(s)\|_\infty \les \|\nabla \Uzer(s)\|_{4+q} + \|\Uzer(s)\|_2\les \delta s^{-\frac34}+\delta. \]
Since $s^{-\frac38}\in L^{\frac85}(0,T)$ we are done. The conclusion is that
 $U\in L^\infty([0,T], H^1(\R^3)\times H^1(\R^3))$.
The continuity in $t$ relative to the $H^1$ norm is
implicit in the above argument, and we skip it. Finally,  time-stepping extends the $H^1$-statement
to all times.

Finally, if $h\in\R$ and $a_j(h)\in\R$ are as above,
then the initial condition~\eqref{eq:Unewinit} is $\calJ$-invariant
by Remark~\ref{rem:Jinv}. Also, we assume that $\pi(0)\in\R^8$.
It remains to derive the system of equations which $(\bar{\pi}, \calJ U)(t)$ obey.
By the assumption that $\pizer(t)\in\R^8$ and $\calJ \Uzer(t)=\Uzer(t)$ for all $t\ge0$, one checks that
\eqref{eq:Unew} implies that
\beeq
\label{eq:calJU}
i\partial_t \calJ U - \Hil(\azer_\infty)\calJ U = \dot{\bar{\pi}} \partial_\pi W(\pizer) + N(\Uzer,\pizer)+V(\pizer)\Uzer,
\eneq
see the proof of Lemma~\ref{lem:UPDE} for more details here.
On the other hand, as in Lemma~\ref{lem:modul}, one obtains the following system which is
equivalent to~\eqref{eq:pinew}, with $E(t)$ as in~\eqref{eq:Edef} and with $\phi=\phi(\cdot,\azer(t))$:
\begin{align}
\dot{\alpha}\alpha^{-1}\|\phi\|_2^2 &= \la U,\dot{\xi}_1(\pizer) \ra -i\la U, E(\pizer)\xi_1(\pizer) \ra
-i \la N(\Uzer,\pizer),\xi_1(\pizer) \ra \label{eq:modul1}\\
\dot{\gatil}\alpha^{-1}\|\phi\|_2^2 &= \la U,\dot{\xi}_2(\pizer) \ra -i \la U, E(\pizer)\xi_2(\pizer) \ra
-i \la N(\Uzer,\pizer),\xi_2(\pizer) \ra \nn \\
2\dot{D_\ell} \|\phi\|_2^2 &= \la U,\dot{\xi}_{2+\ell}(\pizer) \ra -i \la U, E(\pizer)\xi_{2+\ell}(\pizer) \ra
-i \la N(\Uzer,\pizer),\xi_{2+\ell}(\pizer) \ra \nn \\
2\dot{v}_\ell \|\phi\|_2^2 &= \la U,\dot{\xi}_{5+\ell}(\pizer) \ra - i \la U, E(\pizer)\xi_{5+\ell}(\pizer) \ra
-i \la N(\Uzer,\pizer),\xi_{5+\ell}(\pizer) \ra \nn
\end{align}
for all $1\le\ell\le 3$.
This is based on the observation of Lemma~\ref{lem:sigma_tild}, namely that
\[
\dot{\pi}\partial_\pi W(\pizer)= i\Big[  \sum_{\ell=1}^3 (\dot{D_\ell}\,\eta_{5+\ell}(\pizer)-
\dot{v}_\ell\,\eta_{2+\ell}(\pizer)) + \dot{\alpha}\,\eta_2(\pizer)
- \dot{\gatil}\,\eta_1(\pizer) \Big].
\]
Note that $\overline{JE(\pizer)J}=-E(\pizer)$, see \eqref{eq:Edef}.
Taking complex conjugates of the $\dot{\alpha}$ equation~\eqref{eq:modul1} yields
\begin{align*}
\dot{\bar{\alpha}}\bar{\alpha}^{-1}\|\phi\|_2^2 &= \overline{\la JU,J\dot{\xi}_1(\pizer) \ra} +
i\; \overline{\la JU, JE(\pizer)JJ\xi_1(\pizer) \ra}
+ i\; \overline{\la J N(\Uzer,\pizer), J\xi_1(\pizer) \ra} \\
&= \la \overline{JU},\overline{J\dot{\xi}_1(\pizer)} \ra +
i \la \overline{JU}, \overline{JE(\pizer)J}\overline{J\xi_1(\pizer)} \ra
+ i \la \overline{J N(\Uzer,\pizer)},\overline{ J\xi_1(\pizer)} \ra \\
&= \la \calJ U,\dot{\xi}_1(\pizer) \ra -i\la \calJ U, E(\pizer)\xi_1(\pizer) \ra
-i \la N(\Uzer,\pizer),\xi_1(\pizer) \ra
\end{align*}
Taking complex conjugates one therefore derives the following system from the preceding one, see~\eqref{eq:modul1},
\begin{align*}
\dot{\bar{\alpha}}\bar{\alpha}^{-1}\|\phi\|_2^2 &= \la \calJ U,\dot{\xi}_1(\pizer) \ra -i\la \calJ U, E(\pizer)\xi_1(\pizer) \ra
-i \la N(\Uzer,\pizer),\xi_1(\pizer) \ra \\
\dot{\tilde{\bar{\gamma}}}\bar{\alpha}^{-1}\|\phi\|_2^2 &= \la \calJ U,\dot{\xi}_2(\pizer) \ra -i \la \calJ U, E(\pizer)\xi_2(\pizer) \ra
-i \la N(\Uzer,\pizer),\xi_2(\pizer) \ra \\
2\dot{\bar{D}}_\ell \|\phi\|_2^2 &= \la\calJ U,\dot{\xi}_{2+\ell}(\pizer) \ra -i \la\calJ U, E(\pizer)\xi_{2+\ell}(\pizer) \ra
-i \la N(\Uzer,\pizer),\xi_{2+\ell}(\pizer) \ra \\
2\dot{\bar{v}}_\ell \|\phi\|_2^2 &= \la \calJ U,\dot{\xi}_{5+\ell}(\pizer) \ra - i \la\calJ U, E(\pizer)\xi_{5+\ell}(\pizer) \ra
-i \la N(\Uzer,\pizer),\xi_{5+\ell}(\pizer) \ra
\end{align*}
for all $1\le\ell\le 3$. Combining this system with~\eqref{eq:calJU} shows that $(\bar{\pi},\calJ U)$
solves the same equations  as $(\pi, U)$ namely~\eqref{eq:Unew}, \eqref{eq:pinew}. Since their
initial conditions agree if $h\in\R$, we conclude that they agree for all times.
\end{proof}

Next we present a rather simple lemma about bounded solutions to
hyperbolic ODE. This will be the mechanism to determine the unique value of $h$
in~\eqref{eq:Unewinit} so that the solution $U(t)$ constructed in Lemma~\ref{lem:linexist}
remains bounded in $L^2$ for all times.

\begin{lemma}
\label{lem:ODE_stable}
Consider the two-dimensional ODE
\[ \dot{x}(t)-A_0 x(t) = f(t),\qquad x(0)=\binom{x_1(0)}{x_2(0)}\]
where $f=\binom{f_1}{f_2}\in L^\infty([0,\infty),\Compl^2)$ and
$ A_0 = \bm \sigma & 0 \\ 0 & -\sigma \endm $
where $\sigma>0$.
Then $x(t)=\binom{x_1(t)}{x_2(t)}$ remains bounded for all times iff
\beeq
\label{eq:stable}
0 = x_1(0)+\int_0^\infty e^{-\sigma t} f_1(t)\, dt.
\eneq
Moreover, in that case
\beeq
 x_1(t) = - \int_t^\infty e^{-(s-t)\sigma} f_1(s)\, ds, \quad \label{eq:Duh_stab} 
 x_2(t) = e^{-t\sigma}x_2(0)+ \int_0^t e^{-(t-s)\sigma} f_2(s)\, ds. 
\eneq
for all $t\ge0$.
\end{lemma}
\begin{proof}
Clearly,
$x_1(t) = e^{t\sigma}x_1(0) + \int_0^t e^{(t-s)\sigma} f_1(s)\, ds $ and
$ x_2(t) = e^{-t\sigma}x_2(0) + \int_0^t e^{-(t-s)\sigma} f_2(s)\, ds$.
If $\lim_{t\to\infty} e^{-t\sigma}x_1(t)=0$, then
$ 0 = x_1(0) + \int_0^\infty e^{-s\sigma} f_1(s)\, ds,$
which is \eqref{eq:stable}. Conversely, if this holds, then
$ x_1(t) = -e^{t\sigma} \int_t^\infty e^{-s\sigma } f_1(s)\, ds,$
and the lemma is proved.
 \end{proof}

A similar statement can be made for non-autonomous ODEs, see~\cite{schlagODE}.
Recall that we are making all spectral assumptions on the linear operator $\Hil(\alpha)$
as described in Definition~\ref{def:spec_ass}.

\begin{lemma}
\label{lem:U2}
There exists a small constant $c_1$ (depending on the constant $c_0$ in~\eqref{eq:boot2}, \eqref{eq:boot3})
so that the following holds: With
 $\delta>0$  small,
let $R_0\in W^{1,1}(\R^3)\cap W^{1,2}(\R^3)$ satisfy
$P_u^+(\alpha_0)\binom{R_0}{\bar{R}_0}=0$ and
\beeq
\label{eq:R0small}
 \|R_0\|_{W^{1,1}\cap W^{1,2}}\le c_1\delta.
\eneq
Furthermore, let $(\pi,Z)$ be the solution from Lemma~\ref{lem:linexist} for a given
$(\pizer,\Zzer)\in X_\delta$ and $h\in \Compl$, $a_j=a_j(h)$.
Then there exists a unique value of $h\in\R$  so that $(\pi,Z)\in X_\delta$. In other words, under the assumption~\eqref{eq:R0small}, the map $\Psi:X_\delta\to X_\delta$,
see~\eqref{eq:Psidef}.  Moreover,  as a function of $(R_0,\pizer,\Zzer)$, $h=h(R_0,\pizer,\Zzer)$ satisfies
\beeq
\label{eq:1hbd}
|h(R_0,\pizer,\Zzer)| \le  C_0\, \|R_0\|_{W^{1,1}\cap W^{1,2}}^2
\eneq
with a universal constant $C_0$ as well as
\beeq
 |h(R_0,\pizer,\Zzer)-h(R_1,\pi^{(0)},\Zzer)| \le \|R_0-R_1\|_2
\label{eq:hlip}
\eneq
for any $R_0,R_1$ as above.
\end{lemma}
\begin{proof}
Let $(\pizer,\Uzer)\in X_\delta$ be fixed and let $(\pi,U)$ be the solutions constructed in Lemma~\ref{lem:linexist},
with $h\in\R$ arbitrary and $\{a_j(h)\}_{j=1}^8\in\R^8$ the unique choice that guarantees the
orthogonality condition~\eqref{eq:UnewOC}. Moreover,  $\pi$ is real-valued,
and $\calJ U=U$.
We start by decomposing the function $U(t)$ into three pieces $U(t)=\Udis(t)+\Uroot(t)+\Uhyp(t)$
where
\[
\Udis(t)= P_s(\azer_\infty)U(t),\qquad \Uroot(t)=\Proot(\azer_\infty)U(t),\qquad
\Uhyp(t) = \Pim(\azer_\infty)U(t).
\]
Here $\Proot(\alpha)$ and $\Pim(\alpha)$ are the Riesz projections
corresponding to the spectrum at $\{0\}$, and $\{\pm i\sigma\}$ of
$\Hil(\alpha)$, respectively. For ease of notation, let the
elements of the rootspaces of $\Hil(\azer_\infty)$ and
$\Hil(\azer_\infty)^*$ be
\[
 \calN(\azer_\infty) = \{\etazer_j\}_{j=1}^8, \qquad \calN(\azer_\infty)^* = \{\xizer_j\}_{j=1}^8,
\]
respectively, and write accordingly
\beeq
\Uroot(t)= \sum_{j=1}^8 \tilde{a}_j(t) \etazer_j, \qquad \Uhyp(t)= b^+(t) f^+(\azer_\infty) + b^{-}(t) f^{-}(\azer_\infty). \label{eq:U2piece}
\eneq
Since $\Uroot$ and $\Uhyp$ are $\calJ$-invariant, see Remark~\ref{rem:Jinv}, it follows that
$\{\tilde{a}_j\}_{j=1}^8$ and~$b^+,b^-$ are real. Moreover, because of the orthogonality
condition~\eqref{eq:UnewOC}, for all $1\le k\le 8$,
\begin{align}
 0 & =\sum_{j=1}^8 \tilde{a}_j(t) \la \etazer_j,\xi_k(\pizer)(t)\ra + b^+(t) \la f^+(\azer_\infty),\xi_k(\pizer)(t)\ra \label{eq:root_det}\\
& \qquad + b^{-}(t) \la f^{-}(\azer_\infty),\xi_k(\pizer)(t)\ra + \la \Udis(t),\xi_k(\pizer)(t) \ra \nn
\end{align}
for all times $t\ge0$. For small $\delta$ this allows one to solve for $\tilde{a}_j(t)$.
Indeed, by Definition~\ref{def:rootspace} and  Lemma~\ref{lem:rho_infty2},
\[  \sup_{t\ge0}\max_{1\le k\le 8}\|\xizer_k-\xi_k(\pizer)(t)\|_2\les  \delta^2 \la t\ra^{-1}. \]
Also, by Lemma~\ref{lem:orth}, for each $j$ there is $k(j)$ so
that $|\la \eta_j(\pizer)(t),\xi_k(\pizer)(t) \ra|\asymp1$ if
$k=k(j)$ and $=0$ else.
Hence, $|\la \etazer_j,\xi_k(\pizer)(t) \ra|\asymp1$ if $k=k(j)$, but $|\la \etazer_j,\xi_k(\pizer)(t) \ra|\les \delta^2$ if $k\ne k(j)$. Since $j\to k(j)$ is a permutation, it follows that the matrix
$\{ \la \etazer_j,\xi_k(\pizer)(t) \ra \}_{j,k=1}^8$ is invertible with norm of the inverse $\les 1$.
Consequently, there exist uniformly bounded functions $c^{\pm}(t)$, $c_{jk}(t)$ and $d_{jk}(t)$ so that for all $t\ge0$,
\begin{align}
\label{eq:ajsolv}
\tilde{a}_j(t)
&= b^+(t)c^+_j(t)+b^-(t)c^{-}_j(t)+ \sum_{k=1}^8 d_{jk}(t)\la \Udis(t),\xi_k(\pizer)(t) \ra \\
&= \sum_{k=1}^8 c_{jk}(t)\la \Uhyp(t),\xi_k(\pizer)(t) \ra + \sum_{k=1}^8 d_{jk}(t)\la \Udis(t),\xi_k(\pizer)(t) \ra \nn
\end{align}
and therefore, in particular,
\beeq
\label{eq:Uroot_est}
 \|\Uroot(t)\|_{1\cap\infty} \le C(\|\Udis(t)\|_{1+\infty} + \|\Uhyp(t)\|_{1+\infty}),
\eneq
for all $t\ge0$ with a constant $C$ that does not depend on time $t$.
Hence, the solution $U(t)$ is completely determined by $\Udis(t)$ and~$\Uhyp(t)$, and in fact,
 for all $t\ge0$,
\beeq
\label{eq:U_est}
 \|U\|_2 \le C(\|\Udis(t)\|_2 + \|\Uhyp(t)\|_2),
\eneq
with a constant $C$ that does not depend on time $t$.
Clearly,~\eqref{eq:Uroot_est} remains correct with derivatives on the left-hand side (without derivatives
on the right-hand side), and~\eqref{eq:U_est} therefore remains true with derivatives and/or other
$L^p$-norms. For example, it follows that
\beeq
\label{eq:U_est2infty}
 \|U\|_{2+\infty} \le C(\|\Udis(t)\|_{2+\infty} + \|\Uhyp(t)\|_{2+\infty}).
\eneq
In Lemma~\ref{lem:linexist}
we showed that the system \eqref{eq:Unew}, \eqref{eq:pinew} is equivalent to the single equation
\[
i\partial_t U - \Hil(\azer_\infty)U = \Lin(U,\pizer) + \tilde{N}(\Uzer,\pizer)+V(\pizer)U,
\]
see~\eqref{eq:Ualone}.
This equation is $\calJ$-invariant in the sense that $\calJ U$ satisfies the identical equation.
We now split this equation into two equations, one for $\Udis$ and the other for $\Uhyp$.
This yields (with $P_s=P_s(\azer_\infty)$ and $\Pim=\Pim(\azer_\infty)$),
\begin{align}
& i\partial_t \Udis-\Hil(\azer_\infty)\Udis  \nn \\
&= P_s\big[ \Lin_1(\Udis,\pizer)+\Lin_2(\Uhyp,\pizer)+
\tilde{N}(\Uzer,\pizer)+V(\pizer)\Udis+V(\pizer)\Uhyp\big] \label{eq:Udis_eq} \\
& i\partial_t \Uhyp-\Hil(\azer_\infty)\Uhyp \nn \\
&= \Pim\big[ \Lin_1(\Udis,\pizer)+\Lin_2(\Uhyp,\pizer)+
\tilde{N}(\Uzer,\pizer)+V(\pizer)\Udis+V(\pizer)\Udis\big], \nn 
\end{align}
with initial conditions $\Udis(0)=P_s(\azer_\infty) U(0)$ and
$\Uhyp(0)=P_{\rm im}(\azer_\infty)U(0)$, see~\eqref{eq:Unewinit}.
Here the linear terms $\Lin_1$ and $\Lin_2$ are derived from $\Lin$ via expressing
$\Uroot$ as a linear combination of (projections of) $\Udis$ and $\Uhyp$, see~\eqref{eq:root_det}.
More precisely, write
\begin{align*}
& \Lin(U,\pizer)+V(\pizer)U \\
&= \Lin(\Udis,\pizer)+\Lin(\Uhyp,\pizer)+ \sum_{j=1}^8 \tilde{a}_j(t) [\Lin(\etazer_j,\pizer) +
V(\pizer)\etazer_j] + V(\pizer)\Udis+ V(\pizer) \Uhyp\\
                &=: \Lin_1(\Udis,\pizer)+\Lin_2(\Uhyp,\pizer)+ V(\pizer)\Udis+ V(\pizer) \Uhyp ,
\end{align*}
where the second line follows from the first by means of~\eqref{eq:ajsolv}.
Since the functions $\tilde{a}_j(t)$ have the explicit expression in~\eqref{eq:ajsolv},
$\Lin_1,\Lin_2$  satisfy the following estimates as linear operators in the variable~$U$,
\beeq
\label{eq:L1L2}
  \|\Lin_1(U,\pizer)\|_2+\|\Lin_2(U,\pizer)\|_2 \les \delta^2 \la t\ra^{-1}\|U\|_2.
\eneq
see~\eqref{eq:Linbd}. Moreover, they are small as well as of finite rank and
co-rank with  ranges spanned by smooth, exponentially decreasing functions. Hence,
the estimate~\eqref{eq:L1L2} holds with any number derivatives. In particular, we record
the estimate
\beeq
\label{eq:L1L2H1}
  \|\nabla\Lin_1(U,\pizer)\|_{1\cap\infty}+\|\nabla\Lin_2(U,\pizer)\|_{1\cap\infty} \les \delta^2 \la t\ra^{-1}\|U\|_{1+\infty}
\eneq
for future use.
Because of the small parameter $\delta^2$ in~\eqref{eq:L1L2}, we shall solve for $\Udis, \Uhyp$ by means of
a contraction. However, recall that we yet need to determine the value of~$h$.
Thus fix $\Utildis, \Utilhyp\in C([0,\infty),L^2(\R^3)+L^\infty(\R^3))$, with
\beeq
\label{eq:UtilL2} \sup_{t\ge0}\,\la t\ra^{\frac32}\big [\|\Utildis(t)\|_{2+\infty}
+\|\Utilhyp(t)\|_{2+\infty}\big]\le \delta
\eneq
and so that $\calJ \Utildis=\Utildis$ and $\calJ\Utilhyp=\Utilhyp$, and set
\begin{align}
F_1(\Utildis,\Utilhyp) &:= P_s\big[ \Lin_1(\Utildis,\pizer)+\Lin_2(\Utilhyp,\pizer)+
\tilde{N}(\Uzer,\pizer)+V(\pizer)\Utildis+V(\pizer)\Utilhyp\big] \label{eq:F1_def} \\
F_2(\Utildis,\Utilhyp) &:= \Pim\big[ \Lin_1(\Utildis,\pizer)+\Lin_2(\Utilhyp,\pizer)+
\tilde{N}(\Uzer,\pizer)+V(\pizer)\Utildis+V(\pizer)\Utilhyp\big]. \label{eq:F2_def}
\end{align}
In view of the definition~\eqref{eq:NUpi}, \eqref{eq:NUtilbd} and Lemma~\ref{lem:NUbd},
the assumptions on $(\Uzer,\pizer)$ in Definition~\ref{def:boot}, as well as
Lemma~\ref{lem:rho_infty2},  the following bounds hold:
If $t>1$, then
\begin{align}
& \|\tilde{N}(\Uzer,\pizer)(t)\|_{1\cap2}
\les \|\Uzer(t)\|_\infty^2+\|\,|\Uzer|^2\Uzer(t)\,\|_{1\cap 2} \nn \\
&\les \delta^2\la t\ra^{-3} + \delta^2 \la t\ra^{-\frac32} \|\Uzer(t)\|_2
\les \delta^2 \la t\ra^{-\frac32}. \label{eq:nablaU}
\end{align}
On the other hand, if $0<t<1$, then by Sobolev embedding,
\begin{align}
& \|\tilde{N}(\Uzer,\pizer)(t)\|_{1\cap2}
 \les \|\Uzer(t)\|_4^2+\|\,|\Uzer|^2\Uzer(t)\,\|_2 \nn \\
&\les \|\Uzer(t)\|^2_{H^1} +  \|\Uzer(t)\|_{H^1}^3 
\les \delta^2, \nn
\end{align}
so that the bound in~\eqref{eq:nablaU} holds for all $t\ge 0$.
We therefore conclude from~\eqref{eq:UtilL2}, \eqref{eq:L1L2} that for all $t\ge0$
\begin{align}
\label{eq:F12bd}
& \max_{j=1,2}\|F_j(\Utildis,\Utilhyp)(t)\|_{1\cap2} \les \delta^2\la t\ra^{-\frac32}+
\delta^2\la t\ra^{-1}\big[ \|\Utildis(t)\|_{2+\infty}
+\|\Utilhyp(t)\|_{2+\infty}\big], \\
& \max_{j=1,2}\|F_j(\Utildis^{(1)},\Utilhyp^{(1)})(t)-F_j(\Utildis^{(2)},\Utilhyp^{(2)})(t)\|_{1\cap2} \nn \\
&\les \delta^2\la t\ra^{-1}(\|(\Utildis^{(1)}-\Utildis^{(2)})(t)\|_{2+\infty} +
\|(\Utilhyp^{(1)}-\Utilhyp^{(2)})(t)\|_{2+\infty}). \label{eq:F12lip}
\end{align}
Since the  system~\eqref{eq:Udis_eq} is $\calJ$-invariant in the usual sense, it follows that
\[
\calJ F_1(\Utildis,\Utilhyp) = -F_1(\Utildis,\Utilhyp), \; \calJ F_2(\Utildis,\Utilhyp) = -F_2(\Utildis,\Utilhyp).
\]
We now solve
\begin{align}
i\partial_t \Udis-\Hil(\azer_\infty)\Udis &= F_1(\Utildis,\Utilhyp), \qquad \Udis(0)=P_s(\azer_\infty) U(0) \label{eq:til1_eq} \\
i\partial_t \Uhyp-\Hil(\azer_\infty)\Uhyp &= F_2(\Utildis,\Utilhyp),\qquad  \Uhyp(0)=P_{\rm im}(\azer_\infty)U(0). \label{eq:til2_eq}
\end{align}
We can rewrite \eqref{eq:til2_eq} in the equivalent form (using the basis $f^\pm(\azer_\infty)$)
\[ \frac{d}{dt}\binom{b^+}{b^-} + \bm -\sigma(\azer_\infty) & 0 \\ 0 & \sigma(\azer_\infty) \endm
\binom{b^+}{b^-} = \binom{g^+}{g^-}
\]
where $g^{\pm}\in\R$ satisfy
\[ \sup_{t\ge0}\;\la t\ra^{\frac32}\,|g^{\pm}(t)|\les \delta^2, \]
see \eqref{eq:F12bd} and~\eqref{eq:UtilL2}.
We impose the {\em stability condition} from Lemma~\ref{lem:ODE_stable}, i.e.,
\beeq
\label{eq:h_det}
0=b^+(0) + \int_0^\infty e^{-\sigma(\azer_\infty)s} g^+(s)\,ds.
\eneq
We conclude from  the bound on $g^+$ and \eqref{eq:h_det} that
\beeq
\label{eq:b^+bd}
|b^{+}(0)|\les \delta^2.
\eneq
Recall that $b^+(0)$ is the coefficient of $f^\pm(\azer_\infty)$ in~\eqref{eq:U2piece}.
Hence, in view of~\eqref{eq:Unewinit}, we need to choose $h=h(\Utildis,\Utilhyp)$ so that
\begin{align}
\label{eq:hdet}
 b^+(0)f^\pm(\azer_\infty) &= \Pim^+(\azer_\infty)U(0) \\
& = \Pim^+(\azer_\infty)\calG_\infty(\pizer)(0)\Big[\binom{R_0}{\bar{R}_0}
+ hf^+(\azer_\infty) + \sum_{j=1}^8 a_j \eta_j(\azer_\infty)\Big].
\end{align}
We claim that \eqref{eq:b^+bd} implies that $|h|\les \delta^2$. To see this, we of course need
to use the assumption that $\Pim^+(\alpha_0) \binom{R_0}{\bar{R}_0} =0$. Thus, using the notation
and estimates of Corollary~\ref{cor:imag_proj} we conclude that
\begin{align*}
 & \Big\|\Pim^+(\azer_\infty)\calG_\infty(\pizer)(0)\binom{R_0}{\bar{R}_0} \Big\|_{1\cap2} \\
&=  \Big\|\Pim^+(\alpha_0)\binom{R_0}{\bar{R}_0}
-\Pim^+(\azer_\infty)\calG_\infty(\pizer)(0)\binom{R_0}{\bar{R}_0}\Big\|_{1\cap2} \\
&= \Big  \|f^{+}(\alpha_0) \Big\la \binom{R_0}{\bar{R}_0}, \tilde{f}^{+}(\alpha_0)\Big\ra
- f^{+}(\azer_\infty)  \Big\la \binom{R_0}{\bar{R}_0}, \calG_\infty(\pizer)(0)^*\tilde{f}^{+}(\azer_\infty)\Big\ra  \Big\|_{1\cap2} \\
&\les \|f^{+}(\alpha_0)-f^{+}(\azer_\infty)\|_2 \|R_0\|_{2} +  \|R_0\|_2\big\| \calG_\infty(\pizer)(0)^*\tilde{f}^{+}(\alpha_0) - \tilde{f}^{+}(\alpha_0) \big\|_{2}
\les \delta^3.
\end{align*}
To pass to the final inequality, we invoke the bound
\[ \|\calG_\infty(\pizer)(0)^* f -f\|_2 \les \delta^2(\|f\|_{H^1}+\|\la x\ra f\|_2).\]
The appearance of the weight here is the reason we did not estimate
the difference between $\calG_\infty(\pizer)(0)\binom{R_0}{\bar{R}_0}$ and $\binom{R_0}{\bar{R}_0}$.
We conclude from \eqref{eq:hdet}, \eqref{eq:b^+bd}, and~\eqref{eq:aj_bd} that
\beeq
\label{eq:hbd}
|h| \les \delta^2, \quad \sum_{j=1}^8 |a_j(h)|\les \delta^4.
\eneq
Note that this estimate requires the full strength of the assumption $P_u^+(\alpha_0)\binom{R_0}{\bar{R}_0}=0$.  In particular, \eqref{eq:1hbd} holds. It is now easy to prove the Lipschitz bound~\eqref{eq:hlip}.
Indeed, if $h^{(1)}$, $h^{(0)}$ are associated with $R_0, R_1$, respectively, then
\[ (h^{(1)}-h^{(0)})f^+(\azer_\infty) = \Pim^+(\azer_\infty)\Big[\binom{R_0}{\bar{R}_0}-\binom{R_1}{\bar{R}_1}\Big],\]
and \eqref{eq:hlip} follows by taking $L^2$-norms.
For simplicity, let $a_j(\Utildis,\Utilhyp) := a_j(h(\Utildis,\Utilhyp))$.
Define the map $\Psi_0:(\Utildis,\Utilhyp)\mapsto (\Udis,\Uhyp)$ by means of
\begin{align}
 \Udis(t) &= e^{-it\Hil(\azer_\infty)}\Udis(0)-i\int_0^t e^{-i(t-s)\Hil(\azer_\infty)}
F_1(\Utildis,\Utilhyp)(s)\, ds \label{eq:Udissys} \\
 \Udis(0) &= P_s(\azer_\infty)\calG_\infty(\pizer)(0)\Big[ \binom{R_0}{\bar{R}_0}+ h(\Utildis,\Utilhyp) f^+(\azer_\infty)+\sum_{j=1}^8 a_j(\Utildis,\Utilhyp)\eta(\azer_\infty)\Big]   \nn \\
 \Uhyp(t) &= e^{-it\Hil(\azer_\infty)}\Uhyp(0)-i\int_0^t e^{-i(t-s)\Hil(\azer_\infty)}
F_2(\Utildis,\Utilhyp)(s)\, ds \label{eq:Uhypsys}\\
 \Uhyp(0) &= \Pim(\azer_\infty)\calG_\infty(\pizer)(0)\Big[ \binom{R_0}{\bar{R}_0}+ h(\Utildis,\Utilhyp)
f^+(\azer_\infty)+\sum_{j=1}^8 a_j(\Utildis,\Utilhyp)\eta(\azer_\infty)\Big]. \nn
\end{align}
By \eqref{eq:hbd} and \eqref{eq:R0small} one has
\[ \|\Udis(0)\|_{1\cap2}+\|\Uhyp(0)\|_{1\cap2} \les \delta_0+\delta^2, \]
where $\delta_0:=c_1\delta$.
We claim that, with $c_0$ being the small constant from~\eqref{eq:boot2}, \eqref{eq:boot3},
\beeq
\label{eq:Uclaim}
\sup_{t\ge0}\la t\ra^{\frac32}\big[\|\Udis(t)\|_{2+\infty} + \|\Uhyp(t)\|_{2+\infty}\big] \le c_0 \delta.
\eneq
To verify this claim, we use the linear bound of Theorem~\ref{thm:L2stable}  and~\ref{thm:disp}
on~$\Udis$. Because of~\eqref{eq:F12bd} this leads to
\begin{align*}
\|\Udis(t)\|_{2+\infty} &\les \la t\ra^{-\frac32}(\delta_0+\delta^2)
+ \int_0^t \delta^2\la t-s\ra^{-\frac32}\la s\ra^{-\frac32}\, ds \\
& \les \la t\ra^{-\frac32}(c_1\delta+\delta^2) \le c_0\frac{\delta}{2} \la t\ra^{-\frac32}
\end{align*}
for all $t\ge0$ provided $c_1$ was chosen small enough.
Similarly, because of our choice of $h$, see~\eqref{eq:Duh_stab}, we obtain that
 for all $t\ge0$
\[
\|\Uhyp(t)\|_{2} \les \int_t^\infty e^{-\sigma(\azer_\infty)(s-t)} \delta^2\la s\ra^{-\frac32}\, ds +
e^{-\sigma(\azer_\infty) t}(\delta_0+\delta^2) + \int_0^t e^{-\sigma(\azer_\infty)(s-t)}\delta^2
\la s\ra^{-\frac32}\, ds \le c_0\frac{\delta}{2} ,
\]
and \eqref{eq:Uclaim} follows.
Next, we claim that the map $\Psi_0$ is a contraction in the space of $\calJ$-invariant
functions satisfying~\eqref{eq:UtilL2}. To see this, we first remark that there is
the Lipschitz bound
\beeq
\label{eq:hLip}
 |h(\Utildis^{(1)},\Utilhyp^{(1)})-h(\Utildis^{(2)},\Utilhyp^{(2)})| \les \delta^2 \,
\sup_{t\ge0} \, \la t\ra^{\frac32}
(\|(\Utildis^{(1)}-\Utildis^{(2)})(t)\|_{2+\infty} + \|(\Utilhyp^{(1)}-\Utilhyp^{(2)})(t)\|_{2+\infty}).
\eneq
This is a consequence of~\eqref{eq:F12lip} and the explicit expressions for $b^+(0)$ and~$h$
in~\eqref{eq:h_det} and~\eqref{eq:hdet}. Since the coefficients $a_j$ are linear
in~$h$, they satisfy the exact same bounds.
Let $(\Udis^{(j)},\Uhyp^{(j)})=\Psi_0(\Utildis^{(j)},\Utilhyp^{(j)})$ for $j=1,2$.
Subtracting the two equations in~\eqref{eq:Udissys} for $j=1,2$ with the corresponding
difference of initial conditions, and applying Theorems~\ref{thm:L2stable}, \ref{thm:disp} leads to
\[
\sup_{t\ge0}\la t\ra^{\frac32}\|\Udis^{(1)}(t)-\Udis^{(2)}(t)\|_{2+\infty} \les \delta^2 \sup_{s\ge0}
\la s\ra^{\frac32}
(\|(\Utildis^{(1)}-\Utildis^{(2)})(s)\|_{2+\infty} + \|(\Utilhyp^{(1)}-\Utilhyp^{(2)})(s)\|_{2+\infty}).
\]
Note that the difference $(\Uhyp^{(1)}-\Uhyp^{(2)})(t)$ is potentially dangerous, since
we cannot adjust the initial condition to make sure that the stability condition holds.
The point is, however, that this condition holds automatically since
\[ \sup_{t\ge0}\|(\Uhyp^{(1)}-\Uhyp^{(2)})(t)\|_2 < \infty. \]
Lemma~\ref{lem:ODE_stable} therefore guarantees that both~\eqref{eq:stable} and
 \eqref{eq:Duh_stab} hold for $\Uhyp^{(1)}-\Uhyp^{(2)}$.
In particular, one concludes that in this case as well
\[
\sup_{t\ge0}\la t\ra^{\frac32}\|(\Uhyp^{(1)}-\Uhyp^{(2)})(t)\|_2 \les \delta^2
\sup_{s\ge0} \la s\ra^{\frac32}
(\|(\Utildis^{(1)}-\Utildis^{(2)})(s)\|_{2+\infty} + \|(\Utilhyp^{(1)}-\Utilhyp^{(2)})(s)\|_{2+\infty}),
\]
and we have shown that $\Psi_0$ is indeed a contraction.
The conclusion is that there exist $\calJ$-invariant functions $(\Udis,\Uhyp)$ satisfying
\eqref{eq:Uclaim} as well the system~\eqref{eq:Udis_eq}. In addition, there exist $h,a_j(h)\in\R$
as in~\eqref{eq:hbd} determining the initial conditions~\eqref{eq:Unewinit}.

Next, observe that the solution $(\Udis,\Uhyp)$ which we just constructed also satisfies the bound
\beeq
\label{eq:Lzwei_U}
\sup_{t\ge0} (\|\Udis(t)\|_2 + \|\Uhyp(t)\|_2) \le c_0\,\delta.
\eneq
To see this, it suffices to deal with $\Udis(t)$.
Applying Theorem~\ref{thm:L2stable} to~\eqref{eq:Udissys} yields
\begin{align*}
\sup_{t\ge0}\|\Udis(t)\|_2 &\les \|\Udis(0)\|_2 + \int_0^\infty \big[\delta^2\la s\ra^{-1}(\|\Udis(s)\|_{2+\infty}
+\|\Uhyp(s)\|_{2+\infty}) + \delta^2 \la s\ra^{-\frac32} \\
&\qquad + \|V(\pizer)(s)\|_{2\cap\infty}\|(\Udis+\Uhyp)(s)\|_{2+\infty}\big]\, ds
\\
&\les (\delta_0+\delta^2) + \int_0^t \big[\delta^3\la s\ra^{-\frac52}
+ \delta^2 \la s\ra^{-\frac32} + \delta^2 \la s\ra^{-\frac52}\big]\, ds
\\
&\ll c_0\delta,
\end{align*}
as desired.
Retracing our steps we now reintroduce $\Uroot$ via~\eqref{eq:root_det} which leads to
a (weak, i.e., Duhamel) solution $(\pi(t),U(t))$  of the system~\eqref{eq:Unew},
\eqref{eq:pinew} with initial conditions~\eqref{eq:Unewinit}, \eqref{eq:pinewinit}.
Moreover, $U(t)$ is $\calJ$-invariant, and $\pi(t)\in\R^8$ for all $t\ge0$, and the
orthogonality condition~\eqref{eq:UnewOC} holds.

Note that \eqref{eq:Uclaim} insures that
\[ \sup_{t\ge0}\, \la t\ra^{-\frac32} \|U(t)\|_{2+\infty} \le c_0\delta, \;
\sup_{t\ge0} \|U(t)\|\le c_0\delta. \]
Estimating the two terms involving $U$ on the right-hand side of~\eqref{eq:pinew}
by means of this bound and the bounds from Lemma~\ref{lem:rho_infty2} leads
to the estimate
\beeq
\nn
|\dot{\alpha}(t)|+|\dot{v}(t)|+|\dot{\gatil}(t)|+|\dot{D}(t)| \le \delta^2\la t\ra^{-3}
\eneq
for all $t\ge0$ (this is where we need to use the small $c_0$ in~\eqref{eq:boot2}, \eqref{eq:boot3}).
This is precisely~\eqref{eq:boot1}.
Strictly speaking, \eqref{eq:boot1} can be improved by a small factor of $\asymp c_0^2$ on
the right-hand side. However, here and in what follows we ignore this improvement.

It remains to show that our solution $U(t)$ satisfies the other bounds in
\eqref{eq:boot2}--\eqref{eq:boot4}.
Moreover, we have only shown that $(\Udis,\Uhyp)$ satisfies the system~\eqref{eq:Udis_eq}
in the weak (i.e., Duhamel) sense. However, once we prove
\begin{align}
\label{eq:Uclaim2}
\sup_{t\ge0} \Big[\|\nabla \Udis(t)\|_{2} + \|\nabla\Uhyp(t)\|_{2}\Big]
& \leq c_0\delta
\end{align}
it will follow that \eqref{eq:boot2} holds and that
$(\Udis,\Uhyp)$ solves~\eqref{eq:Udis_eq} in the strong sense, i.e., in
\beeq
\label{eq:sol_space}
C([0,\infty),H^1(\R^3)\times H^1(\R^3))\cap C^1([0,\infty),H^{-1}(\R^3)\times H^{-1}(\R^3)).
\eneq
The details of~\eqref{eq:Uclaim2} are as follows:
Clearly, we need to show that the conditions \eqref{eq:boot2}--\eqref{eq:boot4} are consistent
with our contraction scheme. Thus, in addition to~\eqref{eq:UtilL2}
we assume that $\Utildis,\Utilhyp$ satisfy these assumptions
and then check that $\Udis,\Uhyp$ satisfy them as well
First, by the nature of $\Ran(\Pim)$, see
Proposition~\ref{prop:imag},
\[ \|\Uhyp(t)\|_{H^1} \les \|\Uhyp(t)\|_{L^2} \text{\ \ for all\ \ }t\ge0, \]
so that it suffices to deal with $\Udis$. Second, by the Strichartz estimates of
Corollary~\ref{cor:strich_der}, as well as~\eqref{eq:NUtilbd} and Lemma~\ref{lem:NUbd} we
obtain
\begin{align}
& \sup_{0\le t} \Big\|\nabla \int_0^t e^{-i(t-s)\Hil(\azer_\infty)}\, F_1(\Utildis,\Utilhyp)(s) \, ds \Big\|_2  \nn \\
&\les \int_0^\infty \big( \|\Lin_1(\Utildis,\pizer)(s)\|_{H^1}+\|\Lin_2(\Utilhyp,\pizer)(s)\|_{H^1}
+ \|V(\pizer)\Utilhyp(s)\|_{H^1})\, ds  \label{eq:Lin12H1}\\
& \quad + \int_0^\infty \big( \|N(\Uzer,\pizer)(s)\|_1+ \|\Uzer(s)\|_4^2+ \|\Uzer\nabla \Uzer(s)\|_2+
\|V(\pizer)\Utildis(s)\|_{H^1} \big)\, ds  \label{eq:Nuzer}\\
& \quad + \Big( \int_0^\infty \| \, |\Uzer|^2\Uzer(s) \|_{L^{\frac43}(\R^3)}^{\frac85}\, ds \Big)^{\frac58}
  +\Big( \int_0^\infty \|\,  |\Uzer|^2 \nabla \Uzer(s) \|_{L^{\frac43}(\R^3)}^{\frac85}\, ds \Big)^{\frac58}. \label{eq:NUstrich}
\end{align}
In view of~\eqref{eq:L1L2H1} and~\eqref{eq:UtilL2}, the contribution of~\eqref{eq:Lin12H1} is
\[ \les \int_0^\infty \delta^2 \la s\ra^{-1}(\|\Utildis(s)\|_{2+\infty}+\|\Utilhyp(s)\|_{2+\infty})\, ds \les \delta^3. \]
By~\eqref{eq:NUbd1}, \eqref{eq:boot2} and~\eqref{eq:boot3},
\begin{align*}
\| N(\Uzer,\pizer)(s)\|_1 &\les
\min(\|\Uzer(s)\|_\infty^2,\|\Uzer(s)\|_2^2)+\|\Uzer(s)\|_3^3 \\
& \les \delta^2 \la s\ra^{-3} + \delta^2 \la s\ra^{-\frac32} \les \delta^2 \la s\ra^{-\frac32}.
\end{align*}
Furthermore, if $0<s<1$, then we estimate
\begin{align*}
 \|\Uzer(s)\|_4^2 + \|\Uzer\nabla \Uzer(s)\|_2 &\les
\|\Uzer(s)\|_{H^1}^2+\|\Uzer(s)\|_{L^4\cap L^{\frac{2q}{q-2}}} \|\nabla \Uzer(s)\|_{4+q} \\
& \les \delta^2  +\delta^2 s^{-\frac34} \les \delta^2 s^{-\frac34},
\end{align*}
whereas for $s>1$ we have
\begin{align*}
 \|\Uzer(s)\|_4^2 + \|\Uzer\nabla \Uzer(s)\|_2 & \les
\|\Uzer(s)\|_2 \|\Uzer(s)\|_\infty +\|\Uzer(s)\|_{\infty} \|\nabla \Uzer(s)\|_{2} \\
& \les \delta^2 s^{-\frac32} +\delta^2 s^{-\frac32} \les \delta^2 s^{-\frac32}.
\end{align*}
Hence the first three terms in~\eqref{eq:Nuzer} are integrable and their  contribution is $\les \delta^2$.
As far as the final term in~\eqref{eq:Nuzer} is concerned, note that
\begin{align*}
 \|V(\pizer)\Utildis(s)\|_{H^1} & \les \|V(\pizer)\Utildis(s)\|_{2}+\|(\nabla V)(\pizer)\Utildis(s)\|_{2}
+\|V(\pizer)\nabla \Utildis(s)\|_{2} \\
&\les  \delta^3 \la s\ra^{-\frac52} +
\|V(\pizer)\|_{L^4\cap L^{\frac{2q}{q-2}}} \|\nabla \Utildis(s)\|_{4+q} \\
&\les \delta^3 \la s\ra^{-\frac52} + \delta^2 \la s\ra^{-1} \delta s^{-\frac34},
\end{align*}
which  contributes $\les \delta^2$ to~\eqref{eq:Nuzer}. Here we used~\eqref{eq:boot2}-\eqref{eq:boot4},
as well as the bound from Lemma~\ref{lem:rho_infty2}.
Previously, we derived the bound
\[ \|\,  |\Uzer|^2 \nabla \Uzer(s) \|_{L^{\frac43}(\R^3)} \les \delta^3 s^{-\frac38},\]
for $0<s<1$, see \eqref{eq:38}.
On the other hand, if $s>1$, then
\[ \|\,  |\Uzer|^2 \nabla \Uzer(s) \|_{L^{\frac43}(\R^3)} \les \|\nabla \Uzer(s)\|_2 \|\Uzer(s)\|_8^2
\les \delta \|\Uzer(s)\|_\infty^{\frac32}\|\Uzer(s)\|_2^{\frac12} \les \delta^3 s^{-\frac94}.
\]
Hence,
\[ \Big( \int_0^\infty \|\,  |\Uzer|^2 \nabla \Uzer(s) \|_{L^{\frac43}(\R^3)}^{\frac85}\, ds \Big)^{\frac58} \les \delta^3,
\]
and similarly for the term without a gradient in~\eqref{eq:NUstrich}. We have proved~\eqref{eq:Uclaim2}
and therefore also the gradient part of~\eqref{eq:boot2}.

For the remainder of the proof we will subscribe to the somewhat imprecise practice of replacing
the term $V(\pizer)U(t)$ with $V(\pizer)\Uzer(t)$. This will allow us to avoid working with
$\Utildis, \Utilhyp$ and instead make it possible to estimate $U(t)$ directly. The logic here is
that we will only use the bounds \eqref{eq:boot2}--\eqref{eq:boot4} to estimate $\Uzer$, just as
we would in order to show that the contraction scheme is consistent with the remaining conditions
\eqref{eq:boot3}, \eqref{eq:boot4}.

Thus we turn to proving $\|P_s U(t)\|_\infty \le \delta t^{-\frac32}$ for $t>0$.
It will be necessary to bound various terms in $L^1(\R^3)$. One such term is, see~\eqref{eq:NUbd1}
and~\eqref{eq:boot1},
\begin{align*}
\|\dot{\pi}\partial_\pi W(\pizer)+\tilde{N}(\Uzer,\pizer)\|_1 &\les |\dot{\alpha}(s)|+|\dot{v}(s)|+|\dot{\gatil}(s)|+|\dot{D}(s)| + \|\Uzer(s)\|_\infty^2 + \|\Uzer(s)\|_3^3 \\
&\les \delta^2 s^{-3} + \|\Uzer(s)\|_\infty^2+\|\Uzer(s)\|_2^2 \|\Uzer(s)\|_\infty
\les \delta^2 s^{-\frac32},
\end{align*}
provided $s\ge1$. If $0\le s\le1$, then one argues similarly. More precisely, using
\eqref{eq:boot2} and Sobolev embedding instead of~\eqref{eq:boot3}, we obtain
\[ \|\dot{\pi}\partial_\pi W(\pizer)(s)+\tilde{N}(\Uzer,\pizer)(s)\|_1 \les \delta^2.\]
Another term is
\[ \|V(\pizer)\Uzer(s)\|_1 \les \delta^2 \la s\ra^{-\frac52}, \]
valid for all $s\ge0$. This follows from Lemma~\ref{lem:rho_infty2}, and
\eqref{eq:boot2}, \eqref{eq:boot3}.
We now rewrite~\eqref{eq:Unew} via the Duhamel formula.
Let us first consider the case $t\ge1$. Then, by the embedding $W^{1,4}\hookrightarrow L^\infty$,
\begin{align}
\nn &\|P_s U(t)\|_\infty \\
\nn &\le \|e^{-it\Hil(\azer_\infty)}P_s U(0)\|_\infty + \int_0^t \|e^{-i(t-s)\Hil(\azer_\infty)} P_s (\dot{\pi}(s)\partial_\pi W(\pizer)(s)+\tilde{N}(\Uzer,\pizer)(s)+V(\pizer)\Uzer(s))\|_\infty\,ds\\
\nn &\les t^{-\frac32}\|U(0)\|_1 + \!\!\int_0^{t-\half} (t-s)^{-\frac32} (\|\dot{\pi}(s)\partial_\pi W(\pizer)(s)+\tilde{N}(\Uzer,\pizer)(s)\|_1 + \|V(\pizer)\Uzer(s)\|_1)\, ds  \nn \\
& \quad + \int_{t-\half}^t \| e^{-i(t-s)\Hil(\azer_\infty)} P_s (\dot{\pi}(s)\partial_\pi W(\pizer)(s)+\tilde{N}(\Uzer,\pizer)(s)+V(\pizer)\Uzer)(s))\|_{W^{1,4}}\,ds \nn
\end{align}
Invoking the $L^1$ bounds which we just derived on the right-hand side yields
\begin{align}
 &\|P_s U(t)\|_\infty \nn \\
&\les  t^{-\frac32}\|U(0)\|_1 + \int_0^{t-\half} (t-s)^{-\frac32} \delta^2 \la s\ra^{-\frac32}\, ds \nn \\
& \quad + \int_{t-\half}^t (t-s)^{-\frac34}\|\dot{\pi}(s)\partial_\pi W(\pizer)(s)+\tilde{N}(\Uzer,\pizer)(s)+V(\pizer)\Uzer )(s)\|_{W^{1,\frac43}}\,ds \nn \\
&\les t^{-\frac32}\|U(0)\|_1 + \delta^2 t^{-\frac32}
+ \int_{t-\half}^t (t-s)^{-\frac34}  \big[\|\dot{\pi}\partial_\pi W(\pizer)\|_{W^{1,\frac43}} + \||\Uzer|^2\phi\|_{\frac43}
+ \|\Uzer\nabla \Uzer \phi\|_{\frac43} + \||\Uzer|^2\nabla \phi\|_{\frac43}
  \nn \\
&\!\! +\|\,|\Uzer|^2\Uzer\|_{\frac43} +\|\,|\Uzer|^2\nabla \Uzer\|_{\frac43} + \|V(\pizer)\Uzer(s)\|_{\frac43}+
\|\Uzer\nabla V(\pizer)\|_{\frac43}+
\|V(\pizer)\nabla \Uzer\|_{\frac43}\big](s)\,ds. \label{eq:inf_1}
\end{align}
Here we have used the slightly formal notation $|\Uzer|^2\phi$ for the quadratic part of the nonlinearity $N(\Uzer,\pizer)$.
In view of our assumptions \eqref{eq:boot1}-\eqref{eq:boot3} on~$\Uzer$,
\begin{align*}
& \|\dot{\pi}\partial_\pi W(\pizer) (s)\|_{W^{1,\frac43}} + \||\Uzer|^2(s)\phi\|_{\frac43} + \||\Uzer|^2(s)\nabla \phi\|_{\frac43} \les \delta^2 \la s\ra^{-3},  \\
&\|\Uzer\nabla \Uzer(s) \phi\|_{\frac43} \les \|\Uzer(s)\|_\infty \|\nabla \Uzer(s)\|_2 \les \delta^2 s^{-\frac32}
\text{\ \ \ provided\ \ }s\ge1, \\
& \|\Uzer\nabla \Uzer(s) \phi\|_{\frac43} \les \|\Uzer(s)\|_{4} \|\nabla \Uzer(s)\|_2 \les \delta^2
\text{\ \ \ provided\ \ }0<s\le 1, \\
& \|\,|\Uzer|^2\Uzer(s)\|_{\frac43} \les \|\Uzer(s)\|_2^{\frac32}\|\Uzer(s)\|_\infty^{\frac32}\les \delta^3 s^{-\frac94} \text{\ \ \ if\ \ }s\ge1,\\
& \|\,|\Uzer|^2 \Uzer(s)\|_{\frac43} = \|\Uzer(s)\|_4^{3} \les \delta^3 \text{\ \ \ if\ \ }0<s<1.
\end{align*}
Furthermore,
\begin{align*}
& \|\,|\Uzer|^2\nabla \Uzer(s)\|_{\frac43} \le \|\Uzer(s)\|_8^2 \|\nabla \Uzer(s)\|_2 \les \delta^3 s^{-\frac94}
\text{\ \ \ if\ \ }s\ge\half,\\
& \|\,|\Uzer|^2\nabla \Uzer(s)\|_{\frac43} \le \|\Uzer(s)\|_6^{\frac32}\|\Uzer(s)\|_\infty^{\frac12} \|\nabla \Uzer(s)\|_2 \les \delta^3 s^{-\frac34} \text{\ \ if\ \ }0<s<\half,\\
& \|V(\pizer)\Uzer(s)\|_{\frac43} + \|\Uzer(s)\nabla V(\pizer)(s)\|_{\frac43} \les \delta^3 \la s\ra^{-\frac52}.
\end{align*}
It remains to consider the bounds on  $\|V(\pizer)\nabla \Uzer(s)\|_{\frac43}$ for all $s>0$. Note that the latter is always $\les \delta^3\la s\ra^{-1}$,
by~\eqref{eq:boot1} and~\eqref{eq:boot2}, but this is insufficient.
At this point we need to use~\eqref{eq:boot4} to generate more decay. Indeed,
\[ \|V(\pizer)\nabla \Uzer(s)\|_{\frac43} \les \|V(\pizer)\|_{L^{2}\cap L^{\frac{4q}{3q-4}}}
\|\nabla \Uzer(s)\|_{L^4+L^q}
\les \delta^3 \la s\ra^{-1}s^{-\frac34}. \]
Inserting these bounds into~\eqref{eq:inf_1} leads to (recall $t\ge1$),
\[
 \|P_s U(t)\|_\infty \les (\delta_0+\delta^2) t^{-\frac32} + \delta^2 t^{-\frac32} + \int_{t-\half}^t
(t-s)^{-\frac34} \delta^2 s^{-\frac32}\, ds
 \le c_0 \frac{\delta}{2} t^{-\frac32},
\]
provided $c_1\ll c_0$ and $\delta$ are sufficiently small.
To deal with the range $0<t<1$, we perform a similar estimate, using now the small time cases of
the previous bounds:
\begin{align}
& \|P_s U(t)\|_\infty  \nn \\
&\les \|e^{-it\Hil(\azer_\infty)}P_s U(0)\|_\infty \nn \\
&  + \int_0^t \|\nabla e^{-i(t-s)\Hil(\azer_\infty)} P_s
 (\dot{\pi}\partial_\pi W(\pizer)+\tilde{N}(\Uzer,\pizer)+V(\pizer)\Uzer)(s)\|_4\,ds \nn\\
&  + \int_0^t \| e^{-i(t-s)\Hil(\azer_\infty)} P_s
 (\dot{\pi}\partial_\pi W(\pizer)+\tilde{N}(\Uzer,\pizer)+V(\pizer)\Uzer)(s)\|_2\,ds \label{eq:L2stk}\\
\nn &\les t^{-\frac32}\|U(0)\|_1 +  \delta^2+\int_0^{t} (t-s)^{-\frac34} \big[\|\dot{\pi}\partial_\pi W(\pizer)\|_{W^{1,\frac43}} + \||\Uzer|^2\phi\|_{\frac43}
+ \|\Uzer\nabla \Uzer \phi\|_{\frac43} \nn \\
& +\|\,|\Uzer|^2\Uzer\|_{\frac43} + \|\,|\Uzer|^2\nabla \Uzer\|_{\frac43} + \|V(\pizer)\Uzer(s)\|_{\frac43}+ \|\Uzer\nabla V(\pizer)\|_{\frac43}+ \|V(\pizer)\nabla \Uzer\|_{\frac43}\big](s)  \,ds \nn \\
&\les t^{-\frac32}\|U(0)\|_1 + \delta^2 \int_0^t (t-s)^{-\frac34} s^{-\frac34}\,ds
\les t^{-\frac32}(\delta_0+\delta^2) + \delta^2 t^{-\half} \le c_0\frac{\delta}{2} t^{-\frac32}, \label{eq:inf2}
\end{align}
provided $c_1,\delta$ are small. Here~\eqref{eq:L2stk} comes about because of the Sobolev embedding bound
\[ \|f\|_{L^\infty}\les \|\nabla f\|_{4} +\|f\|_2.\]
Since $t<1$ it makes the harmless contribution $\delta^2$ to the following line.

The only remaining bound on the infinite dimensional evolution $P_s U(t)$ is~\eqref{eq:boot4}.
Here $q$ is chosen very large so that the dispersive $L^{q'}(\R^3)\to L^q(\R^3)$ decay is $t^{-\frac32+}$.
The reason we do not take $q=\infty$ can be found in Corollary~\ref{cor:inter} below.
 Thus, with $(t-\half)_+=\max(t-\half,0)$,
\begin{align*}
& \|\nabla P_s U(t)\|_{L^4+L^q} \\
&\le \|\nabla e^{-it\Hil(\azer_\infty)}P_s U(0)\|_{L^4} \\
&\quad  + \int_0^{(t-\half)_+}
\|\nabla e^{-i(t-s)\Hil(\azer_\infty)} P_s (\dot{\pi}\partial_\pi W(\pizer)+\tilde{N}(\Uzer,\pizer)+V(\pizer)\Uzer)(s)\|_q\,ds\\
&\qquad + \int_{t-\half}^{t} \|\nabla e^{-i(t-s)\Hil(\azer_\infty)} P_s (\dot{\pi}\partial_\pi W(\pizer)+\tilde{N}(\Uzer,\pizer)+V(\pizer)\Uzer)(s)\|_4\,ds\;\chi_{[t\ge1]} \\
&\qquad + \int_{0}^{t} \|\nabla e^{-i(t-s)\Hil(\azer_\infty)}P_s (\dot{\pi}\partial_\pi W(\pizer)+\tilde{N}(\Uzer,\pizer)+V(\pizer)\Uzer)(s) \|_4\,ds\;\chi_{[0<t<1]}
\end{align*}
Interchanging $\nabla$ with the evolution as before, and invoking the dispersive estimate yields
\begin{align}
& \|\nabla P_s U(t)\|_{L^4+L^q} \nn \\
& \les t^{-\frac34}\| U(0)\|_{W^{1,\frac43}}
 + \int_0^{(t-\half)_+} (t-s)^{-\frac32+}
\|(\dot{\pi}\partial_\pi W(\pizer)+\tilde{N}(\Uzer,\pizer)+V(\pizer)\Uzer)(s)\|_{W^{1,q'}}\,ds \nn \\
&\qquad + \int_{t-\half}^{t} (t-s)^{-\frac34}\|(\dot{\pi}\partial_\pi W(\pizer)+\tilde{N}(\Uzer,\pizer)+V(\pizer)\Uzer)(s)\|_{W^{1,\frac43}}\,ds\;\chi_{[t\ge1]} \nn \\
&\qquad + \int_{0}^{t} (t-s)^{-\frac34}\|(\dot{\pi}\partial_\pi W(\pizer)+\tilde{N}(\Uzer,\pizer)+V(\pizer)\Uzer)(s)\|_{W^{1,\frac43}}\,ds\;\chi_{[0<t<1]} \label{eq:W1q}
\end{align}
The terms involving the $W^{1,\frac43}$ arose already in~\eqref{eq:inf_1} and~\eqref{eq:inf2} above.
Invoking the same we used there shows that the two final integrals in~\eqref{eq:W1q} contribute
\begin{align}
&\les  \int_{t-\half}^{t} (t-s)^{-\frac34}\delta^2 \la s\ra^{-\frac32}\,ds\;\chi_{[t\ge1]}
 + \int_{0}^{t} (t-s)^{-\frac34}\delta^2 s^{-\frac34}\,ds\;\chi_{[0<t<1]} \nn \\
&\les \delta^2 t^{-\frac32}\,\chi_{[t\ge1]}+\delta^2 t^{-\frac12} \, \chi_{[0<t<1]}
\les \delta^2 t^{-\frac34},\label{eq:lasttwo}
\end{align}
as desired. It remains to bound the integral involving the $W^{1,q'}$ norm in~\eqref{eq:W1q}.
First, we have
\[
\|(\dot{\pi}\partial_\pi W(\pizer)+\tilde{N}(\Uzer,\pizer)+V(\pizer)\Uzer)(s)\|_{L^{q'}}
\les \delta^2\la s\ra^{-3}+\|\Uzer\|_{3q'}^3 \les \delta^2\la s\ra^{-\frac32}.
\]
This bound is a small variation of previous ones, and we skip the details.
Second, we derive the following variant of the $L^{\frac43}$ bounds obtained above:
In view of our assumptions \eqref{eq:boot1}-\eqref{eq:boot4} on~$\Uzer$,
\begin{align*}
& \|\dot{\pi}\partial_\pi W(\pizer) (s)\|_{W^{1,1}} + \||\Uzer|^2(s)\phi\|_{1} + \||\Uzer|^2(s)\nabla \phi\|_{1} \les \delta^2 \la s\ra^{-3},  \\
&\|\Uzer\nabla \Uzer(s) \phi\|_{1} \les \|\Uzer(s)\|_\infty \|\nabla \Uzer(s)\|_2 \les \delta^2 s^{-\frac32}
\text{\ \ \ provided\ \ }s\ge1, \\
& \|\Uzer\nabla \Uzer(s) \phi\|_{1} \les \|\Uzer(s)\|_{2} \|\nabla \Uzer(s)\|_2 \les \delta^2
\text{\ \ \ provided\ \ }0<s\le 1, \\
& \|\,|\Uzer|^2\Uzer(s)\|_{1} \les \|\Uzer(s)\|_2^{2}\|\Uzer(s)\|_\infty\les \delta^3 s^{-\frac32} \text{\ \ \ if\ \ }s\ge1,\\
& \|\,|\Uzer|^2 \Uzer(s)\|_{1} = \|\Uzer(s)\|_3^{3} \les \delta^3 \text{\ \ \ if\ \ }0<s<1.
\end{align*}
Furthermore,
\begin{align*}
& \|\,|\Uzer|^2\nabla \Uzer(s)\|_{1} \le \|\Uzer(s)\|_4^2 \|\nabla \Uzer(s)\|_2 \les \delta^3 s^{-\frac32}
\text{\ \ \ if\ \ }s\ge\half,\\
& \|\,|\Uzer|^2\nabla \Uzer(s)\|_{1} \le \|\Uzer(s)\|_4^{2}\|\nabla \Uzer(s)\|_2 \les \delta^3  \text{\ \ if\ \ }0<s<\half,\\
& \|V(\pizer)\Uzer(s)\|_{1} + \|\Uzer(s)\nabla V(\pizer)(s)\|_{1} \les \delta^3 \la s\ra^{-\frac52}, \\
& \|V(\pizer)\nabla \Uzer(s)\|_{1} \les \|V(\pizer)\|_{L^{\frac43}\cap L^{q'}}
 \|\nabla \Uzer(s)\|_{L^4+L^q}
\les \delta^3 \la s\ra^{-1}s^{-\frac34}.
\end{align*}
We performed these estimates on $L^1$ rather than $L^{q'}$ for simplicity. However, the $L^{q'}$
case is an interpolation of the $L^{\frac43}$ bounds above and the $L^1$ bounds which we just derived.
Thus, the first integral in~\eqref{eq:W1q} which involves the $W^{1,q'}$ norm is
no larger than
\[
\les \int_0^{(t-\half)_+} (t-s)^{-\frac32+} \delta^2\, \la s\ra^{-\frac32}\, ds.
\]
In conjunction with~\eqref{eq:lasttwo} we finally arrive at
\begin{align*}
& \|\nabla P_s U(t)\|_{L^4+L^q} \\
&\les t^{-\frac34}\| U(0)\|_{W^{1,\frac43}} +  \int_0^{(t-\half)_+} (t-s)^{-\frac32+}\delta^2 \la s\ra^{-\frac32+} \, ds + \int_{t-\half}^{t} (t-s)^{-\frac34}\delta^2 \la s\ra^{-\frac32}\,ds\;\chi_{[t\ge1]} \\
& \qquad + \int_{0}^{t} (t-s)^{-\frac34}\delta^2 s^{-\frac34}\,ds\;\chi_{[0<t<1]} \\
&\les t^{-\frac34}(\delta_0+\delta^2) \le \frac{\delta}{2}t^{-\frac34},
\end{align*}
provided $c_1,\delta$ are small.

The conclusion is that $P_s U(t)$ satisfies \eqref{eq:boot3} and~\eqref{eq:boot4}.
As far as  $\Pim U(t)$ is concerned, we claim that it satisfies the stronger
estimate
\beeq
\label{eq:Uimclaim}
\|\Pim U(t)\|_\infty + \|\nabla \Pim U(t)\|_\infty \ll \delta \la t\ra^{-\frac52}
\eneq
for all $t\ge0$. To see this, return to the equation
\[ i\partial_t \Uhyp-\Hil(\azer_\infty)\Uhyp = F_2(\Udis,\Uhyp) \]
see~\eqref{eq:Udis_eq}, which governs the evolution of $\Pim U(t)$.
Here $F_2$ satisfies~\eqref{eq:F12bd}, i.e.,
\[ \|F_2(\Udis,\Uhyp)(t)\|_2\les \delta^2\la t\ra^{-1}(\|\Udis(t)\|_{2+\infty}+\|\Uhyp(t)\|_{2+\infty})\les \delta^3\la t\ra^{-\frac52}.
\]
Writing $\Uhyp(t)= b^+(t) f^+(\azer_\infty) + b^{-}(t) f^{-}(\azer_\infty)$, see~\eqref{eq:U2piece},
we conclude from Lemma~\ref{lem:ODE_stable} that in fact
\[ |b^+(t)|\les \delta^{3}\la t\ra^{-\frac52}, \quad |b^{-}(t)|\les e^{-\sigma(\azer_\infty)t}
\delta_0+ \delta^{3} \la t\ra^{-\frac52},
\]
which implies that
\[ \|\Uhyp(t)\|_{W^{k,p}} \ll \delta \la t\ra^{-\frac52} \]
since the functions $f^{\pm}(\azer_\infty)$ are smooth and decay at infinity.
In particular, \eqref{eq:Uimclaim} holds.  Finally, inserting the bounds on $\Udis$ and $\Uhyp$
into~\eqref{eq:ajsolv} yields
\[ \|\Proot U(t)\|_\infty + \|\nabla \Proot U(t)\|_\infty \ll \delta \la t\ra^{-\frac32}, \]
and the lemma is proved.
\end{proof}

\section{The contraction scheme: part II}
\label{sec:contract2}

It remains to check that $\Psi$ is a contraction.
One guess would be to prove this property in the norm
\beeq
\|(\pi,U)\| = \|\pi\| + \|U\|_{L^\infty([0,\infty),L^2)}, \nn
\eneq
where
\beeq \nn
 \|\pi\|:= \sup_{t\ge0} \la t\ra^{3} (|\dot{\alpha}(t)|+|\dot{v}(t)|+|\dot{\gatil}(t)|+|\dot{D}(t)|).
\eneq
Since the paths are all required to start at the same point $(\alpha_0,0,0,0)$, it suffices to
control their derivatives, which is what this norm does.
Moreover, it is easy to check that the set $X_\delta$ is a complete metric space in this norm.
Unfortunately,  $\Psi$ does not contract in this norm. To see this, suppose we are given
two different data $(\pizer,\Zzer)\in X_\delta$ and $(\pione,\Zone)\in X_\delta$.
Set $(\pitwo,\Ztwo):=\Psi(\pizer,\Zzer)$, $(\pithree,\Zthree):=\Psi(\pione,\Zone)$.
Then the evolutions of $\Ztwo$ and $\Zthree$ are governed by the reference Hamiltonians
$\Hil(\azer_\infty)$ and $\Hil(\aone_\infty)$, respectively.
These Hamiltonians have different spectra, namely their thresholds are $\pm(\azer_\infty)^2$ and
$\pm(\aone_\infty)^2$, respectively.
For this reason one cannot hope to obtain a favorable estimate for $\|\Ztwo(t)-\Zthree(t)\|_2$,
at least for long times. As a model problem, consider the ODEs
\[ i\dot{u}-\alpha_1^2 u =0,\quad i\dot{v}-\alpha_2^2 v=0,\quad u(0)=v(0)\ne0\]
with $\alpha_1\ne\alpha_2$. Evidently, $|u(t)-v(t)|$ will be as large as $|u(0)|$
infinitely often for large~$t$. In contrast to this example, our solutions do disperse
at the rate~$t^{-\frac32}$. Hence, we need to contract in a dispersive norm and
the best decay we can hope for is $t^{-\frac12}$, as can be seen from
\[ |e^{it\alpha_1^2} - e^{it\alpha_2^2}|t^{-\frac32} \les t^{-\half}|\alpha_1-\alpha_2|.\]
Since we incur this loss of $t$ in the $Z$-norm, we also end up
losing $t$ over the decay of $\dot{\pi}$.
  The actual norm is a bit technical, and we introduce it now.

\begin{defi}
\label{def:contr_norm}
For any $(\pi,Z)\in X_\delta$ set $ \dot{\gamma}_*(t) := \dot{\gamma}(t)+\dot{v}\cdot y(t)$
and
\beeq
\label{eq:Xdeltanorm}
\|(\pi,Z)\| := \sup_{0<t\le1}t^{\frac34}|\dot{\pi}_*(t)| + \sup_{t\ge1}t^{2}|\dot{\pi}_*(t)|
+ \sup_{0<t\le 1}t^{\frac34}\|Z(t)\|_{L^4(\R^3)} +
\sup_{t\ge1}t^{\frac12}\|Z(t)\|_{L^4+L^\infty},
\eneq
where ${\pi}_*=(\alpha,v,D,\gamma_*)$. The suprema here are
essential suprema.
\end{defi}

The appearance of $L^4$ rather than the perhaps more natural $L^2$  has to do with
the cubic nonlinearity.
We first make a routine check that $X_\delta$ is indeed complete in this metric.

\begin{lemma}
\label{lem:complete}
If $(\pi,U)\in X_\delta$, then $\|(\pi,U)\|\les \delta$.
For any fixed $\delta>0$ the space $X_\delta$ is a complete metric space in the norm~\eqref{eq:Xdeltanorm}.
\end{lemma}
\begin{proof}
Suppose $\|(\pi_n,U_n)-(\pi_m,U_m)\|\to0$ as $n,m\to\infty$ where $(\pi_n,U_n)\in X_\delta$.
Recall that we are requiring that $\pi_n(0)=(\alpha_0,0,0,0)$.
Thus,
\begin{align*}
& \sup_{t\ge0}[|\alpha_n(t)- \alpha_m(t)|+|v_n(t)-v_m(t)|+|D_n(t)-D_m(t)|] \\
&\le \|(\pi_n,U_n)-(\pi_m,U_m)\| \int_0^\infty(s^{-\frac34}\chi_{[0<s<1]} + s^{-2} \chi_{[s>1]})\,ds\\
&\le C\,\|(\pi_n,U_n)-(\pi_m,U_m)\|.
\end{align*}
Define $(\alpha,D,v):= \lim_{n\to\infty} (\alpha_n,D_n,v_n)$ in the uniform sense. Then by~\eqref{eq:boot1}
\[ |(\alpha,D,v)(t_1)-(\alpha,D,v)(t_2)|=\lim_{n\to\infty}|(\alpha_n,D_n,v_n)(t_1)-(\alpha_n,D_n,v_n)(t_2)|\le \int_{t_1}^{t_2} \delta^2 \la s\ra^{-3}\,ds\]
for all $0\le t_1<t_2$.  This implies that $(\alpha,D,v)\in{\rm Lip}([0,\infty),\R^5)$
and that~\eqref{eq:boot1} also holds for~$(\alpha,D,v)$. Now define $y_\infty(t)=2tv_\infty+D_\infty$
as given by Definition~\ref{def:adm} and set $\dot{\gamma}_*:=\lim_{n\to\infty} (\dot{\gamma}_n)_*$.
Similarly, let $\dot\gamma$ be the limit of $\dot{\gamma}_n=(\dot{\gamma}_n)_* - \dot{v}_n\cdot y_n$.
Then define $\dot{\tilde{\gamma}} := \dot{\gamma}+\dot{v}\cdot y_\infty=\lim_{n\to\infty}\dot{\tilde{\gamma}}_n$. Since each
$\dot{\tilde{\gamma}}_n$ satisfies~\eqref{eq:boot1}, the same argument as before shows that
$\dot{\tilde{\gamma}}$ does, too.
Since
\[ |\dot{\gamma}(t)|\le |\dot{\gatil}(t)|+|\dot{v}(t)|C(1+t)\les \la t\ra^{-2},\]
it follows that $\gamma$ is also Lipschitz and hence $\pi\in{\rm Lip}([0,\infty),\R^8)$,
as required in Definition~\ref{def:boot}.

For a.e.~$t>0$ let $U(t):=\lim_{n\to\infty} U_n(t)$, where the convergence takes place in~$L^4+L^\infty$.
Since $U_n$ satisfy~\eqref{eq:boot2}, for any Schwartz function $\psi$
\begin{align*}
 |\la U(t),\psi\ra|+|\la U(t),\nabla\psi\ra| &= \lim_{n\to\infty}( |\la U_n(t),\psi\ra|+|\la U_n(t),\nabla\psi\ra|) \\
&\le \lim_{n\to\infty} (|\la U_n(t),\psi\ra|+|\la \nabla U_n(t),\psi\ra|) \\
&\le c_0\delta\|\psi\|_2.
\end{align*}
It follows by the usual Hahn-Banach, Riesz-Fischer argument that~\eqref{eq:boot2}
holds for~$U(t)$ and a.e.~$t>0$. For the same reason, the other estimates~\eqref{eq:boot3}, \eqref{eq:boot4} also persist in the limit.
Finally, the $\calJ$ invariance clearly survives in the limit, and we are done.
\end{proof}

Next, we prove a simple technical lemma which will control the variation of various
quantities in the path.

\begin{lemma}
\label{lem:path_var}
Define $\calP:=\{ \pi\in C^1([0,\infty),\R^8)\:|\: \pi(0)=(\alpha_0,0,0,0),\; \|\pi\|\le 1\}$,
where $\|\pi\|$ is the $\pi$-part of~\eqref{eq:Xdeltanorm}, i.e.,
\beeq
\label{eq:pinorm}
 \|\pi\| = \sup_{0<t\le1}t^{\frac34}|\dot{\pi}_*(t)| + \sup_{t\ge1}t^{2}|\dot{\pi}_*(t)|.
\eneq
Then, with $y^{(0)}=y(\pizer)$ and $\theta^{(0)}=\theta(\pizer)$ as in~\eqref{eq:y}, \eqref{eq:theta},
and similarly for $y^{(1)}, \theta^{(1)}$,
\begin{align}
|\azer_\infty-\aone_\infty|+|v^{(0)}_\infty-v^{(1)}_\infty| &\les \|\pizer-\pione\| \label{eq:inf_der} \\
\|e^{i\theta^{(0)}(t)}\phi(\cdot-y^{(0)}(t),\azer(t))-e^{i\theta^{(1)}(t)}\phi(\cdot-y^{(1)}(t),\aone(t)) \|_{L^1\cap L^\infty} &\les \la t\ra\|\pione-\pizer\| \label{eq:fund_der} \\
\|\xitil_j(\pizer)(t)-\xitil_j(\pione)(t)\|_{L^1\cap L^\infty} &\les \la t\ra\|\pizer-\pione\|  \label{eq:xi_der} \\
\|\etatil_j(\pizer)(t)-\etatil_j(\pione)(t)\|_{L^1\cap L^\infty} &\les \la t\ra\|\pizer-\pione\|  \label{eq:eta_der} \\
\|\calS_j(\pizer)(t)-\calS_j(\pione)(t)\|_{L^1\cap L^\infty} &\les \la t\ra\|\pizer-\pione\|
\label{eq:Sj_der}
\end{align}
for all $\pizer,\pione\in\calP$.
For the definitions of the various quantities on the left-hand side see
Lemma~\ref{lem:xiZ}.
The implicit constants here depend on~$\alpha_0$ but
are otherwise absolute.
\end{lemma}
\begin{proof}
In view of the definition of $\pizer_\infty,\pione_\infty$ in Definition~\ref{def:adm}
we have the following bounds:
\begin{align*}
& \sup_{t\ge0}[|\azer(t)-\aone(t)|+|v^{(0)}(t)-v^{(1)}(t)|+|D^{(0)}(t)-D^{(1)}(t)|] \\
& \le
\int_0^\infty (|\dot{\alpha}^{(0)}(s)-\dot{\alpha}^{(1)}(s)|+ |\dot{v}^{(0)}(s)-\dot{v}^{(1)}(s)|+|\dot{D}^{(0)}(s)-\dot{D}^{(1)}(s)|)\, ds\\
&\les \|\pizer-\pione\| \int_0^\infty (s^{-\frac34}\chi_{[0<s<1]} + s^{-2} \chi_{[s>1]})\, ds \les \|\pizer-\pione\|.
\end{align*}
In particular,
\[
|\azer_\infty-\aone_\infty| + |v^{(0)}_\infty - v^{(1)}_\infty| \les \|\pione-\pizer\|,
\]
which is \eqref{eq:inf_der}.
Moreover, (recall that $\pizer(0)=\pione(0)=(\alpha_0,0,0,0)$)
\begin{align*}
|y^{(0)}(t)-y^{(1)}(t)| &\les \int_0^t |v^{(0)}(s)-v^{(1)}(s)|\, ds + |D^{(0)}(t)-D^{(1)}(t)| \\
&\les \Big[\int_0^t \int_0^s \la \sigma\ra^{-2}\, d\sigma ds + \int_0^t \la s\ra^{-2}\,ds\Big]
\|\pizer-\pione\| \\
& \les \la t\ra \|\pizer-\pione\|.
\end{align*}
and
\begin{align*}
|\gamma^{(0)}(t)-\gamma^{(1)}(t)| &\le \int_0^t |\dot{\gamma}^{(0)}_*(s)-\dot{\gamma}^{(1)}_*(s)|\, ds
+ \int_0^t |\dot{v}^{(0)}(s)\cdot y^{(0)}(s)-\dot{v}^{(1)}(s)\cdot y^{(1)}(s)|\, ds \\
& \les  \int_0^t \la s\ra^{-2}\, ds\,\|\pizer-\pione\| + \int_0^t |\dot{v}^{(0)}(s)-\dot{v}^{(1)}(s)|\la s\ra\, ds
\\
&\les \log(2+t)\|\pizer-\pione\| \les \la t\ra\|\pizer-\pione\|.
\end{align*}
Let $\theta$ be as in~\eqref{eq:theta}. Then
\begin{align*}
|\theta^{(0)}(t,x)-\theta^{(1)}(t,x)| &\les |v^{(0)}(t)-v^{(1)}(t)||x| + \int_0^t (|v^{(0)}(s)-v^{(1)}(s)|
+ |\azer(s)-\aone(s)|)\, ds \\
& \quad + |\gamma^{(0)}(t)-\gamma^{(1)}(t)|
 \les (|x|+\la t\ra) \|\pizer-\pione\|.
\end{align*}
The estimate \eqref{eq:fund_der} now follows easily. Indeed, observe that $|x|$ behaves like $t$
in this context. The other estimates \eqref{eq:xi_der}, \eqref{eq:eta_der}, and~\eqref{eq:Sj_der} are easily
deduced from~\eqref{eq:fund_der}.
\end{proof}

We will use the following simple extension of the
contraction principle. Of course if is well-known, but we still record it here.

\begin{lemma}
\label{lem:param_contr}
Let $S\subset X$ be a closed subset of a Banach space $X$ and $T\subset Y$
an arbitrary subset of some normed space~$Y$.
Suppose that $A:S\times T \to S$ so that  with some $0<\gamma<1$
\begin{align*}
 \sup_{t\in T}\|A(x,t)-A(y,t)\|_X &\le \gamma \|x-y\|_X \text{\ \ for all \ \ }x,y\in S,\\
  \sup_{x\in S}\|A(x,t_1)-A(x,t_2)\| &\le C_0 \|t_1-t_2\|_Y \text{\ \ for all \ \ }t_1,t_2\in T.
\end{align*}
Then for every $t\in T$ there exists a unique fixed-point $x(t)\in S$ such that
$A(x(t),t)=x(t)$. Moreover, these points satisfy the bounds
\[ \|x(t_1)-x(t_2)\|_X \le \frac{C_0}{1-\gamma}\|t_1-t_2\|_Y\]
for all $t_1,t_2\in T$.
\end{lemma}
\begin{proof}
Clearly,
$x(t)= \lim_{n\to\infty} A(x_n(t),t)$
where for some fixed (i.e., independent of $t$) $x_0$
\[ x_0(t):= x_0, \quad x_{n+1}(t)= A(x_n(t),t).\]
Then inductively,
\begin{align*}
\|x_{n+1}(t_1)-x_{n+1}(t_2)\|_X & \le \| A(x_n(t_1),t_1) - A(x_n(t_2),t_1)\|_X + \|A(x_n(t_2),t_1)
- A(x_n(t_2),t_2)\|_X \\
&\le \gamma \|x_n(t_1)-x_n(t_2)\|_X + C_0 \|t_1-t_2\|_Y \\
&\le C_0 \sum_{k=0}^n\gamma^{k}\;\|t_1-t_2\|_Y
\end{align*}
for all $n\ge0$. Passing to the limit $n\to\infty$ proves the lemma.
\end{proof}

We are now ready to state the contraction property of $\Psi$.

\begin{lemma}
\label{lem:Psi_contract}
Under the hypotheses of Theorem~\ref{thm:main} the map $\Psi:X_\delta\to X_\delta$
is a contraction in the norm~\eqref{eq:Xdeltanorm}.
Thus $\Psi$ has a fixed point $(\pi,Z)\in X_\delta$, which
is completely determined by~$R_0$.
Hence, the function $h(R_0,\pizer,\Zzer)$ now becomes a
function $h=h(R_0)$ of~$R_0$ alone. It satisfies~\eqref{eq:1hbd} as well as
the Lipschitz bound
\beeq
\label{eq:hwaschno}
|h(R_0)-h(R_1)|\les \delta\trip R_0-R_1\trip
\eneq
for any $R_0,R_1$ satisfying $P_u^+(\alpha_0)\binom{R_j}{\bar{R_j}}=0$, and
$\trip R_j\trip \le c_1\delta$, $j=1,2$.
\end{lemma}
\begin{proof}
Let $(\pizer,\Zzer), (\pione,\Zone)\in X_\delta$ and set
\[ (\pitwo,\Ztwo):=\Psi(\pizer,\Zzer), \quad (\pithree,\Zthree):=\Psi(\pione,\Zone), \]
as well as
\begin{align*}
 \Uzer(t) &:= M(\pizer)(t)\calG_\infty(\pizer)(t)\Zzer(t), \;
 \Utwo(t) := M(\pizer)(t)\calG_\infty(\pizer)(t)\Ztwo(t) \\
\Uone(t) &:= M(\pione)(t)\calG_\infty(\pione)(t)\Zone(t), \;
\Uthree(t) := M(\pione)(t)\calG_\infty(\pione)(t)\Zthree(t).
\end{align*}
Hence, by definition of $\Psi$ we have the linear problems
\begin{align}
i\partial_t \Utwo - \Hil(\azer_\infty)\Utwo &= \dot{\pi}^{(2)} \partial_\pi W(\pizer) + N(\Uzer,\pizer)+V(\pizer)\Uzer
\label{eq:U2new} \\
\la \dot{\pi}^{(2)} \partial_\pi W(\pizer),\xi_j(\pizer) \ra &= -i\la \Utwo,\dot{\xi}_j(\pizer) \ra
- \la \Utwo, E(\pizer)\xi_j(\pizer) \ra - \la N(\Uzer,\pizer),\xi_j(\pizer) \ra \nn
\end{align}
for $1\le j\le8$ and
\begin{align}
i\partial_t \Uthree - \Hil(\aone_\infty)\Uthree &= \dot{\pi}^{(3)} \partial_\pi W(\pione) +
N(\Uone,\pione)+V(\pione)\Uone
\label{eq:U3new} \\
\la \dot{\pi}^{(3)} \partial_\pi W(\pione),\xi_j(\pione) \ra &= -i\la \Uthree,\dot{\xi}_j(\pione) \ra
- \la \Uthree, E(\pione)\xi_j(\pione) \ra - \la N(\Uone,\pione),\xi_j(\pione) \ra \nn
\end{align}
for $1\le j\le8$.
 The initial conditions are
\begin{align}
\Utwo(0) &= \calG_\infty(\pizer)(0)\Big[\binom{R_0}{\bar{R}_0}+h^{(0)}f^+(\azer_\infty)+\sum_{j=1}^8 a^{(0)}_j \eta_j(\azer_\infty)\Big] \label{eq:U2newinit} \\
\Uthree(0) &= \calG_\infty(\pione)(0)\Big[\binom{R_0}{\bar{R}_0}+h^{(1)}f^+(\aone_\infty)+\sum_{j=1}^8 a^{(1)}_j \eta_j(\aone_\infty)\Big]\label{eq:U3newinit} \\
\pitwo(0) &= \pithree(0) = (\alpha_0,0,0,0), \label{eq:pi23newinit}
\end{align}
where we have set
\[ h^{(0)}:=h(R_0,\pizer,\Zzer), \quad h^{(1)}:=h(R_0,\pione,\Zone)\]
for simplicity, and similarly for $a_j=a_j(h)$.
We {\em cannot} compare $\Utwo$ and $\Uthree$ because they are given in terms of reference
Hamiltonians which involve the vectors $\pizer_\infty$ and~$\pione_\infty$ and the latter cannot
be compared. Indeed, since we
only know that $|\dot{\pi}^{(0)}(t)-\dot{\pi}^{(1)}(t)|\le \la t\ra^{-2}\|\pizer-\pione\|$,
the best estimate on the ``terminal translation'' $D_\infty$ here would be
\[ |D^{(0)}_\infty-D^{(1)}_\infty|\les \|\pizer-\pione\|\log \|\pizer-\pione\|^{-1},\]
which is too weak for the contraction.
Therefore, we return to the system~\eqref{eq:Zsys}, \eqref{eq:Zsyszer}. More precisely, with
$\Zzer=\binom{R^{(0)} }{\bar{R}^{(0)}}$ and $\Zone=\binom{R^{(1)} }{\bar{R}^{(1)}}$
one has the systems
\begin{align*}
& i\partial_t \Ztwo(t) + \bm \Laplace + 2|W(\pizer)|^2 & W^2(\pizer) \\ -\bar{W}^2(\pizer) & -\Laplace -2|W(\pizer)|^2 \endm \Ztwo(t)  \\
 & = \dot{v}^{(2)} \binom{-xe^{i\theta(\pizer)(t)} \phi(\cdot-y(\pizer)(t),\azer(t))}
{xe^{-i\theta(\pizer)(t)} \phi(\cdot-y(\pizer)(t),\azer(t))} +
\dot{\gamma}^{(2)} \binom{-e^{i\theta(\pizer)(t)} \phi(\cdot-y(\pizer)(t),\azer(t))}{e^{-i\theta(\pizer)(t)} \phi(\cdot-y(\pizer)(t),\azer(t))}\nn \\
&\quad + i\dot{\alpha}^{(2)}
\binom{e^{i\theta(\pizer)(t)}\partial_\alpha \phi(\cdot-y(\pizer)(t),\azer(t))}{e^{-i\theta(\pizer)(t)}\partial_\alpha \phi(\cdot-y(\pizer)(t),\azer(t))}
 + i\dot{D}^{(2)}
\binom{-e^{i\theta(\pizer)(t)}\nabla \phi(\cdot-y(\pizer)(t),\azer(t))}{-e^{-i\theta(\pizer)(t)}\nabla \phi(\cdot-y(\pizer)(t),\azer(t))} \nn \\
&\quad + \binom{-2|R^{(0)}|^2W(\pizer)(t)-\bar{W}(\pizer)(t)(R^{(0)})^2-|R^{(0)}|^2R^{(0)}}{2|R^{(0)}|^2\bar{W}(\pizer)(t)+W(\pizer)(t)(\bar{R}^{(0)})^2+|R^{(0)}|^2\bar{R}^{(0)}}.
\end{align*}
and
\begin{align*}
& i\partial_t \Zthree(t) + \bm \Laplace + 2|W(\pione)|^2 & W^2(\pione) \\ -\bar{W}^2(\pione) & -\Laplace -2|W(\pione)|^2 \endm \Zthree(t)  \\
 & = \dot{v}^{(3)} \binom{-xe^{i\theta(\pione)(t)} \phi(\cdot-y(\pione)(t),\aone(t))}
{xe^{-i\theta(\pione)(t)} \phi(\cdot-y(\pione)(t),\aone(t))} +
\dot{\gamma}^{(3)} \binom{-e^{i\theta(\pione)(t)} \phi(\cdot-y(\pione)(t),\aone(t))}{e^{-i\theta(\pione)(t)} \phi(\cdot-y(\pione)(t),\aone(t))}\nn \\
&\quad + i\dot{\alpha}^{(3)}
\binom{e^{i\theta(\pione)(t)}\partial_\alpha \phi(\cdot-y(\pione)(t),\aone(t))}{e^{-i\theta(\pione)(t)}\partial_\alpha \phi(\cdot-y(\pione)(t),\aone(t))}
 + i\dot{D}^{(3)}
\binom{-e^{i\theta(\pione)(t)}\nabla \phi(\cdot-y(\pione)(t),\aone(t))}{-e^{-i\theta(\pione)(t)}\nabla \phi(\cdot-y(\pione)(t),\aone(t))} \nn \\
&\quad + \binom{-2|R^{(1)}|^2W(\pione)(t)-\bar{W}(\pione)(t)(R^{(1)})^2-|R^{(1)}|^2R^{(1)}}{2|R^{(1)}|^2\bar{W}(\pione)(t)+W(\pione)(t)(\bar{R}^{(1)})^2+|R^{(1)}|^2\bar{R}^{(1)}}.
\end{align*}
Using the notations from Lemma~\ref{lem:xiZ} we rewrite these systems in the form
\begin{align}
i\partial_t \Ztwo - \Hil(\pizer(t))\Ztwo &=
 i\Big[  \sum_{\ell=1}^3 (\dot{D}^{(2)}_\ell\etatil^{(0)}_{5+\ell}-\dot{v}^{(2)}_\ell\etatil^{(0)}_{2+\ell})
+ \dot{\alpha}^{(2)}\etatil^{(0)}_2 - \dot{\gamma}_*^{(2)}\etatil^{(0)}_1 \Big]  + N_*(\Zzer,\pizer) \nn\\
&=: \dot{\pi}_*^{(2)} \partial_\pi \tilde{W}(\pizer) + N_*(\Zzer,\pizer)
\label{eq:Ztwoeq} \\
i\partial_t \Zthree - \Hil(\pione(t))\Zthree &=
 i\Big[  \sum_{\ell=1}^3 (\dot{D}^{(3)}_\ell\etatil^{(1)}_{5+\ell}-\dot{v}^{(3)}_\ell\etatil^{(1)}_{2+\ell})
+ \dot{\alpha}^{(3)}\etatil^{(1)}_2 - \dot{\gamma}_*^{(3)}\etatil^{(1)}_1 \Big]  + N_*(\Zone,\pione) \nn\\
&=: \dot{\pi}_*^{(3)} \partial_\pi \tilde{W}(\pione) + N_*(\Zone,\pione),
\label{eq:Zthreeeq}
\end{align}
where $\etatil_j^{(0)}:=\etatil_j(\pizer)$, $\etatil_j^{(1)}:=\etatil_j(\pione)$,
$ \dot{\gamma}_*(t):= \dot{\gamma}+ \dot{v}\cdot y(t)$,
and $N_*$ is defined in the obvious way.
By construction, $\Uthree, \Utwo$ satisfy the orthogonality conditions
\[ \la \Utwo(t), \xi_j(\pizer)(t) \ra = \la \Uthree,\xi_j(\pione)(t) \ra =0 \]
for all $1\le j\le8 $ and $t\ge0$. By Lemma~\ref{lem:xiZ} these are equivalent to
\[ \la \Ztwo(t), \xitil_j(\pizer)(t) \ra = \la \Zthree,\xitil_j(\pione)(t)\ra =0.\]
Taking the scalar products of \eqref{eq:Ztwoeq} and~\eqref{eq:Zthreeeq} with $\xitil_j(\pizer)$
and $\xitil_j(\pione)$, respectively, leads to the following modulation equations on the
level of~$Z(t)$:
\begin{align}
\la \dot{\pi}_*^{(2)} \partial_\pi \tilde{W}(\pizer),\xitil_j(\pizer) \ra &=
\la \Ztwo, \dot{\pi}_*^{(0)} \calS_j(\pizer) \ra  - \la N_*(\Zzer,\pizer), \xitil_j(\pizer) \ra
\label{eq:modul2} \\
 \la \dot{\pi}_*^{(3)} \partial_\pi \tilde{W}(\pione),\xitil_j(\pione) \ra &=
\la \Zthree, \dot{\pi}_*^{(1)} \calS_j(\pione) \ra  - \la N_*(\Zone,\pione), \xitil_j(\pione) \ra.
\label{eq:modul3}
\end{align}
Here we used the notation from~\eqref{eq:calSj}.
Subtracting \eqref{eq:Ztwoeq}, \eqref{eq:Zthreeeq}, and \eqref{eq:modul2}, \eqref{eq:modul3},
respectively, we obtain the equations that will provide the estimates for the contraction
step:
\begin{align}
& i\partial_t (\Zthree-\Ztwo) - \Hil(\pizer(t))(\Zthree-\Ztwo)
= (\dot{\pi}_*^{(3)}-\dot{\pi}_*^{(2)})
\partial_\pi \tilde{W}(\pizer) + V(\pizer,\pione)\Zthree \nn \\
& \quad + N_*(\Zone,\pione)- N_*(\Zzer,\pizer)
 + \dot{\pi}_*^{(3)} (\partial_\pi \tilde{W}(\pione) - \partial_\pi \tilde{W}(\pizer)) \label{eq:diffeq1} \\
& (\Zthree-\Ztwo)(0) =  h^{(1)}f^+(\aone_\infty)+\sum_{j=1}^8 a^{(1)}_j \eta_j(\aone_\infty)
-[h^{(0)}f^+(\azer_\infty)+\sum_{j=1}^8 a^{(0)}_j \eta_j(\azer_\infty)] \label{eq:diff_init} \\
& \la (\dot{\pi}_*^{(3)}-\dot{\pi}_*^{(2)}) \partial_\pi \tilde{W}(\pizer),\xitil_j(\pizer) \ra
= \la \Zthree-\Ztwo, \dot{\pi}_*^{(0)} \calS_j(\pizer) \ra + \la \Zthree, \dot{\pi}_*^{(1)} \calS_j(\pione)-\dot{\pi}_*^{(0)} \calS_j(\pizer) \ra \nn \\
& \quad +
\la N_*(\Zone,\pione)- N_*(\Zzer,\pizer) , \xitil_j(\pione) \ra +
\la N_*(\Zzer,\pizer) , \xitil_j(\pione)- \xitil_j(\pizer) \ra  \nn \\
& \quad +
\la \dot{\pi}_*^{(3)} ( \partial_\pi \tilde{W}(\pizer)-\partial_\pi \tilde{W}(\pione)),\xitil_j(\pizer) \ra + \la \dot{\pi}_*^{(3)} \partial_\pi \tilde{W}(\pizer),\xitil_j(\pizer)-\xitil_j(\pione) \ra
\label{eq:diffeq2}
\end{align}
Here
\begin{align*}
&\tilde{V}(\pizer,\pione) := \Hil(\pione)-\Hil(\pizer) \nn \\
& = \bm 2(\phi^2(\cdot-y^{(0)},\azer)-\phi^2(\cdot-y^{(1)},\aone)) &
e^{2i\theta^{(0)}} \phi(\cdot-y^{(0)},\azer)-e^{2i\theta^{(1)}} \phi(\cdot-y^{(1)},\aone) \\
-e^{-2i\theta^{(0)}} \phi(\cdot-y^{(0)},\azer)+e^{-2i\theta^{(1)}} \phi(\cdot-y^{(1)},\aone)
&  -2(\phi^2(\cdot-y^{(0)},\azer)-\phi^2(\cdot-y^{(1)},\aone))
 \endm.
\end{align*}
In view of~\eqref{eq:fund_der},
\beeq
\label{eq:V01_est}
\|\tilde{V}(\pizer,\pione)\|_{L^1\cap L^\infty} \les \la t\ra \|\pizer-\pione\|.
\eneq
Set $\calTzer(t) := M(\pizer)(t)\calG_\infty(\pizer)(t)$ and define $\Util3:=\calTzer(t)\Zthree(t)$.
Hence,
\[ (\Util3-\Utwo)(t) = \calTzer(t) (\Zthree(t)-\Ztwo(t)) \]
satisfies the transformed equation
\begin{align}
& i\partial_t (\Util3-\Utwo) - \Hil(\azer_\infty)(\Util3-\Utwo)
= \calTzer(t)[(\dot{\pi}_*^{(3)}-\dot{\pi}_*^{(2)})
\partial_\pi \tilde{W}(\pizer) + V(\pizer,\pione)\Zthree] \nn \\
&+\calTzer(t)[ N_*(\Zone,\pione)- N_*(\Zzer,\pizer)
 + \dot{\pi}_*^{(3)} (\partial_\pi \tilde{W}(\pione) - \partial_\pi \tilde{W}(\pizer))]+
V(\pizer)(\Util3-\Utwo). \label{eq:diffeqU}
\end{align}
As before, we need to split the evolution into the three pieces
\[ \Util3-\Utwo = P_s(\azer_\infty)(\Util3-\Utwo) + \Proot(\azer_\infty)(\Util3-\Utwo) + \Pim(\azer_\infty)(\Util3-\Utwo).
\]
In view of Lemma~\ref{lem:xiZ} and Lemma~\ref{lem:path_var},
\begin{align*}
|\la \Util3-\Utwo,\xi_j(\pizer)\ra| &= | \la\Zthree-\Ztwo,\xitil_j(\pizer)\ra| =
|\la \Zthree, \xitil_j(\pizer)-\xitil_j(\pione) \ra | \\
&\le \delta^2\|\Zthree\|_{2+\infty} \la t\ra\| \, \|\pizer-\pione\| \les \delta^3\la t\ra^{-\half} \|\pizer-\pione\|,
\end{align*}
provided $j\ne3,4,5$. On the other hand, by Lemma~\ref{lem:rho_infty2},
\begin{align*}
&|\la \Util3-\Utwo,\xi_{2+\ell}(\pizer)\ra| \\
&\le | \la\Zthree-\Ztwo,\xitil_{2+\ell}(\pizer)\ra| +
|y_\ell(\pizer)-(y_\infty)_\ell(\pizer)| \, |\la \Util3-\Utwo,\xi_{1}(\pizer)\ra| \\
&\le |\la \Zthree, \xitil_{2+\ell}(\pizer)-\xitil_{2+\ell}(\pione) \ra | +
|y_\ell(\pizer)-(y_\infty)_\ell(\pizer)| \, |\la \Util3-\Utwo,\xi_{1}(\pizer)\ra|\\
&\les \delta^3 \la t\ra^{-\half} \|\pizer-\pione\| + \delta^2 \|\Util3(t)-\Utwo(t)\|_{4+\infty}
\la t\ra^{-1}.
\end{align*}
The conclusion is that
\beeq
\label{eq:diff_root}
 \|\Proot(\azer_\infty)(\Util3-\Utwo)(t)\|_{4+\infty} \les \delta \la t\ra^{-\half}\|\pizer-\pione\|+\delta^2 \|\Util3(t)-\Utwo(t)\|_{4+\infty}
\la t\ra^{-1}.
\eneq
Next, we turn to the dispersive piece
$P_s(\azer_\infty)(\Util3-\Utwo)$. This requires estimating each of the expressions on
the right-hand side of~\eqref{eq:diffeqU} in the appropriate norms.
In will be convenient to use the notations
\[ \mu_Z(t):= t^{\half}\chi_{[t\ge1]}+t^{\frac34}\chi_{[0<t<1]},\; \mu_\pi(t):= t^{2}\chi_{[t\ge1]}+t^{\frac34}\chi_{[0<t<1]}.\]
Then,
\begin{align}
& \|\calTzer(t)[(\dot{\pi}_*^{(3)}-\dot{\pi}_*^{(2)})
\partial_\pi \tilde{W}(\pizer) + V(\pizer,\pione)\Zthree] \|_{1\cap 2} \nn \\
& \les \|\pi^{(3)}-\pi^{(2)}\| \mu_\pi(t)^{-1} + \delta\la t\ra^{-\half}\|\pizer-\pione\|. \label{eq:stk1}
\end{align}
Moreover, by~\eqref{eq:V01_est},
\begin{align}
&\|V(\pizer) (\Util3(t)-\Utwo(t)) \|_{1\cap2} \les
\|V(\pizer)\|_{1\cap4} \|\Util3(t)-\Utwo(t)\|_{4+\infty} \nn \\
& \les\delta^2 \|\Util3(t)-\Utwo(t)\|_{4+\infty}. \label{eq:stk2}
\end{align}
Another easy term is
\beeq
 \|\calTzer(t) \dot{\pi}_*^{(3)}
(\partial_\pi \tilde{W}(\pione)-\partial_\pi \tilde{W}(\pizer))\|_{1\cap\infty}
 \les \delta^2\la t\ra^{-2}\|\pizer-\pione\|. \label{eq:stk3}
\eneq
Next, we turn to the nonlinear terms $N_*(\Zone,\pione)- N_*(\Zzer,\pizer)$. Recall that
\beeq
\nn
N_*(\Zone,\pione)=\binom{-2|R^{(0)}|^2W(\pizer)(t)-\bar{W}(\pizer)(t)(R^{(0)})^2-|R^{(0)}|^2R^{(0)}}{2|R^{(0)}|^2\bar{W}(\pizer)(t)+W(\pizer)(t)(\bar{R}^{(0)})^2+|R^{(0)}|^2\bar{R}^{(0)}}
\eneq
The right-hand side here naturally divides into three columns, which
we formally write as
\[ |Z|^2\, W, \; Z^2\, \bar{W},\;  |Z|^2 Z,\]
respectively. Let us start with the third column (we suppress $t$ for the most part):
\begin{align*}
& \| |\Zzer|^2\Zzer -|\Zone|^2\Zone\|_{1\cap\frac43} \\
 &\les \| \Zzer-\Zone\|_{4+\infty}
   (\|\,|\Zzer|^2 + |\Zone|^2\|_{\frac43\cap 1} + \|\,|\Zzer|^2 + |\Zone|^2\|_{2\cap \frac43}) \\
 &\les \delta^2 \| \Zzer-\Zone\|_{4+\infty}.
\end{align*}
This estimate is the reason we do not work on $L^2+L^\infty$. Indeed, in the latter
case we would be faced with $\|U(t)\|_\infty^2$, which we can only bound by $t^{-\frac32}$
for small~$t$, see~\eqref{eq:boot4}. This bound is nonintegrable at~$t=0$.
The first column satisfies
\begin{align*}
& \| |\Zzer|^2 W(\pizer) -|\Zone|^2 W(\pione)\|_{1\cap\frac43} \\
&\les \| \Zzer-\Zone\|_{4+\infty}
   \Big(\|\,|\Zzer W(\pizer)| + |\Zone W(\pione)|\,\|_{\frac43\cap 1} +
\|\, |\Zzer W(\pizer)| + |\Zone W(\pione)|\,\|_{2\cap \frac43}\Big) \\
& \qquad + \|W(\pizer)-W(\pione)\|_{1\cap\infty}\|\Zone\|^2_{2+\infty} \\
&\les \delta \| \Zzer-\Zone\|_{4+\infty} + \delta^2\la t\ra^{-2}\|\pizer-\pione\|.
\end{align*}
An analogous bound holds for the second column.
Collecting these bounds yields
\beeq
\label{eq:Ndiff}
\|\calTzer[N_*(\Zone,\pione)- N_*(\Zzer,\pizer)]\|_{1\cap\frac43}
 \les \delta \| \Zzer-\Zone\|_{4+\infty} + \delta^2\la t\ra^{-2}\|\pizer-\pione\|.
\eneq
Combining \eqref{eq:stk1}, \eqref{eq:stk2}, \eqref{eq:stk3}, and~\eqref{eq:Ndiff} leads to
\begin{align}
 \|\text{right-hand side of \eqref{eq:diffeqU}} \|_{1\cap\frac43}
& \les \delta \la t\ra^{-\half} \|\pione-\pizer\| + \delta \|\Zone(t)-\Zzer(t)\|_{4+\infty} \nn \\
& \qquad + \delta^2 \|\Util3(t)-\Utwo(t)\|_{4+\infty} + \|\pi^{(3)}-\pi^{(2)}\|\mu_\pi(t)^{-1}
\label{eq:rhsF}
\end{align}
Denote the right-hand side of \eqref{eq:diffeqU} by $F$.
Estimating the Duhamel version of \eqref{eq:diffeqU} therefore leads to the conclusion that
\begin{align}
& \|P_s(\azer_\infty) (\Util3-\Utwo)(t)\|_{4+\infty} \label{eq:Psdiff1}\\
&\le \|e^{-it\Hil(\azer_\infty)}P_s(\azer_\infty)(\Util3(0)-\Utwo(0))\|_{4+\infty} \nn \\
& \quad + \int_0^{(t-1)_+} \|e^{-i(t-s)\Hil(\azer_\infty)}P_s(\azer_\infty) F(s)\|_{\infty}\, ds
+ \int_{(t-1)_+}^t \|e^{-i(t-s)\Hil(\azer_\infty)}P_s(\azer_\infty) F(s)\|_{4}\, ds \nn \\
& \les t^{-\frac32} \|P_s(\azer_\infty)[\Util3(0)-\Utwo(0)]\|_1 \chi_{[t\ge1]} +
t^{-\frac34} \|P_s(\azer_\infty)[\Util3(0)-\Utwo(0)]\|_{\frac43} \chi_{[0<t<1]}\nn \\
& \quad + \int_0^t \big(\la t-s\ra^{-\frac32}+(t-s)^{-\frac34}\chi_{[(t-1)_+<s<t]}\big)
\Big(\delta \la s\ra^{-\frac12} \|\pione-\pizer\| + \delta^2 \|(\Zone-\Zzer)(s)\|_{4+\infty} \nn \\
&\qquad + \delta^2  \|\Util3(s)-\Utwo(s)\|_{4+\infty} + \|\pi^{(3)}-\pi^{(2)}\|\mu_\pi(s)^{-1}\Big)\, ds\nn
\end{align}
As far as the initial conditions are concerned, we infer from~\eqref{eq:U2newinit}, \eqref{eq:U3newinit},
as well as~\eqref{eq:inf_der}
that
\begin{align}
& \|P_s(\azer_\infty)[\Uthree(0)-\Utwo(0)]\|_{1\cap \frac43} \nn \\
& \les   |h^{(1)}| \|[P_s(\azer_\infty)-P_s(\aone_\infty)]f^+(\azer_\infty)\|_2 + \sum_{k=1}^8
|a_j^{(1)}| \|[P_s(\azer_\infty)-P_s(\aone_\infty)]\eta_k(\azer_\infty)\|_2 \nn \\
&   \label{eq:init_diff}
\les \delta^2 \|\pizer-\pione\|.
\end{align}
Further simplification of~\eqref{eq:Psdiff1} therefore leads to
\begin{align}
& \|P_s(\azer_\infty) (\Util3-\Utwo)(t)\|_{4+\infty}
\les \delta\mu_Z(t)^{-1}\, \|\pizer-\pione\|
+ \delta^2 t^{-\frac12} \sup_{s\ge0} \mu_Z(s)\|\Util3(s)-\Utwo(s)\|_{4+\infty}  \nn \\
& + \delta\, t^{-\frac12} \big(\|\pizer-\pione\|+\sup_{s\ge0} \mu_Z(s) \|\Zzer(s)-\Zone(s)\|_{4+\infty}\big) +  t^{-\half}  \|\pi^{(3)}-\pi^{(2)}\|,  \label{eq:Psdiff2}
\end{align}
where we used the elementary estimate $\int_0^t \la t-s\ra^{-\frac32} \la s\ra^{-\half}\, ds\les \la t\ra^{-\half}$.
It remains to bound
\[ \|\Pim(\azer_\infty)[\Util3(t)-\Utwo(t)]\|_{4+\infty}. \]
To this end write
\[ \Pim(\azer_\infty)[\Util3(t)-\Utwo(t)] = b^+(t) f^+(\azer_\infty) + b^-(t) f^-(\azer_\infty), \]
with coefficients that are governed by the hyperbolic ODE
\beeq
\label{eq:bode}
 \frac{d}{dt}\binom{b^+(t)}{b^-(t)} -
\bm \sigma(\azer_\infty) & 0 \\ 0 & -\sigma(\azer_\infty) \endm \binom{b^+(t)}{b^-(t)} =
\binom{g^+}{g^-}.
\eneq
Here
\[ \Pim(\azer_\infty) F(t) = g^+(t) f^+(\azer_\infty) + g^-(t) f^-(\azer_\infty) \]
where $F$ stands for the right-hand side of~\eqref{eq:diffeqU}.
Clearly, $g^{\pm}(t)$ satisfy the bound from~\eqref{eq:rhsF}.
We need to find $b^\pm(0)$.
To this end compute
\begin{align*}
& \Pim(\azer_\infty)[\Uthree(0)-\Utwo(0)] =  b^+(0) f^+(\azer_\infty) + b^-(0) f^-(\azer_\infty) \\
&= \Pim(\azer_\infty)\calG_\infty(\pizer)(0)\big[(h^{(1)}-h^{(0)}) f^+(\azer_\infty)
- h^{(1)} [f^+(\azer_\infty)  -  f^+(\aone_\infty)] \nn \\
& \quad + \sum_{j=1}^8 [a_j^{(1)} \eta_j(\aone_\infty) - a_j^{(1)} \eta_j(\aone_\infty)] \big].
\end{align*}
Thus,
\begin{align*}
|b^+(0)-(h^{(1)}-h^{(0)})| &\les \delta^2|h^{(1)}-h^{(0)}|+\delta^2 \|\pione-\pizer\|
+\delta^2 \sum_{j=1}^8 |a_j^{(1)}-a_j^{(0)}| \\
&\les \delta^2|h^{(1)}-h^{(0)}|+\delta^2 \|\pione-\pizer\|,
\end{align*}
where we used \eqref{eq:ajdt} in the final inequality. Moreover,
\[
|b^-(0)| \les |h^{(1)}-h^{(0)}|+\delta^2 \|\pione-\pizer\|.
\]
Since $b^{\pm}(t)$ is a bounded solution of the ODE~\eqref{eq:bode}, it follows from Lemma~\ref{lem:ODE_stable} and~\eqref{eq:rhsF} that
\begin{align*}
 |b^+(0)| &\les \int_0^\infty e^{-\sigma(\azer_\infty) t} (|g^+(t)|+|g^-(t)|)\, dt \\
&\les \int_0^\infty e^{-\sigma(\azer_\infty) t} \Big[ \delta \la t\ra^{-\half} \|\pione-\pizer\| + \delta \|\Zone(t)-\Zzer(t)\|_{4+\infty} \nn \\
& \qquad + \delta^2 \|\Util3(t)-\Utwo(t)\|_{4+\infty} + \|\pi^{(3)}-\pi^{(2)}\|\mu_\pi(t)^{-1} \Big]\, dt\\
&\les
\delta \|(\pizer-\pione,\Zzer-\Zone)\|  +
\delta^2 \|(\pi^{(3)}-\pi^{(2)}, \Zthree-\Ztwo)\| + \|\pi^{(3)}-\pi^{(2)}\|
\end{align*}
and thus also
\begin{align}
& |h^{(1)}-h^{(0)}| + |b^{-}(0)| \nn \\
& \les \delta \|(\pizer-\pione,\Zzer-\Zone)\|  +
\delta^2 \|(\pi^{(3)}-\pi^{(2)}, \Zthree-\Ztwo)\| + \|\pi^{(3)}-\pi^{(2)}\|.
\label{eq:h1h0diff}
\end{align}
Furthermore, in view of \eqref{eq:Duh_stab},
\begin{align}
\nn
& \|\Pim(\azer_\infty)[\Util3(t)-\Utwo(t)]\|_{4+\infty} \les |b^+(t)|+|b^-(t)| \\
& \les \int_t^\infty e^{-\sigma(\azer_\infty)(s-t)}  \Big[ \delta \la s\ra^{-\half} \|\pione-\pizer\| + \delta \|\Zone(s)-\Zzer(s)\|_{4+\infty} \nn \\
& \qquad + \delta^2 \|\Util3(s)-\Utwo(s)\|_{4+\infty} + \|\pi^{(3)}-\pi^{(2)}\|\mu_\pi(s)^{-1} \Big]\, ds\nn \\
& + e^{-\sigma(\azer_\infty)t} (\delta \|(\pizer-\pione,\Zzer-\Zone)\|  +
\delta^2 \|(\pi^{(3)}-\pi^{(2)}, \Zthree-\Ztwo)\| + \|\pi^{(3)}-\pi^{(2)}\|) \nn \\
& + \int_0^t e^{-\sigma(\azer_\infty)(t-s)} \Big[ \delta \la s\ra^{-\half} \|\pione-\pizer\| + \delta \|\Zone(s)-\Zzer(s)\|_{4+\infty} \nn \\
& \qquad + \delta^2 \|\Util3(s)-\Utwo(s)\|_{4+\infty} + \|\pi^{(3)}-\pi^{(2)}\|\mu_\pi(s)^{-1} \Big]\, ds
\nn \\
& \les \la t\ra^{-\half} \big[\delta \|(\pizer-\pione,\Zzer-\Zone)\|  +
\delta^2 \|(\pi^{(3)}-\pi^{(2)}, \Zthree-\Ztwo)\| + \|\pi^{(3)}-\pi^{(2)}\|].
\label{eq:Udiff_hyp}
\end{align}
Now set
\[ \eps_0:= \|(\pizer-\pione,\Zzer-\Zone)\|, \quad \eps_2 := \|(\pi^{(3)}-\pi^{(2)}, \Zthree-\Ztwo)\|. \]
Combining \eqref{eq:diff_root}, \eqref{eq:Psdiff2}, and \eqref{eq:Udiff_hyp} leads
to the bound
\beeq
\label{eq:Zpart}
\sup_{t\ge0} \mu_Z(t) \|Z^{(3)}(t)-Z^{(2)}(t)\|_{4+\infty} \les \delta\eps_0+\delta^2\eps_2 +
\|\pi^{(3)}-\pi^{(2)}\|
\eneq
and thus also, in view  \eqref{eq:h1h0diff},
\beeq
|h^{(1)}-h^{(0)}| \les \delta\eps_0+\delta^2\eps_2 + \|\pi^{(3)}-\pi^{(2)}\|.
\label{eq:hpart}
\eneq
We now turn to estimating the difference of the paths $\pi^{(3)}, \pi^{(2)}$.
Indeed, inserting some of the bounds we derived into~\eqref{eq:diffeq2} yields
\begin{align*}
| \dot{\pi}^{(3)}(t)-\dot{\pi}^{(2)}(t)| &\les \delta^2 \la t\ra^{-3} \|(\Zthree-\Ztwo)(t)\|_{4+\infty}
+ \delta \la t\ra^{-\frac32} (\la t\ra^{-\half} \|\pizer-\pione\| +
\delta^2 \la t\ra^{-2}\| \pione-\pizer\|) \\
&\quad + \delta\la t\ra^{-\frac32} \|\Zzer(t)-\Zone(t)\|_{4+\infty} +
\delta^2\la t\ra^{-2}\|\pione-\pizer\| +\delta^2 \la t\ra^{-2} \|\pione-\pizer\| \\
&\les \delta^2 \mu_\pi(t)^{-1} \eps_2 + \delta \mu_\pi(t)^{-1} \eps_0,
\end{align*}
which implies that
\[ \|\pi^{(3)}-\pi^{(2)}\| = \sup_{t\ge0}\mu_\pi(t) | \dot{\pi}^{(3)}(t)-\dot{\pi}^{(2)}(t)|
\les \delta^2\eps_2 + \delta\eps_0.\]
Combining this bound with \eqref{eq:Zpart} yields
that $\eps_2\les \delta\eps_0$, which is the same as
\[ \|\Psi(R_0,\pi^{(0)},\Zzer)-\Psi(R_1,\pi^{(1)},\Zone)\|
\les \delta \|(\pizer,\Zzer)-(\pione,\Zone)\|,\]
where we have included the initial condition $R_0$ in the notation.
We have shown that $\Psi$ is a contraction in~$X_\delta$.
Denote the unique fixed-point in $X_\delta$ by $(\pi(R_0),Z(R_0))$. We claim
that this fixed-point is Lipschitz in $R_0$ in the following sense:
\beeq
\label{eq:lip_fp}
\| (\pi(R_0),Z(R_0)) - (\pi(R_1),Z(R_1)) \| \les \trip R_0-R_1 \trip.
\eneq
In view of Lemma~\ref{lem:param_contr} it suffices to show that
\beeq
\label{eq:Psi_lip}
\| \Psi(R_0,\pizer,\Zzer)-\Psi(R_1,\pizer,\Zzer)\| \les \trip R_0-R_1 \trip.
\eneq
To prove this, set
\[ (\pithree,\Zthree) = \Psi(R_1,\pizer,\Zzer), \; (\pitwo,\Ztwo) = \Psi(R_0,\pizer,\Zzer).\]
The difference of these functions is controlled by
the equations \eqref{eq:diffeq1}, \eqref{eq:diffeq2} with $\pizer=\pione, \Zzer=\Zone$.
Hence,
\begin{align*}
i\partial_t (\Zthree-\Ztwo) - \Hil(\pizer(t))(\Zthree-\Ztwo)
& = (\dot{\pi}_*^{(3)}-\dot{\pi}_*^{(2)})
\partial_\pi \tilde{W}(\pizer) \\
\la (\dot{\pi}_*^{(3)}-\dot{\pi}_*^{(2)}) \partial_\pi \tilde{W}(\pizer),\xitil_j(\pizer) \ra
& = \la \Zthree-\Ztwo, \dot{\pi}_*^{(0)} \calS_j(\pizer) \ra
\end{align*}
with initial conditions
\[
(\Zthree-\Ztwo)(0) =  \binom{R_1}{\bar{R}_1}-\binom{R_0}{\bar{R}_0}+ (h^{(1)}-h^{(0)})f^+(\azer_\infty)
+\sum_{j=1}^8 (a^{(1)}_j - a^{(0)}_j) \eta_j(\azer_\infty),
\]
cf.~\eqref{eq:diff_init}. The orthogonality conditions
\[ \la \Zthree(t),\xitil_j(\pizer)(t)\ra = \la \Ztwo(t),\xitil_j(\pizer)(t)\ra =  0\]
hold for all $t\ge0$ by construction. Setting
\[ \Util3:=\calTzer \Zthree,\; \Utwo:=\calTzer \Ztwo\]
as before, we obtain the transformed equations
\begin{align*}
& i\partial_t (\Util3-\Utwo) - \Hil(\azer_\infty)(\Util3-\Utwo) \\
& = \calTzer[(\dot{\pi}_*^{(3)}-\dot{\pi}_*^{(2)})
\partial_\pi \tilde{W}(\pizer)] + V(\pizer)(\Util3-\Utwo) \\
(\Util3-\Utwo)(0) &=  \calG_\infty(0)[\binom{R_1}{\bar{R}_1}-\binom{R_0}{\bar{R}_0}+ (h^{(1)}-h^{(0)})f^+(\azer_\infty)
+\sum_{j=1}^8 (a^{(1)}_j - a^{(0)}_j) \eta_j(\azer_\infty)].
\end{align*}
The orthogonality conditions are $\la \Util3-\Utwo,\xi_j(\pizer)\ra=0$.
The estimate~\eqref{eq:Psi_lip} now follows by
using the same techniques that we have employed repeatedly in order to control the
solution~$\Util3-\Utwo$. We skip the details.

Combining~\eqref{eq:lip_fp} with~\eqref{eq:hpart} leads to the statement that
\[ |h(R_0)-h(R_1)|\les \delta \trip R_0-R_1\trip,\]
as claimed.
\end{proof}

\begin{proof}[Proof of Theorem~\ref{thm:main}] Given $R_0$, the previous lemma establishes
the existence of $h=h(R_0)\in\R$ as well as $(\pi,Z)=(\pi(R_0),Z(R_0))\in X_\delta$ where $\delta=C_0\trip R_0\trip$, which solve
\begin{align*}
& i\partial_t Z(t) + \bm \Laplace + 2|W(\pi)|^2 & W^2(\pi) \\ -\bar{W}^2(\pi) & -\Laplace -2|W(\pi)|^2 \endm Z(t)  \\
 & = \dot{v} \binom{-xe^{i\theta(\pi)(t)} \phi(\cdot-y(\pi)(t),\alpha(t))}
{xe^{-i\theta(\pi)(t)} \phi(\cdot-y(\pi)(t),\alpha(t))} +
\dot{\gamma} \binom{-e^{i\theta(\pi)(t)} \phi(\cdot-y(\pi)(t),\alpha(t))}{e^{-i\theta(\pi)(t)}
\phi(\cdot-y(\pi)(t),\alpha(t))}\nn \\
&\quad + i\dot{\alpha}
\binom{e^{i\theta(\pi)(t)}\partial_\alpha \phi(\cdot-y(\pi)(t),\alpha(t))}{e^{-i\theta(\pi)(t)}\partial_\alpha \phi(\cdot-y(\pi)(t),\alpha(t))}
 + i\dot{D}
\binom{-e^{i\theta(\pi)(t)}\nabla \phi(\cdot-y(\pi)(t),\alpha(t))}{-e^{-i\theta(\pi)(t)}\nabla
\phi(\cdot-y(\pi)(t),\alpha(t))} \nn \\
&\quad + \binom{-2|R|^2W(\pi)(t)-\bar{W}(\pi)(t)R^2-|R|^2R}{2|R|^2\bar{W}(\pi)(t)+W(\pi)(t)\bar{R}^2+
|R|^2\bar{R}}
\end{align*}
with initial conditions
\begin{align*}
Z(0) &= \binom{R_0}{\bar{R}_0} + h(R_0)f^{+}(\alpha_\infty) + \sum_{j=1}^8 a_j(R_0)\eta_j(\alpha_\infty)\\
\pi(0) &= (\alpha_0,0,0,0).
\end{align*}
Here $a_j=a_j(h(R_0))\in\R$.
Define
\[
\Phi(R_0) :=  h(R_0)f^{+}(\alpha_\infty) + \sum_{j=1}^8 a_j(R_0)\eta_j(\alpha_\infty).
\]
Since
\[ |h(R_0)|+\sum_{j=1}^8 |a_j(h(R_0))| \les \delta^2 \les \trip R_0\trip^2,\]
the estimate \eqref{eq:r0til} follows. Moreover, \eqref{eq:lip} follows from~\eqref{eq:hwaschno}.
Since $Z=\binom{R(t)}{\bar{R}(t)}$ is $\calJ$-invariant, it follows from Lemma~\ref{lem:Zeq} that
\[\psi(t):= W(\pi(t))+R(t)\]
is an $H^1$-solution  of the NLS~\eqref{eq:NLS}. Finally,
\[ \|R(t)\|_{W^{1,2}} \les \delta, \quad \|R(t)\|_\infty \les \delta t^{-\frac32} \]
follows from \eqref{eq:boot2} and \eqref{eq:boot3}, whereas \eqref{eq:boot1} ensures that
the path is admissible and therefore converges to $\pi(\infty)$ with
\[ \sup_{t\ge0} |\pi(t)-\pi(\infty)|\les \delta^2.\]
Finally, we turn to the scattering statement. According to Lemma~\ref{lem:UPDE},
\begin{align}
i\partial_t U - \Hil(\alpha_\infty) U &= \dot\pi \partial_\pi W(\pi) + N(U,\pi) + V(\pi) U
\label{eq:Ufin}\\
U(0) &= \calG_\infty(\pi)(0)\Big[\binom{R_0}{\bar{R}_0} + \binom{\Phi(R_0)}{\overline{\Phi(R_0)}}\Big].
\nn
\end{align}
Denoting the right-hand side \eqref{eq:Ufin} by $F(t)$, we have
\[ U(t) = e^{-it\Hil(\alpha_\infty)} U(0) -i \int_0^t e^{-i(t-s)\Hil(\alpha_\infty)} F(s)\, ds.\]
The estimates \eqref{eq:boot1}--\eqref{eq:boot3} imply that
\[ \|F(s)\|_{\frac32\cap 2}\les \la s\ra^{-\frac52},\quad\|F(s)\|_{1\cap2}\les \la s\ra^{-\frac32}, \quad \int_0^\infty \|F(s)\|_2\, ds< \infty.\]
This allows us to define
\[ U_1 := P U(0) -i\int_0^\infty e^{is\Hil(\alpha_\infty)} P F(s)\, ds \in L^2(\R^3)\]
where we have set $P:=P_s(\alpha_\infty)+\Proot(\alpha_\infty)=1-\Pim(\alpha_\infty)$.
We are using here that
\[ \|e^{is\Hil(\alpha_\infty)} P F(s)\|_2 \les \la s\ra^{-\frac32}, \]
which follows from the fact that growth of $e^{is\Hil(\alpha_\infty)}$ on the root-space
can be at most~$s$.
Clearly, $U_1$ was defined so that
\[ P U(t)-e^{-it\Hil(\alpha_\infty)} U_1 = i \int_t^\infty e^{-i(t-s)\Hil(\alpha_\infty)} PF(s)\, ds\]
which implies that
\[ \|P U(t)-e^{-it\Hil(\alpha_\infty)} U_1\|_2 \les \la t\ra^{-\half} \to 0\]
as $t\to\infty$. As far as the hyperbolic part is concerned, we define
\[ U_2 := \Pim^{-}(\alpha_\infty) U(0) - i\int_0^\infty e^{-s\sigma(\alpha_\infty)}
\Pim^{-}(\alpha_\infty) F(s) \, ds.\]
Because of Lemma~\ref{lem:ODE_stable},
\[
 \Pim(\alpha_\infty) U(t)-e^{-it\Hil(\alpha_\infty)} U_2 = i \int_t^\infty e^{(t-s)\sigma(\alpha_\infty)} \Pim(\alpha_\infty)F(s)\, ds.
\]
In conjunction with the $P$-part this shows that
\beeq
\label{eq:scat_diff}
U(t) - e^{-it\Hil(\alpha_\infty)} (U_1+U_2) = i\int_t^\infty e^{-i(t-s)\Hil(\alpha_\infty)}P F(s)\,ds
+ i \int_t^\infty e^{(t-s)\sigma(\alpha_\infty)} \Pim(\alpha_\infty)F(s)\, ds.
\eneq
Therefore, as $t\to\infty$,
\beeq
 U(t) = e^{-it\Hil(\alpha_\infty)} (U_1+U_2) + o_{L^2}(1).
\label{eq:scat_1}
\eneq
Another consequence of \eqref{eq:scat_diff} is the estimate
\beeq
\|e^{-it\Hil(\alpha_\infty)} (U_1+U_2) \|_{3}  \les \|U(t)\|_{3}
+ \int_t^\infty (t-s)^{-\frac12} \|F(s)\|_{\frac32}\, ds
\les \la t\ra^{-\frac12}. \label{eq:slow_dec}
\eneq
This implies that in fact $\Proot(U_1+U_2)=0$. Seeing this requires some care, as
we do not know that $U_1+U_2\in L^1(\R^3)$. However, \eqref{eq:slow_dec} implies that
\[ \|e^{-it\Hil(\alpha_\infty)} (U_1+U_2) \|_{L^4_t(L^3_x(\R^3))} < \infty.\]
On the other hand, by the Strichartz estimate~\eqref{eq:Strich1},
\[ \|e^{-it\Hil(\alpha_\infty)} P_s(U_1+U_2) \|_{L^4_t(L^3_x(\R^3))} \les \|P_s(U_1+U_2)\|_2<\infty.\]
Hence, also
\[ \|e^{-it\Hil(\alpha_\infty)} \Proot(U_1+U_2) \|_{L^4_t(L^3_x(\R^3))} <\infty.\]
However, this is only possible if in fact $\Proot(U_1+U_2)=0$, as claimed.
Therefore,
\[ U_1 := P_s U(0) -i\int_0^\infty e^{is\Hil(\alpha_\infty)} P_s F(s)\, ds\]
which in particular implies  the dispersive bound
\beeq
\|e^{-it\Hil(\alpha_\infty)} (U_1+U_2) \|_{2+\infty} \les \|U(t)\|_{2+\infty}
+ \int_t^\infty \la t-s\ra^{-\frac32} \|F(s)\|_{1\cap2}\, ds
\les \la t\ra^{-\frac32}. \label{eq:U12dec}
\eneq
It remains to show that one has scattering for the evolution of $\Hil(\alpha_\infty)$.
This is a standard Cook's method argument. Indeed, write
\[ \Hil(\alpha_\infty) = \bm -\Laplace+\alpha_\infty^2 & 0 \\ 0 & \Laplace -\alpha_\infty^2 \endm
+ \bm -2\phi_\infty^2 & -\phi_\infty^2 \\ \phi_\infty^2 & 2\phi_\infty^2 \endm =: \Hil_0(\alpha_\infty) + V,
\]
where $\phi_\infty:=\phi(\cdot,\alpha_\infty)$.
Then
\[
e^{-it\Hil(\alpha_\infty)} (U_1+U_2) = e^{-it\Hil_0(\alpha_\infty)} (U_1+U_2)
-i \int_0^t e^{-i(t-s)\Hil_0(\alpha_\infty)} V e^{-is \Hil(\alpha_\infty)} (U_1+U_2)\, ds\]
and thus
\beeq
\label{eq:scat_2}
 U_3:=\lim_{t\to\infty} e^{it\Hil_0(\alpha_\infty)} e^{-it\Hil(\alpha_\infty)} (U_1+U_2)
\eneq
exists as a strong $L^2$ limit. Indeed, this follows from
\[ \int_0^\infty \|e^{is\Hil_0(\alpha_\infty)} V e^{-is \Hil(\alpha_\infty)} (U_1+U_2)\|_2\, ds
\les \int_0^\infty \|e^{-is \Hil(\alpha_\infty)} (U_1+U_2)\|_{2+\infty}\, ds < \infty,
\]
see \eqref{eq:U12dec}. It follows from \eqref{eq:scat_1} and~\eqref{eq:scat_2} that
\[ U(t) = e^{-it\Hil_0(\alpha_\infty)} U_3 + o_{L^2}(1). \]
Finally,
\[ Z(t) = \calG_\infty(t)^{-1} M(t)^{-1}U(t) = e^{-it\Hil_0} \calG_\infty^{-1}(0) U_3 + o_{L^2}(1),\]
where $\Hil_0=\bm -\Laplace & 0\\ 0 & \Laplace\endm$.
Setting $\calG_\infty^{-1}(0) U_3 = \binom{f_0}{\bar{f}_0}$ and $Z(t)= \binom{R(t)}{\bar{R}(t)}$,
we obtain
\[ R(t) = e^{it\Laplace} f_0 +  o_{L^2}(1),\]
and the theorem is proved.
\end{proof}

\begin{proof}[Proof of Theorem~\ref{thm:main2}]
The idea is as follows: Given $\alpha_0$, consider the
NLS~\eqref{eq:NLS} with initial data $\phi(\cdot,\alpha_0)+R_0$.
Applying the usual eight-parameter family of symmetries (Galilei
giving six parameters, modulation one, and scaling also one ---
scaling here is the same as the parameter $\alpha$), we transform
this to $W(0,\cdot)+R_1$ where $W(0,x)$ is a soliton with a
general parameter vector $\pi_0$ which is close to
$(0,0,0,\alpha_0)$. Hence, we can apply Theorem~\ref{thm:main} to
conclude that these initial data will give rise to global
solutions with the desired properties as long as $W(0,x)+R_1$ lies
on the stable manifold associated with~$W(0,x)$. To prove that we
obtain eight dimensions back in this fashion requires checking
that the derivatives of $W(0,x)$ in its parameters are transverse
to the linear space $\calS$ of Theorem~\ref{thm:main}. However,
these derivatives are basically the elements of the root space
$\calN$ of $\Hil(\alpha_0)$, whereas we know that $\calS$ is
perpendicular to the root space $\calN^*$ of $\Hil(\alpha_0)^*$.
But Lemma~\ref{lem:orth} implies that no nonzero vector in $\calN$
is perpendicular to $\calN^*$, which proves that $\calN$ is
transverse to $\calS$, as desired.
\end{proof}

\section{The linear analysis: $L^2$ theory}
\label{sec:lin_L2}

The purpose of this section is to develop the $L^2$-estimates on the evolution $e^{-it\Hil(\alpha)}$
which were used in the proof of Theorem~\ref{thm:main}. Such $L^2$-estimates were already obtained
by Weinstein~\cite{Wei1}, but his results do not cover the unstable operators which result
from linearizing the supercritical equation~\eqref{eq:NLS}. We start with a brief outline of
the linear theory, both $L^2$ as well as dispersive. The latter will be dealt with in
Section~\ref{sec:lin_dis}.  Some of the following statements are well--known and we do not provide
the proofs of such results. For the most part, references are given.

\begin{itemize}
\item Start from a ground state $-\Laplace \phi+\alpha^2\phi=\phi^3$ in~$\R^3$,
$\phi>0$ and $\phi$ smooth, $\phi=\phi(\cdot,\alpha)$.

\item Consider the system ($\mu=\alpha^2$)
\beeq
\label{eq:Hil100}
\Hil=\Hil(\alpha)= \bm -\Laplace+\alpha^2-2\phi^2 & -\phi^2 \\
\phi^2 & \Laplace -\alpha^2 + 2\phi^2
\endm = \bm -\Laplace +\mu - U & -W \\ W & \Laplace -\mu + U
\endm.
\eneq
We will often write $\Hil=\Hil_0+V$, where
\[
\Hil_0=
 \bm -\Laplace +\mu  & 0 \\ 0 & \Laplace -\mu
\endm, \qquad V= \bm - U & -W \\ W &  U
\endm.
\]
Much of the linear theory here does not depend on the special
structure that results from the nonlinear origin of the linear
operator. In fact, it has been our intent to develop most of the
arguments in the generality of a matrix potential $V$ that has
polynomial decay at infinity. In some instances, we have found it
necessary to impose more stringent requirements, like some
regularity or a small amount of exponential decay. In the
forthcoming work with Erdogan~\cite{ErdSch}, we intend to present a
more general and encompassing version of the linear theory
developed here.

\item The system \eqref{eq:Hil100} can also be written as (relative to different coordinates)
\[
\Hil = \bm 0 & -iL_{-} \\ iL_{+} & 0 \endm.
\]
where $L_{+} = -\Laplace+\mu-3\phi^2$, $L_{-} = -\Laplace+\mu-\phi^2$.

\item The essential spectrum is $(-\infty,-\mu]\cup[\mu,\infty)$. The analytic
Fredholm alternative gives that the resolvent is meromorphic outside of the
essential spectrum. I.e., for every point $\sigma_0\not\in(-\infty,-\mu]\cup[\mu,\infty)$
one has $\|(\Hil-z)^{-1}\|\les |z-\sigma_0|^{-\nu}$ for some positive integer $\nu$.
Let $\nu_0=\nu_0(\sigma_0)$ be the smallest such $\nu$. Then $\nu_0$
 is the Jordan index, i.e., $\ker(\Hil-\sigma_0)^{\nu_0}=\ker(\Hil-\sigma_0)^{\nu_0+1}$
and $\nu$ is the smallest integer with this property.

\item The number of isolated points of the spectrum is finite.

\item Consider the Riesz projection around each point of the discrete spectrum.
I.e.,
\[ P_{\sigma_0} = -\frac{1}{2\pi i}\oint_\gamma (\Hil-z)^{-1}\, dz \]
where $\gamma$ is a small loop around $\sigma_0$.  Then
$\Ran(P_{\sigma_0}) = \ker(\Hil-\sigma_0)^{\nu_0}$, where $\nu_0=\nu_0(\sigma_0)$ is as above,
$L^2(\R^3)\times L^2(\R^3)$ is the
direct (but not orthogonal) sum of $\Ran(P_{\sigma_0})$ and $\Ran(I-P_{\sigma_0})$,
which are closed spaces. Moreover, on each of these spaces the spectrum of $\Hil$ is the natural
one, i.e., $\sigma_0$ in the former case, and the spectrum minus $\sigma_0$ in the latter
case. These are general statements that are valid for closed operators and isolated points
of their spectrum. See Grillakis~\cite{Grill}, as
well as Hislop, Sigal~\cite{HS} for the Riesz projections.

\item One has $L_{-}\phi=0$ and thus $L_{-}\ge0$. In the general case of a matrix potential
$V$ we will impose this condition and refer to it as the ``positivity
assumption''. It has several important consequences, see eg.~\cite{BP1} or~\cite{RSS2}:
the spectrum of $\Hil$ is a subset of $\R\cup i\R$, and all discrete
spectrum of $\Hil$ other than zero consists of eigenvalues with Jordan index one.

\item An adaptation of Agmon's argument~\cite{Agm2} shows that all
eigenvalues as well as elements of the generalized eigenspaces
must decay exponentially (and since our $U,W$ are smooth since $\phi$ is
smooth they are also smooth by elliptic regularity - but this is merely
a convenience). In \cite{RSS2} Agmon's argument was generalized to systems. It
was shown there that eigenfunctions with real energy $-\mu<E<\mu$ decay
exponentially. Replacing $E$ with $\Re E$ in that argument allows one
to conclude exponential decay of eigenfunctions with energy $\Re E \in(-\mu,\mu)$.
For this we need to assume that $W$ has (a small amount of) exponential decay.

\item The root space of $\Hil$ at zero was described by Weinstein~\cite{Wei1}.
It is
\[ \Big\{ \binom{\phi}{-\phi}, \binom{\partial_j\phi}{\partial_j\phi},
\binom{\partial_\alpha\phi}{\partial_\alpha\phi},
\binom{x_j\phi}{-x_j\phi}\:\Big|\: 1\le j\le 3\Big\}. \]
Moreover, in the $L^2$-supercritical case like we are considering
here there is spectrum off the real axis, at points $\pm i\sigma$, $\sigma>0$
(in the case at hand there should be exactly one pair of these points, which
are eigenvalues of multiplicity one). See Sections~\ref{sec:linroot} and~\ref{sec:imspec} above.

\item There are no resonances in the essential spectrum. This means the following:
Suppose $|\lambda|>\mu$ and
\[ f+(\Hil_0-(\lambda\pm i0))^{-1}Vf =0 \]
for some $f\in L^{2,-\sigma}\times L^{2,-\sigma}$ with $\sigma>\half$. Then
$f\in L^2\times L^2$, which then implies that $f$ is an eigenvector with eigenvalue $\lambda$.
This is an analogue of Agmon's result about absence of positive energy resonances for
the scalar case, see~\cite{Agm1}.

\item {\bf Let us assume that there is no embedded point spectrum.} Then one has a limiting
absorption principle
\beeq
\label{eq:lam0}
\sup_{|\lambda|\ge\lambda_0,\,0<\eps} |\lambda|^{\half}\,\|(\Hil-(\lambda\pm i\eps))^{-1}\| < \infty
\eneq
for any $\lambda_0>\mu$, in the operator norms $L^{2,\sigma}\times L^{2,\sigma}\to
L^{2,-\sigma}\times L^{2,-\sigma}$ with $\sigma>\half$.

\item We will {\bf assume} that one can take $\lambda_0=\mu$ in~\eqref{eq:lam0}.
This amounts to {\bf assuming that the edges of the continuous spectrum are neither eigenvalues nor resonances}.

\item Under these two assumptions there is a representation
\beeq
\label{eq:repform}
e^{it\Hil}  = \frac{1}{2\pi i}\int_{|\lambda|\ge \mu} e^{it\lambda}\,
[(\Hil-(\lambda+i0))^{-1}-(\Hil-(\lambda-i0))^{-1}]\,d\lambda + \sum_{j} e^{it\Hil}P_{\zeta_j}
\eneq
where the sum runs over the (finite) discrete spectrum $\{\zeta_j\}$ of $\Hil$ and $P_{\zeta_j}$ is the
Riesz projection onto the corresponding generalized eigenspace
\[ \bigcup_{k=1}^\infty \ker\big((\Hil-\zeta_j)^k\big) \]
of the eigenvalue $\zeta_j$.
This representation formula (and the convergence of the integral)
is to be understood in an appropriate weak sense,
see Lemma~\ref{lem:rep} below.
The proof of~\eqref{eq:repform} starts with the definition of $e^{it\Hil}$ for $t\ge0$ by means of
the Hille-Yoshida theorem as an exponentially bounded (semi) group.
Performing an inverse Laplace transform, it can also be written as a contour integral of the
resolvent. Next, we deform this contour and the projections appear as residues.

\item Let $P_s$ and $P_c$ denote the following Riesz projections:  ${\rm Id}-P_s$ is given by
integrating the resolvent $(\Hil-z)^{-1}$
 around curves surrounding the  imaginary spectrum including zero,
whereas ${\rm Id}-P_c$ is given by integration around curves surrounding the
entire discrete spectrum. Hence $P_s\ne P_c$ only if there is discrete spectrum
on the real axis other than zero. Note that we are not excluding this possibility
in the linear theory. The indices ``s'' and ``c'' refer to {\em stable}
and {\em continuous}, although this terminology is used here in a purely formal way.
Using the representation~\eqref{eq:repform}
 one can show that there is {\em stability} in the following
sense
\[ \sup_{t\in\R} \|e^{it\Hil}P_s\|_{2\to 2} \le C,\]
under the assumption that there is no embedded point spectrum and that the thresholds
$\lambda=\pm \mu$ are neither eigenvalues nor resonances.
Weinstein~\cite{Wei1}, \cite{Wei2} derived this bound in the $L^2$ subcritical case. There the only
instability is due to the root space at zero, and the range of $P_s$ is precisely
the orthogonal complement of the root space at zero of $\Hil^*$. His analysis is
variational, and is based on the definite sign of $\partial_\alpha \|\phi(\cdot,\alpha)\|_2^2$.
In particular, his analysis exploits the specific structure of the system as being derived from
a ground state and does not require our assumptions on eigenvalues and resonances.
Note that what is being proposed above is different in several respects. It applies to very general
systems and does not use any variational arguments (naturally, being in the $L^2$ supercritical case,
we have the wrong sign of $\partial_\alpha \|\phi(\cdot,\alpha)\|_2^2$ and $P_s$ is more complicated).

\item Moreover, under the assumption that there is no embedded point spectrum and that the thresholds
$\lambda=\pm \mu$ are neither eigenvalues nor resonances, there is a dispersive bound
\[ \|e^{it\Hil}P_c\|_{1\to \infty} \le C|t|^{-\frac32}. \]
The proof of this relies on the representation~\eqref{eq:repform} and follows the scheme
of proof in the scalar case $H=-\Laplace+V$ which was developed in~\cite{GolSch}.

\item For nonlinear applications it is important to establish the following Strichartz estimates:
\begin{align*}
 \|e^{-it\Hil} P_c f\|_{L^r_t(L^p_x)} &\le C \|f\|_{L^2}  \\
 \Big\|\int_0^t e^{-i(t-s)\Hil} P_c F(s)\,ds \Big\|_{L^r_t(L^p_x)} &\le C \|F\|_{L^{a'}_t(L^{b'}_x)},
\end{align*}
provided $(r,p), (a,b)$ are admissible, i.e., $2<r\le\infty$ and $\frac{2}{r}+\frac{3}{p}=\frac32$ and
the same for $(a,b)$.

\noindent In the scalar case, these follow from the dispersive bound via a $TT^*$ argument,
where $(Tf)(t,x)=(e^{-it\Hil}f)(x)$. But since $\Hil$ is not self-adjoint, this approach does
not work here. We therefore use a different method, which originates in~\cite{RodSch}. This method
is perturbative in nature, and starts from the Strichartz estimates for $\Hil_0$. Since the
perturbing potential $V$ is not small, though, the perturbation theory requires again knowing
that there are no embedded eigenvalues in the essential spectrum of~$\Hil$, and that the resolvent
remains bounded in the aforementioned weighted $L^{2,\sigma}$ spaces at the thresholds~$\pm \mu$.
More technically speaking, we make use of Kato's notion of $\Hil_0$ and $\Hil$-smoothing
operators.
\end{itemize}

We now turn to a more detailed discussion.
We strive to keep it mostly on a general level, i.e., consider
general systems without any mention of ground states.
 Let $H_0=-\Laplace+\mu$ on~$L^2(\R^3)$ where $\mu>0$. Set
\beeq
\label{eq:Hildef}
 \Hil_0 = \bm 0 &i H_0 \\ -i H_0 & 0 \endm,\quad V= \bm 0 & iV_1 \\ -iV_2 & 0
\endm,\quad  \Hil = \Hil_0+V=i\bm 0 & H_0+V_1 \\ -H_0-V_2 & 0 \endm.
\eneq
By means of the matrix $J=\bm 0 & i \\ -i & 0 \endm$ one can also
write
\[ \Hil_0 = \bm H_0 & 0 \\ 0 & H_0 \endm J,\quad \Hil = \bm H_0+V_1 & 0 \\ 0 &
H_0+V_2 \endm J.\]
Clearly, $\Hil$ is a closed operator on the domain $\Dom(\Hil)=W^{2,2}\times W^{2,2}$.
Since $\Hil_0^*=\Hil_0$ it follows that $\spec(\Hil_0)\subset \R$.
One checks that for $\Re z\ne0$
\bear
(\Hil_0-z)^{-1} &=& -(\Hil_0+z)\bm (H_0^2-z^2)^{-1} & 0 \\ 0 & (H_0^2-z^2)^{-1}
\endm  \nn \\
&=& -\bm (H_0^2-z^2)^{-1} & 0 \\ 0 & (H_0^2-z^2)^{-1} \endm (\Hil_0+z)
\label{eq:Bspec} \\
(\Hil-z)^{-1} &=& (\Hil_0-z)^{-1}\!\! -\! (\Hil_0-z)^{-1} W_1 \Bigl[1+W_2 J (\Hil_0-z)^{-1}
\!W_1\Bigr]^{-1}\! W_2 J (\Hil_0-z)^{-1}
\label{eq:grill}
\eear
where \eqref{eq:grill} also requires the expression in brackets to be
invertible, and with
\[
 W_1 = \bm |V_1|^{\half} & 0 \\ 0 & |V_2|^{\half} \endm,\quad W_2 =
\bm |V_1|^{\half}\sign(V_1) & 0 \\ 0 & |V_2|^{\half}\sign(V_2) \endm.
\]
It follows from \eqref{eq:Bspec} that $\spec(\Hil_0) = (-\infty,-\mu]\cup [\mu, \infty)\subset \R$.
If $V_1(x)\to0$ and $V_2(x)\to0$ as $x\to\infty$,
then it follows from Weyl's theorem, 
and the representation~\eqref{eq:grill} for
the resolvent of~$\Hil$, that
$\spec_{ess}(\Hil)=\spec_{ess}(\Hil_0)=(-\infty,-\mu]\cup [\mu,
\infty)\subset \R$. 
Moreover, \eqref{eq:grill} implies via the analytic Fredholm
alternative that $(\Hil-z)^{-1}$
is a meromorphic function in~$\Compl\setminus (-\infty,-\mu]\cup [\mu, \infty)$
and the poles are eigenvalues of~$\Hil$ of finite multiplicity.
Note that since $\Hil$ is not self-adjoint, it can happen that
\[ \ker(\Hil-z)^2 \ne \ker (\Hil-z)\]
for some $z\in\Compl$. In other words, $\Hil$ can possess {\em generalized eigenspaces}.
In the NLS applications this typically does happen at $z=0$, see Section~\ref{sec:linroot}.

\begin{lemma}
\label{lem:sys_equiv}
Let
\[ \binom{\psi_1}{\psi_2} = \bm 1 & i\\ 1 & -i \endm
\binom{f_1}{f_2}.\]
Then, with $H_0=-\Lapl+\mu$,
\begin{equation}
\label{eq:s_1} \partial_t \binom{f_1}{f_2} + \bm 0 & -H_0-V_1 \\ H_0+V_2
& 0 \endm \binom{f_1}{f_2} =0
\end{equation}
holds if and only if
\begin{equation}
\label{eq:s_2}
\frac{1}{i}\partial_t \binom{\psi_1}{\psi_2} + \bm H_0+U & -W \\ W &
-H_0-U \endm \binom{\psi_1}{\psi_2}
=0
\end{equation}
where $U=\half(V_1+V_2)$, $W=\half(V_1-V_2)$.
\end{lemma}

Abusing notation, we will also write
\begin{equation}
\label{eq:Bdef}
 \Hil_0 = \bm H_0 & 0 \\ 0 & -H_0  \endm,\quad V= \bm U & -W \\ W & -U
\endm,\quad  \Hil = \Hil_0+V=\bm H_0+U & -W \\ W & -H_0-U \endm.
\end{equation}
with $U,W$ real-valued. Moreover, they need to satisfy the following requirements, which we
shall use without any further mention.

\begin{defi}
\label{def:main}
$U$, $W$ are such that $-\Laplace+U$ and $-\Laplace+W$ are
self-adjoint on the domain $W^{2,2}(\R^3)$.
Typically, we want to assume $U,W$ real-valued, bounded and polynomial decaying.
Also, $L_1:=-\Laplace+\mu+U+W\ge0$ (positivity assumption).
\end{defi}

The positivity assumption has some important consequences on the
location and the nature of the spectrum.

\begin{lemma}
\label{lem:spec} 
$\spec(\Hil)=-\spec(\Hil)=\overline{\spec(\Hil)}=\spec(\Hil^*)$
and $\spec(\Hil)\subset \R\cup i\R$.
The discrete spectrum of $\spec(\Hil)$ consists of
eigenvalues $\{z_j\}_{j=1}^N$.
Each $z_j$ with $z_j\ne0$ is an eigenvalue of finite multiplicity,
and $\ker(\Hil-z_j)^2=\ker(\Hil-z_j)$ for those $z_j$.
Thus, only $z=0$ can possess a generalized eigenspace $\bigcup_{k=1}^\infty \ker(\Hil^k)$.
There is a finite $m\ge1$ such that $\ker(\Hil^k)=\ker(\Hil^{k+1})$ for all $k\ge m$.
\end{lemma}
\begin{proof}
The symmetries are consequences of the commutation properties of $\Hil$ with the Pauli matrices
\[ \bm 0 & 1\\ 1& 0\endm, \quad \bm 0 & i\\ -i &  0\endm,
\quad \bm 1 & 0\\ 0 & -1\endm. \]
That the spectrum is only real or purely imaginary under the positivity assumption and
that the discrete spectrum consists only of eigenvalues with trivial Jordan form (up to the zero
eigenvalue) is a standard
argument involving the selfadjoint operator  $\sqrt{L_{-}}L_{+}\sqrt{L_{-}}$,
see Lemma~11.10 in~\cite{RSS2}.
Let $P_0$ be the Riesz projection at zero. Then, on the one hand one checks that
\[ \Ran P_0 \supset \ker(\Hil^m) \text{\ \ for all\ \ }m\ge1.\]
On the other hand, if $\|(\Hil-z)^{-1}\|\les |z|^{-\nu}$, then
\[ \Hil^\nu P_0 =0.\]
Thus $\Ran P_0\subset \ker (\Hil^\nu)$. See~\cite{HS} chapter~6 for these general statements
about Riesz projections.
\end{proof}

So the discrete spectrum takes the form (under the positivity assumption $L_1\ge0$):
Eigenvalues $i\sigma_j$ and $-i\sigma_j$ for $1\le j\le N$, with $\sigma_j>0$, as well
as $-\lambda_\ell$ and $\lambda_\ell>0$ for $1\le \ell\le M$; all of these
are true eigenspaces with trivial Jordan form. Then there is $\omega_0=0$ (if present
at all) which may possess a nontrivial generalized eigenspace.
If we assume that $W$ decays exponentially,
then an adaptation of Agmon's argument~\cite{Agm2} implies that any eigenfunction
or any function which belongs to a generalized eigenspace
 is exponentially decaying.  Such an argument already appeared in~\cite{RSS2},
but it was written there only for real
eigenvalues. The following lemma presents the small changes that
need to be made to deal with the situation at hand.
This is one case where we make use of some exponential decay
of the potential (of $W$ to be precise).

\begin{lemma}
\label{lem:agmon}
Let $\Hil$ be as in \eqref{eq:Hildef} with $U,W$ real-valued,
continuous and $W$ exponentially decaying, whereas $U$
is only required to tend to zero.
If $f\in\ker(A-E)^k$ for some $E$ with $-\mu<\Re E<\mu$ and
some positive integer~$k$, then $f$ decays exponentially.
\end{lemma}
\begin{proof} We will use a variant of Agmon's argument~\cite{Agm2}.
More precisely,
suppose that for some $E$, $-\mu < \Re E< \mu$, there are
$\psi_1,\psi_2\in H^2(\R^3)$ so that
\bear
(\Laplace - \mu + U)\psi_1 - W\psi_2 &=& E\psi_1 \nn \\
W\psi_1 + (-\Laplace + \mu - U)\psi_2 &=& E\psi_2. \label{eq:eig_sys}
\eear
Suppose $|W(x)|\les e^{-b|x|}$.  Then define the Agmon metrics
\bear
\rho_E^{\pm}(x) &=& \inf_{\gamma:0\to x} L^\pm_{\rm Ag}(\gamma) \nn \\
L^\pm_{\rm Ag}(\gamma) &=& \int_0^1 \min\Big( \sqrt{(\mu\pm \Re E-U(\gamma(t)))_{+}}\,,\,
b/2 \Big) \|\dot \gamma(t)\|\,dt \label{eq:Ag_def}
\eear
where $\gamma(t)$ is a $C^1$-curve with $t\in[0,1]$,
and the infimum is to be taken over such curves that
connect $0,x$.  These functions satisfy
\beeq
\label{eq:ag_met}
 |\nabla \rho_E^{\pm}(x)|\le \sqrt{(\mu\pm \Re E-U(x))_{+}}.
\eneq
Moreover, one has $\rho_E^{\pm}(x)\le b|x|/2$ by construction.
Now fix some small $\eps>0$ and
set $\omega^\pm(x):= e^{2(1-\eps)\rho^{\pm}_E(x)}$.
Our goal is to show that
\beeq
\label{eq:Agmon}
\int \Big[ \omega^+(x) |\psi_1(x)|^2 + \omega^{-}(x)|\psi_2(x)|^2\Big] \, dx < \infty.
\eneq
Not only does this exponential decay in the mean suffice for our applications,
but it can also be improved
to pointwise decay using regularity estimates for $\psi_1,\psi_2$.
We do not elaborate on this, see for example \cite{Agm2} and
Hislop, Sigal~\cite{HS}.

\noindent Fix $R$ arbitrary and large. For technical reasons, we set
\[ \rho_{E,R}^{\pm}(x):= \min\Big(2(1-\eps)\rho^{\pm}_E(x),R\Big),\quad
\omega_R^\pm(x):= e^{\rho_{E,R}^{\pm}(x)}. \]
Notice that \eqref{eq:ag_met} remains valid in this case, and also that
$\rho_E^{\pm}(x)\le \min(b|x|/2,R)$. Furthermore,
by choice of $E$ there is a smooth functions $\phi$ that is
equal to one for large $x$ so that
\[
\supp(\phi)  \subset \{\mu+\Re E-U>0\} \cap \{\mu-\Re E-U>0\}.
\]
It will therefore suffice to prove the following modified form
of~\eqref{eq:Agmon}:
\beeq
\label{eq:Agmon'}
\sup_{R}\int \Big[ \omega_R^+(x) |\psi_1(x)|^2 +
\omega_R^{-}(x)|\psi_2(x)|^2 \Big]\phi^2(x)\, dx < \infty.
\eneq
All constants in the following argument will be independent of~$R$.
By construction, there is $\delta>0$ such that
\bear
&& \delta \int  \omega_R^+(x) |\psi_1(x)|^2 \phi^2(x)\,dx \le
\int \omega_R^+(x) ( \mu+\Re E-U(x)) |\psi_1(x)|^2 \phi^2(x)\,dx  \label{eq:delta_int} \\
&& =
\Re \int \omega_R^+(x) (\Laplace \psi_1 - W\psi_2)(x) \bar\psi_1(x) \phi^2(x)\,dx  \nn \\
&& = -  \Re \int \nabla(\omega_R^{+}(x) \phi^2(x)) \nabla \psi_1(x) \bar\psi_1(x)\,dx
 -  \Re \int \omega_R^{+}(x) \phi^2(x) |\nabla \psi_1(x)|^2\,dx  \label{eq:nabla_rho} \\
&&  -  \Re \int \omega_R^+(x) W(x) \psi_1(x)\bar\psi_2(x)\phi^2(x) \, dx. \label{eq:cross}
\eear
As far as the final term \eqref{eq:cross} is concerned, notice that
$\sup_{x,R} |\omega_R^+(x) \phi^2(x) W(x)|\les 1$ by construction, whence
$ |\eqref{eq:cross}| \les \|\psi_1\|_2\|\psi_2\|_2$.
Furthermore, by \eqref{eq:ag_met} and Cauchy-Schwarz,
the first integral in~\eqref{eq:nabla_rho} satisfies
\begin{align}
& \left| \int \nabla(\omega_R^{+}(x) \phi^2(x)) \nabla \psi_1(x) \bar\psi_1(x)\,dx
\right| \nn \\
& \le 2(1-\eps) \Bigl( \int \omega_R^+(x)(\mu+\Re E-U(x)) \phi^2(x)|\psi_1(x)|^2 \, dx \Bigr)^{\half}
\Bigl( \int \omega_R^+(x) \phi(x)^2\, |\nabla\psi_1(x)|^2 \, dx \Bigr)^{\half} \label{eq:int_again} \\
& \qquad + 2 \Bigl( \int \omega_R^+(x) \phi^2(x)|\nabla \psi_1(x)|^2 \, dx \Bigr)^{\half}
\Bigl( \int \omega_R^+(x) |\nabla\phi(x)|^2\, |\psi_1(x)|^2 \, dx \Bigr)^{\half}. \nn
\end{align}
Since the first integral in~\eqref{eq:int_again}
is the same as that in~\eqref{eq:delta_int}, inserting~\eqref{eq:int_again}
into~\eqref{eq:nabla_rho} yields after some simple manipulations
\bear
 \eps \int \omega_R^+(x) ( \mu+\Re E-U(x)) |\psi_1(x)|^2 \phi^2(x)\,dx
&\le& \eps^{-1} \int \omega_R^+(x) |\nabla\phi(x)|^2\, |\psi_1(x)|^2\, dx \nn \\
&&  + \int \omega_R^+(x) \phi(x)^2 |W(x)| |\psi_2(x)|\,|\bar\psi_1(x)|\, dx. \nn
\eear
Since $\nabla\phi$ has compact support, and by our previous considerations
involving $\omega_R^+ W$, the entire right-hand side is bounded independently of~$R$,
and thus also~\eqref{eq:delta_int}. A symmetric argument applies to
the integral with $\psi_2$, and~\eqref{eq:Agmon'}, \eqref{eq:Agmon} hold.
This method also shows that functions belonging to generalized eigenspaces
 decay exponentially. Indeed, suppose $(A-E)\vec g=0$ and
$(A-E)\vec f=\vec g$. Then
\bear
(\Laplace - \mu + U)f_1 - W f_2 &=& E f_1 + g_1\nn \\
W f_1 + (-\Laplace + \mu - U) f_2 &=& E f_2 + g_2 \nn
\eear
with $g_1,g_2$ exponentially decaying.
Decreasing the value of $b$ in~\eqref{eq:Ag_def}
if necessary allows one to use the same argument as before
to prove~\eqref{eq:Agmon} for~$\vec f$. By induction, one then
deals with all values of~$k$ as in the statement of the lemma.
\end{proof}

\noindent As far as the essential spectrum is concerned, we recall
some weighted $L^2$ estimates for the free resolvent $(\Hil_0-z)^{-1}$ which go by the
name of limiting absorption principle.
The weighted $L^2$-spaces here are the usual ones $L^{2,\sigma}=\la x\ra^{-\sigma} L^2$.
It will be convenient to introduce the space
\[ X_\sigma := L^{2,\sigma}(\R^3)\times L^{2,\sigma}(\R^3). \]
Clearly, $X_\sigma^* = X_{-\sigma}$. The statement is that
\beeq
\label{eq:H0agm}
 \sup_{|\lambda|\ge\lambda_0,\,0<\eps} |\lambda|^{\half}\,\|(\Hil_0-(\lambda\pm i\eps))^{-1}\|_{X_\sigma\to X_\sigma^*} < \infty
\eneq
provided $\lambda_0>\mu$ and $\sigma>\half$ and was proved in this form by Agmon~\cite{Agm1}.
By the explicit expression for the kernel of the free resolvent in $\R^3$
one obtains the existence of the limit
\[ \lim_{\eps\to0+} \la (\Hil_0-(\lambda\pm i\eps))^{-1} \phi,\psi \ra \]
for any $\lambda\in\R$ and any pair of Schwartz functions $\phi,\psi$, say. Hence
$(\Hil_0-(\lambda\pm i0))^{-1}$  satisfies the same bound as in~\eqref{eq:H0agm}
provided $|\lambda|>\mu$.
We will also need a bound which is valid for all energies. It takes the form
\beeq
\label{eq:HS1}
 \sup_{z\in\Compl}\|(\Hil_0-z)^{-1}\|_{X_\sigma\to X_\sigma^*} < \infty
\eneq
provided $\sigma>1$. This is much easier to obtain than~\eqref{eq:H0agm} and only uses
that convolution with $|x|^{-1}$ is bounded from $L^{2,\sigma}(\R^3)\to L^{2,-\sigma}(\R^3)$
provided $\sigma>1$. In fact, it is Hilbert-Schmidt in these norms. With these preparations, let
us state a lemma about absence of embedded resonances.

\begin{lemma}
\label{lem:invert}
Assume that there are no embedded eigenvalues in the essential spectrum of $\Hil=\Hil_0+V$.
Suppose $|V(x)|\les \la x\ra^{-\beta}$ with $\beta>1$. Then for any $\lambda\in\R$,
$|\lambda|>\mu$ the operator
$(\Hil_0-(\lambda\pm i0))^{-1}V$ is a compact operator on $X_{-\half-}\to X_{-\half-}$
and
\[ I + (\Hil_0-(\lambda\pm i0))^{-1}V \]
is invertible on these spaces.
\end{lemma}
\begin{proof}
The compactness is standard and we refer the reader to~\cite{Agm1} or~\cite{RS4}.
Let $\lambda>\mu$.
By the Fredholm alternative,
the invertibility statement requires excluding solutions $(\psi_1,\psi_2)\in X_{-\half-}$
of the system
\bear
0 &=& \psi_1 + R_0(\lambda-\mu+i0)(U\psi_1+W\psi_2) \nn \\
0 &=& \psi_2 - R_0(-\lambda-\mu)(W\psi_1+U\psi_2), \nn
\eear
where $R_0(z)$ is the free, scalar resolvent $(-\Laplace-z)^{-1}$.
Notice that these equations imply that $\psi_2\in L^2$ and that
\bear
0 &=& \la \psi_1, U\psi_1\ra + \la \psi_1, W\psi_2\ra + \la R_0(\lambda-\mu+i0)(U\psi_1+W\psi_2),
(U\psi_1+W\psi_2) \ra \nn \\
0 &=& \la \psi_2, W\psi_1 \ra - \la R_0(-\lambda-\mu)(W\psi_1+U\psi_2), W\psi_1 \ra \nn \\
0 &=& \la \psi_2, U \psi_2 \ra - \la R_0(-\lambda-\mu) W\psi_1, U\psi_2\ra - \la R_0(-\lambda-\mu)U\psi_2, U\psi_2\ra.
\nn
\eear
Since $U,W$ are real, inspection of these equations reveals that
\[ \Im \la R_0(\lambda-\mu+i0)(U\psi_1+W\psi_2), (U\psi_1+W\psi_2) \ra =0.\]
So Agmon's well-known bootstrap lemma~\cite{Agm1} can be used to conclude that $\psi_1\in L^2$.
But then we have an embedded eigenvalue at $\lambda$, which is impossible.
So one can invert
\[ I+ (\Hil-(\lambda\pm i0))^{-1}V\]
on $X_{-\half-}$ and we are done.
\end{proof}

Next, we need to ensure that the thresholds $\pm\mu$ are neither eigenvalues nor resonances.
This simply means that the following variant of the previous lemma holds also for $\lambda=\pm\mu$.

\begin{defi}
\label{def:edges}
Let $|V(x)|\les \la x\ra^{-\beta}$ with $\beta>2$.
We require that
\beeq
\label{eq:thres}
 I + (\Hil_0-(\lambda \pm i0))^{-1}V: X_{-1-} \to X_{-1-}
\eneq
is invertible for $\lambda=\pm\mu$.
In this case, we refer to the thresholds $\pm\mu$ as being {\em regular}.
\end{defi}

From now on, we will work under the assumptions of this definition.

\begin{prop}
\label{prop:lim_ap}
Assume that there are no embedded eigenvalues in the essential spectrum
and that the thresholds $\pm\mu$ are regular.
Then there is $0<\mu'<\mu$ so that
\beeq
\label{eq:limap}
\sup_{|\lambda|\ge\mu',\,0<\eps} |\lambda|^{\half}\,\|(\Hil-(\lambda\pm i\eps))^{-1}\| < \infty
\eneq
where the norm is the one from  $X_{1+}\to X_{-1-}$. If the supremum in~\eqref{eq:limap}
is only taken over $|\lambda|\ge\lambda_0$ where $\lambda_0>\mu$, then~\eqref{eq:limap}
also holds in the norms of $X_{\half+}\to X_{-\half-}$.
\end{prop}
\begin{proof}
We start by observing that
\[  I + (\Hil_0-(\lambda \pm i0))^{-1}V: X_{-1-} \to X_{-1-} \]
is invertible for all $\lambda$ in neighborhoods of $\pm\mu$.
This follows from~\eqref{eq:thres} and the fact that
\[ \|(\Hil_0-(\lambda \pm i0))^{-1} - (\Hil_0-(\mu\pm i0))^{-1}\| \to 0\]
as $\lambda\to \mu$ in the Hilbert-Schmidt norm of $X_{1+}\to X_{-1-}$
(cf.~the three-dimensional free resolvent close to zero energy).
Let $\mu'\in (0,\mu)$ be such that these neighborhoods contain $\pm\mu'$.
Let $z=\lambda+i\eps$, $\lambda\ge\mu'$, $\eps\ne0$. By the resolvent identity
and the fact that the spectrum of $\Hil$ belongs to $\R\cup i\R$,
\beeq
\label{eq:res}
 (\Hil-z)^{-1} = (I+(\Hil_0-z)^{-1}V)^{-1}(\Hil_0-z)^{-1}
\eneq
as operators on $L^2(\R^3)$. Because of the $|\lambda|^{-\half}$-decay in~\eqref{eq:H0agm},
there exists a positive radius $r_V$ such that
\[ \|(\Hil_0-z)^{-1}V\|<\half \]
for all $|z|>r_V$ in the operator norm of $X_{1+}\to X_{-1-}$. In conjunction with~\eqref{eq:res}
this implies that
\[ \|(\Hil-z)^{-1}\| \les |z|^{-\half} \]
for all $|z|>r_V$ in the operator norm of $X_{-1-}\to X_{-1-}$.
Now suppose~\eqref{eq:limap} fails. It then follows from~\eqref{eq:res} and~\eqref{eq:HS1}
that there exist a sequence $z_n$ with $\Re(z_n)>\mu'$ and functions
$f_n \in X_{-1-}$ with $\|f_n\|_{X_{-1-}}=1$ and such that
\beeq
\label{eq:small_norm}
 \|[I+(\Hil_0-z_n)^{-1}V] f_n \|_{X_{-1-}} \to 0
\eneq
as $n\to\infty$. Necessarily, the $z_n$ accumulate at some point $\lambda\in [\mu',r_V]$.
Without loss of generality, $z_n\to \lambda$ and $\Im(z_n)>0$ for all $n\ge1$.
Next, we claim that~\eqref{eq:small_norm} also holds in the following
form:
\beeq
\label{eq:small_norm2}
 \|[I+(\Hil_0-(\lambda+i0))^{-1}V] f_n \|_{X_{-1-}} \to 0
\eneq
as $n\to\infty$. If so, then it would clearly contradict Lemma~\ref{lem:invert} or
 our discussion involving Definition~\ref{def:edges}. To prove~\eqref{eq:small_norm2},
let
\[ S:= I+(\Hil_0-(\lambda+i0))^{-1}V \]
for simplicity. Then
\begin{align}
I+(\Hil_0-z_n)^{-1}V &= S + ((\Hil_0-z_n)^{-1}-(\Hil_0-(\lambda+i0))^{-1}) V \nn \\
&= \big[I + ((\Hil_0-z_n)^{-1}-(\Hil_0-(\lambda+i0))^{-1}) VS^{-1}\big]S. \label{eq:SHnew}
\end{align}
Our claim now follows from the fact that the expression in brackets
is an invertible operator for large $n$ on $X_{-1-}$. The final statement
concerning $X_{\half+}$ is implicit in the preceding.
\end{proof}

As in the case of the free Hamiltonian $\Hil_0$, it is now possible to
define the boundary values of the resolvent $(\Hil-z)^{-1}$. More
precisely, the following corollary holds.

\begin{cor}
\label{cor:Hbdry}
Assume that there are no embedded eigenvalues in the essential spectrum and that the
thresholds $\pm\mu$ are regular. In that case it is possible to define
\beeq
\label{eq:Hbdrydef}
(\Hil-(\lambda\pm i0))^{-1}:=(I+(\Hil_0-(\lambda\pm i0))^{-1}V)^{-1}(\Hil_0-(\lambda\pm i0))^{-1}
\eneq
for all $|\lambda|>\mu'$ where $\mu'$ is as in~\eqref{eq:limap}.
Then as $\eps\to0+$,
\[ \|(\Hil-(\lambda\pm i\eps))^{-1} - (\Hil-(\lambda\pm i0))^{-1}\| \to0\]
in the norm of $X_{1+}\to X_{-1-}$
and one can extend  \eqref{eq:limap} to $\eps\ge0$. The same statements hold with
$X_{\half+}\to X_{-\half-}$ provided $|\lambda|\ge \lambda_0>\mu$.
\end{cor}
\begin{proof}
$(\Hil-(\lambda\pm i\eps))^{-1}$ is well-defined for $|\lambda|>\mu'$ by
Lemma~\ref{lem:invert} or
 our discussion involving Definition~\ref{def:edges} in the previous proof. Hence
\eqref{eq:Hbdrydef} is legitimate for those~$\lambda$, see also~\eqref{eq:res}.
Moreover, by~\eqref{eq:SHnew},
\begin{align*}
& [I+(\Hil_0-(\lambda+i\eps))^{-1}V]^{-1} - [I+(\Hil_0-(\lambda+i0))^{-1}V]^{-1} \\
& = S^{-1}\big[I + ((\Hil_0-(\lambda+i\eps))^{-1}-(\Hil_0-(\lambda+i0))^{-1}) VS^{-1}\big]^{-1} - S^{-1} \\
& = \sum_{k=1}^\infty S^{-1}
\big[-((\Hil_0-(\lambda+i\eps))^{-1}-(\Hil_0-(\lambda+i0))^{-1}) VS^{-1}\big]^{k}
\end{align*}
tends to zero in the norm of $X_{-1-}$ as $\eps\to0+$. Also,
\[  \|(\Hil_0-(\lambda\pm i\eps))^{-1} - (\Hil_0-(\lambda\pm i0))^{-1}\| \to 0\]
as $\eps\to0$ in the norm of $X_{1+}\to X_{-1-}$. Combining these two convergence
statements finishes the proof of the corollary.
\end{proof}

For later applications to dispersive estimates we will also need to
control derivatives of the resolvents in weighted $L^2$ spaces.

\begin{cor}
\label{cor:der_resolv}
Assume that there are no embedded eigenvalues, and that the matrix potential
$V$ decays like $|V(x)|\les \la x\ra^{-3-}$. Then for every $\lambda_0>\mu$,
\begin{align*}
\sup_{|\lambda|\ge\lambda_0}\big\| \partial_\lambda (\Hil-(\lambda\pm i0))^{-1} \big\|_{X_{\frac32+}\to X_{-\frac32-}}  &\les 1 \\
\sup_{|\lambda|\ge\lambda_0}\big\| \partial^2_\lambda (\Hil-(\lambda\pm i0))^{-1} \big\|_{X_{\frac52+}\to X_{-\frac52-}}  &\les 1.
\end{align*}
\end{cor}
\begin{proof}
We first note that the
 derivatives $\partial_\lambda^j(\Hil_0-(\lambda\pm i0))^{-1}$
of the free (matrix) Hamiltonian $\Hil_0$
satisfy the uniform bounds
$$\sup_{\lambda}\left\|\partial_\lambda^j(\Hil_0-(\lambda\pm i0))^{-1}f
\right\|_{X_{-\sigma}} \les \|f\|_{X_\sigma}$$
for all $\sigma > j + \frac12$ and $j\ge1$. This can be seen from the kernel of
the three-dimensional free scalar resolvent. Next, we transfer this to the perturbed resolvent
by means of
\[
(\Hil-(\lambda\pm i0))^{-1}=(I+(\Hil_0-(\lambda\pm i0))^{-1}V)^{-1}(\Hil_0-(\lambda\pm i0))^{-1}.
\]
Writing $R_0^\pm$ and $R_V^\pm$ for the free and perturbed resolvents, respectively, and
setting $S^\pm(\lambda):=I + R_0^\pm (\lambda)V$, we obtain that
\beeq
\label{eq:der_id}
\partial_\lambda R_V^{\pm}(\lambda)  = -S^\pm(\lambda)^{-1}\partial_\lambda R_0^\pm(\lambda)\,VS^\pm(\lambda)^{-1}
R_0^\pm(\lambda)+ S^\pm(\lambda)^{-1}\partial_\lambda R_0^\pm(\lambda),
\eneq
and since $\sup_{\lambda>\lambda_0}\|S^\pm(\lambda)^{-1}\|_{X_{-\half-}\to X_{-\half-}} < \infty$,
it follows that also
\begin{equation}
\label{eq:rvdev} \sup_{\lambda>\lambda_0}\|\partial_\lambda\,R_V^\pm
(\lambda^2)\|_{X_{\frac32+}\to X_{-\frac32-}} \les 1
\end{equation}
Note from~\eqref{eq:der_id} that one needs to assume the decay $|V(x)|\les (1+|x|)^{-2-\eps}$ for
this to hold. Indeed, $V$ needs to take $X_{-\half-} \to X_{\frac32+}$.
By a similar argument,
\[ \|\partial_\lambda^2 R_V^\pm(\lambda)\|_{X_{\frac52+}\to X_{-\frac52-}} \les 1.\]
This estimate requires the decay $|V(x)|\les (1+|x|)^{-3-}$ by an analogous formula
to~\eqref{eq:der_id}.
\end{proof}

Another easy corollary of Proposition~\ref{prop:lim_ap} is the finiteness of
the discrete spectrum.

\begin{cor}
\label{cor:discfin}
Under the previous assumptions the discrete spectrum of $\Hil$ is finite.
\end{cor}
\begin{proof}
Suppose not. Since  the spectrum of $\Hil$ is a subset of $\R\cup i\R$ and since the
discrete spectrum can only accumulate on the essential spectrum, any
accumulation point would have to be either one of $\pm\mu$.
However, this would contradict~\eqref{eq:limap}.
We are using here that any eigenfunction has to belong to $X_{1+}$ and
not just to $L^2$ because of Agmon's lemma~\ref{lem:agmon}.
\end{proof}

We can now define a projection $P_s$ (the projection onto the {\em stable subspace}).
Indeed, consider the Riesz projection
\beeq
\label{eq:Psdef}
 P_s = I-\frac{1}{2\pi i}\oint_\Gamma (\Hil-z)^{-1} \, dz,
\eneq
where $\Gamma$ is a contour that surrounds $\spec(H)\cap i\R$ but does not contain any other
portion of the spectrum. Similarly, one defines $P_c$ as the projection
corresponding to the essential spectrum (the projection onto the {\em continuous subspace}).
In other words, $I-P_c$ is the Riesz projection given by a contour integral
of the resolvent around a curve which encircles the entire discrete spectrum but is disjoint
from the essential spectrum.
In this connection, one has the following variant of Lemma~7.2 from~\cite{RSS1}.
We give an independent proof here.

\begin{lemma}
\label{lem:L2split}
Let $\{L_j\}_{j=0}^M$ and $\{L_j^*\}_{j=0}^M$ be the generalized
eigenspaces of $\Hil$ and $\Hil^*$, respectively.
Then there is a direct sum decomposition
\begin{equation}
\label{eq:L2split}
L^2(\R^3)\times L^2(\R^3) = \sum_{j=0}^M L_j + \Bigl(\sum_{j=0}^M L_j^*\Bigr)^\perp
\end{equation}
This means that the individual summands are linearly independent,
but not necessarily orthogonal.
The decomposition \eqref{eq:L2split} is invariant under~$\Hil$.
The Riesz projection ${\rm Id}-P_s$, see~\eqref{eq:Psdef}, is precisely the projection onto those
$L_j$ which correspond to eigenvalues on $i\R$ and which
is induced by the splitting~\eqref{eq:L2split} (i.e., it preserves all the other summands in
that direct sum representation).
Similarly, $P_c$ is the projection
onto $\Bigl(\sum_{j=0}^M L_j^*\Bigr)^\perp$ which is induced by the
splitting~\eqref{eq:L2split}.
\end{lemma}
\begin{proof}
This is immediate from the definition of the Riesz projections. First,
\beeq
\label{eq:rp}
 I-P_c = \frac{1}{2\pi i} \oint_\gamma (zI-\Hil)^{-1}\, dz
\eneq where $\gamma$ is a simple closed curve that encloses the
entire discrete spectrum of~$\Hil$ and lies within the resolvent
set. Then, on the one hand,
\[ L^2(\R^3)\times L^2(\R^3) = \ker(P_c) + \Ran(P_c) = \ker(P_c) + \ker (P_c^*)^{\perp}.\]
On the other hand,
\[ \ker(P_c) = \Ran(I-P_c) = \sum_{j=0}^M L_j\]
as well as
\[ \ker (P_c^*) = \sum_{j=0}^M L_j^*.\]
This last equality uses that $P_c^*$ is the same as
the Riesz projection off the discrete spectrum of $\Hil^*$,
as can be seen by taking adjoints of~\eqref{eq:rp}.
The argument for $P_s$ is analogous.
\end{proof}

Clearly,
\beeq
\label{eq:Psmooth} P_s,\;P_c:W^{2,2}\to W^{2,2}.
\eneq
The range of $P_s$ consists of linear combinations of eigenfunctions
and elements of generalized eigenspaces. By Lemma~\ref{lem:agmon}, any
such linear combination has to be exponentially decaying.

We can now state the representation formula for $e^{it\Hil}$ which
is basic to all estimates on this evolution which we prove here.
In~\cite{ErdSch} this same statement is derived without any
assumption on the thresholds.

\begin{lemma}
\label{lem:rep}
Assume that there are no embedded eigenvalues in the essential spectrum
and that the thresholds $\pm\mu$ are regular. Then
there is the representation
\beeq
e^{it\Hil} = \frac{1}{2\pi i}\int_{|\lambda|\ge\mu} e^{it\lambda}\, [(\Hil-(\lambda+i0))^{-1}-(\Hil-(\lambda-i0))^{-1}]\,d\lambda + \sum_{j} e^{it\Hil} P_{\zeta_j},
\label{eq:ac}
\eneq
where the sum runs over the entire discrete spectrum $\{\zeta_j\}_j$  and
 $P_{\zeta_j}$ is the Riesz projection corresponding to the eigenvalue $\zeta_j$.
The formula~\eqref{eq:ac} and the convergence of the integral are
to be understood in the following weak sense: If $\phi,\psi$ belong to
$[W^{2,2}\times W^{2,2}(\R^3)]\cap X_{1+}$, then
\[
\la e^{it\Hil}\phi,\psi\ra = \lim_{R\to\infty} \frac{1}{2\pi i}\int_{R\ge|\lambda|\ge \mu}\!\! e^{it\lambda}
\big\la [(\Hil-(\lambda+i0))^{-1}-(\Hil-(\lambda-i0))^{-1}]\phi,\psi\big\ra\,d\lambda + \sum_{j} \la e^{it\Hil}P_{\zeta_j}\phi,\psi \ra.
\]
for all $t$, where the integrand is well-defined by the limiting absorption principle, see
Proposition~\ref{prop:lim_ap} and Corollary~\ref{cor:Hbdry}.
\end{lemma}
\begin{proof}
The evolution $e^{it\Hil}$ is defined via the Hille-Yoshida theorem. Indeed, let $a>0$ be large.
Then $i\Hil-a$ satisfies (with $\rho$ the resolvent set)
\[
 \rho(i\Hil-a)\supset (0,\infty)\text{\ \ and\ \ }\|(i\Hil-a-\lambda)^{-1}\|\le |\lambda|^{-1}\text{\ \ for all\ \ }\lambda>0.
\]
Hence $\{e^{t(i\Hil-a)}\}_{t\ge0}$ is a contractive semigroup, so that $\|e^{it\Hil}\|_{2\to2}\le e^{|t|a}$
for all $t\in\R$.
If $\Re z>a$, then there is the Laplace transform
\beeq
\label{eq:Lap}
 (i\Hil-z)^{-1}=-\int_0^\infty e^{-tz}\,e^{it\Hil}\,dt
\eneq
as well as its inverse (with $b>a$ and $t>0$)
\beeq
\label{eq:invLap}
e^{it\Hil} = -\frac{1}{2\pi i} \int_{b-i\infty}^{b+i\infty} e^{tz}\, (i\Hil-z)^{-1}\, dz.
\eneq
While \eqref{eq:Lap} converges in the norm sense, defining~\eqref{eq:invLap} requires more care.
The claim is that for any $\phi,\psi\in \Dom(\Hil)=W^{2,2}\times W^{2,2}$,
\beeq
\label{eq:limR}
 \la e^{it\Hil} \phi,\psi \ra =
- \lim_{R\to\infty} \frac{1}{2\pi i} \int_{b-iR}^{b+iR} e^{tz}\, \la (i\Hil-z)^{-1}\phi,\psi\ra\, dz.
\eneq
To verify this, let $t>0$ and use \eqref{eq:Lap} to conclude that
\bear
-\frac{1}{2\pi i} \int_{b-iR}^{b+iR} e^{tz}\, \la (i\Hil-z)^{-1}\phi,\psi\ra\, dz &=&
 \frac{1}{2\pi i} \int_{b-iR}^{b+iR} e^{tz}\,  \int_0^\infty e^{-sz}\,\la e^{is\Hil}\,\phi,\psi\ra\, dsdz \nn \\
&=& \frac{1}{\pi} \int_0^\infty  e^{(t-s)b}\,\frac{\sin((t-s)R)}{t-s}\,\la e^{is\Hil}\,\phi,\psi\ra\, ds.
\label{eq:dirker}
\eear
Since $e^{(t-s)b}\,\la e^{is\Hil}\,\phi,\psi\ra$ is a $C^1$ function in~$s$ (recall $\phi\in \Dom(\Hil)$)
 as well as  exponentially decaying in $s$ (because of $b>a$),
it follows from standard properties of the Dirichlet kernel that
the limit in~\eqref{eq:dirker} exists and equals $\la e^{it\Hil}\phi,\psi\ra$, as claimed.
Note that if $t<0$, then the limit is zero. Therefore, it follows that for any $b>a$,
\begin{align*}
 \la e^{it\Hil} \phi,\psi \ra &=
- \lim_{R\to\infty} \Big\{
\frac{1}{2\pi i} \int_{b-iR}^{b+iR} e^{tz}\, \la (i\Hil-z)^{-1}\phi,\psi\ra\, dz
- \frac{1}{2\pi i} \int_{-b-iR}^{-b+iR} e^{tz}\, \la (i\Hil-z)^{-1}\phi,\psi\ra\, dz \Big\} \\
&= \lim_{R\to\infty}
\frac{1}{2\pi i} \int_{-R}^{R} e^{it\lambda}\, \la [e^{-bt}(\Hil-(\lambda+ib))^{-1}
- e^{bt}(\Hil-(\lambda-ib))^{-1}]\phi,\psi\ra\, d\lambda.
\end{align*}
Next, assume that $\phi,\psi$ are as in the statement of the theorem,
and shift the contour in the previous
integrals by sending $b\to0+$.
More precisely, we apply the residue theorem to the contour integrals over the rectangles with
vertices $\pm R+ib$, $\pm R+i0$ and the reflected one below the real axis.  The horizontal
segments on the real axis need to avoid the poles, which can be achieved by
 surrounding each of the at most finitely many real poles of the resolvent $(\Hil-z)^{-1}$
by a small semi-circle. Combining each such semi-circle with its reflection yields a
small closed loop and the resulting integral is precisely the Riesz projection corresponding to
that real eigenvalue.
The Riesz projections corresponding to eigenvalues on the imaginary axis are obtained as
residues.
On the other hand, we also need to show that the contribution
by the horizontal segments is zero in the limit $R\to\infty$. This, however, follows
from the limiting absorption principle of Proposition~\ref{prop:lim_ap} and Corollary~\ref{cor:Hbdry}.
The lemma follows.
\end{proof}

Note that $\|e^{it\Hil}f\|_2\sim |t|^m\|f\|_2$ if $\Hil^m f\ne0$ but $\Hil^{m+1}f=0$.
This shows that in general $e^{it\Hil}$ is not $L^2$ bounded. However, using the previous
lemma one can show the following.

\begin{theorem}
\label{thm:L2stable}
Assume as before that there are no embedded eigenvalues and that
the thresholds $\pm\mu$ are regular.
Then the following {\em stability} bound holds:
\beeq
\sup_{t\in\R} \|e^{it\Hil} P_s\|_{2\to2} \le C. \label{eq:L2stable}
\eneq
\end{theorem}
\begin{proof}
We need to check via~\eqref{eq:ac} that the $L^2$-norms do not grow, see~\eqref{eq:L2stable}.
To do so, fix again 'nice' $\phi,\psi$ as in the previous theorem and with $P_s\phi=\phi, P_s\psi=\psi$.
Note that this can be done by the aforementioned mapping properties of $P_s$.
Then
\bear
&&  \sup_{t}|\la e^{it\Hil}\phi,\psi \ra| \nn \\
&\le& \limsup_{\eps\to0} \int_{|\lambda|\ge\mu}
\big|\,\la [(\Hil-(\lambda+i\eps))^{-1}-(\Hil-(\lambda-i\eps))^{-1}] \phi,\psi \ra\,\big|\, d\lambda
+\|\phi\|_2\|\psi\|_2 \nn \\
& \le& \limsup_{\eps\to0} 2\eps \int_{|\lambda|\ge\mu} \|(\Hil-(\lambda-i\eps))^{-1}\phi\|_2\,
\|(\Hil^*-(\lambda-i\eps))^{-1}\psi\|_2 \, d\lambda +\|\phi\|_2\|\psi\|_2. \nn
\eear
The terms $\|\phi\|_2\|\psi\|_2$ here are due to possible real eigenvalues which are not removed
by~$P_s$.
It therefore suffices to show that (and similarly for $\Hil^*$)
\[
\limsup_{\eps\to0}\, \eps \int_{|\lambda|\ge\mu} \|(\Hil-(\lambda-i\eps))^{-1}\phi\|_2^2\,d\lambda
\les \|\phi\|_2^2.
\]
To do so, we apply the resolvent identity again:
\bear
&& 2\eps i \int_{|\lambda|\ge\mu} \|(\Hil-(\lambda+i\eps))^{-1}\phi\|_2^2\,d\lambda \nn \\
&=& -\sum_{\zeta_j\in\R,\,\zeta_j\ne0}\|P_{\zeta_j}\phi\|_2^2 + O(\eps) +\int_{-\infty}^\infty \la (\Hil-(\lambda+i\eps))^{-1}\phi, \phi \ra \, d\lambda \label{eq:p1} \\
&& - \int_{-\infty}^\infty \la (\Hil^*-(\lambda-i\eps))^{-1}\phi, \phi \ra \, d\lambda \label{eq:p2} \\
&& + \int_{|\lambda|\ge\mu}  \|M(\Hil-(\lambda+i\eps))^{-1}\phi\|_2^2\,d\lambda \label{eq:p3},
\eear
where
\[ M = \bm 0 & 2 W\\ -2W & 0 \endm. \]
To obtain the two integrals on the right-hand sides of \eqref{eq:p1} and~\eqref{eq:p2}
we needed to add (or subtract) the integrals over the segments $(-\mu,\mu)$ that were not
present on the left-hand side. This is equivalent to adding the sum of the projections
\[ \sum_{\zeta_j\in\R,\,\zeta_j\ne0} \la P_{\zeta_j}\phi,\phi\ra + O(\eps)\]
 where the $O$-term results from
vertical segments of length $\les\eps$.
The integrals on the right-hand sides of \eqref{eq:p1} and~\eqref{eq:p2} do not
present much of a problem: the integrands are either analytic or co-analytic and one
can therefore shift the contour to large $\eps>0$ (to do so requires applying the limiting
absorption estimate in the same way as above). If $\eps>0$ is sufficiently large
(in fact, if $\eps>a$ from above), then one has
\bear
 \lim_{L\to\infty} \int_{-L}^L \la (\Hil-(\lambda-i\eps))^{-1} \phi,\phi \ra \, d\lambda
 &=& -i\lim_{L\to\infty} \int_{-L}^L \int_0^\infty e^{-is(\lambda-i\eps)} \la e^{is\Hil} \phi,\phi \ra \,ds d\lambda  \nn \\
&=& -2i \lim_{L\to\infty} \int_0^\infty \frac{\sin(sL)}{s} e^{-s\eps} \la e^{is\Hil} \phi,\phi \ra  \,ds \nn \\
&=& \const\,\cdot \|\phi\|_2^2. \nn
\eear
The term~\eqref{eq:p3} is controlled by Kato smoothing theory.
We show that it is $\les \|\phi\|_2^2$ uniformly for $0<\eps<1$.
Start from the resolvent identity
\[ (\Hil-(\lambda+i\eps))^{-1}-(\Hil_0-(\lambda+i\eps))^{-1} =
 -(\Hil_0-(\lambda+i\eps))^{-1}V(\Hil-(\lambda+i\eps))^{-1}. \]
Set $\rho(x)=\la x\ra^{-1-}$, say, and define
\beeq
\label{eq:Mtil} \tilde{M}= \bm \rho & 0 \\ 0 & \rho \endm.
\eneq
Then
\[ \tilde{M}(\Hil-(\lambda+i\eps))^{-1}= \tilde{M}(\Hil_0-(\lambda+i\eps))^{-1}
 - \tilde{M}(\Hil_0-(\lambda+i\eps))^{-1}V\tilde{M}^{-1}\tilde{M}(\Hil-(\lambda+i\eps))^{-1}, \]
which implies that
\[ \tilde{M}(\Hil-(\lambda+i0))^{-1}= \big[1+\tilde{M}(\Hil_0-(\lambda+i0))^{-1}V\tilde{M}^{-1}\big]^{-1}
\tilde{M}(\Hil_0-(\lambda+i0))^{-1}, \]
provided the term in brackets is actually invertible.
However, since $\tilde{M}(\Hil_0-(\lambda+i0))^{-1}V\tilde{M}^{-1}$ is compact on $L^2$,
this inverse exists iff there is no nonzero $f\in L^2(\R^3)$ with the property that
\[ [1+\tilde{M}(\Hil_0-(\lambda+i0))^{-1}V\tilde{M}^{-1}]f = 0 \]
or equivalently,
\[ [1+(\Hil_0-(\lambda+i0))^{-1}V] \tilde{M}^{-1}f = 0 \]
where now $\tilde{M}^{-1}f\in X_{-1-}$.
However, the existence of such an $f\ne0$
would contradict Lemma~\ref{lem:invert} or Definition~\ref{def:edges}.
It follows that the aforementioned inverses exist.
In fact, not only do the inverses of the expressions in brackets exist, but the norms of
these inverses are uniformly bounded in $|\lambda|\ge\mu$.
Thus, we have
\[ \|\tilde{M}(\Hil-(\lambda+i0))^{-1}\phi\| \le C\, \|\tilde{M}(\Hil_0-(\lambda+i0))^{-1}\phi\|\]
for all $|\lambda|\ge\mu$. But since on the one hand $\|M\tilde{M}^{-1}\|_\infty <\infty$ (assuming that
$U,W$ decay like $(1+|x|)^{-1-}$),
and on the other hand
\[ \int_{|\lambda|\ge\mu} \|\tilde{M}(\Hil_0-(\lambda+i0))^{-1}\phi\|_2^2 \,d\lambda \le C\|\phi\|_2^2,\]
(this is the statement that $\rho(x)$ is $-\Laplace$-smooth which is standard) we conclude that
\beeq
\label{eq:Hilsmooth} \int_{|\lambda|\ge\mu}  \|M(\Hil-(\lambda+i0))^{-1}\phi\|_2^2\,d\lambda \le C\|\phi\|_2^2,
\eneq
as desired.
\end{proof}

\section{The linear analysis: dispersive theory}
\label{sec:lin_dis}

We now turn to the dispersive bound.

\begin{theorem}
\label{thm:disp}
Let $|V(x)|\les \la x\ra^{-\beta}$ for all $x\in\R^3$ with some $\beta>3$.
Assume again that there are no embedded eigenvalues and that the thresholds
$\pm\mu$ are regular.
Then, with ${\rm Id}-P_c$ being the Riesz projection corresponding to the
entire discrete spectrum, there is the dispersive bound
\beeq
\|e^{it\Hil} P_c\|_{1\to\infty} \les |t|^{-\frac32}. \label{eq:dec}
\eneq
\end{theorem}
\begin{proof}
We will use the method of proof from \cite{GolSch}.
We start from Lemma~\ref{lem:rep}, i.e.,
\beeq
\label{eq:scal_flow}
e^{it\Hil}P_{c} = \frac{1}{2\pi i}\int_{|\lambda|\ge\mu}
 e^{it\lambda} [(\Hil-(\lambda+i0))^{-1} - (\Hil-(\lambda-i0))^{-1}]\, d\lambda.
\eneq
We distinguish between energies close to the thresholds $\pm\mu$ and those separated from these points.
Thus let $\chi_+(\lambda)=1$ if $\lambda-\mu>2\lambda_1$ and $=0$ if $\lambda-\mu\le\lambda_1$
where $\lambda_1>0$ is some small number.
Similarly, $\chi_-(\lambda)=1$ if $\lambda+\mu<-2\lambda_1$ and $=0$ if $\lambda+\mu\ge-\lambda_1$.
We will use the notation $\chi_+(\Hil)$ and $\chi_-(\Hil)$ formally with the obvious meaning.
Let $R_0^{\pm}(\lambda^2)$ and $R_V^{\pm}(\lambda^2)$ be the resolvents of $\Hil_0$, and
$\Hil=\Hil_0+V$, respectively.
Then, by a finite resolvent expansion and a change of variables $\lambda\to\lambda^2+\mu$,
\begin{align}
& \la e^{it\Hil}\chi_+(\Hil)P_{c}f,g\ra = \frac{e^{it\mu}}{\pi i}\int_0^\infty \lambda e^{it\lambda^2}
 \chi_+(\lambda^2+\mu) \big\la [R_V^{+}(\lambda^2+\mu)-R_V^{-}(\lambda^2+\mu)] f,g \big\ra \, d\lambda \nn
\\
&= e^{it\mu}\sum_{\ell=0}^{2m-1} \frac{(-1)^\ell}{\pi i}\int_0^\infty \lambda  e^{it\lambda^2}
\chi_+(\lambda^2+\mu) \big\la [R_0^+(\lambda^2+\mu)(VR_0^+(\lambda^2+\mu))^\ell \label{eq:2mborn}\\
& \qquad\qquad\qquad - R_0^{-}(\lambda^2+\mu)(VR_0^{-}(\lambda^2+\mu))^\ell]f,g \big\ra \, d\lambda \nn\\
& +\frac{e^{it\mu}}{\pi i}\int_0^\infty \lambda  e^{it\lambda^2} \chi_+(\lambda^2+\mu)
\big\la [(R_0^+(\lambda^2+\mu)V)^m R_V^+(\lambda^2+\mu)(VR_0^+(\lambda^2+\mu))^m \nn\\
& \qquad - (R_0^-(\lambda^2+\mu)V)^m R_V^{-}(\lambda^2+\mu) (V R_0^{-}(\lambda^2+\mu))^m]f,g \big\ra \, d\lambda. \label{eq:born_ser}
\end{align}
We need to show that each of the $2m$ terms in the finite (Born) sum is in absolute value $\le C(\ell,V)\, |t|^{-\frac32}\|f\|_1\|g\|_1$, and similarly for the remaining term containing $R_V$.

Each of the first $2m$ terms of the Born series is written out explicitly
using the free scalar resolvent ($\Im z>0, \Im\sqrt{z}>0$)
\[ (-\Laplace-z)^{-1}(x,y) = \frac{e^{i\sqrt{z}|x-y|}}{4\pi|x-y|},\]
which implies for the matrix case
\beeq
\label{eq:R0}
 R_0^{\pm}(\lambda^2+\mu)(x,y) = \bm \frac{e^{\pm i\lambda|x-y|}}{4\pi|x-y|} & 0 \\
0 & \frac{e^{-\sqrt{2\mu+\lambda^2}|x-y|}}{4\pi|x-y|} \endm.
\eneq
Consider the case $\ell=0$ in~\eqref{eq:born_ser}.
Upon recombining the two $\pm$ parts the lower right-hand corner of~\eqref{eq:R0} drops
out, and one is lead to proving an oscillatory integral bound of the form
\beeq
\label{eq:oscint}
\left|\int_0^\infty e^{it\lambda^2}\lambda\chi_+(\lambda^2+\mu) \sin(\lambda |x-y|)\,d\lambda \right|
\les t^{-\frac32} |x-y|,
\eneq
To prove~\eqref{eq:oscint}, we argue as follows:
\begin{align}
& \left|\int_0^\infty e^{it\lambda^2}\lambda\chi_+(\lambda^2+\mu) \sin(\lambda |x-y|)\,d\lambda \right|
 = \half \left|\int_{-\infty}^\infty e^{it\lambda^2}\lambda\chi_+(\lambda^2+\mu) \sin(\lambda |x-y|)\,d\lambda \right| \nn \\
&\les t^{-1}|x-y| \left|\int_{-\infty}^\infty e^{it\lambda^2} \chi_+(\lambda^2+\mu) \cos(\lambda |x-y|)\,d\lambda \right|+
t^{-1} \left|\int_{-\infty}^\infty e^{it\lambda^2} \lambda\chi'_+(\lambda^2+\mu) \sin(\lambda |x-y|)\,d\lambda \right| \nn \\
&\les t^{-\frac32} |x-y| \big\| [\chi_+(\lambda^2+\mu) \cos(\lambda |x-y|)]^{\vee} \big\|_{\calM}
+ t^{-\frac32} \big\| [\lambda\chi_+'(\lambda^2+\mu) \sin(\lambda |x-y|)]^{\vee} \big\|_{\calM} \nn \\
&\les t^{-\frac32} |x-y|. \label{eq:oscint2}
\end{align}
Here we used the $L^1\to L^\infty$ estimate for the one-dimensional Schr\"odinger equation,
as well as the elementary facts
\begin{align*}
\sup_{a\in\R}\big\| [\chi_+(\lambda^2+\mu) \cos(\lambda a)]^{\vee} \big\|_{\calM} & \le C \\
\sup_{a\in\R}|a|^{-1}\big\| [\lambda\chi'_+(\lambda^2+\mu) \sin(\lambda a)]^{\vee} \big\|_{\calM} & \le C,
\end{align*}
where $\|\cdot\|_\calM$ stands for the total variation norm of measures.
The first is proved by writing it as the convolution of two measure of mass $\les 1$
uniformly in~$a$. The second is done similarly, but first write
\beeq
\label{eq:sintocos}
\sin(\lambda a) = \lambda \int_0^a \cos(\lambda \alpha)\, d\alpha.
\eneq
This yields that the
$\ell=0$ in~\eqref{eq:2mborn}  contributes $\les t^{-\frac32} \|f\|_1\|g\|_1$, as desired.
Next, we sketch the argument for the case $\ell=1$. The argument for larger $\ell$ is similar, and
we will discuss it later.
 Writing $f=\binom{f_1}{f_2}, g=\binom{g_1}{g_2}$ this term becomes
(we ignore the factor $e^{it\mu}$ as well as other constants and we write $dx=dx_0dx_1dx_2$ for simplicity)
\begin{align}
& \int_{\R^9}\int_0^\infty e^{it\lambda^2}\lambda \chi_+(\lambda^2+\mu) \sin(\lambda(|x_0-x_1|+|x_1-x_2|))\,d\lambda \, \frac{U(x_1)f_1(x_0) \bar{g}_1(x_2)}{|x_0-x_1||x_1-x_2|}\, dx\label{eq:PU1} \\
& +\int_{\R^9}\int_0^\infty e^{it\lambda^2}\lambda \chi_+(\lambda^2+\mu) \sin(\lambda |x_0-x_1|)
e^{-\sqrt{2\mu+\lambda^2}|x_2-x_1|}   \,d\lambda \,  \frac{W(x_1)f_1(x_0) \bar{g}_2(x_2)}{|x_0-x_1||x_1-x_2|}\,dx \label{eq:PW1} \\
& -\int_{\R^9}\int_0^\infty e^{it\lambda^2}\lambda \chi_+(\lambda^2+\mu) \sin(\lambda |x_2-x_1|)
e^{-\sqrt{2\mu+\lambda^2}|x_1-x_0|} \,d\lambda \, \frac{W(x_1)f_2(x_0) \bar{g}_1(x_2)}{|x_0-x_1||x_1-x_2|}\,dx \label{eq:PW2}
\end{align}
The term~\eqref{eq:PU1} can be treated by means of~\eqref{eq:oscint}. Indeed, using this bound
it reduces to
\[
 \les t^{-\frac32} \sup_{x\in\R^3}\int_{\R^3} \frac{|U(y)|}{|x-y|}\, dy \|f\|_1\|g\|_1.
\]
Hence it is enough to assume that the so-called {\em Kato norm}
\[ \|U\|_{\kato}:=\sup_{x\in\R^3}\int_{\R^3} \frac{|U(y)|}{|x-y|}\, dy < \infty\]
in order to obtain the desired decay for that term. Since we are assuming the pointwise bound
$|U(x)|\les \la x\ra^{-3-}$, the Kato norm is indeed finite.
Now consider the $\lambda$-integral in~\eqref{eq:PW1}. Extending the integral to $(-\infty,\infty)$ and
integrating by parts yields
\begin{align}
& 2i t \int_{-\infty}^\infty e^{it\lambda^2} \lambda \chi_+(\lambda^2+\mu) \sin(\lambda |x_0-x_1|)
e^{-\sqrt{2\mu+\lambda^2}|x_2-x_1|}   \,d\lambda \nn\\
& = -\int_{-\infty}^\infty e^{it\lambda^2} 2\lambda \chi_+'(\lambda^2+\mu) \sin(\lambda |x_0-x_1|)
e^{-\sqrt{2\mu+\lambda^2}|x_2-x_1|}   \,d\lambda \label{eq:I1} \\
&  -\int_{-\infty}^\infty e^{it\lambda^2} \chi_+(\lambda^2+\mu) \cos(\lambda |x_0-x_1|)
e^{-\sqrt{2\mu+\lambda^2}|x_2-x_1|}   \,d\lambda \; |x_0-x_1| \label{eq:I2} \\
&  +\int_{-\infty}^\infty e^{it\lambda^2} \chi_+(\lambda^2+\mu) \sin(\lambda |x_0-x_1|)
e^{-\sqrt{2\mu+\lambda^2}|x_2-x_1|}  \frac{2\lambda}{\sqrt{\mu+\lambda^2}} \,d\lambda \;|x_1-x_2|
\label{eq:I3}
\end{align}
The integrals in~\eqref{eq:I1} and~\eqref{eq:I2} can be treated by the same type of arguments
which lead up to~\eqref{eq:oscint2} provided we show that
\begin{align}
&\sup_{b\ge0} \big\| \int_{-\infty}^\infty e^{-b\sqrt{2\mu+\lambda^2}} e^{-i\lambda u}\, d\lambda \big\|_{\calM_u} \nn \\
&= \sup_{\mu\ge0} \big\| \int_{-\infty}^\infty e^{-\sqrt{\mu+\lambda^2}} e^{-i\lambda u}\, d\lambda \big\|_{\calM_u} <\infty. \label{eq:L1mu}
\end{align}
Now
\begin{align}
\partial_\lambda e^{-\sqrt{\mu+\lambda^2}} &= -\frac{\lambda}{\sqrt{\mu+\lambda^2}}  e^{-\sqrt{\mu+\lambda^2}} \nn \\
\partial^2_\lambda e^{-\sqrt{\mu+\lambda^2}} &= \Big(-\frac{\mu}{(\mu+\lambda^2)^{_\frac32}} +\frac{\lambda^2}{\mu+\lambda^2}\Big)e^{-\sqrt{\mu+\lambda^2}} \label{eq:der_1}
\end{align}
are both in $L^1(\R)$, and their $L^1$ norms are uniformly bounded in~$\mu>0$.
It follows that
\[ \sup_{\mu\ge0}\; (1+u^2)
\Big| \int_{-\infty}^\infty e^{-\sqrt{\mu+\lambda^2}} e^{-i\lambda u}\, d\lambda \Big|
\les 1
\]
and \eqref{eq:L1mu} holds. Therefore, arguing as in~\eqref{eq:oscint2},
\[ |\eqref{eq:I1}|+|\eqref{eq:I2}|\les t^{-\half}|x_0-x_1|.\]
To deal with~\eqref{eq:I3}, note that because of \eqref{eq:sintocos}, the same type of
arguments as before will yield
\[ |\eqref{eq:I3}| \les t^{-\frac12}|x_0-x_1| \]
provided we can show that
\begin{align}
&\sup_{b>0} \big\| \int_{-\infty}^\infty  e^{-i\lambda u} \lambda \partial_\lambda e^{-b\sqrt{2\mu+\lambda^2}}\, d\lambda \big\|_{\calM_u} \nn \\
&= \sup_{\mu>0} \big\| \int_{-\infty}^\infty  e^{-i\lambda u}\lambda \partial_\lambda e^{-\sqrt{\mu+\lambda^2}}\, d\lambda \big\|_{\calM_u} <\infty.  \label{eq:L1mu2}
\end{align}
We leave it to reader to check that
\beeq
\sup_{\mu>0} \Big[ \|\lambda \partial_\lambda e^{-\sqrt{\mu+\lambda^2}}\|_1
+ \|\partial_\lambda \lambda \partial_\lambda e^{-\sqrt{\mu+\lambda^2}}\|_1
+ \|\partial_\lambda^2 \lambda \partial_\lambda e^{-\sqrt{\mu+\lambda^2}}\|_1\Big] <\infty,
\label{eq:der_2}
\eneq
which implies that
\[ \sup_{\mu\ge0}\; (1+u^2)
\Big| \int_{-\infty}^\infty e^{-i\lambda u}\lambda \partial_\lambda e^{-\sqrt{\mu+\lambda^2}}\, d\lambda \Big|
\les 1
\]
and \eqref{eq:L1mu2} holds.
As as side remark, let us note the difference between \eqref{eq:der_1} and~\eqref{eq:der_2}.
If $\mu=0$, then the former holds because $\partial_\lambda^2 e^{-|\lambda|}$ contains
a $\delta$-measure at the origin. Hence it is not possible to increase this to three derivatives.
On the other hand, $\partial_\lambda^3 \lambda e^{-|\lambda|}$ is again a measure, which makes~\eqref{eq:der_2} hold.
Hence, we conclude that  for all $t>0$
\[
|\eqref{eq:PW1}|+|\eqref{eq:PW2}|\les t^{-\frac32} \sup_{x\in\R^3}\int_{\R^3} \frac{|W(y)|}{|x-y|}\, dy \|f\|_1\|g\|_1 \les t^{-\frac32} \|f\|_1\|g\|_1.
\]
Recall that this leads to the desired dispersive bound for the term $\ell=1$
in~\eqref{eq:2mborn}. The cases $\ell>1$ are similar. Indeed, the reader will easily check that
in the general case one arrives at oscillatory integrals of the form, cf.~\eqref{eq:PU1}, \eqref{eq:PW1},
\eqref{eq:PW2},
\[
\int_{-\infty}^\infty e^{it\lambda^2} \lambda \chi_+(\lambda^2+\mu) \sin\big(\lambda\sum_{j\in\calJ}
|x_{j+1}-x_j|\big) \exp\big( -\sqrt{2\mu+\lambda^2} \sum_{k\in\calJ^*} |x_{k+1}-x_k|\big)\,d\lambda
\]
where $\calJ\cup\calJ^* = \{0,1,\ldots,\ell\}$ is a disjoint partition with $\calJ\ne\emptyset$.
This integral is exactly of the type that we have just dealt with. Therefore, it is bounded by
\[ \les t^{-\frac32} \sum_{j\in\calJ} |x_{j+1}-x_j|. \]
Combining the oscillatory integral with the potentials that accompany it, we are lead to
estimating
\begin{align*}
&\int_{\R^{3(\ell+2)}} \sum_{j\in\calJ} |x_{j+1}-x_j| \prod_{k=1}^n \frac{|V(x_k)|}{|x_{k+1}-x_k|}\, \frac{|f(x_0)||g(x_{\ell+1})|}{|x_0-x_1|} dx \\
& \les (\ell+1) \|V\|_{\kato}^{\ell} \|f\|_1\|g\|_1.
\end{align*}
To pass to the final inequality we invoked a simple lemma from~\cite{RodSch} which says
that for any positive integer $\ell$
\begin{equation}
\nn
\sup_{x_0,x_{\ell+1}\in\R^3}\int_{\R^{3\ell}} \frac{\prod_{j=1}^\ell |V(x_j)|}{\prod_{j=0}^\ell|x_j-x_{j+1}|}\sum_{\ell=0}^\ell |x_\ell-x_{\ell+1}|\; dx_1\ldots\,dx_\ell \le (\ell+1) \|V\|_{\kato}^\ell.
\end{equation}
See Section~2 of \cite{RodSch} for the proof of this. It follows that each of the first~$2m$
terms in~\eqref{eq:2mborn} satisfy the desired dispersive bound.

In order to bound the ``remainder'' in~\eqref{eq:born_ser}, which is the
final summand containing the perturbed resolvents $R_V^{\pm}(\lambda^2+\mu)$,
we need to regard the resolvents as operators $L^{2,\sigma}\to L^{2,-\sigma}$
with $\sigma>\half$ (this is the  limiting absorption principle from Proposition~\ref{prop:lim_ap}.
Note that we only need $\sigma>\half$ rather than $\sigma>1$ since the energies are
separated from the thresholds, although this is not too important).
Moreover, not only are the resolvents bounded $L^{2,\half+}\to L^{2,-\half-}$,
but their operator norms decay like $\lambda^{-\half}$.
Note that this makes the composition of
resolvents and $V$, which appears in~\eqref{eq:born_ser}, well-defined
provided $|V(x)|\les (1+|x|)^{-1-}$ (recall that we are assuming $-3-$ decay).
Set
\[
G_{\pm,x}(\lambda^2)(x_1):= \bm e^{\mp i\lambda|x|} & 0 \\ 0 & 1 \endm R_0^{\pm}(\lambda^2+\mu)(x_1,x)
= \bm \frac{e^{\pm i\lambda(|x_1-x|-|x|)}}{4\pi|x_1-x|} & 0 \\ 0 & \frac{e^{-\sqrt{2\mu+\lambda^2}|x-x_1|}}{4\pi|x-x_1|} \endm.
\]
Let $e_1=\binom{1}{0}$ and $e_2=\binom{0}{1}$.
Removing $f,g$ from~\eqref{eq:born_ser}, we are led to proving that
\begin{align}
& \left|
\int_0^\infty e^{it\lambda^2}e^{\pm i\lambda(|x|+|y|)}\,\chi(\lambda)
\lambda \Big\la VR^{\pm}_V(\lambda^2)V (R_0^{\pm}(\lambda^2)V)^m G_{\pm,y}(\lambda^2)e_1, (R_0^{\mp}(\lambda^2)V)^m G_{\pm,x}^*(\lambda^2)e_1 \Big\ra \, d\lambda
\right|  \label{eq:dt} \\
& \les |t|^{-\frac32}, \nn
\end{align}
uniformly in $x,y\in\R^3$ as well as
\begin{align*}
& \left|
\int_0^\infty e^{it\lambda^2}e^{\pm i\lambda |x|}\,\chi(\lambda)
\lambda \Big\la VR^{\pm}_V(\lambda^2)V (R_0^{\pm}(\lambda^2)V)^m G_{\pm,y}(\lambda^2)e_2, (R_0^{\mp}(\lambda^2)V)^m G_{\pm,x}^*(\lambda^2)e_1 \Big\ra \, d\lambda
\right|   \\
& + \left|
\int_0^\infty e^{it\lambda^2}e^{\pm i\lambda |y|}\,\chi(\lambda)
\lambda \Big\la VR^{\pm}_V(\lambda^2)V (R_0^{\pm}(\lambda^2)V)^m G_{\pm,y}(\lambda^2)e_1, (R_0^{\mp}(\lambda^2)V)^m G_{\pm,x}^*(\lambda^2)e_2 \Big\ra \, d\lambda
\right|   \\
& + \left|
\int_0^\infty e^{it\lambda^2}\,\chi(\lambda)
\lambda \Big\la VR^{\pm}_V(\lambda^2)V (R_0^{\pm}(\lambda^2)V)^m G_{\pm,y}(\lambda^2)e_2, (R_0^{\mp}(\lambda^2)V)^m G_{\pm,x}^*(\lambda^2)e_2 \Big\ra \, d\lambda
\right|
 \les |t|^{-\frac32}, \nn
\end{align*}
uniformly in $x,y\in\R^3$. We first verify~\eqref{eq:dt}.
It is a simple matter to check that
the derivatives of $G_{+,x}(\lambda^2)$ satisfy the estimates
\begin{align}
\label{eq:Gest}
\sup_{x\in\R^3} \Big\| \frac{d^j}{d\lambda^j} G_{+,x}(\lambda^2)e_k\Big\|_{L^{2,-\sigma}} &\les {\la x\ra}^{-1}
\text{\ \ provided\ \ } \sigma>\frac32+j \\
\sup_{x\in\R^3} \Big\| \frac{d^j}{d\lambda^j} G_{+,x}(\lambda^2)e_k\Big\|_{L^{2,-\sigma}} &\les 1
\text{\ \ provided\ \ } \sigma>\frac12+j \nn
\end{align}
for all $j\ge 0$ and $k=1,2$.
Rewrite the integral in~\eqref{eq:dt} in the form
\begin{equation}
\nn     
I^{\pm}(t,x,y):=\int_0^\infty e^{it\lambda^2 \pm i\lambda(|x|+|y|)} a^{\pm}_{x,y}(\lambda)\, d\lambda.
\end{equation}
Then in view of the limiting absorption principle of Corollaries~\ref{cor:Hbdry}, \ref{cor:der_resolv}
and the estimate~\eqref{eq:Gest}
one concludes that $a^{\pm}_{x,y}(\lambda)$ has two
derivatives in~$\lambda$ and
\begin{align}
\label{eq:adec}
\Big|\frac{d^j}{d\lambda^j} a^{\pm}_{x,y}(\lambda)\Big| &\les (1+\lambda)^{-2+} (\la x\ra\la y\ra)^{-1}
\text{\ \ for\ \ } j = 0,1, \text{\ \ and all\ \ }\lambda>1 \\
\Big|\frac{d^2}{d\lambda^2} a^{\pm}_{x,y}(\lambda)\Big| &\les (1+\lambda)^{-2+}
\text{\ \ for all\ \ }\lambda>1 \nn
\end{align}
which in particular shows that the integral in~\eqref{eq:dt} is absolutely convergent.
This requires that one takes $m$ sufficiently large  and
that $|V(x)|\les (1+|x|)^{-\beta}$ for some $\beta>3$. The latter condition
arises as follows: Consider, for example, the case where two derivatives fall
one of the $G$-terms at the ends. Then $V$ has to compensate for $\frac52+$ powers
because of~\eqref{eq:Gest}, and also a $\frac12+$ power from
\[ \|R_0^{\pm}(\lambda^2) f\|_{X_{-\half-}} \les \lambda^{-1+} \|f\|_{X_{\half+}}. \]
Similarly with the other terms.

As far as estimating $I^{+}(t,x,y)$ is concerned, note that
on the support of~$a^{\pm}_{x,y}(\lambda)$ the phase $t\lambda^2+\lambda(|x|+|y|)$ has
no critical point.  Two integrations by parts yield the bound
$|I^{+}(t,x,y)|\les t^{-2}$.
In the case of $I^{-}(t,x,y)$ the phase $t\lambda^2-\lambda(|x|+|y|)$
has a unique critical point at $\lambda_0=(|x|+|y|)/(2t)$. If $\lambda_0\ll \lambda_1$, then
two integration by parts again yield a bound of~$t^{-2}$.
If $\lambda_0\gtrsim \lambda_1$ then the bound
$\max(|x|,|y|)\gtrsim t$ is also true, and stationary phase contributes
$t^{-\half}(\la x\ra \la y\ra)^{-1}\les t^{-\frac32}$,
as desired.  Strictly speaking, these estimates are only useful when $t>1$.
On the other hand, when $0<t<1$ there is nothing to prove since
$I^{\pm}(t,x,y) \les 1$ by~\eqref{eq:adec}.

Now consider the other three terms following~\eqref{eq:dt} which involve one or more~$e_2$.
The two integrals involving exactly one $e_2$ can be handled by the exact same argument as~\eqref{eq:dt},
the only difference being that the critical point is at~$\frac{|x|}{2t}$ or~$\frac{|y|}{2t}$.
But since~\eqref{eq:Gest} takes the same form for $e_2$ (actually a better estimate holds here,
but we ignore that since it is of no use), no other changes are needed.
Finally, concerning the integral involving two $e_2$'s: It is estimated by two integration
by parts if $t>1$, and by putting absolute values inside it $0<t<1$. Indeed, in this case
the critical point is at $\lambda=0$, which falls outside the support of the integrand.
Hence, two integration by parts give a decay of~$t^{-2}$.

The conclusion of the preceding is
that \eqref{eq:2mborn} and~\eqref{eq:born_ser} satisfy the desired dispersive bounds. Therefore,
\[ |\la e^{it\Hil}\chi_+(\Hil)P_{c}f,g\ra| \les t^{-\frac32} \|f\|_1\|g\|_1,\]
and the same bound holds for $e^{it\Hil}\chi_-(\Hil)P_{c}$.

We now deal with the contribution by those $\lambda$ which are close to $\pm\mu$.
 This requires showing that
\beeq
 \la e^{it\Hil}(1-\chi_+(\Hil))P_{c}f,g\ra = \frac{e^{it\mu}}{\pi i}\int_0^\infty \lambda e^{it\lambda^2} (1-\chi_+)(\lambda^2+\mu) \big\la [R_V^{+}(\lambda^2+\mu)-R_V^{-}(\lambda^2+\mu)] f,g \big\ra \, d\lambda \label{eq:scal_low}
\eneq
is $\les t^{-\frac32}\|f\|_1\|g\|_1$ in absolute value, and similarly for $\chi_-$.
We use the resolvent identity in the form
\begin{equation}
\label{eq:res_id}R_V^\pm(\lambda^2+\mu) = R_0^\pm(\lambda^2+\mu) - R_0^\pm(\lambda^2+\mu)V
 (I + R_0^\pm(\lambda^2+\mu)V)^{-1} R_0^\pm(\lambda^2+\mu)
\end{equation}
and write $R_0^\pm(\lambda^2+\mu) = R_0^{\pm}(\mu) + B^\pm(\lambda)$.  Then
\beeq
\label{eq:scal_exp}
[I+R_0^\pm(\lambda^2+\mu)V]^{-1} = S_0^{-1}[I+B^\pm(\lambda)VS_0^{-1}]^{-1},
\eneq
where $S_0=I+R^\pm_0(\mu)V$. In view of~\eqref{eq:R0}
\[
 R_0^{\pm}(\mu)(x,y) = \bm \frac{1}{4\pi|x-y|} & 0 \\
0 & \frac{e^{-\sqrt{2\mu}|x-y|}}{4\pi|x-y|} \endm.
\]
As far as the invertibility of $S_0$ is concerned, we note the following:
First, if $\sigma, \alpha > \frac12$, and $\sigma + \alpha > 2$, then one checks from
the explicit form of the scalar, free resolvent that
$$\sup_{\lambda}\norm[R_0^\pm(\lambda^2)][HS(\sigma,-\alpha)] \le C_{\sigma,\alpha}$$
where $HS(\sigma,-\alpha)$ refers to the Hilbert-Schmidt norm of $X_\sigma\to X_{-\alpha}$.
Hence, if $|V(\x)| \les \langle\x\rangle^{-\beta}$ for some $\beta > 3$,
it follows that the operator $R_0^\pm(\lambda)V$ is compact on the weighted space $X_\sigma(\R^3)$
for all choices of $-\frac52 \le \sigma < -\frac12$. Thus,
the invertibility of $S_0$ depends only on whether
a solution exists in $X_\sigma$ to the equation $\psi = -R_0(\mu)V\psi$.
However, if such a solution $\psi$ satisfies  $\psi\in X_\sigma$ for some
$\sigma \ge -\frac52$, then
$\psi = -R_0(\mu)V\psi \in X_\alpha$ for any choice of $\alpha < -\frac32$.
Applying this bootstrapping process again, we see that the solution $\psi$
must lie in $X_\alpha$ for all $\alpha < -\frac12$. Evidently, this would
contradict the requirement of Definition~\ref{def:edges}.

Returning to~\eqref{eq:scal_exp}, a simple estimation of the explicit kernel
\beeq
\label{eq:Bpm} B^\pm(\lambda)(\x,\y) =
\bm
\frac{e^{\pm i\lambda|\x-\y|} - 1}{4\pi |\x-\y|} & 0 \\
0 & \frac{e^{-\sqrt{2\mu+\lambda^2}|x-y|}-e^{-\sqrt{2\mu}|x-y|}}{4\pi|x-y|}
 \endm
\eneq
shows that if $|V(x)| \les \langle x\rangle^{-\beta}$ for some choice of $\beta > 3$, then
$$\lim_{\lambda\to 0} \norm[B^\pm(\lambda)VS_0^{-1}][HS(\sigma,\sigma)] = 0 $$
for all $\sigma \in (-\frac52, -\frac12)$.
For sufficiently small $\lambda^2 < \lambda_1$, it is then possible to expand
$$\tilde{B}^\pm(\lambda) := [I + B^\pm(\lambda)VS_0^{-1}]^{-1}$$
as a Neumann series in the norm $\norm[\cdot][HS(\sigma,\sigma)]$ for
all values $-\frac52 < \sigma < -\frac12$.
Moreover, the symmetry $\tilde{B}^-(\lambda) = \Btp(-\lambda)$ holds.
For ease of notation, define $\chi_0(\lambda) = (1-\chi_+)(\lambda^2+\mu)$ and
extend it as an even function of $\lambda$.
In view of \eqref{eq:scal_low} and~\eqref{eq:res_id} we wish to control the size of
\begin{align*}
& \sup_{\x,\y\in\R^3}\bigg|\int_0^\infty e^{it\lambda^2} \lambda \chi_0(\lambda)
\Big[ \big[R_0^+(\lambda^2+\mu) - R_0^-(\lambda^2+\mu)\big] \\
&\qquad  -\big[R_0^+(\lambda^2+\mu) VS_0^{-1}\Btp(\lambda)R_0^+(\lambda^2+\mu) -
R_0^-(\lambda^2+\mu)VS_0^{-1}
   \tilde{B}^-(\lambda)R_0^-(\lambda^2+\mu)\big]\Big](\x,\y)\, d\lambda \bigg|
\end{align*}
which is
\begin{align}
& \les  \sup_{\x,\y\in\R^3}
\Big|\int_{-\infty}^\infty e^{it\lambda^2}\lambda \chi_0(\lambda)
\frac{e^{i\lambda|\x-\y|}}{4\pi |\x-\y|} d\lambda \Big| \label{eq:freie}\\
& + \sup_{\x,\y\in\R^3} \Big|\int_{-\infty}^\infty e^{it\lambda^2}\lambda
   \iint_{\R^6} \frac{U(x_4)e^{i\lambda|\y-x_4|}}{|\y-x_4|}
  \la e_1,\big(S_0^{-1} (\chi_0\Btp)(\lambda)(x_4,x_1)\big) e_1\ra
  \frac{e^{i\lambda|\x-x_1|}}{|\x-x_1|}\, dx_1 dx_4 d\lambda \Big| \label{eq:V11} \\
& + \sup_{\x,\y\in\R^3} \Big|\int_{-\infty}^\infty e^{it\lambda^2}\lambda
   \iint_{\R^6} \frac{W(x_4)e^{i\lambda|\y-x_4|}}{|\y-x_4|}
  \la e_1,\big(S_0^{-1} (\chi_0\Btp)(\lambda)(x_4,x_1)\big)e_2\ra
  \frac{e^{-{\sqrt{2\mu+\lambda^2}}|\x-x_1|}}{|\x-x_1|}\, dx_1 dx_4 d\lambda \Big| \label{eq:V12} \\
& + \sup_{\x,\y\in\R^3} \Big|\int_{-\infty}^\infty e^{it\lambda^2}\lambda
   \iint_{\R^6} \frac{W(x_4)e^{-\sqrt{2\mu+\lambda^2}|\y-x_4|}}{|\y-x_4|}
  \la e_2\big(S_0^{-1} (\chi_0\Btp)(\lambda)(x_4,x_1)\big), e_1\ra
  \frac{e^{i\lambda|\x-x_1|}}{|\x-x_1|}\, dx_1 dx_4 d\lambda \Big| \label{eq:V21} \\
& + \sup_{\x,\y\in\R^3} \Big|\int_{-\infty}^\infty e^{it\lambda^2}\lambda
   \iint_{\R^6} \frac{U(x_4)e^{-{\sqrt{2\mu+\lambda^2}}|\y-x_4|}}{|\y-x_4|}
  \la e_2,\big(S_0^{-1} (\chi_0\Btp)(\lambda)(x_4,x_1)\big)e_2\ra
  \frac{e^{-\sqrt{2\mu+\lambda^2}|\x-x_1|}}{|\x-x_1|}\, dx_1 dx_4 d\lambda \Big| \label{eq:V22}
\end{align}
The first term~\eqref{eq:freie} is simply the low-energy part of the free Schr\"odinger
evolution, which is known to be dispersive.
The second term~\eqref{eq:V11} can be integrated by parts once, leaving
\begin{equation} \label{eq:IBP}
\sup_{\x,\y \in\R^3} \frac1{2t} \Big|
 \int_{-\infty}^\infty e^{it\lambda^2} \iint_{\R^6}\dfrac{d}{d\lambda}\Big[
 \frac{U(x_4)e^{i\lambda|\y-x_4|}}{|\y-x_4|}
  \big(S_0^{-1} (\chi_0 \Btp)(\lambda)(x_4,x_1)\big)
 \frac{e^{i\lambda|\x-x_1|}}{|\x-x_1|} \Big] \,dx_1 dx_4 d\lambda  \Big|
\end{equation}
to be controlled. Note that we have dropped  $e_1$ on both sides
of the matrix operator in the middle. This does no harm, as long
as the absolute value on the outside is interpreted entry-wise.
The same comment is in effect for the remainder of the proof. The
other terms~\eqref{eq:V12}, \eqref{eq:V21}, and~\eqref{eq:V22} are
treated similarly to~\eqref{eq:V11}. In fact, we verified
in~\eqref{eq:L1mu} that for $a>0$
\[ \int_{-\infty}^\infty e^{i\tau\lambda} e^{-a\,\sqrt{2\mu+\lambda^2}} \, d\lambda =: \nu_a(d\tau)\]
is a measure with mass $\sup_{a>0} \|\nu_a\| < \infty$. This simple fact allows
one to use the same argument which is sketched here for~\eqref{eq:V11}
in the other three cases as well up to some obvious modifications.
We now return to~\eqref{eq:IBP}, which is essentially identical with the analogous term
arising in the scalar case treated in~\cite{GolSch}. Since we see no reason to
repeat the details verbatim, we will provide a sketch and refer the reader to~\cite{GolSch}
for more details.
  Consider the term where $\frac{d}{d\lambda}$ falls on
$\Btp(\lambda)$.  The others will be similar.
Using Parseval's identity, and the fact that
$\norm[(e^{it(\cdot)^2})^\wedge(u)][L^\infty(u)] = C t^{-1/2}$,
this is less than
$$\sup_{\x,\y\in\R^3} \frac1{t^{3/2}} \int_{-\infty}^\infty
  \Big|  \iint_{\R^6} \frac{U(x_4)}{|\y-x_4|} S_0^{-1}
  \big[\chi_0(\Btp)'\big]^\vee
  \big(u+|\y-x_4| + |\x-x_1|\big)(x_4,x_1) \frac{1}{|\x-x_1|}\, dx_1 dx_4
  \Big| \, du.$$
If the absolute value is taken inside the inner integral, then Fubini's
theorem may be used to exchange the order of integration to obtain
\begin{equation*}
\begin{aligned}
\sup_{\x,\y\in\R^3} \frac1{t^{3/2}} \iint_{\R^6} \int_{-\infty}^\infty
\frac{|U(x_4)|}{|\y-x_4|}  \Big| S_0^{-1}\big[\chi_0(\Btp)'\big]^\vee
\big(u + |\y-x_4| + |\x- x_1|\big)
{\scriptstyle (x_4, x_1)} \Big|
\frac{1}{|\x-x_1|} \, du\, dx_1 dx_4
\\
\le \sup_{\x,\y\in\R^3} \frac1{t^{3/2}}
\Big\|\frac{|U(\cdot)|}{|\y- \cdot|}\Big\|_{L^{2,2+}}
 \ \big\|{\textstyle \int} |S_0^{-1} [\chi_0(\Btp)']
  ^\vee(u)|du \big\|_{OP(-1-,-2-)} \ \big\||\x-\cdot|^{-1}\big\|_{L^{2,-1-}},
\end{aligned}
\end{equation*}
where $OP(-1-,-2-)$ stands for the operator norm from $X_{-1-}\to X_{-2-}$.
The two norms at the ends of the last line are easily seen to be uniformly bounded
in $\x,\y \in\R^3$.
It therefore only remains
to control the size of
$$ \big\|{\textstyle \int} |S_0^{-1} [\chi_0(\Btp)']^\vee(u)|du
  \big\|_{OP(-1-,-2-)}.$$
Minkowski's Inequality allows us to bring the norm inside the integral.
Recall that $S_0^{-1}$ is a bounded operator on $L^{2,-2-}$. Furthermore, it
is an easy matter to check that the operator $|S_0^{-1}|$ whose kernel is
the absolute value of the kernel of $S_0^{-1}$, is also a bounded operator on $L^{2,-2-}$.
The problem then reduces to showing that
\begin{equation} \label{eq:desired}
 \int_{-\infty}^\infty
 \norm[[\chi_0(\Btp)']^\vee(u)][OP(-1-,-2-)]\, du  < \infty
\end{equation}
provided the support of $\chi_0$ is sufficiently small. The
operators $\Btp(\lambda)$ are defined by the convergent Neuman
series
$$\Btp(\lambda) \ = \ [I + \Bp(\lambda)VS_0^{-1}]^{-1}\  =\  \sum_{n=0}^\infty
 \big(- \Bp(\lambda)VS_0^{-1}\big)^n.$$
Exploiting the explicit form of the kernel of $\Bp$, see~\eqref{eq:Bpm},
it is possible to control the Fourier transform in $\lambda$ of each term
in this Neuman series in the appropriate weighted $L^2$ spaces, leading to~\eqref{eq:desired}
upon summation.
For these details we refer the reader to the end of the paper~\cite{GolSch}.
\end{proof}

Finally, we discuss Strichartz estimates. The usual derivation for Strichartz estimates
involves $TT^*$ arguments where $(Tf)(t,x)=(e^{-itH}f)(x)$. This relies on the unitarity of
the evolution, since one wants $TT^*F(t,x)=\int_{-\infty}^\infty (e^{-i(t-s)H} F(s,\cdot))(x)\,ds$.
In the system case, this cannot be done. We therefore rely on a different approach which
is perturbative in nature. It uses Kato's notion of an $\Hil_0$-smooth and $\Hil$-smooth
operator, and originates in~\cite{RodSch}. In addition, we use the following lemma,
which is due to Christ-Kiselev~\cite{CrKi}. See also Sogge, Smith~\cite{SoSm}.

\begin{lemma}
Let $X, Y$ be Banach spaces and let $K(t,s)$ be the kernel of the
operator \[ K: L^p([0,T]; X)\to L^q([0,T]; Y).\] Denote by $\|K\|$ the
operator norm of $K$. Define the lower diagonal operator
\[ \tilde K:\,L^p([0,T]; X)\to L^q([0,T]; Y) \]
to be
$$
\tilde Kf(t) = \il_0^t K(t,s) f(s)\,ds
$$
Then the operator $\tilde K$ is bounded from $L^p([0,T]; X)$ to $L^q([0,T]; Y)$
and its norm $\|\tilde K\|\le c \|K\|$, provided that $p<q$.
\end{lemma}

Now we can state the Strichartz estimates.

\begin{cor}
\label{cor:strich}
Under the same assumptions as in Theorem~\ref{thm:disp}, one has the
Strichartz estimates
\begin{align}
 \|e^{-it\Hil} P_c f\|_{L^r_t(L^p_x)} &\le C \|f\|_{L^2} \label{eq:Strich1} \\
 \Big\|\int_0^t e^{-i(t-s)\Hil} P_c F(s)\,ds \Big\|_{L^r_t(L^p_x)} &\le C \|F\|_{L^{a'}_t(L^{b'}_x)} \label{eq:Strich2},
\end{align}
provided $(r,p), (a,b)$ are admissible, i.e., $2<r\le\infty$ and $\frac{2}{r}+\frac{3}{p}=\frac32$ and
the same for $(a,b)$.
\end{cor}
\begin{proof}
Let ($\calS$ for ``Strichartz'')
\[ (\calS F)(t,x)=\int_0^t (e^{-i(t-s)\Hil}P_c\,F(s,\cdot))(x)\, ds.\]
In this proof it will be understood that all times are $\ge0$. Then by~\eqref{eq:L2stable},
\[ \|\calS F\|_{L^\infty_t(L^2_x)} \les \|F\|_{L^1_t(L^2_x)}, \]
and more generally, by the usual fractional integration argument based on Theorem~\ref{thm:disp},
\beeq
\label{eq:S1} \|\calS F\|_{L^r_t(L^p_x)} \les \|F\|_{L^{r'}_t(L^{p'}_x)}
\eneq
for any admissible pair $(r,p)$. In the unitary case this implies~\eqref{eq:Strich1} via
a $TT^*$ argument, but this reasoning does not apply here. Instead, we rely on a Kato theory
type approach as in \cite{RodSch}, Section~4. Since $\Hil=\Hil_0+V$, Duhamel's formula yields
\beeq
\label{eq:du1}
 e^{-it\Hil}P_c = e^{-it\Hil_0}P_c -i\int_0^t e^{-i(t-s)\Hil_0}Ve^{-is\Hil}P_c \, ds.
\eneq
Writing $V=\tilde{M}\tilde{M}^{-1}V$, where $\tilde{M}$ is as in~\eqref{eq:Mtil}, observe firstly that
\[
 \Big\| \int_0^\infty e^{-i(t-s)\Hil_0} \tilde{M} g(s)\, ds\Big\|_{L^r_t(L^p_x)} \les
\Big\| \int_0^\infty e^{is\Hil_0} \tilde{M} g(s)\Big\|_{L^2}
\les \|g\|_{L^2_s(L^2_x)},
\]
where the last inequality is the dual of the smoothing bound
\[ \int_0^\infty \Big\| \tilde{M}e^{-is\Hil_0^*} \psi\Big\|_2^2 \, ds \les \|\psi\|_2^2. \]
Now one applies the Christ-Kiselev lemma to conclude that
\[
\Big\| \int_0^t e^{-i(t-s)\Hil_0} \tilde{M} g(s)\, ds\Big\|_{L^r_t(L^p_x)} \les \|g\|_{L^2_s(L^2_x)}
\]
for any admissible pair $(r,p)$. Hence, continuing in~\eqref{eq:du1}, one obtains (using
that $\|P_c f\|_2\les \|f\|_2$)
\[
\|e^{-it\Hil}P_c f \|_{L^r_t(L^p_x)} \les \|f\|_2 +
\Big\| \tilde{M}^{-1}V e^{-is\Hil}P_c f\Big\|_{L^2_s(L^2_x)}.
\]
It remains to show that $\tilde{M}^{-1}V$ is $\Hil P_c$-smoothing, i.e.,
\beeq
\label{eq:smft}
\Big\| \tilde{M}^{-1}V e^{-is\Hil}P_c f\Big\|_{L^2_s(L^2_x)} \les \|f\|_2.
\eneq
Taking the Fourier transform in $s$, shows that \eqref{eq:smft} is equivalent with
\beeq
\label{eq:Hsmooth2}
 \int_{-\infty}^\infty \| \tilde{M}^{-1}V [P_c(\Hil-\lambda-i0)P_c]^{-1}P_cf\|_2^2 \, d\lambda
\les \|f\|_2^2.
\eneq
However, this was already shown in~\eqref{eq:Hilsmooth}. Indeed, setting $\phi=P_c f$ and restricting
$\lambda$ to $|\lambda|\ge\mu$ leads to the same expression as in~\eqref{eq:Hilsmooth}.
On the other hand, if $|\lambda|\le \mu$, then one simply notes that
\[ \sup_{|\lambda|\le\mu} \|[P_c(\Hil-\lambda-i0)P_c]^{-1}\|_{2\to2} \les 1,\]
so that the entire integral in~\eqref{eq:Hsmooth2} is controlled. The conclusion is that
\[  \|e^{-it\Hil}P_c f \|_{L^r_t(L^p_x)} \les \|f\|_2 \]
for any admissible $(r,p)$, which is~\eqref{eq:Strich1}.
The proof of \eqref{eq:Strich2} is now the usual interpolation argument. Indeed,
in view of the preceding one has the following bounds on~$\calS$ for any admissible
pair $(r,p)$:
\begin{align}
\calS: &L^1_t(L^2_x) \to L^r_t(L^p_x) \label{eq:St1}\\
\calS: &L^{r'}_t(L^{p'}_x) \to L^r_t(L^p_x) \label{eq:St2}\\
\calS: &L^{r'}_t(L^{p'}_x) \to L^\infty_t(L^2_x).  \label{eq:St3}
\end{align}
These estimates arise as follows: \eqref{eq:St2} is exactly~\eqref{eq:S1}, whereas~\eqref{eq:St1}
follows from~\eqref{eq:Strich1} by means of Minkowski's inequality. Finally, \eqref{eq:St3}
is dual to the bound
\beeq
\label{eq:again} \Big\|\int_t^\infty e^{i(t-s)\Hil^*}\tilde{P}_c G(s)\,ds \Big\|_{ L^r_t(L^p_x)} \les
\|G\|_{L^1_t(L^2_x)}.
\eneq
Here $\tilde{P}_c$ corresponds to $\Hil^*$ in the same way that $P_c$ does to $\Hil$.
In particular, one has
\[ \|e^{-it\Hil^*}\tilde{P}_c\|_{1\to\infty}\les t^{-\frac32}\]
and therefore, \eqref{eq:again} is derived be the same methods as~\eqref{eq:St1}.
It is important to notice  that $P_c^*=\tilde{P}_c$ which is essential for the duality argument here.
This can be seen, for example, by writing the Riesz projections onto (generalized) eigenspaces
as contour integrals around circles surrounding the eigenvalues. Since the (complex) eigenvalues
always come in pairs, the adjoints have the desired property.
Interpolating between \eqref{eq:St1} and \eqref{eq:St2} yields~\eqref{eq:Strich2}
for the range $a'\le r'$ or $a\ge r$, whereas interpolating between \eqref{eq:St1} and
\eqref{eq:St2} yields~\eqref{eq:Strich2} in the range $a\le r$.
\end{proof}

Finally, we introduce derivatives into the  estimates
of Theorems~\ref{thm:L2stable}, \ref{thm:disp} and Corollary~\ref{cor:strich}.

\begin{cor}
\label{cor:inter}
Under the same assumptions as in Theorem~\ref{thm:disp},
\[ \big\|e^{it\Hil} P_c f\big\|_{W^{k,p'}(\R^3)} \les t^{-\frac32(\frac{1}{p}-\frac{1}{p'})} \|f\|_{W^{k,p}(\R^3)}\]
for  $0\le k\le 2$ and $1<p\le2$.
\end{cor}
\begin{proof}
The case $k=0$ is obtained by interpolating between Theorems~\ref{thm:L2stable} and~\ref{thm:disp}
and holds for the entire range $1\le p\le2$. We need to require $p>1$ only for the derivatives.
If $a$ is sufficiently large, then
\[ (\Hil-ia)^{-1}: L^2\times L^2 \to W^{2,2}\times W^{2,2} \]
is an isomorphism. More generally,
\[ (\Hil-ia)^{-\half}: L^p\times L^p \to W^{2,p}\times W^{2,p}\]
is an isomorphism for $1<p<\infty$. This can be seen from the resolvent
identity
\[ (\Hil-ia)^{-1} = (\Hil_0-ia)^{-1}[1+V(\Hil_0-ia)^{-1}]^{-1},\]
since $\|V\|_\infty <\infty$ implies that
\[ \|V(\Hil_0-ia)^{-1}\|_{p\to p}<\half\]
if $a$ is large enough, and because
\[ (\Hil_0-ia)^{-\half}: L^p\times L^p \to W^{2,p}\times W^{2,p}\]
for any $a\ne0$ as an isomorphism. Hence,
\begin{align*}
\|\Laplace e^{it\Hil}P_c f\|_{p'} &\les \|(\Hil-ia) e^{it\Hil} f\|_{p'} =  \|e^{it\Hil}(\Hil-ia)f\|_{p'}\\
& \les t^{-\frac32(\frac{1}{p}-\frac{1}{p'})} \|(\Hil-ia)f\|_p
\les t^{-\frac32(\frac{1}{p}-\frac{1}{p'})} \|f\|_{W^{2,p}(\R^3)}.
\end{align*}
This gives the case $k=2$ of the lemma, whereas $k=1$ follows
by interpolating between $k=0$ and $k=2$.
\end{proof}

And now the same for the Strichartz estimates.

\begin{cor}
\label{cor:strich_der}
Under the same assumptions as in Corollary~\ref{cor:strich}, one has the
Strichartz estimates
\begin{align}
 \|e^{-it\Hil} P_c f\|_{L^r_t(W^{k,p}_x)} &\le C \|f\|_{W^{k,2}} \label{eq:Strich1_der} \\
 \Big\|\int_0^t e^{-i(t-s)\Hil} P_c F(s)\,ds \Big\|_{L^r_t(W^{k,p}_x)} &\le C \|F\|_{L^{a'}_t(W^{k,b'}_x)} \label{eq:Strich2_der},
\end{align}
provided $(r,p), (a,b)$ are admissible, i.e., $2<r\le\infty$ and $\frac{2}{r}+\frac{3}{p}=\frac32$ and
the same for $(a,b)$. Here $k$ is an integer, $0\le k\le2$.
\end{cor}
\begin{proof}
The case $k=0$ is just Corollary~\ref{cor:strich}.
As in the previous proof, we rely on the fact that (because of $\|V\|_\infty<\infty$),
\[ \|\Laplace f\|_q\les \|(\Hil-ia) f\|_q \]
for any $1<q<\infty$. Hence,
\begin{align*}
\|e^{-it\Hil} P_c f\|_{L^r_t(W^{2,p}_x)} &\les  \|(\Hil-ia) e^{-it\Hil} P_c f\|_{L^r_t(L^p_x)}
= \| e^{-it\Hil} P_c (\Hil-ia) f\|_{L^r_t(L^p_x)}  \\
&\les \|(\Hil-ia)f\|_2 \les \|f\|_{W^{2,2}},
\end{align*}
which is \eqref{eq:Strich1_der} for $k=2$. Similarly, one proves~\eqref{eq:Strich2_der}
for $k=1$. The case $k=1$ is then obtained by interpolation.
\end{proof}

{\em Acknowledgments:} The author is grateful to Avy Soffer, Igor
Rodnianski, and Joachim Krieger for their interest as well as
useful comments on a preliminary version of this work. In
particular, he wishes to thank Igor Rodnianski for pointing out
that in~\cite{RSS2} the contraction can be carried out with
$t^{-2}$ control on the derivatives of the paths. He is also
grateful to Israel Michael Sigal for an invitation to give a talk
at Notre Dame University, and he thanks him, as well as Avy
Soffer, and Francois Ledrappier for helpful criticism.

\bibliographystyle{amsplain}

\begin{thebibliography}{99}

\bibitem[Agm1]{Agm1} Agmon, S. {\em Spectral properties of Schr\"odinger
operators and scattering theory.}
Ann.\ Scuola Norm.\ Sup.\ Pisa Cl.\ Sci. (4) 2 (1975), no.~2, 151--218.

\bibitem[Agm2]{Agm2} Agmon, S. {\em  Lectures on exponential decay of solutions of second-order
elliptic equations: bounds on eigenfunctions of $N$-body Schr\"odinger operators.}
Mathematical Notes, 29. Princeton University Press, Princeton, NJ;
University of Tokyo Press, Tokyo, 1982.

\bibitem[BerCaz]{BerCaz} Berestycki, H., Cazenave, T.
{\em Instabilit\'e des \'etats stationnaires dans les \'equations de Schr\"odinger et de Klein-Gordon non lin\'eaires.}
C.\ R.\ Acad.\ Sci.\ Paris S\'er.~I Math. 293 (1981), no. 9, 489--492.

\bibitem[BerLio]{BerLio} H. Berestycki, P.L. Lions, {\em Existence d'oudes
solitaires daus les problemes nonlineares du type Klein-Gordon.}
C.R.\ Acad.\ Sci.~288 (1979), no.~7, 395--398.

\bibitem[BusPer1]{BP1} Buslaev, V.\ S., Perelman, G.\ S. {\em Scattering for the
nonlinear Schr\"odinger equation: states that are close to a soliton.} (Russian)
Algebra i Analiz  4  (1992),  no.~6, 63--102;  translation in  St.\ Petersburg Math.\ J.~4
(1993),  no.~6, 1111--1142.

\bibitem[BusPer2]{BP2} Buslaev, V.\ S., Perelman, G.\ S. {\em On the stability of solitary waves
for nonlinear Schr\"odinger equations.  Nonlinear evolution equations,}  75--98,
Amer.\ Math.\ Soc.\ Transl.\ Ser.~2, 164, Amer.\ Math.\ Soc., Providence, RI, 1995.

\bibitem[CazLio]{CazLio} Cazenave, T., Lions, P.-L. {\em Orbital stability
of standing waves for some nonlinear Schr\"odinger equations.} Comm.\ Math.\ Phys.~85 (1982), 549--561

\bibitem[CriKis]{CrKi} Christ, M., Kiselev, A. {\em Maximal functions associated with
filtrations,} J.\ Funct.\ Anal.\ 179 (2001), 409-425.

\bibitem[Cof]{Cof} Coffman, C. V. {\em Uniqueness of positive solutions of
$\Laplace u - u+u^{3}=0$ and a variational characterization of other
solutions.} Arch.\ Rat.\ Mech.\ Anal.~46 (1972), 81--95.

\bibitem[CosSof]{CosSof} Costin, O., Soffer, A.{\em Resonance theory for
Schrödinger operators. } Comm.\ Math.\ Phys.~224 (2001), no.~1,
133--152.

\bibitem[Cuc]{Cuc} Cuccagna, S.  {\em Stabilization of solutions to nonlinear Schr\"odinger equations.}
Comm.\ Pure Appl.\ Math.~54  (2001),  no.~9, 1110--1145.

\bibitem[CucPelVou]{CPV}  Cuccagna, S., Pelinovsky, D., Vougalter,
V. {\em Spectra of positive and negative energies in the
linearized NLS problem}. Preprint, 2003.

\bibitem[CucPel]{CP}  Cuccagna, S., Pelinovsky, D.
V. {\em Bifurcation from the end points of the essential spectrum
in the linearized NLS problem}. Preprint, 2004.


\bibitem[ErdSch]{ErdSch} Erdogan, B., Schlag, W. {\em Dispersive estimates in the presence of
resonances and/or eigenvalues at thresholds}, in preparation.




\bibitem[GolSch]{GolSch} Goldberg, M., Schlag, W. {\em Dispersive estimates
for Schr\"odinger operators in dimensions one and three.} To appear in
Comm.\ Math.\ Phys., preprint 2003.


\bibitem[Gri]{Grill} Grillakis, M. {\em
Analysis of the linearization around a critical point of an infinite dimensional
Hamiltonian system.} Comm.\ Pure Appl.\ Math.~41 (1988), no.~6, 747--774.

\bibitem[GriShaStr1]{GSS1} Grillakis, M., Shatah, J., Strauss, W.
{\em Stability theory of
solitary waves in the presence of symmetry. I.}
J.\ Funct.\ Anal.~74  (1987),  no.~1, 160--197.

\bibitem[GriShaStr2]{GSS2} Grillakis, M., Shatah, J., Strauss, W.
{\em Stability theory of
solitary waves in the presence of symmetry. II.}
J.\ Funct.\ Anal.~94  (1990), 308--348.

\bibitem[HisSig]{HS} Hislop, P.\ D., Sigal, I.\ M.
{\em Introduction to spectral theory. With applications to Schr\"odinger operators.}
Applied Mathematical Sciences, 113.\ Springer-Verlag, New York, 1996.


\bibitem[KriSch]{KS} Krieger, J., Schlag, W. {\em Stable manifolds
for supercritical focusing NLS in one dimension}, in preparation.

\bibitem[Kwo]{Kwo} Kwong, M.\ K. {\em Uniqueness of positive solutions
of $\Laplace u - u +u^{p}=0$ in $\R^{n}$.} Arch.\ Rat.\ Mech.\ Anal.~65
(1989), 243--266.



\bibitem[Per1]{Pe1} Perelman, G. {\em Some Results on the Scattering of
Weakly Interacting Solitons for Nonlinear Schr\"odinger Equations}
in "Spectral theory, microlocal analysis, singular manifolds",
Akad.\ Verlag (1997), 78--137.

\bibitem[Per2]{Pe2} Perelman, G. {\em On the formation of singularities in solutions
of the critical nonlinear Schrödinger equation.} Ann.\ Henri
Poincar\'e 2 (2001), no.~4, 605--673.

\bibitem[Per3]{Pe3} Perelman, G. {\em Asymptotic Stability of
multi-soliton solutions for nonlinear Schr\"odinger equations.}
preprint 2003.

\bibitem[PilWay]{PW} Pillet, C.\ A., Wayne, C.\ E. {\em Invariant
manifolds for a class of dispersive, Hamiltonian, partial differential
equations.} J. Diff. Eq. 141 (1997), no. 2, 310--326


\bibitem[ReeSim4]{RS4}  Reed, M., Simon, B. {\em Methods of modern mathematical
physics. IV.} Academic Press [Harcourt Brace Jovanovich, Publishers],
New York-London, 1979.

\bibitem[RodSch]{RodSch} Rodnianski, I.,  Schlag, W. {\em Time decay for solutions of
Schr\"odinger equations with rough and time-dependent potentials.} to appear in Invent.\ Math.

\bibitem[RodSchSof1]{RSS1} Rodnianski, I., Schlag, W., Soffer, A.
{\em Dispersive Analysis of Charge Transfer Models}, preprint 2002, to appear in CPAM

\bibitem[RodSchSof2]{RSS2} Rodnianski, I., Schlag, W., Soffer, A.
{\em Asymptotic stability of $N$-soliton states of NLS}, preprint 2003, submitted to CPAM



\bibitem[Sch]{schlagODE} Schlag, W. {\em A remark on linear perturbations of hyperbolic ODEs in the plane.} preprint, 2004.
available online at www.its.caltech.edu/\~{}schlag

\bibitem[Sha]{Sha} Shatah, J. {\em Stable standing waves of nonlinear
Klein-Gordon equations.} Comm.\ Math.\ Phys.~91 (1983), no.~3, 313--327

\bibitem[ShaStr]{ShaStr} Shatah, J., Strauss, W.
{\em Instability of nonlinear bound states.}
  Comm.\ Math.\ Phys.~100  (1985),  no.~2, 173--190.

\bibitem[SogSmi]{SoSm} Smith, H., Sogge, C. {\em Global Strichartz estimates for nontrapping
perturbations of the Laplacian.}
Comm.\ Partial Differential Equations  25  (2000),  no.~11-12, 2171--2183.

\bibitem[SofWei1]{SofWei1} Soffer, A., Weinstein, M.
{\em Multichannel nonlinear scattering for nonintegrable equations.}
Comm.\ Math.\ Phys.~133 (1990), 119--146

\bibitem[SofWei2]{SofWei2} Soffer, A., Weinstein, M.
{\em Multichannel nonlinear scattering, II. The case of anysotropic
potentials and data.} J.\ Diff.\ Eq.~98 (1992), 376--390



\bibitem[SulSul]{SulSul} Sulem, C., Sulem, P.-L.
{\em The nonlinear Schr\"odinger equation. Self-focusing and wave collapse.}
 Applied Mathematical Sciences, 139. Springer-Verlag, New York, 1999.


\bibitem[TsaYau]{TY} Tsai, T.-P., Yau, H.-T. {\em Stable directions for excited
states of nonlinear Schr\"odinger equations.} Comm.\ Partial
Differential Equations 27 (2002), no.~11-12, 2363--2402.

\bibitem[Wei1]{Wei1} Weinstein, Michael I.
{\em Modulational stability of ground states of nonlinear Schr\"odinger equations.}
  SIAM J.\ Math.\ Anal.~16  (1985),  no.~3, 472--491.

\bibitem[Wei2]{Wei2} Weinstein, Michael I.
{\em Lyapunov stability of ground states of nonlinear dispersive evolution equations.}
Comm.\ Pure Appl.\ Math.~39  (1986),  no.~1, 51--67.

\end{thebibliography}

\noindent

\medskip\noindent
\textsc{Division of Astronomy, Mathematics, and Physics,
253-37 Caltech,\\ Pasadena, CA 91125, U.S.A.}\\
{\em email: }\textsf{\bf schlag@its.caltech.edu}

\end{document}